\documentclass[letterpaper, 11pt,  reqno]{amsart}
\usepackage{graphicx} 
\usepackage{amsmath,amssymb,amscd,amsthm,amsxtra,esint}
\usepackage{color}
\usepackage[margin=1.1in,marginparwidth=1.5cm, marginparsep=0.5cm]{geometry} 
\usepackage{tikz-cd}
\tikzset{smalltext/.style={"\textup{\small #1}" description}}

\usepackage[markup=default]{changes}
\definechangesauthor[name={Reviewer 1}, color=orange]{R1} 
\definechangesauthor[name={Reviewer 2}, color=blue]{R2} 
\definechangesauthor[name={Reviewer 3}, color=green]{R3}

\usepackage[implicit=true]{hyperref}

\allowdisplaybreaks[4]

\definecolor{gr}{rgb}   {0.,   0.69,   0.23 }
\definecolor{bl}{rgb}   {0.,   0.5,   1. }
\definecolor{mg}{rgb}   {0.85,  0.,    0.85}
\definecolor{yl}{rgb}   {0.8,  0.7,   0.}
\definecolor{or}{rgb}  {0.7,0.2,0.2}

\DeclareMathOperator*{\supp}{supp}

\newcommand{\noi}{\noindent}
\newcommand{\R}{\mathbb{R}}
\newcommand{\T}{\mathbb{T}}
\newcommand{\Z}{\mathbb{Z}}
\newcommand{\N}{\mathbb{N}}

\newcommand{\E}{\mathbb{E}}

\newcommand{\B}{\mathcal{B}}
\newcommand{\C}{\mathcal{C}}
\newcommand{\Id}{\textup{Id}}

\newcommand{\al}{\alpha}
\newcommand{\be}{\beta}
\newcommand{\ta}{\theta}
\newcommand{\s}{\sigma}
\newcommand{\eps}{\varepsilon}
\newcommand{\om}{\omega}
\newcommand{\Om}{\Omega}
\newcommand{\Dl}{\Delta}
\newcommand{\dl}{\delta}
\newcommand{\ld}{\lambda}
\newcommand{\nb}{\nabla}

\newcommand{\dt}{\partial_t}
\newcommand{\ind}{\mathbf 1}
\newcommand{\ft}{\widehat}
\newcommand{\wt}{\widetilde}
\newcommand{\les}{\lesssim}
\newcommand{\ges}{\gtrsim}
\newcommand{\cj}{\overline}

\newcommand{\jb}[1]{\langle #1 \rangle}
\newcommand{\jbb}[1]{[\hspace{-0.6mm}[ #1 ]\hspace{-0.6mm}]}
\newcommand{\wick}[1]{:\!{#1}\!:}
\newcommand{\deff}{\stackrel{\textup{def}}{=}}

\let\Re=\undefined\DeclareMathOperator*{\Re}{Re}
\let\Im=\undefined\DeclareMathOperator*{\Im}{Im}

\newtheorem{theorem}{Theorem}[section]
\newtheorem{lemma}[theorem]{Lemma}
\newtheorem{proposition}[theorem]{Proposition}

\newtheorem{remark}[theorem]{Remark}

\newtheorem*{ackno}{Acknowledgements}

\numberwithin{equation}{section}
\numberwithin{theorem}{section}

\makeatletter
\@namedef{subjclassname@2020}{%
  \textup{2020} Mathematics Subject Classification}
\makeatother

\title[Large torus limit of global dynamics of the 2d dispersive Anderson model]{Large torus limit of global dynamics of the two-dimensional dispersive Anderson model}

\author[R.~Liu and N.~Tzvetkov]{Ruoyuan Liu and Nikolay Tzvetkov}

\address{
Ruoyuan Liu, Mathematical Institute\\
University of Bonn\\
Endenicher Allee 60\\
53115\\
Bonn\\
Germany}

\email{ruoyuanl@math.uni-bonn.de}

\address{
Nikolay Tzvetkov, Ecole Normale Sup\'erieure de Lyon\\ 
UMPA\\ 
UMR CNRS-ENSL 5669, 46, all\'ee d'Italie\\ 69364-Lyon Cedex 07\\
France}

\email{nikolay.tzvetkov@ens-lyon.fr}


\begin{document}

\baselineskip = 14pt


\begin{abstract} 
We continue the study of the two-dimensional dispersive Anderson model (DAM), i.e.~the nonlinear Schr\"odinger equation with a multiplicative spatial white noise. For this model, global well-posedness on the periodic domain was established by Visciglia and the second author (2023), and global well-posedness on the full space was established by Debussche, Visciglia, and the authors (2024). We show that, under suitable initial conditions and suitable periodization procedure of the noise, the periodic global dynamics of the DAM converges in spaces of local domains to that of the DAM on the full space as the period goes to infinity. In order to control the growth of the noise and obtain a priori bounds for solutions independent of the periodicity, we introduce periodic weights and construct weighted function spaces on periodic domains. In Appendix, we also discuss the same problem for the parabolic Anderson model.
\end{abstract}


\maketitle

\tableofcontents

\section{Introduction}

\subsection{Dispersive Anderson model}
\label{SUB:DAM}

We consider the following Cauchy problem for the nonlinear Schr\"odinger equation (NLS) with a multiplicative spatial white noise on $\R^2$:
\begin{align}
\begin{cases}
i \dt u = \Dl u + \xi u - \ld |u|^{p - 1} u \\
u|_{t = 0} = u_0,
\end{cases}
\label{NLSAnd}
\end{align}

\noi
where $p > 1$, $\ld \geq 0$, and $\xi$ denotes the real-valued white noise in space. The equation \eqref{NLSAnd} is also known as the dispersive Anderson model (DAM), which is the dispersive counterpart of the well-studied parabolic Anderson model (with $i \dt u$ replaced by $\dt u$ and $\ld = 0$; see, for example, \cite{HL15, HL18}).

The solution theory for nonlinear Schr\"odinger equations with a power type nonlinearity (i.e.~\eqref{NLSAnd} without the $\xi u$ term) has been studied extensively over the last few decades, with much efforts devoted to the Euclidean setting or the periodic setting. See, for example, \cite{Caz79, Kat87, Tsu, Caz03, BGTz}. Recently, there is a growing interest in the study of nonlinear Schr\"odinger equations with stochastic forcing terms. Among these varieties, a particularly interesting case is the multiplicative spatial white noise, which is useful for the study of solitary waves and Anderson localization.

Let us go over the well-posedness theory for the DAM \eqref{NLSAnd} and keep our discussion at a formal level in this subsection. The equation \eqref{NLSAnd} was first studied by Debussche and Weber in \cite{DW} on the two-dimensional torus $\T^2 = (\R / \Z)^2$. In particular, to deal with the ill-defined nature of $\xi u$, they used the Doss-Sussmann transform \cite{Doss, Suss} by writing $v = e^Y u$ with $Y = \Dl^{-1} \xi$, where the equation for $v$ contains more but smoother terms (see \eqref{vNLSAnd0} below). This transform had previously been used by Hairer and Labb\'e in \cite{HL15} for the parabolic Anderson model. However, unlike the parabolic equations, the Schr\"odinger equations do not have good smoothing properties. Because of this, the authors in \cite{DW} relied instead on the Hamiltonian structure of the DAM \eqref{NLSAnd}. Specifically, by mollifying the noise and establishing $H^2 (\T^2)$ a priori bounds for the solutions to the equations with the mollified noise, they proved the convergence of the approximating solutions and obtained global well-posedness for \eqref{NLSAnd} for a (sub-)cubic nonlinearity (i.e.~$1 \leq p \leq 3$) via a Gronwall argument. Global well-posedness for the DAM on $\T^2$ with an arbitrary power type nonlinearity (i.e.~$p > 3$) was then established by the second author and Visciglia in \cite{TzV23, TzV23-2} using the modified energies technique and also the dispersive smoothing effect.

Compared to the situation on periodic domains, well-poseness for the DAM on the full space presents an additional serious difficulty due to the unboundedness of the white noise $\xi$. On the other hand, the dispersive effect on the full space is stronger, a fact which has not been used so far in the analysis on the DAM. Global well-posedness for the DAM on $\R^2$ was first established by Debussche and Martin in \cite{DM} in the case of a sub-quadratic nonlinearity (i.e.~$1 < p < 2$). To deal with the logarithmic growth of the white noise, they incorporated the weighted Sobolev and Besov spaces in estimating the terms in the conserved quantities of the DAM \eqref{NLSAnd}. Later, Debussche, Visciglia, and the authors showed in \cite{DLTV} global well-posedness for the DAM on $\R^2$ with an arbitrary power type nonlinearity (i.e.~$p \geq 2$) by combining the aforementioned methods.

Given the fact that we have unique global-in-time solutions to the DAM \eqref{NLSAnd} in both the periodic setting and the full space setting, a next natural question is whether there is any relations between the periodic solution and the solution on the full space. In this paper, we aim to investigate this problem.
More precisely, our goal in this paper is to show that the global-in-time solution to \eqref{NLSAnd} on $\R^2$ can be realized as a limit of the following $L$-periodic DAM on $\T_L^2 = (\R / L \Z)^2$ as $L$ tends to infinity:
\begin{align}
\begin{cases}
i \dt u_L = \Dl u_L + \xi_L u_L - \ld |u_L|^{p - 1} u_L \\
u_L|_{t = 0} = u_{0, L}.
\end{cases}
\label{NLSAndL}
\end{align}

\noi
Here, we choose $\xi_L$ to be the $L$-periodized version of $\xi$ defined by
\begin{align}
\xi_L = \sum_{k \in \Z^2} \xi (\cdot + L k) \ind_{[- \frac{L}{2}, \frac{L}{2})^2} (\cdot + Lk),
\label{defxiL}
\end{align}

\noi
and we choose $u_{0, L}$ to be a suitable $L$-periodized version of $u_0$ to be specified in Subsection~\ref{SUB:main} below (see \eqref{defu0L}).

We refer to the above problem as the large torus limit problem. Recently, there has been several works in this direction of research in stochastic parabolic PDEs and random dispersive PDEs. In the former case, the authors in \cite{MW17} and \cite{GH19} studied the stochastic nonlinear heat equation with an additive space-time white noise on $\R^2$ and $\R^3$, respectively, where they both used the large torus limit of periodic solutions to construct a unique strong solution on the full space. See also the work \cite{DGR} on a similar large torus limit construction of the solution to the stochastic Burger's equation on $\R$. We point out that both \cite{MW17} and \cite{GH19} used a periodization of the space-time white noise similar to \eqref{defxiL}. In the latter case, much efforts have been devoted to constructing Gibbs measure dynamics of dispersive PDEs using large volume limit. In the one-dimensional setting, see \cite{Bour00, BS} for the NLS case and \cite{FKV} for the modified KdV case. In the two-dimensional setting, see \cite{OTWZ} for the study of stochastic damped nonlinear wave equation, where the authors used the finite speed of propagation of the wave equation in a crucial manner. We point out that in \cite{BS, FKV, OTWZ}, the source of randomness at each volume of size $L$ is produced via an application of the Skorokhod representation theorem and so these works require a change of probability space. Although such an application can give the strong almost sure convergence of periodic dynamics, the convergence as $L \to \infty$ in these works holds only along a subsequence of $L \in \N$.

In this paper, we study the large torus limit problem for a two-dimensional stochastic NLS. We choose to proceed with a periodization of the spatial white noise $\xi_L$ in \eqref{defxiL}, which has the same spirit as in \cite{MW17, GH19}. We find this periodization a natural procedure because it not only ensures that the only source of randomness in this problem comes from the white noise $\xi$ on $\R^2$, but also is convenient for showing the convergence of $\xi_L$ to $\xi$ and also their related stochastic terms as $L$ goes to infinity, even for $L \in [1, \infty)$ going along a continuum. It would be of interest to see if an analogous procedure can be applied to the setting of the Gibbs measure initial data.

Nevertheless, even with deterministic initial data with high regularity, establishing the convergence of $L$-periodic solution $u_L$ to \eqref{NLSAndL} to the full space solution $u$ to \eqref{NLSAnd} is highly nontrivial. There are two main challenges. The first challenge is due to the fact that NLS has infinite speed of propagation. Such an issue already appears for the deterministic NLS (i.e.~\eqref{NLSAnd} without the $\xi u$ term) and forces us to study the large torus convergence of solutions only on a bounded local domain. A useful tool for dealing with this problem is an exponential weight with slow decay introduced in \cite{Bour00} (see also \cite{BS}). 
The other challenge was already alluded above: the white noise $\xi$ exhibits logarithmic growth in space. As a result, the $L$-periodic spatial white noise $\xi_L$ in \eqref{defxiL} also has growth which is logarithmic in $L$. Consequently, the exponential stochastic object $e^{Y_L}$ has polynomial growth in $L$, which is further amplified by the Gronwall argument. Such loss in $L$ makes it extremely difficult to establish convergence of solutions as $L$ tends to infinity. To overcome this issue, we need to find a way to effectively control the loss in $L$.

In this paper, we introduce a new notion called {\it the $L$-periodic weights}. These periodic weights are carefully constructed and behave similarly to the polynomial weights $\jb{x}^\mu$ (with $\mu \in \R$) on $\R^2$. We then use these $L$-periodic weights to construct weighted function spaces on the $L$-periodic domain and show some useful properties that are analogous to those of weighted function spaces on $\R^2$, such as the Besov embedding, the duality estimate, and the interpolation inequality. Moreover, our $L$-periodic weights are smooth functions whose first order derivatives can be controlled by lower order weights. It is not yet clear to us if our $L$-periodic weights satisfy the admissible conditions as in Triebel's book \cite{Tri06}, since we do not know how to control higher order derivatives of our weights. Still, we are able to develop tools that are already sufficient for our purposes. More specifically, we are able to establish norm bounds for stochastic objects uniformly-in-$L$. These then allow us to show uniform-in-$L$ bounds for solutions to the $L$-periodic DAM, which are crucial for proving the large torus convergence of solutions to the DAM.

\subsection{Setup and the main result}
\label{SUB:main}

Let us now introduce our models and ansatz in a more rigorous manner.

Given a probability space $(\Omega, \mathcal{F}, \mathbb{P})$, a white noise $\xi$ on $\R^2$ is a linear map $f \mapsto (\xi, f)$ from $L^2 (\R^2)$ to $L^2 (\Om)$ such that for each $f \in L^2 (\R^2)$, $(\xi, f)$ is a centered Gaussian random variable such that $\E [ | (\xi, f) |^2 ] = \| f \|_{L^2 (\R^2)}^2$. We also assume that $\xi$ is real-valued in the sense that, if $f$ is real-valued, then $(\xi, f)$ is a real-valued Gaussian random variable. Given $L \geq 1$, we define $\xi_L$ as the $L$-periodized version of $\xi$ as in \eqref{defxiL}. This means that, given a test function $f \in \mathcal{D} (\R^2)$, we have
\begin{align}
(\xi_L, f) = \big( \xi, f_L \ind_{[ - \frac{L}{2}, \frac{L}{2} )^2} \big) ,
\label{xiL}
\end{align}

\noi
where
\begin{align*}
f_L (\cdot) = \sum_{k \in \Z^2} f (\cdot + Lk) .
\end{align*}

Let $G \in C^\infty (\R^2 \setminus \{0\})$ be a radial truncated Green's function such that $G (x) = - \frac{1}{2 \pi} \log |x|$ for $|x| \leq \frac 18$ and $G(x) = 0$ for $|x| \geq \frac 14$. 
By denoting $*_{\R^2}$ as the convolution on $\R^2$, we define
\begin{align}
Y = G *_{\R^2} \xi \quad \text{and} \quad Y_L = G *_{\R^2} \xi_L ,
\label{defY}
\end{align}

\noi
which solve the equations
\begin{align}
\Dl Y = \xi + \varphi *_{\R^2} \xi \quad \text{and} \quad \Dl Y_L = \xi_L + \varphi *_{\R^2} \xi_L,
\label{eqY}
\end{align}

\noi
respectively, for some $\varphi \in C_c^\infty (\R^2)$ supported on $\{ x \in \R^2: |x| < \frac 14 \}$. Note that we can write
\begin{align*}
Y = \int_{\R^2} G (\cdot - z) \xi (dz)
\end{align*}

\noi
and
\begin{align*}
Y_L = \int_{\R^2} G (\cdot - z) \xi_L (dz) = \int_{[ - \frac{L}{2}, \frac{L}{2} )^2} \sum_{k \in \Z^2} G (\cdot - z + Lk) \xi (dz),
\end{align*}

\noi
where we abused notations by writing $\xi (dz)$ (and $\xi_L (dz)$) as the Gaussian stochastic measure on $\R^2$ induced by the white noise $\xi$ (and the $L$-periodized noise $\xi_L$, respectively).

We now proceed as in \cite{HL15} by writing
\begin{align}
v = e^{Y} u \quad \text{and} \quad v_L = e^{Y_L} u_L,
\label{gauge}
\end{align}

\noi
which transform the DAM \eqref{NLSAnd} into
\begin{align}
\begin{cases}
i \dt v = \Dl v - 2 \nb Y \cdot \nb v + ( |\nb Y|^2 - \varphi *_{\R^2} \xi) v - \ld e^{- (p - 1) Y} |v|^{p - 1} v \\
v|_{t = 0} = e^{-Y} u_0
\end{cases}
\label{vNLSAnd0}
\end{align}

\noi
and transform the $L$-periodic DAM \eqref{NLSAndL} into
\begin{align}
\begin{cases}
i \dt v_L = \Dl v_L - 2 \nb Y_L \cdot \nb v_L + ( |\nb Y_L|^2 - \varphi *_{\R^2} \xi_L ) v_L - \ld e^{- (p - 1) Y_L} |v_L|^{p - 1} v_L \\
v_L|_{t = 0} = e^{-Y_L} u_{0, L} .
\end{cases}
\label{vNLSAndL0}
\end{align}

\noi
From \eqref{defY} and the fact that the spatial white noises $\xi$ and $\xi_L$ have regularity $-1 -$, we see that $\nb Y$ and $\nb Y_L$ have regularity $0-$, so that $\nb Y$ and $\nb Y_L$ are merely distributions and the square terms $|\nb Y|^2$ in \eqref{vNLSAnd0} and $|\nb Y_L|^2$ in \eqref{vNLSAndL0} are ill-defined. This requires us to renormalized the equations \eqref{vNLSAnd0} and \eqref{vNLSAndL0} by replacing the terms $|\nb Y|^2$ and $|\nb Y_L|^2$ with the Wick orderings $\wick{|\nb Y|^2}$ and $\wick{|\nb Y_L|^2}$ formally given by the following multiple Wiener-Ito integrals:
\begin{align*}
\wick{|\nb Y|^2} \, = \int_{\R^2} \int_{\R^2} \nb G (\cdot - z_1) \cdot \nb G (\cdot - z_2) \xi (d z_1) \xi (d z_2)
\end{align*}

\noi
and
\begin{align*}
\wick{|\nb Y_L|^2} \, = \int_{\R^2} \int_{\R^2} \nb G (\cdot - z_1) \cdot \nb G (\cdot - z_2) \xi_L (d z_1) \xi_L (d z_2) ,
\end{align*}

\noi
which are read in the following distributional sense:
\begin{align}
( \wick{|\nb Y|^2} \, , f ) = \int_{\R^2} \int_{\R^2} \bigg( \int_{\R^2} \nb G (x - z_1) \cdot \nb G (x - z_2) \cj{f (x)} dx \bigg) \xi (d z_1) \xi (d z_2)
\label{Y2_sto}
\end{align}

\noi
and
\begin{align}
( \wick{|\nb Y_L|^2} \, , f ) = \int_{\R^2} \int_{\R^2} \bigg( \int_{\R^2} \nb G (x - z_1) \cdot \nb G (x - z_2) \cj{f (x)} dx \bigg) \xi_L (d z_1) \xi_L (d z_2) 
\label{YL2_sto}
\end{align}

\noi
for any $f \in \mathcal{D} (\R^2)$. In fact, for $\wick{|\nb Y_L|^2}$\,, we will first define it using Fourier series and then show that it satisfies \eqref{YL2_sto} when we need this form; see \eqref{YL2_fs} and Lemma~\ref{LEM:YL2_sto} below. Such a renormalization procedure has already become a standard technique for solving stochastic PDEs containing singular terms; see also \cite{HL15, DW, DM, TzV23, TzV23-2, DLTV}.

From now on, we focus on the following renormalized DAM
\begin{align}
\begin{cases}
i \dt v = \Dl v - 2 \nb Y \cdot \nb v + \wt{\wick{|\nb Y|^2}} \, v - \ld e^{- (p - 1) Y} |v|^{p - 1} v \\
v|_{t = 0} = v_0 \deff e^{-Y} u_0
\end{cases}
\label{vNLSAnd}
\end{align}

\noi
and the renormalized $L$-periodic DAM
\begin{align}
\begin{cases}
i \dt v_L = \Dl v_L - 2 \nb Y_L \cdot \nb v_L + \wt{\wick{|\nb Y_L|^2}} \, v_L - \ld e^{- (p - 1) Y_L} |v_L|^{p - 1} v_L \\
v_L|_{t = 0} = v_{0, L} \deff e^{-Y_L} u_{0, L} ,
\end{cases}
\label{vNLSAndL}
\end{align}

\noi
where we used the notations
\begin{align}
\wt{\wick{|\nb Y|^2}} \deff \, \wick{|\nb Y|^2} - \, \varphi *_{\R^2} \xi \quad \text{and} \quad \wt{\wick{|\nb Y_L|^2}} \deff \, \wick{|\nb Y_L|^2} - \, \varphi *_{\R^2} \xi_L .
\label{wY2t}
\end{align}

\noi
For later convenience, we simply refer to the equation \eqref{vNLSAnd} as DAM and the equation \eqref{vNLSAndL} as the $L$-periodic DAM.

We also need to specify the initial conditions for $v_0$ and $v_{0, L}$ for the equations \eqref{vNLSAnd} and \eqref{vNLSAndL}. For this purpose, we need the weighted Sobolev spaces. The weighted Sobolev norms on $\R^2$ are now classical, but for the moment let us take the following equivalence as the definition for simplicity:
\begin{align*}
\| f \|_{H^2_\mu (\R^2)} \sim \| \jb{\cdot}^\mu f \|_{L^2 (\R^2)} + \| \jb{\cdot}^\mu \Dl f \|_{L^2 (\R^2)},
\end{align*}

\noi
where $\mu \in \R$; see Subsection~\ref{SUB:poly} below for more details. The weighted Sobolev norms on the $L$-periodic domain $\T_L^2$ are constructed later in this paper as one of our main novelties. Again, for now let us take the following equivalence as the definition for simplicity:
\begin{align*}
\| f_L \|_{H^2_\mu (\T_L^2)} \sim \| \jb{\cdot}_L^\mu f \|_{L^2 (\T_L^2)} + \| \jb{\cdot}_L^\mu \Dl f \|_{L^2 (\T_L^2)},
\end{align*}

\noi
where $\mu \in \R$ and $\jb{\cdot}_L$ is an $L$-periodic function such that $\jb{x}_L = \jb{x}$ for $x \in [- \frac{L}{2}, \frac{L}{2})^2$.
Let us assume that $v_0 \in H^2_{\mu_0} (\R^2)$ with $\mu_0 > 1$ and define
\begin{align}
v_{0, L} \deff \sum_{k \in \Z^2} v_0 (\cdot + L k).
\label{v0L}
\end{align}

\noi
Thanks to Lemma~\ref{LEM:per} below, we see that $v_{0, L}$ makes sense as a function in $H_{\mu}^2 (\T_L^2)$ for any $0 < \mu < \mu_0 - 1$. This provides the choice of the initial data $u_{0, L}$ for the original (formal) equation \eqref{NLSAndL}:
\begin{align}
u_{0, L} = e^{Y_L} v_{0, L} = e^{Y_L} \sum_{k \in \Z^2} e^{- Y (\cdot + Lk)} u_0 (\cdot + Lk).
\label{defu0L}
\end{align}

Let us first state the following result on global well-posedness of the $L$-periodic DAM \eqref{vNLSAndL} and also a good bound for the solution. 

\begin{theorem}
\label{THM:GWP}
Let $p \geq 2$, $\ld \geq 0$, $L \geq 1$, $\mu_0 > 0$, and $v_{0, L} \in H_{\mu_0}^2 (\T_L^2)$ be a periodic function. Then, for any $1 < s < 2$, there exists a unique global-in-time solution $v_L$ to the $L$-periodic DAM \eqref{vNLSAndL} with $v_L |_{t = 0} = v_{0, L}$ in the space $C (\R_+; H^s (\T_L^2))$ almost surely. Moreover, for any $T \geq 1$, we have the bound
\begin{align*}
\| v_L \|_{C ([0, T]; H^s (\T_L^2))} \leq C \big( \om, T, \| v_{0, L} \|_{H_{\mu_0}^2 (\T_L^2)} \big) ,
\end{align*}

\noi
where $C (\om, T, \| v_{0, L} \|_{H_{\mu_0}^2 (\T_L^2)}) > 0$ is almost surely finite.
\end{theorem}

In the above statement, we state the result for general periodic initial data $v_{0, L}$ with some regularity assumptions. If $v_{0, L}$ is given by \eqref{v0L} with $v_0 \in H_{\mu_0}^2 (\R^2)$ for $\mu_0 > 1$, then $v_{0, L}$ satisfies the assumption in Theorem~\ref{THM:GWP} and the bound for the solution $v_L$ in Theorem~\ref{THM:GWP} is uniform in $L$ thanks to Lemma~\ref{LEM:per} below.

To prove Theorem~\ref{THM:GWP}, we show the convergence of a sequence of approximating solutions $\{ v_{L, \eps} \}_{\eps > 0}$ to equations with the mollified noise. However, in order to obtain a uniform-in-$L$ bound for the solution, we cannot proceed as in \cite{DW, TzV23, TzV23-2} on the usual periodic domain $\T^2$ due to the growth of the noise in $L$ as mentioned above. Instead, we introduce $L$-periodic weights and construct weighted Besov and Sobolev spaces on the $L$-periodic domain $\T_L^2$; see Section~\ref{SEC:BesL} below. With these tools in hand, the estimates for the stochastic terms $\xi_L$, $Y_L$, $\wick{|\nb Y_L|^2}$\,, and $e^{- Y_L}$ can be made independent of $L$; see Subsection~\ref{SUB:sto} below. We then adapt the weighted function spaces to the framework of \cite{DLTV} (with slightly simplified steps) on the full space $\R^2$, where all the deterministic estimates such as the Besov embedding, the Strichartz estimates, and the product estimates on $\T_L^2$ are also uniform in $L$.

\medskip
We are now ready to state the main result in this paper.

\begin{theorem}
\label{THM:convL}
Let $p \geq 2$, $\ld \geq 0$, $\mu_0 > 1$, $v_0 \in H_{\mu_0}^2 (\R^2)$, and $v_{0, L} \in H_{\mu_1}^2 (\T_L^2)$ with $0 < \mu_1 < \mu_0 - 1$ be as defined in \eqref{v0L} given any $L \geq 1$.  Let $v$ be the global-in-time solution to DAM \eqref{vNLSAnd} with initial data $v_0$ guaranteed by \cite[Theorem~1.2]{DLTV} and let $v_L$ be the global-in-time solution to the $L$-periodic DAM \eqref{vNLSAndL} with initial data $v_{0, L}$ guaranteed by Theorem~\ref{THM:GWP}. Then, for any $0 \leq s < 2$ and bounded open set $U \subset \R^2$, $v_L$ converges to $v$ in $C(\R_+; H^s (U))$ in probability as $L \to \infty$.
\end{theorem}

In the above statement, the space $H^s (U)$ given an open set $U \subset \R^2$ denotes the restriction of the $H^s (\R^2)$ space on $U$ (see \eqref{defXK} below), and the space $C(\R_+; H^s (U))$ is endowed with the compact-open topology in time. In particular, we have from the statement of Theorem~\ref{THM:convL} and also Sobolev's embedding that $v_L$ converges in probability to $v$ uniformly on compact subsets of $\R_+ \times \R^2$ as $L$ goes to infinity. We also mention that the condition on the weight parameter $\mu_0$ in Theorem~\ref{THM:convL} is sharp up to the endpoint in the sense that if $\mu_0 < 1$, then there exists $v_0 \in H_{\mu_0}^2 (\R^2)$ such that the periodic function $v_{0, L}$ may not be constructed via the form \eqref{v0L}; see Remark~\ref{RMK:mu0} below.

Given that we have a uniform-in-$L$ bound for the solution $v_L$ to the $L$-periodic DAM \eqref{vNLSAndL}, in order to prove Theorem~\ref{THM:convL}, we need to establish the convergence of the stochastic terms $\xi_L$, $Y_L$, $\wick{|\nb Y_L|^2}$\,, and $e^{- Y_L}$ to their full space limits $\xi$, $Y$, $\wick{|\nb Y|^2}$\,, and $e^{- Y}$. Nevertheless, at the same time, we need to address the issue that the NLS has infinite speed of propagation. Here, we follow \cite{Bour00, BS} by incorporating an exponential weight with slow decay to help us establish the large torus convergence; see the function $\s_R$ \eqref{sigR} below, but we often use it as $\s_R^{-1}$. In Appendix~\ref{SEC:NLS}, we separately present such an argument for the large torus convergence problem for deterministic NLS, since we are not aware of such a result in existing literature and this might be of general interest.

The incorporation of the weight $\s_R^{-1}$ in the estimate of the aforementioned stochastic terms is not immediate. The reason is that $\nb Y_L$, $\wick{|\nb Y_L|^2}$\,, $\nb Y$, and $\wick{|\nb Y|^2}$ all live in weighted Besov spaces with negative regularity. It is possible to solve this issue by viewing $\s_R^{-1}$ essentially as the classical $\jb{\cdot}^{-1}$ with some loss due to the slowly decaying nature of $\s_R^{-1}$, but one has to carefully perform the analysis with fractional order derivatives on these weights. This procedure should work eventually, but in this paper, we proceed in another way which we believe is easier. Instead of directly showing the convergence of the solutions $v_L$ to $v$, we choose to first show the convergence of the approximating solutions $v_{L, \eps}$ to $v_\eps$, which correspond to solutions to equations with mollified noise. The advantage of this approach is that we do not need to worry about regularities in estimating the difference of solutions, since the smoothed stochastic terms $\nb Y_{L, \eps}$, $\wick{|\nb Y_{L, \eps}|^2}$\,, $\nb Y_\eps$, and $\wick{|\nb Y_\eps|^2}$ all live in Lebesgue spaces with logarithmic loss in $\eps$. See Subsection~\ref{SUB:sto_conv} for the large torus convergence of the smoothed stochastic terms.


After obtaining good convergence of stochastic terms, we then proceed with a Gronwall argument. Namely, our goal is to obtain an estimate that looks like
\begin{align*}
\frac{d}{dt} \int_{\R^2} \s_R^{-1} e^{- 2 Y_{L, \eps}} |v_\eps (t) - v_{L, \eps} (t)|^2 dx \les C_1 L^{- \dl} + C_2 \int_{\R^2} \s_R^{-1} e^{- 2 Y_{L, \eps}} |v_\eps (t) - v_{L, \eps} (t)|^2 dx ,
\end{align*}

\noi
where $\dl > 0$ and $C_1, C_2 > 0$ are constants independent of $L$. The factor $L^{- \dl}$ implies the desired convergence and it comes from both the large torus convergence of stochastic terms and also the slowly decaying nature of the weight $\s_R^{-1}$.
For neatness and readability of this paper, we choose to split the proof of Theorem~\ref{THM:GWP} and Theorem~\ref{THM:convL} into the linear case ($\lambda = 0$) (see Section~\ref{SEC:convL1}) and the nonlinear case ($\lambda > 0$) (see Section~\ref{SEC:convL2}).

Before moving on, we would like to mention that, to the best of our knowledge, the large torus limit problem for the parabolic Anderson model is not yet explored in existing literature. In Appendix~\ref{SEC:heat}, we state and prove a large torus limit result for the linear parabolic Anderson model in the two-dimensional setting, combining the tools from \cite{HL15} and those in our paper. We remark that it is also possible to consider the same problem for a (power-type) nonlinear parabolic equation with a multiplicative spatial white noise, but one may need some additional techniques to deal with the nonlinearity, such as the energy estimates in \cite{MW17}.

\medskip
We close the discussion by stating several remarks.

\begin{remark} \rm
In this paper, we mainly consider the defocusing case $\ld > 0$ of DAM \eqref{vNLSAnd}. Nevertheless, the same analysis for proving the large torus limit in Theorem~\ref{THM:convL} holds for the focusing case $\ld < 0$, as long as we impose an additional smallness assumption on the initial data $\| v_0 \|_{H_{\mu_0}^1 (\R^2)}$. This smallness assumption is needed for obtaining a weighted $H^1$ a priori bound for the soltion $v$ of the DAM \eqref{vNLSAnd} and also, along with Lemma~\ref{LEM:per} below, a weighted $H^1$ a priori bound for the soltion $v_L$ of the periodic DAM \eqref{vNLSAndL}; see \cite[Remark~5.2]{DLTV}. Once we have these a priori bound for solutions of the full space DAM and the periodic DAM, the proof of the large torus limit in Theorem~\ref{THM:convL} is then irrelevant of the sign of $\ld$.
\end{remark}

\begin{remark} \rm
In Theorem~\ref{THM:convL}, we prove the large torus convergence of solutions to DAM in probability as $L \to \infty$. At this point, it is not clear if one can prove such convergence in an almost sure manner. Even for the parabolic Anderson model, we need to add the restriction that $L$ goes along positive integers in order to prove the large torus convergence almost surely; see Remark~\ref{RMK:as2} in Appendix~\ref{SEC:heat}. Unfortunately, such result is not directly applicable to DAM due to the use of different function spaces. It would be of interest to investigate almost sure large torus convergence for DAM, but we choose not to pursue this issue in this work.
\end{remark}

\begin{remark} \rm
In this paper, we consider $p \geq 2$ for the large torus limit problem of DAM \eqref{NLSAnd}. In particular, we cover the cubic nonlinearity which corresponds to $p = 3$. The condition $p \geq 2$ is required for the estimates of modified energies in Subsection~\ref{SUB:mod} for getting the $H^2$-a priori bounds for the solutions. As for the lower order nonlinearity with $1 < p < 2$, such $H^2$-a priori bounds can be obtained via a slightly easier approach as in \cite{DM}. We choose not to add such an argument to further complicate the presentation in this paper.
\end{remark}

\begin{remark} \rm
Besides the gauge transform $v = e^Y u$ for solving the DAM \eqref{NLSAnd}, there is another approach based on realizing the Anderson Hamiltonian $H = \Dl + \xi$ as a self-adjoint operator on $L^2$; see \cite{GUZ, Ugur, MZ, DVJZ}. It would be of interest to explore the large torus limit problem for DAM from the viewpoint of the Anderson Hamiltonian.
\end{remark}

\subsection{Notations}
For two quantities $a, b > 0$, we write $a \les b$ if $a \leq C b$ for some constant $C > 0$ independent with the set where $a$ and $b$ are allowed to vary. We also use subscripts such as ``$\les_a$'' to emphasize dependence on external parameters. We also write $a \ll b$ if $a < cb$ for some small constant $c > 0$ independent with the set where $a$ and $b$ are allowed to vary. Unless we make explicit explanation, the underlying constants in this paper depend on all parameters except for the dyadic frequency scale $N$, the size $L$ of the torus, the spatial localization parameter $R$, and the rate $\eps$ of the regularization. In fact, at the end of proving the convergence, the main large parameter is $L$, and the small $\eps$ and the large $R$ are chosen as functions of $L$.

Throughout the paper, we use the notation $\jb{\cdot} = (1 + |\cdot|^2)^{\frac 12}$. 
We distinguish the convolutions on different domains by writing $*_{\R^2}$ as the convolution on $\R^2$ and $*_{\T_L^2}$ as the convolution on the torus $\T_L^2$ of size $L \geq 1$. We also write $\mathcal{F}_{\R^2}$ (or $\mathcal{F}^{-1}_{\R^2}$) as the Fourier transform (or the inverse Fourier transform, respectively) on $\R^2$ and $\mathcal{F}_{\T_L^2}$ (or $\mathcal{F}^{-1}_{\T_L^2}$) as the Fourier transform (or the inverse Fourier transform, respectively) on the torus $\T_L^2$ of size $L \geq 1$. Given $T > 0$, $1 \leq p \leq \infty$, and a Banach space $X$, we use the shorthand notations $L_T^p X$ for the space $L^p ([0, T]; X)$ and $C_T X$ for the space $C ([0, T]; X)$.

\subsection{Organization of the paper}

This paper is organized as follows. In Section~\ref{SEC:lem}, we introduce some basic tools and present some preliminary lemmas. In Section~\ref{SEC:BesL}, we introduce $L$-periodic weights and weighted function spaces on periodic domains and prove some useful estimates. In Section~\ref{SEC:sto}, we establish regularity and convergence properties of the stochastic objects mentioned in Subsection~\ref{SUB:main}. In Section~\ref{SEC:convL1}, we prove  Theorem~\ref{THM:GWP} and Theorem~\ref{THM:convL} in the linear case $\ld = 0$. In Section~\ref{SEC:convL2}, we prove Theorem~\ref{THM:GWP} and Theorem~\ref{THM:convL} in the nonlinear case $\ld > 0$. In Appendix~\ref{SEC:NLS}, we show the large torus limit for the deterministic NLS. Lastly, in Appendix~\ref{SEC:heat}, we show the large torus limit for the parabolic Anderson model.

\section{Basic tools and preliminary lemmas}
\label{SEC:lem}

In this section, we introduce various function spaces and present some preliminary lemmas.

\subsection{Periodic distributions}

In this subsection, we clarify the definition and properties of periodic distributions. 
To avoid confusion, we use $(\cdot, \cdot)_{\R^2}$ to denote the duality pairing on $\R^2$ and use $(\cdot, \cdot)_{\T_L^2}$ to denote the duality pairing on $\T_L^2$.

Let $\mathcal{D} (\R^2)$ be the space of test functions on $\R^2$, i.e.~the space of infinitely differentiable functions of compact support, and let $\mathcal{D}' (\R^2)$ be the space of distributions, i.e.~the space of continuous linear functionals on $\mathcal{D} (\R^2)$, whose topology of locally convex topological vector space is given by a suitable family of semi-norms which is essentially characterized by convergence of all partial derivatives on compact sets. Given $L \geq 1$, a distribution $S \in \mathcal{D}' (\R^2)$ is said to be $L$-periodic if $(S, f)_{\R^2} = (S, f_{L, k})_{\R^2}$ with $f_{L, k} (\cdot) = f (\cdot - Lk)$ for all $f \in \mathcal{D} (\R^2)$ and $k \in \Z^2$. Given $S$ an $L$-periodic distribution and $f_L \in C^\infty (\T_L^2)$ viewed as an $L$-periodic function, we set
\begin{align}
(S, f_L)_{\T_L^2} \deff (S, \phi_L f_L)_{\R^2},
\label{Sper}
\end{align}

\noi
where $\phi_L \in \mathcal{D} (\R^2)$ satisfies
\begin{align}
\sum_{k \in \Z^2} \phi_L (x + L k) = 1
\label{phiL_cond}
\end{align}

\noi
for any $x \in \R^2$. 
To show that \eqref{Sper} is well-defined, i.e. the definition is independent of the choice of $\phi_L$, we need the following lemma. 

\begin{lemma}
\label{LEM:Sfull}
Let $S \in \mathcal{D}' (\R^2)$ be an $L$-periodic distribution and $\phi_L \in \mathcal{D} (\R^2)$ satisfying \eqref{phiL_cond}. Then, for any $f \in \mathcal{D} (\R^2)$, we have
\begin{align*}
(S, f)_{\R^2} = \Big( S, \phi_L \sum_{k \in \Z^2} f (\cdot + L k) \Big)_{\R^2}.
\end{align*}
\end{lemma}

\begin{proof}
From \eqref{phiL_cond} and the fact that $S$ is an $L$-periodic distribution, we get
\begin{align*}
(S, f)_{\R^2} &= \Big( S, \sum_{k \in \Z^2} \phi_L (\cdot - L k) f  \Big)_{\R^2} = \sum_{k \in \Z^2} (S, \phi_L (\cdot - L k) f)_{\R^2} \\
&= \sum_{k \in \Z^2} (S, \phi_L f (\cdot + L k))_{\R^2} = \Big( S, \phi_L \sum_{k \in \Z^2} f (\cdot + L k) \Big)_{\R^2} ,
\end{align*}

\noi
where the summation is a finite sum due to the compact support of $\phi_L$ and $f$. 
\end{proof}

Thanks to Lemma~\ref{LEM:Sfull}, given $f_L \in C^\infty (\T_L^2)$ viewed as an $L$-periodic function and $\phi_{1, L}, \phi_{2, L} \in \mathcal{D} (\R^2)$ both satisfying \eqref{phiL_cond}, we have
\begin{align*}
(S, \phi_{1, L} f_L)_{\R^2} &= \Big( S, \phi_{2, L} \sum_{k \in \Z^2} \phi_{1, L} (\cdot + L k) f_L (\cdot + L k) \Big)_{\R^2} \\
&= \Big( S, \phi_{2, L} f_L \sum_{k \in \Z^2} \phi_{1, L} (\cdot + L k) \Big)_{\R^2} \\
&= (S, \phi_{2, L} f_L)_{\R^2} .
\end{align*}

\noi
This shows that \eqref{Sper} is well-defined.

\medskip
Given an $L$-periodic distribution $S$ and a function $g_L \in L^1 (\T_L^2)$, we define the convolution $g_L *_{\T_L^2} S$ as the $L$-periodic distribution via
\begin{align}
(g_L *_{\T_L^2} S, f_L)_{\T_L^2} \deff (S, \wt{g_L} *_{\T_L^2} f_L)_{\T_L^2} ,
\label{def_convL}
\end{align}

\noi
where $\wt{g_L} (\cdot) = g_L (- \cdot)$. If $g_L \in C^\infty (\T_L^2)$, a standard argument (see \cite[Theorem~2.3.20]{Gra}) shows that $g_L *_{\T_L^2} S \in C^\infty (\T_L^2)$ and
\begin{align*}
(g_L *_{\T_L^2} S) (x) = (S, \wt{g_L} (\cdot - x)),
\end{align*}

\noi
which we will use frequently throughout this paper.

\begin{lemma}
\label{LEM:GL}
Let $L \geq 1$, $f_L$ be an $L$-periodic distribution, and $g \in L^1 (\R^2)$ supported on $\{ x \in \R^2 : |x| < \frac{L}{4} \}$. Then, $g *_{\R^2} f_L$ is an $L$-periodic distributions and
\begin{align*}
g *_{\R^2} f_L = g_L *_{\T_L^2} f_L,
\end{align*}

\noi
where
\begin{align*}
g_L (\cdot) = \sum_{k \in \Z^2} g (\cdot + Lk).
\end{align*}
\end{lemma}

\begin{proof}
We write $\wt g (\cdot) = g (- \cdot)$ and $\wt{g_L} (\cdot) = g_L (- \cdot)$. Given any $\varphi \in \mathcal{D} (\R^2)$, $k \in \Z^2$, and $\varphi_{L, k}$ defined by $\varphi_{L, k} (\cdot) = \varphi (\cdot - Lk)$, we have $(\wt g *_{\R^2} \varphi) (x) = (\wt g *_{\R^2} \varphi_{L, k}) (x + Lk)$ for any $x \in \R^2$ and $k \in \Z^2$. Thus, from the fact that $f_L$ is an $L$-periodic distribution, we have
\begin{align*}
(g *_{\R^2} f_L, \varphi)_{\R^2} = (f_L, \wt g *_{\R^2} \varphi)_{\R^2} = (f_L, \wt g *_{\R^2} \varphi_{L, k})_{\R^2} = (g *_{\R^2} f_L, \varphi_{L, k})_{\R^2} ,
\end{align*}

\noi
which shows that $g *_{\R^2} f_L$ is an $L$-periodic distribution.

For any $\varphi_L \in C^\infty (\T_L^2)$, from the definition in \eqref{Sper} and Lemma~\ref{LEM:Sfull} with $\phi_L \in \mathcal{D}(\R^2)$ satisfying \eqref{phiL_cond}, we have
\begin{align}
\begin{split}
(g *_{\R^2} f_L, \varphi_L)_{\T_L^2} &= (g *_{\R^2} f_L, \phi_L \varphi_L)_{\R^2} \\
&= (f_L, \wt g *_{\R^2} \phi_L \varphi_L)_{\R^2} \\
&= \bigg( f_L, \phi_L (\cdot) \sum_{k \in \Z^2} \int_{\R^2} \wt g (\cdot + Lk - y) \phi_L (y) \varphi_L (y) dy \bigg)_{\R^2} .
\end{split}
\label{gper1}
\end{align}

\noi
By using a change of variable and \eqref{phiL_cond}, we have that for any $x \in \R^2$,
\begin{align}
\begin{split}
\sum_{k \in \Z^2} \int_{\R^2} \wt g (x + Lk - y) \phi_L (y) \varphi_L (y) dy 
&= \sum_{k \in \Z^2} \int_{\R^2} \wt g (x - y) \phi_L (y + Lk) \varphi_L (y) dy \\
&= (\wt g *_{\R^2} \varphi_L) (x) .
\end{split}
\label{gper2}
\end{align}

\noi
Also, for any $\varphi_L \in C^\infty (\T_L^2)$, we easily check that
\begin{align}
\wt g *_{\R^2} \varphi_L = \wt{g_L} *_{\T_L^2} \varphi_L .
\label{GgL}
\end{align}

\noi
Thus, from \eqref{gper1}, \eqref{gper2},  \eqref{GgL}, the definition in \eqref{Sper}, and the definition in \eqref{def_convL}, we obtain
\begin{align*}
(g *_{\R^2} f_L, \varphi_L)_{\T_L^2} 
= \big( f_L, \phi_L (\wt{g_L} *_{\T_L^2} \varphi_L) \big)_{\R^2}
= (f_L, \wt{g_L} *_{\T_L^2} \varphi_L)_{\T_L^2}
= (g_L *_{\T_L^2} f_L, \varphi_L)_{\T_L^2},
\end{align*}

\noi
which gives the desired identity.
\end{proof}

Let us now consider the $L$-periodic white noise $\xi_L$ defined in \eqref{xiL}. 
One can easily see from the definition that $(\xi_L, f)_{\R^2} = (\xi_L, f_{L, k})_{\R^2}$ with $f_{L, k} (\cdot) = f (\cdot - Lk)$ for all $f \in \mathcal{D} (\R^2)$ and $k \in \Z^2$, so that $\xi_L$ can be viewed almost surely as an $L$-periodic distribution.
For any $f_L \in C^\infty (\T_L^2)$ viewed as an $L$-periodic function, we define
\begin{align}
(\xi_L, f_L)_{\T_L^2} \deff ( \xi_L, \phi_L f_L )_{\R^2} , 
\label{defxiL2}
\end{align}

\noi
where $\phi_L \in \mathcal{D} (\R^2)$ satisfies \eqref{phiL_cond}. Note that by the definition \eqref{xiL} (on $\R^2$) and the fact that $f_L$ is $L$-periodic, we have
\begin{align*}
\begin{split}
( \xi_L, f_L )_{\T_L^2} &= ( \xi_L, \phi_L f_L )_{\R^2} \\
&= \Big( \xi, \sum_{k \in \Z^2} \phi_L (\cdot + L k) f_L (\cdot + L k) \ind_{[- \frac{L}{2}, \frac{L}{2})^2} \Big)_{\R^2} \\
&= \Big( \xi, \sum_{k \in \Z^2} \phi_L (\cdot + L k) f_L \ind_{[- \frac{L}{2}, \frac{L}{2})^2} \Big)_{\R^2} ,
\end{split}
\end{align*}

\noi
so that from \eqref{phiL_cond}, we get
\begin{align}
( \xi_L, f_L )_{\T_L^2} = \big( \xi, f_L \ind_{[- \frac{L}{2}, \frac{L}{2})^2} \big)_{\R^2} .
\label{xiL_fL}
\end{align}

\noi
The identity \eqref{xiL_fL} aligns with the definition \eqref{xiL} and shows that $(\xi_L, f_L)_{\T_L^2}$ is a centered Gaussian random variable such that 
\begin{align*}
\E \big[ | (\xi_L, f_L)_{\T_L^2} |^2 \big] = \| f_L \|_{L^2 (\T_L^2)}^2,
\end{align*}

\noi
so that the definition \eqref{defxiL2} can be extended to $L$-periodic functions in $L^2 (\T_L^2)$. Moreover, if $f_L$ is real-valued, then $(\xi_L, f_L)_{\T_L^2}$ is a real-valued Gaussian random variable.

\subsection{Weighted function spaces on the plane}
\label{SUB:poly}


Given $\mu \in \R$ and $1 \leq p \leq \infty$, we define the $L^p$ space with a polynomial weight $\jb{x}^\mu$ as
\begin{align*}
\| f \|_{L_\mu^p (\R^2)} \deff \bigg( \int_{\R^2} |f (x)|^p \jb{x}^{\mu p} dx \bigg)^{\frac{1}{p}} = \| \jb{\cdot}^\mu f \|_{L^p (\R^2)} ,
\end{align*}

\noi
with the obvious interpretation if $p = \infty$. When $\mu = 0$, we write $L^p = L^p_0$.

We now define weighted Besov spaces on $\R^2$. Let $\phi \in C_c^\infty (\R^2)$ be supported on $\{ \xi \in \R^2: |\xi| \leq \frac 85 \}$ and $\phi \equiv 1$ on $\{ \xi \in \R^2: |\xi| \leq \frac 54 \}$. Given a dyadic number $N \geq 1$, we define 
\begin{align}
\psi_N (\xi) = 
\begin{cases}
\phi (\xi) & \text{if } N = 1 \\
\phi (\tfrac{\xi}{N}) - \phi (\tfrac{2\xi}{N}) & \text{if } N \geq 2.
\end{cases}
\label{defpsi}
\end{align}

\noi
We define the Littlewood-Paley projector $\Dl_N$ as the Fourier multiplier operator with symbol $\psi_N$, so that
\begin{align}
\Id = \sum_{\substack{N \geq 1 \\ \text{dyadic}}} \Dl_N .
\label{DN}
\end{align}

\noi
We also define $\Dl_{\frac 12} = 0$. Given a dyadic $N \geq 1$, we will frequently use the identity
\begin{align*}
\Dl_N = \Dl_N \big( \Dl_{\frac{N}{2}} + \Dl_N + \Dl_{2N} \big) .
\end{align*}

Given $s, \mu \in \R$ and $1 \leq p, q \leq \infty$, we define the weighted inhomogeneous Besov space $\B_{p, q, \mu}^s (\R^2)$ as
\begin{align*}
\| f \|_{\B_{p, q, \mu}^s (\R^2)} \deff \bigg( \sum_{\substack{N \geq 1 \\ \text{dyadic}}} N^{s q} \| \Dl_N f \|_{L_\mu^p (\R^2)}^q  \bigg)^{\frac{1}{q}} = \bigg( \sum_{\substack{N \geq 1 \\ \text{dyadic}}} N^{s q} \| \jb{\cdot}^\mu \Dl_N f \|_{L^p (\R^2)}^q  \bigg)^{\frac{1}{q}} .
\end{align*}

\noi
When $\mu = 0$, we omit the weight parameter and write $\B^s_{p, q} (\R^2) = \B^s_{p, q, 0} (\R^2)$, which is the usual Besov space on $\R^2$.
We have the following convenient property of the weighted Besov spaces (see \cite[Theorem~6.5]{Tri06}):
\begin{align}
\| f \|_{\B_{p, q, \mu}^s (\R^2)} \sim \| \jb{\cdot}^\mu f \|_{\B_{p, q}^s (\R^2)} ,
\label{Bpq_pull}
\end{align}

\noi
which can be used to translate results from the unweighted spaces to the weighted analogues.
When $p = q = 2$, we denote
\begin{align*}
H^s_\mu (\R^2) \deff \B_{2, 2, \mu}^s (\R^2).
\end{align*}

\noi
From \cite[Theorem~6.9]{Tri06}, for $n \in \N$ and $\mu \in \R$, we have the following equivalence:
\begin{align}
\| f \|_{H^n_{\mu}} \sim \| f \|_{L^2_\mu} + \sum_{\substack{k \in \Z_{\geq 0}^2 \\ 0 < |k|_1 \leq n}} \| \partial^k f \|_{L^2_\mu}
\label{Hn_equi}
\end{align}

\noi
with $\Z_{\geq 0}$ denoting the set of nonnegative integers and $|k|_1 = |(k_1, k_2)|_1 = |k_1| + |k_2|$. When $p = q = \infty$, we denote
\begin{align*}
\C^s_\mu (\R^2) \deff \B^s_{\infty, \infty, \mu} (\R^2) .
\end{align*}

\noi
If $s > 0$ is not an integer, the space $\C^s_\mu (\R^2)$ coincides with the classical H\"older space with the equivalent norm (see \cite[Section~6.2]{Tri06}):
\begin{align}
\| f \|_{\C^s_\mu (\R^2)} \sim \sum_{\substack{k \in \Z_{\geq 0}^2 \\ |k|_1 \leq \lfloor s \rfloor}} \sup_{x \in \R^2} \jb{x}^\mu |\partial^k f (x)| + \sum_{\substack{k \in \Z_{\geq 0}^2 \\ |k|_1 = \lfloor s \rfloor}} \sup_{0 < |x - y| \leq 1} \jb{x}^\mu \frac{| \partial^k f (x) - \partial^k f (y) |}{|x - y|^{s - \lfloor s \rfloor}} .
\label{Ca-norm}
\end{align}


Let $X$ be any norm defined above. For any open set $U \subset \R^2$, we define the localized norm $X (U)$ as
\begin{align}
\| f \|_{X (U)} \deff \inf \big\{ \| g \|_X : f = g \text{ on } U \big\} .
\label{defXK}
\end{align} 

\noi
When $X = H^s$ with $s \in \N$ and $U$ is a Lipschitz domain, we have the following equivalence from \cite[Theorem~3.30~(ii)]{McL}:
\begin{align}
\| f \|_{H^s (U)} \sim \| f \|_{L^2 (U)} + \sum_{\substack{k \in \Z_{\geq 0}^2 \\ 0 < |k|_1 \leq s}} \| \partial^k f \|_{L^2 (U)} .
\label{loc_equi}
\end{align}

The embeddings in the following lemma follow either from the definition directly or from the unweighted estimates in \cite{BCD} thanks to the equivalence \eqref{Bpq_pull}.

\begin{lemma}
\label{LEM:emb}
\textup{(i)} Let $1 \leq p \leq \infty$, $1 \leq q_1, q_2 \leq \infty$, $s_1 < s_2$, and $\mu \in \R$. Then, we have
\begin{align*}
\| f \|_{\B_{p, q_1, \mu}^{s_1} (\R^2)} \les \| f \|_{\B_{p, \infty, \mu}^{s_2} (\R^2)} \les \| f \|_{\B_{p, q_2, \mu}^{s_2} (\R^2)} .
\end{align*}

\smallskip \noi
\textup{(ii)} Let $1 \leq p \leq \infty$ and $\mu \in \R$. Then, we have
\begin{align*}
\| f \|_{\B_{p, \infty, \mu}^0 (\R^2)} \les \| f \|_{L^p_\mu (\R^2)} \les \| f \|_{\B_{p, 1, \mu}^0 (\R^2)} .
\end{align*}
\end{lemma}

We also record the following product estimates; see \cite{BCD, PT16}, which are for unweighted spaces but can be easily adapted to weighted spaces.
\begin{lemma}
\label{LEM:prod}
Let $s, s_1, s_2 \in \R$ be such that $s_1 + s_2 > 0$ and $s = \min (s_1 + s_2, s_1, s_2)$, $s' > 0$, $1 \leq p, p_1, p_2 \leq \infty$ be such that $\frac{1}{p} = \frac{1}{p_1} + \frac{1}{p_2}$, $1 \leq q \leq \infty$, and $\mu, \mu_1, \mu_2 \in \R$ be such that $\mu = \mu_1 + \mu_2$. Then, for any $\kappa > 0$, we have
\begin{align*}
\| fg \|_{\B_{p, p, \mu}^{s - \kappa} (\R^2)} &\les \| f \|_{\B_{p_1, p_1, \mu_1}^{s_1} (\R^2)} \| g \|_{\B_{p_2, p_2, \mu_2}^{s_2} (\R^2)} , \\
\| fg \|_{\C_\mu^s (\R^2)} &\les \| f \|_{\C_{\mu_1}^{s_1} (\R^2)} \| g \|_{\C_{\mu_2}^{s_2} (\R^2)} .
\end{align*}
\end{lemma}

We recall that $G \in C^\infty (\R^2 \setminus \{0\})$ is a radial truncated Green's function such that $G (x) = - \frac{1}{2 \pi} \log |x|$ for $|x| \leq \frac 18$ and $G(x) = 0$ for $|x| \geq \frac 14$. We show the following smoothing properties of $G$ and $\nb G$.

\begin{lemma}
\label{LEM:G}
Let $1 \leq p, q \leq \infty$, $s \in \R$, and $\mu \in \R$. Then, we have
\begin{align*}
\| G *_{\R^2} f \|_{\B_{p, q, \mu}^s (\R^2)} &\les \| f \|_{\B_{p, q, \mu}^{s - 2} (\R^2)} , \\
\| \nb G *_{\R^2} f \|_{\B_{p, q, \mu}^s (\R^2)} &\les \| f \|_{\B_{p, q, \mu}^{s - 1} (\R^2)}.
\end{align*}
\end{lemma}

\begin{proof}
We only show the estimate for $\nb G$, since the estimate for $G$ is similar and simpler. Given a dyadic $N \geq 1$, we have
\begin{align}
\| \Dl_N (\nb G *_{\R^2} f) \|_{L^p_\mu (\R^2)} \leq \sum_{M = \frac{N}{2}, N, 2N} \| \Dl_N \nb G *_{\R^2} \Dl_M f \|_{L^p_\mu (\R^2)} .
\label{G0}
\end{align}

\noi
A direct computation shows that $|\nb G (x)| \les |x|^{-1} \ind_{|x| \leq \frac 14}$, and so $\nb G \in L^1_{|\mu|} (\R^2)$. We recall the definition of $\psi_N$ in \eqref{defpsi} and write $\eta_N = \mathcal{F}_{\R^2}^{-1} (\psi_N)$. When $N = 1$, we use Young's convolution inequalities to obtain
\begin{align}
\begin{split}
\| \Dl_1 \nb G *_{\R^2} \Dl_M f \|_{L^p_\mu (\R^2)} &\les \| \Dl_1 \nb G \|_{L^1_{|\mu|} (\R^2)} \| \Dl_M f \|_{L^p_\mu (\R^2)} \\
&\les \| \eta_1 \|_{L^1_{|\mu|} (\R^2)} \| \nb G \|_{L_{|\mu|}^1 (\R^2)} \| \Dl_M f \|_{L^p_\mu (\R^2)} \\
&\les \| \Dl_M f \|_{L^p_\mu (\R^2)} .
\end{split}
\label{G1}
\end{align}

\noi
When $N \geq 2$, we note that $F = G *_{\R^2} f$ satisfies $\Dl F = f + \varphi *_{\R^2} f$ for some $\varphi \in C_c^\infty (\R^2)$, so that
\begin{align}
\Dl_N \nb G *_{\R^2} \Dl_M f = - \Dl_N \nb (- \Dl)^{-1} \Dl_M f - \Dl_N \nb (- \Dl)^{-1} \varphi *_{\R^2} \Dl_M f.
\label{G2}
\end{align}

\noi
By defining
\begin{align*}
\wt \psi_N (\xi) = \frac{N \xi}{|\xi|^{2}} \psi_N (\xi) \quad \text{and} \quad \wt \eta_N = \mathcal{F}_{\R^2}^{-1} (\wt \psi_N),
\end{align*}

\noi
we use \eqref{G2} and Young's convolution inequalities to obtain
\begin{align*}
\begin{split}
\| &\Dl_N \nb G *_{\R^2} \Dl_M f \|_{L^p_\mu (\R^2)} \\
&\les N^{-1} \| \wt \eta_N *_{\R^2} \Dl_M f \|_{L^p_\mu (\R^2)} + N^{-1} \| \wt \eta_N *_{\R^2} \varphi *_{\R^2} \Dl_M f \|_{L^p_\mu (\R^2)} \\
&\les N^{-1} \| \wt \eta_N \|_{L^1_{|\mu|} (\R^2)} \| \Dl_M f \|_{L^p_\mu (\R^2)} + N^{-1} \| \wt \eta_N \|_{L^1_{|\mu|} (\R^2)} \| \varphi \|_{L^1_{|\mu|} (\R^2)} \| \Dl_M f \|_{L^p_\mu (\R^2)} ,
\end{split}
\end{align*}

\noi
so that by using a change of variable to deal with $\| \wt \eta_N \|_{L^1_{|\mu|} (\R^2)}$, we get
\begin{align}
\| \Dl_N \nb G *_{\R^2} \Dl_M f \|_{L^p_\mu (\R^2)} \les N^{-1} \| \Dl_M f \|_{L^p_\mu (\R^2)} .
\label{G3}
\end{align}

\noi
The desired estimate then follows easily from \eqref{G0}, \eqref{G1}, and \eqref{G3}.
\end{proof}

Before moving on, we record the following estimate on convolution with a smooth and compactly supported function. This estimate is a slight extension of that in \cite[Lemma~2.5]{DLTV} and the proof follows with minor modifications.
\begin{lemma}
\label{LEM:varphi}
For any $1 \leq p, q \leq \infty$, $s_1, s_2 \in \R$, $\mu \in \R$, and $\varphi \in C_c^\infty (\R^2)$, we have
\begin{align*}
\| \varphi *_{\R^2} f \|_{\B_{p, q, \mu}^{s_1} (\R^2)} \les \| f \|_{\B_{p, q, \mu}^{s_2} (\R^2)} .
\end{align*}
\end{lemma}

\subsection{Some useful estimates from counting and probability theory}

In this subsection, we present some more useful estimates.

We first show the following convolution lemma. 
\begin{lemma}
\label{LEM:convc}
Let $0 < \al, \be < 2$ be such that $\al + \be > 2$. Then, we have
\begin{align*}
\int_{\R^2} \frac{1}{\jb{x - a}^\al \jb{x}^\be} dx \les \jb{a}^{2 - \al - \be}.
\end{align*}
\end{lemma}

\begin{proof}
We write
\begin{align*}
\int_{\R^2} \frac{1}{\jb{x - a}^\al \jb{x}^\be} dx = \textup{I}_1 + \textup{I}_2 + \textup{I}_3 ,
\end{align*}

\noi
where $\textup{I}_1$ denotes the contribution from $\{ |x - a| \leq \frac{|a|}{2} \}$, $\textup{I}_2$ denotes the contribution from $\{ |x| \leq \frac{|a|}{2} \}$, and $\textup{I}_3$ denotes the contribution from $\{ |x - a| > \frac{|a|}{2}, |x| > \frac{|a|}{2} \}$. For $\textup{I}_1$, by the triangle inequality, we have $|x| \geq |a| - |x - a| \geq \frac{|a|}{2}$, so that
\begin{align*}
\textup{I}_1 \les \jb{a}^{- \be} \int_{\{ |x - a| \leq \frac{|a|}{2} \}} \frac{1}{\jb{x - a}^\al} dx \les \jb{a}^{2 - \al - \be} .
\end{align*}

\noi
For $\textup{I}_2$, by the triangle inequality, we have $|x - a| \geq |a| - |x| \geq \frac{|a|}{2}$, so that
\begin{align*}
\textup{I}_2 \les \jb{a}^{- \al} \int_{\{ |x| \leq \frac{|a|}{2} \}} \frac{1}{\jb{x}^\be} dx \les \jb{a}^{2 - \al - \be} .
\end{align*}

\noi
For $\textup{I}_3$, by the condition $|x - a| > \frac{|a|}{2}$ and the triangle inequality, we have
\begin{align*}
|x - a| > \frac{3}{8}|a| + \frac 14 |x - a| \geq \frac 14 |x| + \frac 18 |a| > \frac 14 |x| .
\end{align*}

\noi
Thus, we have
\begin{align*}
\textup{I}_3 \les \int_{\{ |x| > \frac{|a|}{2} \}} \frac{1}{\jb{x}^{\al + \be}} dx \les \jb{a}^{2 - \al - \be}.
\end{align*}

\noi
Thus, we have finished the proof.
\end{proof}

Lastly, we record the following lemma from probability theory, which is a special case of \cite[Proposition~3.1]{HV} (see also \cite[Proposition~2.3]{Ver}).
\begin{lemma}
\label{LEM:Ver}
Let $(X, \| \cdot \|_X)$ be a separable Banach space and let $X'$ be its dual space. Let $\{ \zeta_j \}_{j \in \N}$ be a sequence of $X$-valued random variables. Suppose that there exists $0 < \s < 1$ such that for all $f \in X'$, we have
\begin{align*}
\E \big[ | (\zeta_j, f) |^2 \big] \leq \s^{2 j} \| f \|_{X'}^2 .
\end{align*}

\noi
Then, we have
\begin{align*}
\E \Big[ \sup_{j \in \N} \| \zeta_j \|_X \Big] \leq \sup_{j \in \N} \E [ \| \zeta_j \|_X ] + C(\s)
\end{align*}

\noi
for some constant $C(\s) > 0$.
\end{lemma}

\section{Weighted function spaces on torus with arbitrary size}
\label{SEC:BesL}

In this section, we introduce weighted function spaces on periodic domains. 

\subsection{Periodic weights}
Let us first introduce weights on periodic domains. Given $L \geq 1$, let
\begin{align*}
\jb{x}_{L} = \jb{x - Lk} \quad \text{if } x \in Lk + [ - \tfrac{L}{2}, \tfrac{L}{2} )^2 \text{ with } k \in \Z^2,
\end{align*}

\noi
which is continuous and $L$-periodic. We now perform a regularization procedure. Let $\rho$ be a smooth nonnegative function on $\R^2$ supported on $\{|x| \leq \frac 14\}$ such that $\int_{\R^2} \rho = 1$. Given $\mu \in \R$, let
\begin{align*}
\jbb{x}_{L, \mu} = \big( \rho *_{\R_2} \jb{\cdot}_L^\mu \big) (x),
\end{align*}

\noi
which is smooth and $L$-periodic. Note that $\jbb{x}_{L, 0} = 1$. 

Let us show the following lemma on various properties of of the $L$-periodic weights $\jbb{x}_{L, \mu}$.

\begin{lemma}
\label{LEM:Lwei}
Let $L \geq 1$.

\smallskip \noi
\textup{(i)} For any $\mu \geq 0$ and $x \in \R^2$, we have
\begin{align*}
\jbb{x}_{L, \mu} > 2^{- \mu}.
\end{align*}

\smallskip \noi
\textup{(ii)} For any $\mu \in \R$ and $x \in \R^2$, we have
\begin{align*}
2^{- |\mu|} \jb{x}_L^{\mu} < \jbb{x}_{L, \mu} < 2^{|\mu|} \jb{x}_L^\mu .
\end{align*}

\smallskip \noi
\textup{(iii)} For any $\mu_1, \mu_2 \in \R$ and $x \in \R^2$, we have
\begin{align*}
2^{- \max (|\mu_1|, |\mu_2|)} \jbb{x}_{L, \mu_1 + \mu_2} < \jbb{x}_{L, \mu_1} \jbb{x}_{L, \mu_2} < 2^{\min (|\mu_1|, |\mu_2|)} \jbb{x}_{L, \mu_1 + \mu_2} .
\end{align*}

\smallskip \noi
\textup{(iv)} For any $\mu \in \R$, $a > 0$, and $x \in \R^2$, we have
\begin{align*}
2^{- 2 a |\mu|} \jbb{x}_{L, a \mu} < \jbb{x}_{L, \mu}^a < 2^{2 a |\mu|} \jbb{x}_{L, a \mu}.
\end{align*}

\smallskip \noi
\textup{(v)} For any $\mu \in \R$ and $x, y \in \R^2$, we have
\begin{align*}
\jbb{x + y}_{L, \mu} \leq 2^{4 |\mu|} \jbb{x}_{L, |\mu|} \jbb{y}_{L, \mu}.
\end{align*}

\smallskip \noi
\textup{(vi)} For any $\mu \in \R$ and $x \in \R^2$, we have
\begin{align*}
\nb \jbb{x}_{L, \mu} = \mu \int_{\R^2} \rho (x - y) \wt{\jb{y}}_{L, \mu - 1} dy ,
\end{align*}

\noi
where
\begin{align*}
\wt{\jb{y}}_{L, \mu - 1} = \jb{y - Lk}^{\mu - 2} (y - Lk) \quad \text{if } y \in Lk + [-\tfrac{L}{2}, \tfrac{L}{2})^2 \text{ with } k \in \Z^2.
\end{align*}

\noi
Consequently, we have
\begin{align*}
\big| \nb \jbb{x}_{L, \mu} \big| \leq |\mu| \jbb{x}_{L, \mu - 1} .
\end{align*}
\end{lemma}

\begin{proof}
(i) We first show that if $x, y \in \R^2$ and $|x - y| \leq \frac 12$, then
\begin{align}
\frac 14 \jb{y}_L^2 < \jb{x}_L^2 < 4 \jb{y}_L^2.
\label{xyL}
\end{align}

\noi
We may assume without loss of generality that $x \in [-\frac{L}{2}, \frac{L}{2})^2$. If $y \in [-\frac{L}{2}, \frac{L}{2})^2$, then a straightforward computation gives 
\begin{align*}
\frac 12 \jb{y}^2 \leq \jb{x}^2 \leq 2 \jb{y}^2,
\end{align*} 

\noi
which implies \eqref{xyL}. If $y \in Lk + [-\frac{L}{2}, \frac{L}{2})^2$ for some $k \neq 0$, then since $|x - y| \leq \frac 12$, we must have $|x| \geq \frac{L}{2} - \frac 12$ and $|y - Lk| \geq \frac{L}{2} - \frac 12$. Thus, we have
\begin{align*}
\jb{y}_L^2 \leq 1 + \frac{L^2}{2} \leq 2 + 4 |x|^2 < 4 \jb{x}^2 = 4 \jb{x}_L^2
\end{align*}

\noi
and
\begin{align*}
\jb{x}_L^2 \leq 1 + \frac{L^2}{2} \leq 2 + 4 |y - Lk|^2 < 4 \jb{y - Lk}^2 = 4 \jb{y}_L^2,
\end{align*}

\noi
so that \eqref{xyL} holds.

From the support property of $\rho$ and \eqref{xyL}, we have
\begin{align*}
\jbb{x}_{L, \mu} = \int_{\R^2} \rho (x - y) \jb{y}_L^{\mu} dy > 2^{- \mu} \jb{x}_L^\mu \geq 2^{- \mu},
\end{align*}

\noi
as desired.

\smallskip \noi
(ii) Using the support property of $\rho$ and \eqref{xyL}, we have
\begin{align*}
\jbb{x}_{L, \mu} = \int_{\R^2} \rho (x - y) \jb{y}_L^{\mu} dy < 2^{|\mu|} \jb{x}_L^{\mu} \int_{\R^2} \rho (x - y) dy = 2^{|\mu|} \jb{x}_L^{\mu} 
\end{align*}

\noi
and
\begin{align*}
\jbb{x}_{L, \mu} = \int_{\R^2} \rho (x - y) \jb{y}_L^{\mu} dy > 2^{-|\mu|} \jb{x}_L^{\mu} \int_{\R^2} \rho (x - y) dy = 2^{-|\mu|} \jb{x}_L^{\mu} ,
\end{align*}

\noi
as desired.

\smallskip \noi
(iii) We may assume without loss of generality that $|\mu_1| \leq |\mu_2|$. Using the support property of $\rho$ and \eqref{xyL}, we have
\begin{align*}
\jbb{x}_{L, \mu_1} \jbb{x}_{L, \mu_2} &= \int_{\R^2} \int_{\R^2} \rho (x - y_1) \jb{y_1}_L^{\mu_1} \rho (x - y_2) \jb{y_2}_L^{\mu_2} dy_1 dy_2 \\
&< 2^{|\mu_1|} \int_{\R^2} \rho(x - y_1) \rho(x - y_2) \jb{y_2}_L^{\mu_1 + \mu_2} dy_1 dy_2 \\
&= 2^{|\mu_1|} \jbb{x}_{L, \mu_1 + \mu_2}
\end{align*}

\noi
and
\begin{align*}
\jbb{x}_{L, \mu_1} \jbb{x}_{L, \mu_2} &= \int_{\R^2} \int_{\R^2} \rho (x - y_1) \jb{y_1}_L^{\mu_1} \rho (x - y_2) \jb{y_2}_L^{\mu_2} dy_1 dy_2 \\
&> 2^{- |\mu_2|} \int_{\R^2} \rho(x - y_1) \rho(x - y_2) \jb{y_1}_L^{\mu_1 + \mu_2} dy_1 dy_2 \\
&= 2^{- |\mu_2|} \jbb{x}_{L, \mu_1 + \mu_2},
\end{align*}

\noi
as desired.

\smallskip \noi
(iv) Using the support property of $\rho$, \eqref{xyL}, and part (ii), we have
\begin{align*}
\jbb{x}_{L, a \mu} = \int_{\R^2} \rho (x - y) \jb{y}_L^{a \mu} dy < 2^{a |\mu|} \jb{x}_L^{a \mu} < 2^{2 a |\mu|} \jbb{x}_{L, \mu}^a
\end{align*}

\noi
and
\begin{align*}
\jbb{x}_{L, a \mu} = \int_{\R^2} \rho (x - y) \jb{y}_L^{a \mu} dy > 2^{- a |\mu|} \jb{x}_L^{a \mu} > 2^{- 2 a |\mu|} \jbb{x}_{L, \mu}^a,
\end{align*}

\noi
as desired.

\smallskip \noi
(v) From part (ii), we only need to show
\begin{align*}
\jb{x + y}_L \leq 2 \jb{x}_L \jb{y}_L.
\end{align*}

\noi
Due to $L$-periodicity, we only need to show for $x, y \in [- \frac{L}{2}, \frac{L}{2})^2$,
\begin{align*}
\jb{x + y}_L \leq 2 \jb{x} \jb{y}.
\end{align*}

\noi
This is obvious since $\jb{x + y}_L \leq \jb{x + y} \leq 2 \jb{x} \jb{y}$.

\smallskip \noi
(vi) Using the dominated convergence theorem, we have
\begin{align*}
\nb \jbb{x}_{L, \mu} &= \int_{\R^2} \nb_x \rho (x - y) \jb{y}_L^\mu dy = - \int_{\R^2} \nb_y \rho (x - y) \jb{y}_L^\mu dy.
\end{align*}

\noi
The identity and the estimate then follow from an integration by parts and the fact that $| \wt{\jb{y}}_{L, \mu - 1} | \leq \jb{y}_L^{\mu - 1}$. However, the use of the integration by parts needs to be justified since $\jb{\cdot}_L^\mu$ is not a $C^1$ function. Let us write $y = (y_1, y_2) \in \R^2$. By Fubini's theorem, we only need to show
\begin{align}
- \int_\R \partial_{y_1} \rho (x - y) \jb{y}_L^\mu dy_1 = \int_\R \rho (x - y) \partial_{y_1} \jb{y}_L^\mu dy_1
\label{ibp}
\end{align}

\noi
for any $y_2 \notin \{ \frac{L}{2} + Lk' : k' \in \Z \}$. Since $\rho$ is a smooth and compactly supported function and $y_1 \mapsto \jb{y}_L^\mu$ is an absolutely continuous function, we can apply \cite[Theorem~3.35 and Theorem~3.36]{Fol} to obtain the desired identity \eqref{ibp}.
\end{proof}

\begin{remark} \rm
Part (vi) of Lemma~\ref{LEM:Lwei} allows us to differentiate the $L$-periodic weights and estimate the derivatives in a similar manner to the polynomial weights $\jb{\cdot}^\mu$ (with $\mu \in \R$) on $\R^2$. This is a crucial point for our analysis in this paper; see Lemma~\ref{LEM:St_H1} below. At this point, we do not know how to obtain a similar estimate for the second or higher derivatives of the $L$-periodic weights, but fortunately we do not need it in this paper.
\end{remark}

\begin{remark} \rm
One may also proceed with the weights $\jb{\cdot}_L^\mu$ (with $\mu \in \R$) on $\T_L^2$ without regularization. Note that the function $\jb{\cdot}_L^\mu$ is not $C^1$ since it is not differentiable on $\{ (\frac{L}{2} + L k'_1, \frac{L}{2} + L k'_2) \in \R^2 : k'_1, k'_2 \in \Z \}$, and so some separate analysis is required on this set if one wants to differentiate the weight. This is relevant below in the commutator estimate in Lemma~\ref{LEM:Lcomm} via the mean value theorem and the energy estimate in Lemma~\ref{LEM:St_H1} via an integration by parts. In this paper, however, we choose to regularize our weights on $\T_L^2$ so as to avoid a separate discussion on the set where $\jb{\cdot}_L^\mu$ is not differentiable. Moreover, it might be of interest to investigate other properties of the weights such as their higher order derivatives, and so the smoothness of the weights guaranteed by the regularization procedure makes it more convenient.
\end{remark}

Below and throughout the rest of this paper, we will suppress the constant dependence on the weight parameter $\mu$. There is only one place where such dependence matters (Lemma~\ref{LEM:YregLe} below), but we will make the assumption $|\mu| \leq 1$ so that the constants can be made independent of $\mu$.

\subsection{Periodic weighted function spaces}
We now discuss weighted function spaces on $\T_L^2$ given $L \geq 1$. Given $L \geq 1$, $\mu \in \R$, and $1 \leq p \leq \infty$, we define
\begin{align*}
\| f_L \|_{L^p_{\mu} (\T_L^2)} \deff \bigg( \int_{\T_L^2} |f_L (x)|^p \jbb{x}_{L, \mu}^p dx \bigg)^{\frac{1}{p}} = \big\| \jbb{\cdot}_{L, \mu} f_L \big\|_{L^p (\T_L^2)}.
\end{align*}

\noi
with the obvious interpretation if $p = \infty$.
Since $\jbb{\cdot}_{L, 0} \equiv 1$, we have $\| f_L \|_{L_0^p (\T_L^2)} = \| f_L \|_{L^p (\T_L^2)}$, and so we write $L^p (\T_L^2) = L_0^p (\T_L^2)$. 
When $1 \leq p < \infty$, by using Lemma~\ref{LEM:Lwei}~(iii), we easily see the duality relation $(L_{\mu}^p (\T_L^2))' = L_{- \mu}^{p'} (\T_L^2)$ with $\frac{1}{p} + \frac{1}{p'} = 1$.

\medskip
We also need to define weighted Besov spaces. For this purpose, let us first clarify the Fourier transform and define $L^2$-based Sobolev spaces on $\T_L^2$.
Let $L \geq 1$. Given $n \in \Z^2$, we define the Fourier transform of a function $f_L$ on $\T_L^2$ at frequency $\frac{n}{L}$ by
\begin{align*}
\mathcal{F}_{\T_L^2} (f_L) (\tfrac{n}{L}) = \ft {f_L} ( \tfrac{n}{L} ) \deff \frac{1}{L} \int_{\T_L^2} f_L (x) e^{- 2 \pi i \frac{n}{L} \cdot x} dx.
\end{align*}

\noi
The inverse Fourier transform is then given by
\begin{align*}
f_L (x) = \frac{1}{L} \sum_{n \in \Z^2} \ft{f_L} ( \tfrac{n}{L} ) e^{2 \pi i \frac{n}{L} \cdot x} .
\end{align*}

\noi
We have the following Plancherel's (or Parseval's) identity on $\T_L^2$:
\begin{align}
\int_{\T_L^2} f_L (x) \cj{g_L (x)} dx = \sum_{n \in \Z^2} \ft{f_L} ( \tfrac{n}{L} ) \cj{\ft{g_L} ( \tfrac{n}{L} )}.
\label{planL}
\end{align}

\noi
For any $L$-periodic functions $f_L$ and $g_L$ on $\T_L^2$, we have
\begin{align}
\frac{1}{L} (f_L *_{\T_L^2} g_L)^\wedge (\tfrac{n}{L}) = \ft{f_L} (\tfrac{n}{L}) \ft{g_L} (\tfrac{n}{L}).
\label{convo}
\end{align}

\noi
For any Schwartz function $f$ on $\R^2$, we have the following Poisson summation formula (see \cite[Theorem~3.2.8]{Gra}):
\begin{align}
L \sum_{k \in \Z^2} f (x + L k) = \frac{1}{L} \sum_{n \in \Z^2} \mathcal{F}_{\R^2} (f) (\tfrac{n}{L}) e^{2 \pi i \frac{n}{L} \cdot x} ,
\label{Poi}
\end{align}

\noi
where $\mathcal{F}_{\R^2}$ denotes the Fourier transform on $\R^2$.

We denote by $\jb{\nb_{\T_L^2}}$ the Fourier multiplier operator with symbol $\jb{\frac{n}{L}}$. Given $s \in \R$, we define the $L^2$-based Sobolev space $H^s (\T_L^2)$ via the norm
\begin{align*}
\| f_L \|_{H^s (\T_L^2)} = \| \jb{\nb_{\T_L^2}}^s f_L \|_{L^2 (\T_L^2)} = \bigg( \sum_{n \in \Z^2} \jb{\tfrac{n}{L}}^{2s} | \ft{f_L} ( \tfrac{n}{L} ) |^2 \bigg)^{\frac 12},
\end{align*}

\noi
where the second equality follows from Plancherel's identity \eqref{planL}. It is not hard to see that
\begin{align*}
\| f_L \|_{H^1 (\T_L^2)} &\sim \| f_L \|_{L^2 (\T_L^2)} + \| \nb f_L \|_{L^2 (\T_L^2)}, \\
\| f_L \|_{H^2 (\T_L^2)} &\sim \| f_L \|_{L^2 (\T_L^2)} + \| \Dl f_L \|_{L^2 (\T_L^2)}
\end{align*}

\noi
with the underlying constants independent of $L$. Given $\mu \in \R$, we define the weighted $H^1 (\T_L^2)$-norm and the weighted $H^2 (\T_L^2)$-norm as
\begin{align}
\begin{split}
\| f_L \|_{H_\mu^1 (\T_L^2)} &\deff \| f_L \|_{L_\mu^2 (\T_L^2)} + \| \nb f_L \|_{L_\mu^2 (\T_L^2)}, \\
\| f_L \|_{H_\mu^2 (\T_L^2)} &\deff \| f_L \|_{L_\mu^2 (\T_L^2)} + \| \Dl f_L \|_{L_\mu^2 (\T_L^2)} .
\end{split}
\label{H2mu}
\end{align}

We now define weighted Besov spaces on $\T_L^2$. 
Given $L \geq 1$ and a dyadic number $N \geq 1$, we define $\Dl_N^L$ as the Fourier multiplier operator with symbol $\psi_N$ defined in \eqref{defpsi}. This means that for any $L$-periodic function $f_L$, we have
\begin{align*}
\Dl_N^L f_L = \frac{1}{L} \sum_{n \in \Z^2} \psi_N (\tfrac{n}{L}) \ft{f_L} (\tfrac{n}{L}) e^{2 \pi i \frac{n}{L} \cdot x} .
\end{align*}

\noi
Note that from \eqref{defpsi}, for any $L \geq 1$, we have
\begin{align}
\Id = \sum_{\substack{N \geq 1 \\ \text{dyadic}}} \Dl_N^L.
\label{idDNL}
\end{align}

\noi
We also denote $\Dl_{\frac 12}^L = 0$. Given a dyadic $N \geq 1$, we will frequently use the identity
\begin{align*}
\Dl_N^L = \Dl_N^L \big( \Dl_{\frac{N}{2}}^L + \Dl_N^L + \Dl_{2N}^L \big) .
\end{align*}

Given $1 \leq p, q \leq \infty$, $s \in \R$, and $\mu \in \R$, we define the weighted Besov space $\B_{p, q, \mu}^s (\T_L^2)$ as
\begin{align*}
\| f_L \|_{\B_{p, q, \mu}^s (\T_L^2)} \deff \bigg( \sum_{\substack{N \geq 1 \\ \text{dyadic}}} N^{s q} \| \Dl_N^L f_L \|_{L_\mu^p (\T_L^2)}^q \bigg)^{\frac 1q} = \bigg( \sum_{\substack{N \geq 1 \\ \text{dyadic}}} N^{s q} \big\| \jbb{\cdot}_{L, \mu} \Dl_N^L f_L  \big\|_{L^p (\T_L^2)}^q \bigg)^{\frac 1q} .
\end{align*}

\noi
When $p = q = \infty$, we denote
\begin{align*}
\C_\mu^s (\T_L^2) \deff \B_{\infty, \infty, \mu}^s (\T_L^2).
\end{align*}


\noi
As usual, we drop the weight parameter if $\mu = 0$. Note that when $p = q = 2$ and $\mu = 0$, we have the equivalence
\begin{align}
\| f_L \|_{H^s (\T_L^2)} \sim \| f_L \|_{\B^s_{2, 2} (\T_L^2)}
\label{Hs_equi}
\end{align}

\noi
as a direct consequence of Plancherel's identity \eqref{planL} and the support property of $\psi_N$ defined in \eqref{defpsi}, where the underlying constant is independent of $L$.

We also define
\begin{align}
\eta_N \deff \mathcal{F}_{\R^2}^{-1} (\psi_N) \quad \text{and} \quad \eta_N^L \deff \mathcal{F}_{\T_L^2}^{-1} (\psi_N) = \frac 1L \sum_{n \in \Z^2} \psi_N (\tfrac{n}{L}) e^{2 \pi i \frac{n}{L} \cdot x}.
\label{defeta}
\end{align}

\noi
From \eqref{convo}, we have
\begin{align}
\Dl_N^L f_L = \frac{1}{L} (\eta_N^L *_{\T_L^2} f_L) .
\label{convo2}
\end{align}

\noi
By the Poisson summation formula \eqref{Poi}, we get that for any $x \in [-\frac{L}{2}, \frac{L}{2})^2$, 
\begin{align}
\eta_N^L (x) = \frac{1}{L} \sum_{n \in \Z^2} \psi_N (\tfrac{n}{L}) e^{2 \pi i \frac{n}{L} \cdot x} = L \sum_{k \in \Z^2} \eta_N (x + Lk).
\label{Poi_eta}
\end{align}

Let us now show the following embedding results uniformly in the size of the torus.
\begin{lemma}
\label{LEM:embL}
Let $L \geq 1$ and let $f_L$ be an $L$-periodic distribution. 

\smallskip \noi
\textup{(i)} Let $1 \leq p \leq \infty$, $1 \leq q_1, q_2 \leq \infty$, $s_1 < s_2$, and $\mu \in \R$. Then, we have
\begin{align*}
\| f_L \|_{\B_{p, q_1, \mu}^{s_1} (\T_L^2)} \les \| f_L \|_{\B_{p, \infty, \mu}^{s_2} (\T_L^2)} \les \| f_L \|_{\B_{p, q_2, \mu}^{s_2} (\T_L^2)}.
\end{align*}

\smallskip \noi
\textup{(ii)} Let $1 \leq p \leq \infty$ and $\mu \in \R$. Then, we have
\begin{align*}
\| f_L \|_{\B_{p, \infty, \mu}^0 (\T_L^2)} \les \| f_L \|_{L_\mu^p (\T_L^2)} \les \| f_L \|_{\B_{p, 1, \mu}^0 (\T_L^2)}.
\end{align*}

\smallskip \noi
\textup{(iii)} Let $1 \leq p_2 \leq p_1 \leq \infty$ and $\mu \in \R$. Then, for any dyadic $N \geq 1$, we have
\begin{align*}
\| \Dl_N^L f_L \|_{L_\mu^{p_1} (\T_L^2)} \les N^{2 (\frac{1}{p_2} - \frac{1}{p_1})} \| \Dl_N^L f_L \|_{L_\mu^{p_2} (\T_L^2)} ,
\end{align*}

\noi
where the underlying constant is independent of $p_1$ and $p_2$. Consequently, with $1 \leq q_2 \leq q_1 \leq \infty$ and $s_2 \geq s_1 + 2 (\frac{1}{p_2} - \frac{1}{p_1})$, we have
\begin{align*}
\| f_L \|_{\B_{p_1, q_1, \mu}^{s_1} (\T_L^2)} \les \| f_L \|_{\B_{p_2, q_2, \mu}^{s_2} (\T_L^2)} ,
\end{align*}

\noi
where the underlying constant is independent of $p_1$ and $p_2$.

\smallskip \noi
\textup{(iv)} Let $1 \leq p, q \leq \infty$, $s \in \R$, and $\mu_1, \mu_2 \in \R$ with $\mu_1 \leq \mu_2$. Then, we have
\begin{align*}
\| f_L \|_{L^p_{\mu_1} (\T_L^2)} &\les \| f_L \|_{L^p_{\mu_2} (\T_L^2)} , \\
\| f_L \|_{\B^s_{p, q, \mu_1} (\T_L^2)} &\les \| f_L \|_{\B^s_{p, q, \mu_2} (\T_L^2)}.
\end{align*}

\smallskip \noi
\textup{(v)} Let $1 \leq p \leq \infty$ and $\mu \in \R$. Then, for any dyadic $N \geq 1$, we have
\begin{align*}
\| \Dl_N^L \nb f_L \|_{L_\mu^p (\T_L^2)} &\les \sum_{M = \frac{N}{2}, N, 2N} M \| \Dl_M^L f_L \|_{L_\mu^p (\T_L^2)} , \\
\| \Dl_N^L \Dl f_L \|_{L_\mu^p (\T_L^2)} &\les \sum_{M = \frac{N}{2}, N, 2N} M^2 \| \Dl_M^L f_L \|_{L_\mu^p (\T_L^2)} .
\end{align*}

\noi
Consequently, with $1 \leq q \leq \infty$ and $s > 0$, we have
\begin{align*}
\| f_L \|_{\B_{p, q, \mu}^{s} (\T_L^2)} &\sim \| f_L \|_{L^p_\mu (\T_L^2)} + \| \nb f_L \|_{\B_{p, q, \mu}^{s - 1} (\T_L^2)} , \\
\| f_L \|_{\B_{p, q, \mu}^{s} (\T_L^2)} &\sim \| f_L \|_{L^p_\mu (\T_L^2)} + \| \Dl f_L \|_{\B_{p, q, \mu}^{s - 2} (\T_L^2)} .
\end{align*}
\end{lemma}

\begin{proof}
(i) The first embedding follows directly from the $\ell^{q_1}$-summability of $\{ N^{s_1 - s_2} \}_{N \in 2^{\N}}$. The second embedding is immediate due to $\ell^{q_2} \hookrightarrow \ell^\infty$.

\smallskip \noi
(ii) The second embedding follows directly from the definition and the triangle inequality. To prove the first embedding,
we note that for any $L \geq 1$ and dyadic $N \geq 1$, we have from \eqref{convo2}, Young's convolution inequality along with Lemma~\ref{LEM:Lwei}~(v), \eqref{Poi_eta}, and Lemma~\ref{LEM:Lwei}~(ii) that
\begin{align*}
\| \Dl_N^L f_L \|_{L^p_\mu (\T_L^2)} &= \frac{1}{L} \| \eta_N^L *_{\T_L^2} f_L \|_{L^p_{\mu} (\T_L^2)} \\
&\leq \frac{1}{L} \big\| \jbb{\cdot}_{L, |\mu|} \eta_N^L \big\|_{L^1 (\T_L^2)} \| f_L \|_{L^p_{\mu} (\T_L^2)} \\
&\les \| \jb{\cdot}^{|\mu|} \eta_N \|_{L^1 (\R^2)} \| f_L \|_{L^p_{\mu} (\T_L^2)} ,
\end{align*}

\noi
so that by using a change of variable to deal with $\| \jb{\cdot}^{|\mu|} \eta_N \|_{L^1 (\R^2)}$, we get
\begin{align*}
\| \Dl_N^L f_L \|_{L^p_\mu (\T_L^2)} \les \| f_L \|_{L^p_{\mu} (\T_L^2)} .
\end{align*}

\noi
The first embedding then follows from taking the supremum over dyadic $N \geq 1$. 

\smallskip \noi
(iii) Given $L \geq 1$ and a dyadic number $N \geq 1$, we use \eqref{convo2} and Young's convolution inequality along with Lemma~\ref{LEM:Lwei}~(v) to get
\begin{align*}
\| \Dl_N^L f_L \|_{L_\mu^{p_1} (\T_L^2)} &= \big\| \big( \Dl_{\frac{N}{2}}^L + \Dl_N^L + \Dl_{2N}^L \big) \Dl_N^L f_L \big\|_{L_\mu^{p_1} (\T_L^2)} \\
&\leq \frac{1}{L} \sum_{M = \frac{N}{2}, N, 2N} \| \eta_M^L *_{\T_L^2} \Dl_N^L f_L \|_{L_\mu^{p_1} (\T_L^2)} \\
&\leq \frac{1}{L} \sum_{M = \frac{N}{2}, N, 2N} \big\| \jbb{\cdot}_{L, |\mu|} \eta_M^L \big\|_{L^r (\T_L^2)} \| \Dl_N^L f_L \|_{L_\mu^{p_2} (\T_L^2)} ,
\end{align*}

\noi
where $\frac 1r = 1 + \frac{1}{p_1} - \frac{1}{p_2}$. By \eqref{Poi_eta}, Lemma~\ref{LEM:Lwei}~(ii), H\"older's inequality, and a change of variable, we get
\begin{align*}
\big\| \jbb{\cdot}_{L, |\mu|} \eta_M^L \big\|_{L^r ([-\frac{L}{2}, \frac{L}{2})^2)} 
&= L \Big\| \jb{x}^{|\mu|} \sum_{k \in \Z^2} \eta_M (x + L k) \Big\|_{L_x^r (\T_L^2)} \\
&\leq L \bigg\| \jb{x}^{|\mu|} \Big( \sum_{k \in \Z^2} \jb{k}^{- 10 r'} \Big)^{\frac{1}{r'}} \| \jb{k}^{10} \eta_M (x + L k) \|_{\ell_k^r (\Z^2)} \bigg\|_{L_x^r ([-\frac{L}{2}, \frac{L}{2})^2)} \\
&\les L \| \jb{\cdot}^{10 + |\mu|} \eta_M (\cdot) \|_{L^r (\R^2)} .
\end{align*}

\noi
The desired embedding then follows from a change of variable
\begin{align*}
\| \jb{\cdot}^{10 + |\mu|} \eta_M \|_{L^r (\R^2)} \leq M^{2 - \frac{2}{r}} \| \jb{\cdot}^{10 + |\mu|} \eta_1 \|_{L^r (\R^2)} \les M^{2 (\frac{1}{p_2} - \frac{1}{p_1})}
\end{align*}
 
\noi 
and a summation over the dyadic pieces, where the underlying constant can be made independent of $r$ thanks to the fact that $\eta_1$ is a Schwartz function.

\smallskip \noi
(iv) The proof follows directly from the definition and Lemma~\ref{LEM:Lwei}~(iii) and (i).

\smallskip \noi
(v) We focus on the estimates for $\nb f_L$, since the estimates for $\Dl f_L$ follow from a similar manner. Note that for any dyadic number $N \geq 1$, we have
\begin{align*}
\Dl_N^L \nb f_L = \frac{2 \pi i}{L} \sum_{n \in \Z^2} \frac{n}{L} \psi_N (\tfrac{n}{L}) \ft{f_L} (\tfrac{n}{L}) e^{2 \pi i \frac{n}{L} \cdot x} = \frac{2 \pi i N}{L} \sum_{n \in \Z^2} \frac{n}{NL} \psi_N (\tfrac{n}{L}) \ft{f_L} (\tfrac{n}{L}) e^{2 \pi i \frac{n}{L} \cdot x} .
\end{align*}

\noi
We define
\begin{align*}
\wt \psi_{1, N} (\xi) = \frac{\xi}{N} \psi_N (\xi)
\end{align*}

\noi
and denote by $\wt \Dl_{1, N}^L$ the Fourier multiplier operator by $\wt \psi_{1, N}$ on $\T_L^2$, which by the same proof in part (ii) is bounded in $L^p_\mu (\T_L^2)$ uniformly in $L$ and $N$. Thus, we have $\Dl_N^L \nb = 2 \pi i N \wt \Dl_{1, N}^L$, so that
\begin{align*}
\| \Dl_N^L \nb f_L \|_{L^p_\mu (\T_L^2)} &\leq \sum_{M = \frac{N}{2}, N, 2N} \| \Dl_N^L \nb \Dl_M^L f_L \|_{L_\mu^p (\T_L^2)} \\
&= \sum_{M = \frac{N}{2}, N, 2N} 2 \pi N \| \wt \Dl_{1, N}^L \Dl_M^L f_L \|_{L^p_\mu (\T_L^2)} \\
&\les \sum_{M = \frac{N}{2}, N, 2N} M \| \Dl_M^L f_L \|_{L^p_\mu (\T_L^2)} ,
\end{align*}

\noi
which gives the desired inequality. After taking the $\ell^q$-summation in $N$ and using the embeddings in part (ii) and (i), we get
\begin{align*}
\| f_L \|_{L_\mu^p (\T_L^2)} + \| \nb f_L \|_{\B_{p, q, \mu}^{s - 1} (\T_L^2)} \les \| f_L \|_{\B_{p, q, \mu}^{s} (\T_L^2)} ,
\end{align*}

\noi
which gives one direction of the desired equivalence.

For the other direction of the equivalence, given $N \geq 2$, we define
\begin{align*}
\wt \psi_{2, N} (\xi) = \frac{N}{\xi} \psi_N (\xi)
\end{align*}

\noi
which is smooth since $\psi_N$ is supported away from the origin when $N \geq 2$. We denote by $\wt \Dl_{2, N}^L$ the Fourier multiplier operator by $\wt \psi_{2, N}$ on $\T_L^2$, which by the same proof in part (ii) is bounded in $L_\mu^p (\T_L^2)$ uniformly in $L$ and $N$. Thus, for $N \geq 2$, we have $N \Dl_N^L = (2 \pi i)^{-1} \wt \Dl_{2, N}^L \nb$, so that
\begin{align*}
N \| \Dl_N^L f_L \|_{L^p_{\mu} (\T_L^2)} &\leq \sum_{M = \frac{N}{2}, N, 2N} N \| \Dl_N^L \Dl_M^L f_L \|_{L^p_\mu (\T_L^2)} \\
&= \sum_{M = \frac{N}{2}, N, 2N} (2 \pi)^{-1} \| \wt \Dl_{2, N}^L \Dl_M^L \nb f_L \|_{L^p_\mu (\T_L^2)} \\
&\les \sum_{M = \frac{N}{2}, N, 2N} \| \Dl_M^L \nb f_L \|_{L^p_\mu (\T_L^2)} , 
\end{align*}

\noi
and so after taking the $\ell^q$-summation in $N$ and using the $L^p_\mu (\T_L^2)$-boundedness of $\Dl_1^L$ from part (ii), we get
\begin{align*}
\| f_L \|_{\B_{p, q, \mu}^{s} (\T_L^2)} &\les \| \Dl_1^L f_L \|_{L^p_\mu (\T_L^2)} + \bigg( \sum_{\substack{N \geq 2 \\ \text{dyadic}}} N^{sq} \| \Dl_N^L f_L \|_{L^p_\mu (\T_L^2)}^q  \bigg)^{\frac{1}{q}} \\
&\les \| f_L \|_{L^p_\mu (\T_L^2)} + \bigg( \sum_{\substack{N \geq 2 \\ \text{dyadic}}} N^{(s - 1)q} \sum_{M = \frac{N}{2}, N, 2N} \| \Dl_M^L \nb f_L \|_{L^p_\mu (\T_L^2)}^q  \bigg)^{\frac{1}{q}} \\
&\les \| f_L \|_{L^p_\mu (\T_L^2)} + \| \nb f_L \|_{\B_{p, q, \mu}^{s - 1} (\T_L^2)} .
\end{align*}

\noi
Thus, we obtain the desired equivalence.
\end{proof}

We also have the following duality estimate uniformly in the size of the torus.

\begin{lemma}
\label{LEM:dualL}
Let $L \geq 1$ and let $f_L$, $g_L$ be two $L$-periodic distributions. Let $s \in \R$, $1 \leq p, p', q, q' \leq \infty$ be such that $\frac{1}{p} + \frac{1}{p'} = 1$ and $\frac{1}{q} + \frac{1}{q'} = 1$, and $\mu \in \R$. Then, we have
\begin{align*}
\bigg| \int_{\T_L^2} f_L \cj{g_L} dx \bigg| \les \| f_L \|_{\B_{p, q, \mu}^s (\T_L^2)} \| g_L \|_{\B_{p', q', -\mu}^{- s} (\T_L^2)} .
\end{align*}
\end{lemma}

\begin{proof}
The estimate follows directly from the Littlewood-Paley decomposition, Lemma~\ref{LEM:Lwei}~(iii), and H\"older's inequalities.
\end{proof}

We then record the following product estimates.

\begin{lemma}
\label{LEM:prodL}
Let $L \geq 1$ and let $f_L$, $g_L$ be two $L$-periodic distributions. 

\smallskip \noi
\textup{(i)} Let $s > 0$, $1 \leq p, p_1, p_2, p_3, p_4 \leq \infty$ be such that $\frac{1}{p} = \frac{1}{p_1} + \frac{1}{p_2} = \frac{1}{p_3} + \frac{1}{p_4}$, and $1 \leq q \leq \infty$. Then, we have
\begin{align*}
\| f_L g_L \|_{\B_{p, q}^s (\T_L^2)} \les \| f_L \|_{\B_{p_1, q}^s (\T_L^2)} \| g_L \|_{L^{p_2} (\T_L^2)} + \| f_L \|_{L^{p_3} (\T_L^2)} \| g_L \|_{\B_{p_4, q}^s (\T_L^2)} .
\end{align*}

\smallskip \noi
\textup{(ii)} Let $s, s_1, s_2 \in \R$ be such that $s_1 + s_2 > 0$ and $s = \min (s_1 + s_2, s_1, s_2)$, $1 \leq p, p_1, p_2 \leq \infty$ be such that $\frac{1}{p} = \frac{1}{p_1} + \frac{1}{p_2}$, and $\mu, \mu_1, \mu_2 \in \R$ be such that $\mu = \mu_1 + \mu_2$. Then, for any $\kappa > 0$, we have
\begin{align*}
\| f_L g_L \|_{\B_{p, p, \mu}^{s - \kappa} (\T_L^2)} \les \| f_L \|_{\B_{p_1, p_1, \mu_1}^{s_1} (\T_L^2)} \| g_L \|_{\B_{p_2, p_2, \mu_2}^{s_2} (\T_L^2)}
\end{align*}

\noi
and
\begin{align*}
\| f_L g_L \|_{\C_\mu^s (\T_L^2)} \les \| f_L \|_{\C_{\mu_1}^{s_1} (\T_L^2)} \| g_L \|_{\C_{\mu_2}^{s_2} (\T_L^2)} .
\end{align*}
\end{lemma}

\begin{proof}
The proof of the product estimates follows from adapting the steps in \cite{BCD, PT16, MW17} to the setting of $\T_L^2$ along with Lemma~\ref{LEM:Lwei}~(iii).
\end{proof}

We now prove the following interpolation inequality.
\begin{lemma}
\label{LEM:interpL}
Let $L \geq 1$ and let $f_L$ be an $L$-periodic distribution. Let $1 \leq p, p_1, p_2 \leq \infty$, $\mu, \mu_1, \mu_2 \in \R$, and $0 < \ta < 1$ be such that $\frac{1}{p} = \frac{\ta}{p_1} + \frac{1 - \ta}{p_2}$ and $\mu = \ta \mu_1 + (1 - \ta) \mu_2$. Then,  we have
\begin{align*}
\| f_L \|_{L_\mu^p (\T_L^2)} \les \| f_L \|_{L_{\mu_1}^{p_1} (\T_L^2)}^\ta \| f_L \|_{L_{\mu_2}^{p_2} (\T_L^2)}^{1 - \ta} .
\end{align*}

\noi
Consequently, if we further let $1 \leq q, q_1, q_2 \leq \infty$ and $s, s_1, s_2 \in \R$ be such that $\frac{1}{q} = \frac{\ta}{q_1} + \frac{1 - \ta}{q_2}$ and $s = \ta s_1 + (1 - \ta) s_2$, we have
\begin{align*}
\| f_L \|_{\B_{p, q, \mu}^s (\T_L^2)} \les \| f_L \|_{\B_{p_1, q_1, \mu_1}^{s_1} (\T_L^2)}^\ta \| f_L \|_{\B_{p_2, q_2, \mu_2}^{s_2} (\T_L^2)}^{1 - \ta} .
\end{align*}
\end{lemma}

\begin{proof}
By using Lemma~\ref{LEM:Lwei}~(iv) and (iii), we have
\begin{align*}
\jbb{x}_{L, \mu}^p \les \jbb{x}_{L, \mu_1}^{\ta p} \jbb{x}_{L, \mu_2}^{(1 - \ta) p} ,
\end{align*}

\noi
so that by H\"older's inequality, we have 
\begin{align*}
\begin{split}
\| f_L \|_{L_\mu^p (\T_L^2)} 
&\les \bigg( \int_{\T_L^2} \jbb{x}_{L, \mu_1}^{\ta p} \jbb{x}_{L, \mu_2}^{(1 - \ta) p} |f_L (x)|^{\ta p} |f_L (x)|^{(1 - \ta) p} dx \bigg)^{\frac{1}{p}} \\
&\leq \bigg( \int_{\T_L^2} \jbb{x}_{L, \mu_1}^{p_1} |f_L (x)|^{p_1} dx \bigg)^{\frac{\ta}{p_1}} \bigg( \int_{\T_L^2} \jbb{x}_{L, \mu_2}^{p_2} |f_L (x)|^{p_2} dx \bigg)^{\frac{1 - \ta}{p_2}} \\
&= \| f_L \|_{L_{\mu_1}^{p_1} (\T_L^2)}^\ta \| f_L \|_{L_{\mu_2}^{p_2} (\T_L^2)}^{1 - \ta} ,
\end{split}
\end{align*}

\noi
which gives the first estimate.

To prove the second estimate for the weight Besov norms, we use the first estimate and H\"older's inequality to obtain
\begin{align*}
\| f_L \|_{\B_{p, q, \mu}^s (\T_L^2)} 
&\les \bigg( \sum_{\substack{N \geq 1 \\ \text{dyadic}}} N^{\ta s_1 q} \| \Dl_N^L f_L \|_{L_{\mu_1}^{p_1} (\T_L^2)}^{\ta q} N^{(1 - \ta) s_2 q} \| \Dl_N^L f_L \|_{L_{\mu_2}^{p_2} (\T_L^2)}^{(1 - \ta) q} \bigg)^{\frac{1}{q}} \\
&\leq \bigg( \sum_{\substack{N \geq 1 \\ \text{dyadic}}} N^{s_1 q_1} \| \Dl_N^L f_L \|_{L_{\mu_1}^{p_1} (\T_L^2)}^{q_1} \bigg)^{\frac{\ta}{q_1}} \bigg( \sum_{\substack{N \geq 1 \\ \text{dyadic}}} N^{s_2 q_2} \| \Dl_N^L f_L \|_{L_{\mu_2}^{p_2} (\T_L^2)}^{q_2} \bigg)^{\frac{1 - \ta}{q_2}} \\
&= \| f_L \|_{\B_{p_1, q_1, \mu_1}^{s_1} (\T_L^2)}^\ta \| f_L \|_{\B_{p_2, q_2, \mu_2}^{s_2} (\T_L^2)}^{1 - \ta},
\end{align*}

\noi
as desired.
\end{proof}

We now show the equivalence the weighted Besov norm $\B_{2, 2, \mu}^0 (\T_L^2)$ and the $L_\mu^2 (\T_L^2)$-norm.
Let us first show the following commutator estimate. 
\begin{lemma}
\label{LEM:Lcomm}
Let $L \geq 1$, $1 \leq p \leq \infty$, $\mu \in \R$, dyadic $N \geq 1$, and $f_L$ be an $L$-periodic function. Then, we have
\begin{align*}
\big\| [\Dl_N^L, \jbb{\cdot}_{L, \mu}] f_L \big\|_{L^p (\T_L^2)} \les N^{-1} \| f_L \|_{L_{\mu - 1}^p (\T_L^2)} ,
\end{align*}

\noi
where 
\begin{align*}
[\Dl_N^L, \jbb{\cdot}_{L, \mu}] = \Dl_N^L \jbb{\cdot}_{L, \mu} - \jbb{\cdot}_{L, \mu} \Dl_N^L .
\end{align*}
\end{lemma}

\begin{proof}
We have
\begin{align*}
[\Dl_N^L, \jbb{\cdot}_{L, \mu}] f_L (x) = \int_{\T_L^2} \big( \jbb{y}_{L, \mu} - \jbb{x}_{L, \mu} \big) \eta_N^L (x - y) f_L (y) dy ,
\end{align*}

\noi
where $\eta_N^L$ is defined in \eqref{defeta}.
By viewing $x, y \in \R^2$ and using the mean value theorem and Lemma~\ref{LEM:Lwei}~(vi) and (v), we have
\begin{align}
\begin{split}
\big| \jbb{y}_{L, \mu} - \jbb{x}_{L, \mu} \big| &\leq \big( \big| \nb \jbb{x}_{L, \mu} \big| + \big| \nb \jbb{y}_{L, \mu} \big| \big) \min_{k \in \Z^2} |x - y + Lk| \\
&\les \big( \jbb{x}_{L, \mu - 1} + \jbb{y}_{L, \mu - 1} \big) \min_{k \in \Z^2} |x - y + Lk| \\
&\les \jbb{y}_{L, \mu - 1} \big( \jbb{x - y}_{L, |\mu - 1|} + 1 \big) \min_{k \in \Z^2} |x - y + Lk| .
\end{split}
\label{yx_dist}
\end{align}

\noi
Thus, from \eqref{Poi_eta}, \eqref{yx_dist}, and Lemma~\ref{LEM:Lwei}~(ii), we have
\begin{align*}
&\sup_{x \in \T_L^2} \int_{\T_L^2} \big| \jbb{y}_{L, \mu} - \jbb{x}_{L, \mu} \big| |\eta_N^L (x - y)| \jbb{y}_{L, \mu - 1}^{-1} dy \\
&\les \sup_{x \in [-\frac{L}{2}, \frac{L}{2})^2} \int_{[-\frac{L}{2}, \frac{L}{2})^2}  \big( \jbb{x - y}_{L, |\mu - 1|} + 1 \big) \sum_{k \in \Z^2} |x - y + Lk| |\eta_N (x - y + Lk)| dy \\
&= \sup_{x \in [-\frac{L}{2}, \frac{L}{2})^2} \int_{\R^2} \jb{x - y}^{|\mu - 1|} |x - y| |\eta_N (x - y)| dy ,
\end{align*}

\noi
so that by using a change of variable on the last expression, we get
\begin{align}
\sup_{x \in \T_L^2} \int_{\T_L^2} \big| \jbb{y}_{L, \mu} - \jbb{x}_{L, \mu} \big| |\eta_N^L (x - y)| \jbb{y}_{L, \mu - 1}^{-1} dy \les N^{-1} .
\label{comm1}
\end{align}

\noi
Similarly, we have
\begin{align}
\sup_{y \in \T_L^2} \int_{\T_L^2} \big| \jbb{y}_{L, \mu} - \jbb{x}_{L, \mu} \big| |\eta_N^L (x - y)| \jbb{y}_{L, \mu - 1}^{-1} dx \les N^{-1} .
\label{comm2}
\end{align}

\noi
Thus, we conclude the desired estimate by Schur's test along with \eqref{comm1} and \eqref{comm2}.
\end{proof}

Using Lemma~\ref{LEM:Lcomm}, we prove the following norm equivalence.

\begin{lemma}
\label{LEM:L2equiv}
Let $L \geq 1$, $\mu \in \R$ with $|\mu| \leq 1$, and $f_L$ be an $L$-periodic function. Then, we have
\begin{align*}
\| f_L \|_{L^2_{\mu} (\T_L^2)} \sim \| f_L \|_{\B^0_{2, 2, \mu} (\T_L^2)}.
\end{align*}
\end{lemma}

\begin{proof}
By using \eqref{Hs_equi} and Lemma~\ref{LEM:Lcomm}, we have
\begin{align*}
\begin{split}
\| f_L \|_{\B^0_{2, 2, \mu}}^2 &= \sum_{\substack{N \geq 1 \\ \text{dyadic}}} \big\| \jbb{\cdot}_{L, \mu} \Dl_N^L f_L \big\|_{L^2 (\T_L^2)}^2 \\
&\les \sum_{\substack{N \geq 1 \\ \text{dyadic}}} \big\| \Dl_N^L \jbb{\cdot}_{L, \mu} f_L \big\|_{L^2 (\T_L^2)}^2 + \sum_{\substack{N \geq 1 \\ \text{dyadic}}} \big\| [\Dl_N^L, \jbb{\cdot}_{L, \mu}] f_L \big\|_{L^2 (\T_L^2)}^2 \\
&\les \| f_L \|_{L_\mu^2 (\T_L^2)}^2 + \| f_L \|_{L_{\mu - 1}^2 (\T_L^2)}^2 ,
\end{split}
\end{align*}

\noi
so that by using the embedding in Lemma~\ref{LEM:embL}~(iv) on the last expression, we get
\begin{align}
\| f_L \|_{\B^0_{2, 2, \mu}} \les \| f_L \|_{L_\mu^2 (\T_L^2)} .
\label{equi_help}
\end{align}

\noi
This gives one direction of the equivalence.

For the other direction, we first consider $0 \leq \mu \leq 1$. Using \eqref{Hs_equi}, Lemma~\ref{LEM:Lcomm}, and the fact that $\jbb{\cdot}_{L, \mu - 1} \les 1$ from Lemma~\ref{LEM:Lwei}~(iii) and (i), we have
\begin{align}
\begin{split}
\| f_L \|_{L^2_{\mu} (\T_L^2)}^2 
&\sim \sum_{\substack{N \geq 1 \\ \text{dyadic}}} \big\| \Dl_N^L \jbb{\cdot}_{L, \mu} f_L \big\|_{L^2 (\T_L^2)}^2 \\
&\les \sum_{\substack{N \geq 1 \\ \text{dyadic}}} \big\| \jbb{\cdot}_{L, \mu} \Dl_N^L f_L \big\|_{L^2 (\T_L^2)}^2 + \sum_{\substack{N \geq 1 \\ \text{dyadic}}} \big\| [\Dl_N^L, \jbb{\cdot}_{L, \mu}] f_L \big\|_{L^2 (\T_L^2)}^2 \\
&\les \| f_L \|_{\B^0_{2, 2, \mu}}^2 + \| f_L \|_{L^2 (\T_L^2)}^2 .
\end{split}
\label{equi1}
\end{align}

\noi
By \eqref{Hs_equi} and the embedding in Lemma~\ref{LEM:embL}~(iv), we have
\begin{align}
\| f_L \|_{L^2 (\T_L^2)} \sim \| f_L \|_{\B_{2, 2}^0 (\T_L^2)} \les \| f_L \|_{\B_{2, 2, \mu}^0 (\T_L^2)} .
\label{equi2}
\end{align}

\noi
Thus, from \eqref{equi1} and \eqref{equi2}, we get the desired
\begin{align*}
\| f_L \|_{L^2_{\mu} (\T_L^2)} \les \| f_L \|_{\B_{2, 2, \mu}^0 (\T_L^2)} .
\end{align*}

\noi
For $-1 \leq \mu < 0$, we use duality, Lemma~\ref{LEM:dualL}, \eqref{equi_help}, and Lemma~\ref{LEM:Lwei}~(iii) to obtain
\begin{align*}
\| f_L \|_{L_\mu^2 (\T_L^2)} &= \sup_{\| g_L \|_{L^2 (\T_L^2)} \leq 1} \bigg| \int_{\T_L^2} \jbb{\cdot}_{L, \mu} f_L \cj{g_L} dx \bigg| \\
&\les \sup_{\| g_L \|_{L^2 (\T_L^2)} \leq 1} \| f_L \|_{\B_{2, 2, \mu}^0 (\T_L^2)} \big\| \jbb{\cdot}_{L, \mu} g_L \big\|_{\B_{2, 2, - \mu}^0 (\T_L^2)} \\
&\sim \| f_L \|_{\B_{2, 2, \mu}^0 (\T_L^2)} .
\end{align*}

\noi
Thus, we get the desired equivalence.
\end{proof}

As a consequence of Lemma~\ref{LEM:L2equiv}, we know from the definitions in \eqref{H2mu} and the equivalence in Lemma~\ref{LEM:embL}~(v) that for any $\mu \in \R$,
\begin{align*}
\begin{split}
\| f_L \|_{H_\mu^1 (\T_L^2)} &\sim \| f_L \|_{L_\mu^2 (\T_L^2)} + \| \nb f_L \|_{\B_{2, 2, \mu}^0 (\T_L^2)} \sim \| f_L \|_{\B_{2, 2, \mu}^1 (\T_L^2)} \\
&\leq \| f_L \|_{\B_{2, 2, \mu}^2 (\T_L^2)} \sim \| f_L \|_{L_\mu^2 (\T_L^2)} + \| \Dl f_L \|_{\B_{2, 2, \mu}^0 (\T_L^2)} \sim \| f_L \|_{H_\mu^2 (\T_L^2)} ,
\end{split}
\end{align*}

\noi
so that
\begin{align}
\| f_L \|_{H_\mu^1 (\T_L^2)} \les \| f_L \|_{H_\mu^2 (\T_L^2)},
\label{H1H2}
\end{align}

\noi
which is not easily seen from the definitions in \eqref{H2mu} directly.


\subsection{Preliminary estimates on function spaces}

In this subsection, we show some useful estimates using the function spaces introduced in previous subsections.

Let us first show the following two lemmas on the periodization of a function with sufficient decay property.
\begin{lemma}
\label{LEM:per}
Let $L \geq 1$, $\mu_0 > 1$, and $f \in L^2_{\mu_0} (\R^2)$. Let $f_L$ be an $L$-periodic function defined by
\begin{align}
f_L \deff \sum_{k \in \Z^2} f (\cdot + L k).
\label{perd}
\end{align}

\noi
Then, for any $0 < \mu < \mu_0 - 1$, $f_L \in L_\mu^2 (\T_L^2)$ and, more precisely,
\begin{align}
\| f_L \|_{L_\mu^2 (\T_L^2)} \les \| f \|_{L_{\mu_0}^2 (\R^2)} .
\label{fLbdd}
\end{align}

\noi
Moreover, we have
\begin{align}
\| f - f_L \|_{L^2 ([- \frac{L}{2}, \frac{L}{2})^2)} \les L^{- \mu_0} \| f \|_{L_{\mu_0}^2 (\R^2)}.
\label{fLdiff}
\end{align}
\end{lemma}

\begin{proof}
We first show \eqref{fLbdd}. Note that by using Lemma~\ref{LEM:Lwei}~(ii), the Cauchy-Schwarz inequality in $k$, the fact that $\jb{x}_L = \jb{x + Lk}_L \leq \jb{x + Lk}$ for any $L \geq 1$ and $k \in \Z^2$, and the fact that $2 (\mu_0 - \mu) > 2$, we have
\begin{align}
\begin{split}
\| f_L \|_{L_\mu^2 ( \T_L^2 )}^2 &= \int_{[- \frac{L}{2}, \frac{L}{2} )^2} \jb{x}_L^{2 \mu} \bigg| \sum_{k \in \Z^2 } f (x + L k) \bigg|^2 dx \\
&\leq \int_{[- \frac{L}{2}, \frac{L}{2} )^2} \Big( \sum_{k \in \Z^2} \frac{1}{\jb{x + L k}^{2 (\mu_0 - \mu)}} \Big) \Big( \sum_{k \in \Z^2} \jb{x + L k}^{2 \mu_0} | f (x + L k) |^2 \Big) d x \\
&\les \int_{[- \frac{L}{2}, \frac{L}{2} )^2} \sum_{k \in \Z^2} \jb{x + L k}^{2 \mu_0} | f (x + L k) |^2 dx .
\end{split}
\label{per1}
\end{align}

\noi
Using a change of variable, we know that the last expression in \eqref{per1} is equal to
\begin{align*}
\sum_{k \in \Z}\int_{[- \frac{L}{2}, \frac{L}{2} )^2 - Lk} \jb{x}^{2 \mu_0} |f (x)|^2 dx = \| f \|_{L_{\mu_0}^2 (\R^2)}^2 ,
\end{align*}

\noi
which gives the desired estimate \eqref{fLbdd}.

We now show the difference estimate \eqref{fLdiff}. By the Cauchy-Schwarz inequality in $k$, we have
\begin{align}
\begin{split}
\| &f - f_L \|_{L^2 ([- \frac{L}{2}, \frac{L}{2} )^2)}^2 \\
&= \int_{[- \frac{L}{2}, \frac{L}{2} )^2} \bigg| \sum_{k \in \Z^2 \setminus \{0\}} f (x + L k) \bigg|^2 dx \\
&\leq \int_{[- \frac{L}{2}, \frac{L}{2} )^2} \Big( \sum_{k \in \Z^2 \setminus \{0\}} \frac{1}{\jb{x + L k}^{2 \mu_0}} \Big) \Big( \sum_{k \in \Z^2 \setminus \{0\}} \jb{x + L k}^{2 \mu_0} | f (x + L k) |^2 \Big) d x .
\end{split}
\label{per2}
\end{align}

\noi
Note that for any $x \in [- \frac{L}{2}, \frac{L}{2})^2$ and $k \in \Z^2 \setminus \{0\}$, we have 
\begin{align}
|x + Lk| \geq L|k| - |x| \geq L|k| - \frac{L}{\sqrt{2}} > \frac{L}{4} |k| .
\label{per_help}
\end{align}

\noi
From \eqref{per_help}, we get
\begin{align}
\sum_{k \in \Z^2 \setminus \{0\}} \frac{1}{\jb{x + L k}^{2 \mu_0}} \les L^{- 2 \mu_0} \sum_{k \in \Z^2 \setminus \{0\}} \frac{1}{|k|} \les L^{- 2 \mu_0} .
\label{per3}
\end{align}

\noi
Thus, from \eqref{per2}, \eqref{per3}, and a change of variable, we have
\begin{align*}
\| f - f_L \|_{L^2 ([- \frac{L}{2}, \frac{L}{2} )^2)}^2 
&\les L^{- 2 \mu_0} \int_{[- \frac{L}{2}, \frac{L}{2} )^2} \sum_{k \in \Z^2 \setminus \{0\}} \jb{x + L k}^{2 \mu_0} | f (x + L k) |^2 dx \\
&= L^{- 2 \mu_0} \sum_{k \in \Z^2 \setminus \{0\}} \int_{[- \frac{L}{2}, \frac{L}{2} )^2 + Lk} \jb{x}^{2 \mu_0} | f (x) |^2 dx \\
&\leq L^{- 2 \mu_0} \| f \|_{L^2_{\mu_0} (\R^2)}^2 .
\end{align*}

\noi
This finishes the proof of the lemma.
\end{proof}

\begin{remark} \rm
\label{RMK:mu0}
In Lemma~\ref{LEM:per}, the condition of the weight parameter $\mu_0$ is sharp up to the endpoint $\mu_0 = 1$. Indeed, suppose that $\mu_0 < 1$ and consider the function $h(x) = \jb{x}^{-2 + \dl}$ for $x \in \R^2$ and $0 < \dl < 1 - \mu_0$. Then, a direct computation shows that $h \in L^2_{\mu_0} (\R^2)$. However, for any $x \in \R^2$, it is not hard to see that
\begin{align*}
\sum_{k \in \Z^2} \jb{x + Lk}^{-2 + \dl} = \infty ,
\end{align*}

\noi
so that the periodization of the form \eqref{perd} for the function $h$ cannot be constructed.
\end{remark}

We then show that an $L$-periodic function can be put in an $L^2$-based localized function space, whose proof is inspired by \cite{BOZ}. We recall the localized function space in \eqref{defXK}.
\begin{lemma}
\label{LEM:Hsloc}
Let $0 \leq s < 2$ and $L \geq 1$. Let $f_L$ be an $L$-periodic function such that $f_L \in H^s (\T_L^2)$. Then, for any bounded and open set $U \subset \R^2$, we know that $f_L \in H^s (U)$ and we have the estimate
\begin{align*}
\| f_L \|_{H^s (U)} \les_U \| f_L \|_{H^s (\T_L^2)}.
\end{align*}
\end{lemma}

\begin{proof}
The case $s = 0$ is trivial, and so we only consider the case $0 < s < 2$. Let $\{ U_j \}_{j = 1}^M$ (with $M \in \N$) be a finite open cover of $\cj{U}$ such that each $U_j$ is a ball of radius $\frac 14$. Let $\{ \phi_j \}_{j = 1}^{M}$ be a smooth partition of unity on $\cj{U}$ subordinated to $\{ U_j \}_{j = 1}^{M}$. Then, we see that $f_L = \sum_{j = 1}^M \phi_j f_L$ on $U$, and so we have
\begin{align*}
\| f_L \|_{H^s (U)} \leq \Big\| \sum_{j = 1}^M \phi_j f_L \Big\|_{H^s (\R^2)} \leq \sum_{j = 1}^M \| \phi_j f_L \|_{H^s (\R^2)} .
\end{align*}

\noi
After using a translation in space, we see that it suffices to show
\begin{align}
\| \phi f_L \|_{H^s (\R^2)} \les \| f_L \|_{H^s (\T_L^2)} ,
\label{loc_goal}
\end{align}

\noi
where $\phi$ is a smooth function supported on $[- \frac{1}{4}, \frac{1}{4}]^2$. 

Let $\mathcal{G}_s$ be the kernel of $\jb{\nb_{\R^2}}^s$, which denotes the Fourier multiplier operator on $\R^2$ with symbol $\jb{\xi}^s$. From \cite[Lemma~2.1]{BOZ}, we know that $\mathcal{G}_s$ coincides with a smooth function on $\R^2 \setminus \{0\}$ and that
\begin{align}
|\mathcal{G}_s (x)| \les |x|^{-2 - s}
\label{Gs}
\end{align}

\noi
for any $|x| \ges 1$. Let $\chi_L \in C_c^\infty (\R^2)$ be supported on $\{ x \in \R^2 : |x| \leq \frac{L}{2} \}$ and $\chi_L \equiv 1$ on $\{ x \in \R^2 : |x| \leq \frac{3}{8} L \}$.
Then, we can write
\begin{align}
\begin{split}
\jb{\nb_{\R^2}}^s (\phi f_L) (x) &= (\mathcal{G}_s *_{\R^2} (\phi f_L)) (x) \\
&= (1 - \chi_L (x)) \int_{[-\frac{1}{4}, \frac{1}{4}]^2} \mathcal{G}_s (x - y) (\phi f_L) (y) dy + \chi_L (x) (\mathcal{G}_s *_{\R^2} (\phi f_L)) (x) \\
&\deff \textup{I}_1 + \textup{I}_2 ,
\end{split}
\label{loc0}
\end{align}

\noi
where the integral in $\textup{I}_1$ makes sense since 
\begin{align*}
|x - y| \geq |x| - |y| \geq \frac{3}{8} L - \frac{1}{2 \sqrt{2}} \ges 1
\end{align*} 

\noi
given $|x| \geq \frac 38 L$ and $y \in [-\frac 14, \frac 14]^2$.

For $\textup{I}_1$, by Young's convolution inequality and \eqref{Gs}, we have
\begin{align}
\| \textup{I}_1 \|_{L^2 (\R^2)} \leq \| \ind_{\{ |\cdot| \ges 1 \}} \mathcal{G}_s \|_{L^1 (\R^2)} \| \phi f_L \|_{L^2 (\R^2)} \les \| f_L \|_{L^2 (\T_L^2)} \leq \| f_L \|_{H^s (\T_L^2)} .
\label{loc1}
\end{align}

\noi
For $\textup{I}_2$, we let $\phi_L$ be the $L$-periodic extension of $\phi$ on $\R^2$ defined by
\begin{align*}
\phi_L = \sum_{k \in \Z^2} \phi (\cdot + Lk).
\end{align*}

\noi
By the Poisson summation formula \eqref{Poi}, we have for $x \in [- \frac{L}{2}, \frac{L}{2})^2$ that
\begin{align*}
\jb{\nb_{\T_L^2}}^s (\phi_L f_L) (x) 
&= \frac{1}{L} \sum_{n \in \Z^2} \jb{\tfrac{n}{L}}^s \mathcal{F}_{\T_L^2} (\phi_L f_L) ( \tfrac{n}{L} ) e^{2 \pi i \frac{n}{L} \cdot x} \\
&= \frac{1}{L^2} \sum_{n \in \Z^2} \jb{\tfrac{n}{L}}^s \mathcal{F}_{\R^2} (\phi f_L) ( \tfrac{n}{L} ) e^{2 \pi i \frac{n}{L} \cdot x} \\
&= \sum_{k \in \Z^2} \jb{\nb_{\R^2}}^s (\phi f_L) (x + L k).
\end{align*}

\noi
This gives
\begin{align}
\begin{split}
\textup{I}_2 (x) &= \chi_L (x) \jb{\nb_{\T_L^2}}^s (\phi_L f_L) (x) - \chi_L (x) \sum_{k \in \Z^2 \setminus \{0\}} \jb{\nb_{\R^2}}^s (\phi f_L) (x + Lk) \\
&= \chi_L (x) \jb{\nb_{\T_L^2}}^s (\phi_L f_L) (x) - \chi_L (x) \sum_{k \in \Z^2 \setminus \{0\}} \int_{[- \frac 14, \frac 14]^2} \mathcal{G}_s (x + Lk - y) (\phi f_L) (y) dy ,
\end{split}
\label{loc2-1}
\end{align}

\noi
where we used the fact that given $x \in [- \frac{L}{2}, \frac{L}{2})^2$, $y \in [-\frac 14, \frac 14]^2$, and $k \in \Z^2 \setminus \{0\}$, we have 
\begin{align*}
|x + Lk - y| \geq |x + Lk| - |y| \geq \frac{L}{2} - \frac{1}{2 \sqrt{2}} \ges 1 .
\end{align*}

\noi
By the product estimate in Lemma~\ref{LEM:prodL}~(i) along with \eqref{Hs_equi}, the embeddings in Lemma~\ref{LEM:embL}~(i) and (ii), and the equivalence in Lemma~\ref{LEM:embL}~(v), we have
\begin{align}
\begin{split}
\| \jb{\nb_{\T_L^2}}^s (\phi_L f_L) \|_{L^2 (\T_L^2)} 
&\les \| \phi_L \|_{L^{\infty} (\T_L^2)} \| f_L \|_{H^s (\T_L^2)} + \| \phi_L \|_{\B_{\infty, 2}^s (\T_L^2)} \| f_L \|_{L^2 (\T_L^2)} \\
&\les \| \phi_L \|_{\C^2 (\T_L^2)} \| f_L \|_{H^s (\T_L^2)} \\
&\les \big( \| \phi_L \|_{L^\infty (\T_L^2)} + \| \Dl \phi_L \|_{L^\infty (\T_L^2)} \big) \| f_L \|_{H^s (\T_L^2)} .
\end{split}
\label{loc2-2}
\end{align}

\noi
Also, by the Cauchy-Schwarz inequality in $y$, we have
\begin{align}
\begin{split}
\bigg\| &\int_{[-\frac{1}{4}, \frac{1}{4}]^2} \sum_{k \in \Z^2 \setminus \{0\}} \mathcal{G}_s (x + Lk - y) (\phi f_L) (y) dy \bigg\|_{L_x^2 ([- \frac{L}{2}, \frac{L}{2})^2)} \\
&\leq \bigg( \int_{[- \frac{L}{2}, \frac{L}{2})^2} \int_{[-\frac 14, \frac 14]^2} \Big| \sum_{k \in \Z^2 \setminus \{0\}} \mathcal{G}_s (x + Lk - y_1) \Big| dy_1 \\
&\quad \times \int_{[-\frac 14, \frac 14]^2} \Big| \sum_{k \in \Z^2 \setminus \{0\}} \mathcal{G}_s (x + Lk - y_2) \Big| |(\phi f_L) (y_2)|^2 dy_2 dx \bigg)^{\frac 12} ,
\end{split}
\label{loc2-3}
\end{align}

\noi
and so we can use H\"older's inequalities in $x$ and in $y_2$, change of variables, and \eqref{Gs} to bound the last expression in \eqref{loc2-3} by
\begin{align}
\begin{split}
&\sup_{x \in [- \frac{L}{2}, \frac{L}{2})^2} \bigg( \sum_{k \in \Z^2 \setminus \{0\}} \int_{[- \frac 14, \frac 14]^2 - Lk} |\mathcal{G}_s (x - y_1)| dy_1 \bigg)^{\frac 12} \\
&\qquad \times \sup_{y_2 \in [- \frac 14, \frac 14]^2} \bigg( \sum_{k \in \Z^2 \setminus \{0\}} \int_{[- \frac{L}{2}, \frac{L}{2})^2 + L k} |\mathcal{G}_s (x - y_2)| dx \bigg)^{\frac 12} \| \phi f_L \|_{L^2 (\T_L^2)} \\
&\quad \les \| \ind_{\{ |\cdot| \ges 1 \}} \mathcal{G}_s \|_{L^1 (\R^2)} \| f_L \|_{L^2 (\T_L^2)} \les \| f_L \|_{L^2 (\T_L^2)} .
\end{split}
\label{loc2-4}
\end{align} 

\noi
Combining \eqref{loc2-1}, \eqref{loc2-2}, \eqref{loc2-3}, and \eqref{loc2-4}, we obtain
\begin{align}
\| \textup{I}_2 \|_{L^2 (\R^2)} \les \| f_L \|_{H^s (\T_L^2)} .
\label{loc2}
\end{align}

The desired estimate \eqref{loc_goal} follows from \eqref{loc0}, \eqref{loc1}, and \eqref{loc2}, and so we have finished the proof of the lemma.
\end{proof}

Given $R \geq 1$, we define the function
\begin{align}
\s_R (x) \deff e^{\jb{\tfrac{x}{R}}} .
\label{sigR}
\end{align}

\noi
One can easily verify that
\begin{align}
|\nb \s_R^{-1} (x)| &\les R^{-1} \s_R^{-1} (x) .
\label{sigR1}
\end{align}

\noi
Also, for any $\mu > 0$, we have
\begin{align}
\jb{x}^{\mu} \leq R^\mu \jb{\tfrac{x}{R}}^\mu \les R^\mu \s_R (x).
\label{sigR2} 
\end{align}

We show the following useful estimate with the weight $\s_R^{-1}$.
\begin{lemma}
\label{LEM:sig}
Let $L \geq 1$, $1 \leq q < \infty$, and $f_L$ be an $L$-periodic function such that $f_L \in L^q (\T_L^2)$. Then, for any $R \geq 1$, we have
\begin{align*}
\| \s_R^{-1} f_L \|_{L^q (\R^2)} \les (1 + R L^{-1})^{\frac{2}{q}} \| f_L \|_{L^q (\T_L^2)}  .
\end{align*}
\end{lemma}

\begin{proof}
We write
\begin{align*}
\int_{\R^2} \s_R^{-q} |f_L|^q dx = \sum_{k \in \Z^2} \int_{[- \frac{L}{2}, \frac{L}{2})^2 + Lk} \s_R^{-q} |f_L|^q dx .
\end{align*}

\noi
Note that for any $x \in [- \frac{L}{2}, \frac{L}{2})^2 + Lk$ with $k \in \Z^2$, we have from \eqref{per_help} that $|x| \geq \frac{L |k|}{4}$. Thus, we have
\begin{align*}
\int_{\R^2} \s_R^{-q} |f_L|^q dx &\leq \Big( \sum_{k \in \Z^2} e^{- q \jb{\frac{L |k|}{4 R}}} \Big) \| f_L \|_{L^q (\T_L^2)}^q \\
&\les \bigg( 1 + \int_{\R^2} e^{- q \jb{\frac{L |x|}{4 R}}} dx \bigg) \| f_L \|_{L^q (\T_L^2)}^q \\
&\les (1 + R^2 L^{-2}) \| f_L \|_{L^q (\T_L^2)}^q ,
\end{align*}

\noi
where we used a change of variable in the last step. 
\end{proof}

We also have the following extension of Lemma~\ref{LEM:per}.
\begin{lemma}
\label{LEM:perw}
Let $L \geq 1$, $\mu_0 > 1$, and $f \in L_{\mu_0}^2 (\R^2)$. Let $f_L$ be an $L$-periodic function
\begin{align*}
f_L \deff \sum_{k \in \Z^2} f (\cdot + L k) .
\end{align*}

\noi
Then, for any $R \geq 1$ and $0 < \mu < \mu_0 - 1$, we have
\begin{align*}
\| \s_R^{-1} (f - f_L) \|_{L_\mu^2 (\R^2)} \les (L^{- \mu_0 + \mu} + R^{1 + 2 \mu} L^{-1 - \mu}) \| f \|_{L_{\mu_0}^2 (\R^2)} .
\end{align*}
\end{lemma}

\begin{proof}
We write
\begin{align}
\| \s_R^{-1} (f - f_L) \|_{L_\mu^2 (\R^2)} \leq \| \s_R^{-1} \jb{\cdot}^\mu (f - f_L) \|_{L^2 ([-\frac{L}{2}, \frac{L}{2})^2)} + \| \s_R^{-1} \jb{\cdot}^\mu (f - f_L) \|_{L^2 (\R^2 \setminus [-\frac{L}{2}, \frac{L}{2})^2)} .
\label{perw1}
\end{align}

\noi
For the first term in \eqref{perw1}, from Lemma~\ref{LEM:per}, we have
\begin{align}
\| \s_R^{-1} \jb{\cdot}^\mu (f - f_L) \|_{L^2 ([-\frac{L}{2}, \frac{L}{2})^2)} \les L^\mu \| f - f_L \|_{L^2 ([-\frac{L}{2}, \frac{L}{2})^2)} \les L^{- \mu_0 + \mu} \| f \|_{L_{\mu_0}^2 (\R^2)} .
\label{perw2}
\end{align}

\noi
For the second term in \eqref{perw1}, we write
\begin{align}
\begin{split}
\| &\s_R^{-1} \jb{\cdot}^\mu (f - f_L) \|_{L^2 (\R^2 \setminus [-\frac{L}{2}, \frac{L}{2})^2)}^2 \\
&= \sum_{k \in \Z^2 \setminus \{0\}} \int_{[-\frac{L}{2}, \frac{L}{2})^2 + L k} \s_R^{-2} (x) \jb{x}^{2 \mu} \Big| \sum_{k' \in \Z^2 \setminus \{0\}} f (x + Lk') \Big|^2 dx .
\end{split}
\label{perw3-1}
\end{align}

\noi
Note that from \eqref{per_help} and the fact that $e^{-y} \les y^{- a}$ for any $y > 0$ and $a > 0$, we have
\begin{align}
\begin{split}
\sum_{k \in \Z^2 \setminus \{0\}} \sup_{x \in [-\frac{L}{2}, \frac{L}{2})^2 + L k} \s_R^{-2} (x) \jb{x}^{2 \mu} 
&\les \sum_{k \in \Z^2 \setminus \{0\}} e^{- 2 \jb{\frac{L |k|}{4 R}}} L^{2 \mu} |k|^{2 \mu} \\
&\les \sum_{k \in \Z^2 \setminus \{0\}} \Big( \frac{R}{L |k|} \Big)^{2 + 4 \mu}  L^{2 \mu} |k|^{2 \mu} \\
&\les R^{2 + 4 \mu} L^{-2 - 2 \mu} .
\end{split}
\label{perw3-2}
\end{align}

\noi
Moreover, for any $k \in \Z^2 \setminus \{0\}$, we use the Cauchy-Schwarz inequality and a change of variable to obtain
\begin{align}
\begin{split}
&\int_{[-\frac{L}{2}, \frac{L}{2})^2 + L k} \Big| \sum_{k' \in \Z^2 \setminus \{0\}} f (x + Lk') \Big|^2 dx \\
&\leq \int_{[-\frac{L}{2}, \frac{L}{2})^2 + L k} \Big( \sum_{k' \in \Z^2} \frac{1}{\jb{x + L k'}^{2 \mu_0}} \Big) \Big( \sum_{k' \in \Z^2} \jb{x + L k'}^{2 \mu_0} |f (x + L k')|^2 \Big) dx \\
&\les \sum_{k' \in \Z^2} 
\int_{[-\frac{L}{2}, \frac{L}{2})^2 + L (k + k')} \jb{x}^{2 \mu_0} |f (x)|^2 dx = \| f \|_{L_{\mu_0}^2 (\R^2)}^2 .
\end{split}
\label{perw3-3}
\end{align}

\noi
Thus, from \eqref{perw3-1} and H\"older's inequality, we have
\begin{align*}
\begin{split}
\| \s_R^{-1} \jb{\cdot}^\mu (f - f_L) \|_{L^2 (\R^2 \setminus [-\frac{L}{2}, \frac{L}{2})^2)}^2 
&\leq \Big( \sum_{k \in \Z^2 \setminus \{0\}} \sup_{x \in [-\frac{L}{2}, \frac{L}{2})^2 + L k} \s_R^{-2} (x) \jb{x}^{2 \mu} \Big) \\
&\quad \times \sup_{k \in \Z^2 \setminus \{0\}} \int_{[-\frac{L}{2}, \frac{L}{2})^2 + L k} \Big| \sum_{k' \in \Z^2 \setminus \{0\}} f (x + Lk') \Big|^2 dx ,
\end{split}
\end{align*}

\noi
so that from \eqref{perw3-2} and \eqref{perw3-3}, we get
\begin{align}
\| \s_R^{-1} \jb{\cdot}^\mu (f - f_L) \|_{L^2 (\R^2 \setminus [-\frac{L}{2}, \frac{L}{2})^2)}^2 \les R^{2 + 4 \mu} L^{-2 - 2 \mu} \| f \|_{L_{\mu_0}^2 (\R^2)}^2 .
\label{perw3}
\end{align}

\noi
The desired estimate follows from \eqref{perw1}, \eqref{perw2}, and \eqref{perw3}.
\end{proof}

We now show some estimates on the approximating identity. Let $\rho \in C_c^\infty (\R^2)$ be such that $\supp \rho \in \{ x \in \R^2: |x| < \frac 14 \}$, $\rho \geq 0$, and $\int_{\R^2} \rho = 1$. Given $0 < \eps < 1$, we define $\rho_\eps (x) = \eps^{-2} \rho (\eps^{-1} x)$.

We first record the following lemma, which follows from similar steps in \cite[Lemma~2.6]{DLTV}.

\begin{lemma}
\label{LEM:rhoB}
For any $0 < \eps < \frac 12$ and $\mu > 0$, we have
\begin{align*}
\| \rho_\eps \|_{\B_{1, 1, \mu}^0 (\R^2)} \les |\log \eps|.
\end{align*}
\end{lemma}

Given $L \geq 1$, we denote by $\rho_{L, \eps}$ the $L$-periodic extension of $\rho_\eps$ on $\R^2$:
\begin{align}
\rho_{L, \eps} (x) \deff \sum_{k \in \Z^2} \rho_\eps (x + Lk) .
\label{def_rhoLe}
\end{align}

\begin{lemma}
\label{LEM:rhoDNL}
For any $L \geq 1$, $0 < \eps < \frac 12$, $\mu \geq 0$, and dyadic $N \geq 1$, we have
\begin{align}
\| \Dl_N^L \rho_{L, \eps} \|_{L_\mu^1 (\T_L^2)} \les \| \Dl_N \rho_\eps \|_{L_\mu^1 (\R^2)} .
\label{rhoL1}
\end{align} 

\noi
Consequently, we have
\begin{align}
\| \rho_{L, \eps} \|_{\B_{1, 1, \mu}^0 (\T_L^2)} \les |\log \eps| .
\label{rhoL2}
\end{align}
\end{lemma}

\begin{proof}
The estimate \eqref{rhoL2} follows directly from \eqref{rhoL1} and Lemma~\ref{LEM:rhoB}, so that we only need to prove \eqref{rhoL1}.

We recall $\eta_N$ and $\eta_N^L$ defined in \eqref{defeta}. By \eqref{convo2}, \eqref{Poi_eta}, \eqref{def_rhoLe}, a change of variable, and the fact that $\jbb{\cdot}_{L, \mu}$ is $L$-periodic, we have
\begin{align*}
\| \Dl_N^L \rho_{L, \eps} \|_{L_{\mu}^1 (\T_L^2)} 
&= \frac{1}{L} \int_{\T_L^2} \jbb{x}_{L, \mu} \bigg| \int_{\T_L^2} \eta_N^L (x - y) \rho_{L, \eps} (y) dy \bigg| dx \\
&= \int_{[ -\frac{L}{2}, \frac{L}{2} )^2} \jbb{x}_{L, \mu} \bigg| \int_{[ -\frac{L}{2}, \frac{L}{2} )^2} \sum_{k \in \Z^2} \eta_N (x + L k - y) \rho_\eps (y) dy \bigg| dx \\
&\leq \sum_{k \in \Z^2} \int_{[ -\frac{L}{2}, \frac{L}{2} )^2 + L k} \jbb{x}_{L, \mu} \bigg| \int_{[ -\frac{L}{2}, \frac{L}{2} )^2} \eta_N (x - y) \rho_\eps (y) dy \bigg| dx ,
\end{align*} 

\noi
so that by the support property of $\rho_\eps$, and Lemma~\ref{LEM:Lwei}~(ii), we get
\begin{align*}
\| \Dl_N^L \rho_{L, \eps} \|_{L_{\mu}^1 (\T_L^2)} 
\les \int_{\R^2} \jb{x}_L^\mu \bigg| \int_{\R^2} \eta_N (x - y) \rho_\eps (y) dy \bigg| dx 
\les \| \Dl_N \rho_\eps \|_{L_\mu^1 (\R^2)} .
\end{align*}

\noi
This gives the desired result.
\end{proof}

\begin{lemma}
\label{LEM:rho_eps}
Let $L \geq 1$ and let $f_L$ be an $L$-periodic distribution. Let $s_1, s_2, s \in \R$ be such that $0 \leq s_1 - s_2 < 1$ and $\mu \geq 0$. Then, for any $0 < \eps < 1$ and $\dl > 0$, we have $\rho_\eps *_{\R^2} f_L = \rho_{L, \eps} *_{\T_L^2} f_L$ and the estimates
\begin{align}
\| \rho_\eps *_{\R^2} f_L \|_{\C_{- \mu}^{s_1} (\T_L^2)} \les \eps^{s_2 - s_1} \| f_L \|_{\C_{- \mu}^{s_2} (\T_L^2)}
\label{rhoe1}
\end{align}

\noi
and
\begin{align}
\| \rho_\eps *_{\R^2} f_L - f_L \|_{\C_{- \mu}^s (\T_L^2)} \les \eps^\dl \| f_L \|_{\C_{- \mu}^{s + \dl} (\T_L^2)} .
\label{rhoe2}
\end{align}
\end{lemma}

\begin{proof}
The fact that $\rho_\eps *_{\R^2} f_L = \rho_{L, \eps} *_{\T_L^2} f_L$ follows from Lemma~\ref{LEM:GL}.

For the estimates, we first consider \eqref{rhoe1}. Given a dyadic number $N \geq 1$, we use Young's convolution inequality along with Lemma~\ref{LEM:Lwei}~(v) to obtain
\begin{align*}
N^{s_1} \| \Dl_N^L (\rho_{L, \eps} *_{\T_L^2} f_L) \|_{L_{- \mu}^\infty (\T_L^2)} &= N^{s_1} \big\| \big( \Dl_{\frac{N}{2}}^L + \Dl_N^L + \Dl_{2 N}^L \big) \rho_{L, \eps} *_{\T_L^2} \Dl_N^L f_L \big\|_{L_{- \mu}^\infty (\T_L^2)} \\
&\les N^{s_2} \| \Dl_N^L f_L \|_{L_{- \mu}^\infty (\T_L^2)} \sum_{M = \frac{N}{2}, N, 2N} M^{s_1 - s_2} \| \Dl_M^L \rho_{L, \eps} \|_{L_{\mu}^1 (\T_L^2)}.
\end{align*}

\noi
Thus, it suffices to show that for any dyadic $N \geq 1$, we have
\begin{align*}
(\eps N)^{s_1 - s_2} \| \Dl_N^L \rho_{L, \eps} \|_{L_\mu^1 (\T_L^2)} \les 1
\end{align*}

\noi
with the underlying constant independent of $L$ and $N$. From \eqref{rhoL1} in Lemma~\ref{LEM:rhoDNL}, we only need to show
\begin{align*}
(\eps N)^{s_1 - s_2} \| \Dl_N \rho_\eps \|_{L_\mu^1 (\R^2)} \les 1 ,
\end{align*}

\noi
which has already been established in \cite[Lemma~2.7]{DLTV}. This implies
\begin{align}
N^{s_1} \| \Dl_N^L (\rho_{L, \eps} *_{\T_L^2} f_L) \|_{L_{- \mu}^\infty (\T_L^2)} \les \eps^{s_2 - s_1} N^{s_2} \| \Dl_N^L f_L \|_{L^\infty_{- \mu} (\T_L^2)},
\label{rhoe_help}
\end{align}

\noi
which implies \eqref{rhoe1}.

We now consider \eqref{rhoe2}. Let $N \geq 1$ be a dyadic number. If $\eps N \geq \frac 12$, we use \eqref{rhoe_help} with $s_1 = s_2 = s$ to obtain
\begin{align}
\begin{split}
N^s &\| \Dl_N^L (\rho_{L, \eps} *_{\T_L^2} f_L) - \Dl_N^L f_L \|_{L_{- \mu}^\infty (\T_L^2)} \\
&\leq N^s \| \rho_{L, \eps} *_{\T_L^2} \Dl_N^L f_L \|_{L_{- \mu}^\infty (\T_L^2)} + N^s \| \Dl_N^L f_L \|_{L_{- \mu}^\infty (\T_L^2)} \\
&\les N^s \| \Dl_N^L f_L \|_{L_{- \mu}^\infty (\T_L^2)} \les \eps^\dl N^{s + \dl} \| \Dl_N^L f_L \|_{L_{- \mu}^\infty (\T_L^2)} .
\end{split}
\label{eN1}
\end{align}

\noi
If $\eps N < \frac 12$, we note that
\begin{align}
\Dl_N^L (\rho_{L, \eps} *_{\T_L^2} f_L) - \Dl_N^L f_L = \sum_{M = \frac{N}{2}, N, 2N} \big( \Dl_M^L (\rho_{L, \eps} *_{\T_L^2} \Dl_N^L f_L) - \Dl_M^L \Dl_N^L f_L \big) .
\label{eN2-0}
\end{align}

\noi
We use $\eta_N$ and $\eta_N^L$ as defined in \eqref{defeta}. By \eqref{convo2}, the fact that $\int_{\T_L^2} \rho_{L, \eps} = 1$, Lemma~\ref{LEM:Lwei}~(iii), H\"older's inequality, and the support property of $\rho_\eps$, we have for any $x \in [- \frac{L}{2}, \frac{L}{2})^2$ that
\begin{align*}
&\big| \Dl_M^L (\rho_{L, \eps} *_{\T_L^2} \Dl_N^L f_L) (x) - \Dl_M^L \Dl_N^L f_L (x) \big| \\
&\les \frac{1}{L} \bigg| \int_{\T_L^2} \int_{\T_L^2} \big( \eta_M^L (x - z) - \eta_M^L (x - y) \big)  \rho_{L, \eps} (z - y) \Dl_N^L f_L (y) dy dz \bigg| \\
&\les \frac{1}{L} \sum_{\substack{k \in \Z^2 \\ |k| \leq 2}} \int_{[- \frac{L}{2}, \frac{L}{2} )^2} \int_{[- \frac{L}{2}, \frac{L}{2} )^2} \big| \eta_M^L (x - z) - \eta_M^L (x - y) \big|  \rho_{\eps} (z - y + L k) \jbb{y}_{L, \mu} dz dy \| \Dl_N^L f_L \|_{L_{- \mu}^\infty (\T_L^2)} ,
\end{align*}

\noi
so that by \eqref{Poi_eta}, we get
\begin{align}
\begin{split}
\big| &\Dl_M^L (\rho_{L, \eps} *_{\T_L^2} \Dl_N^L f_L) (x) - \Dl_M^L \Dl_N^L f_L (x) \big| \\
&\leq \sum_{\substack{k \in \Z^2 \\ |k| \leq 2}} \sum_{k' \in \Z^2} \int_{[- \frac{L}{2}, \frac{L}{2} )^2} \int_{[- \frac{L}{2}, \frac{L}{2} )^2} \big| \eta_M (x - z + L k') - \eta_M (x - y + L k + L k') \big| \\
&\qquad \qquad \qquad \qquad \qquad \qquad \times \rho_{\eps} (z - y + L k) \jbb{y}_{L, \mu} dz dy \| \Dl_N^L f_L \|_{L_{- \mu}^\infty (\T_L^2)} .
\end{split}
\label{eN2-1}
\end{align}

\noi
For any $k \in \Z^2$ with $|k| \leq 2$, we use change of variables, the fact that $\jbb{\cdot}_{L, \mu}$ is $L$-periodic, and Lemma~\ref{LEM:Lwei}~(ii) to obtain
\begin{align*}
&\sum_{k' \in \Z^2} \int_{[- \frac{L}{2}, \frac{L}{2} )^2} \int_{[- \frac{L}{2}, \frac{L}{2} )^2} \big| \eta_M (x - z + L k') - \eta_M (x - y + L k + L k') \big| \rho_\eps (z - y + Lk) \jbb{y}_{L, \mu} dz dy \\
&= \sum_{k' \in \Z^2} \int_{[- \frac{L}{2}, \frac{L}{2} )^2 - L k - L k'} \int_{[- \frac{L}{2}, \frac{L}{2} )^2} \big| \eta_M (x - z + L k') - \eta_M (x - y) \big| \rho_\eps (z - y - L k') \jbb{y}_{L, \mu} dz dy \\
&\les \int_{\R^2} \sum_{k' \in \Z^2} \int_{[ - \frac{L}{2}, \frac{L}{2} )^2 - L k'} \big| \eta_M (x - z) - \eta_M (x - y) \big| \rho_\eps (z - y) \jb{y}_{L}^\mu dz dy , 
\end{align*}

\noi
so that by \eqref{xyL}, we get
\begin{align}
\begin{split}
&\sum_{k' \in \Z^2} \int_{[- \frac{L}{2}, \frac{L}{2} )^2} \int_{[- \frac{L}{2}, \frac{L}{2} )^2}  \big| \eta_M (x - z + L k') - \eta_M (x - y + L k + L k') \big| \\
&\quad \times \rho_\eps (z - y + Lk) \jbb{y}_{L, \mu} dz dy \\
&\les \int_{\R^2} \int_{\R^2} \big| \eta_M (x - z) - \eta_M (x - y) \big| \rho_\eps (z - y) \jb{z}_L^\mu dz dy
\end{split}
\label{eN2-2}
\end{align}

\noi
Thus, from \eqref{eN2-0}, \eqref{eN2-1}, \eqref{eN2-2}, and the fact that $\jbb{x}_{L, - \mu} \les \jbb{x}_{L, \mu}^{-1}$ from Lemma~\ref{LEM:Lwei}~(iii), we get
\begin{align}
\begin{split}
\| &\Dl_N^L (\rho_{L, \eps} *_{\T_L^2} f_L) - \Dl_N^L f_L \|_{L_{- \mu}^\infty (\T_L^2)} \\ 
&\les \sup_{x \in [-\frac{L}{2}, \frac{L}{2})^2} \sum_{M = \frac{N}{2}, N, 2N} \jbb{x}_{L, \mu}^{-1} \\
&\qquad \times \int_{\R^2} \int_{\R^2} \big| \eta_M (x - z) - \eta_M (x - y) \big| \rho_\eps (z - y) \jb{z}_L^\mu dz dy \| \Dl_N^L f_L \|_{L_{- \mu}^\infty (\T_L^2)}.
\end{split}
\label{rhoe_est1}
\end{align}

\noi
Using change of variables and Lemma~\ref{LEM:Lwei}~(ii), we have
\begin{align}
\begin{split}
&\jbb{x}_{L, \mu}^{-1} \int_{\R^2} \int_{\R^2} \big| \eta_M (x - z) - \eta_M (x - y) \big| \rho_\eps (z - y) \jb{z}_L^\mu dz dy \\
&= \eps^{-2} M^2 \jbb{x}_{L, \mu}^{-1} \int_{\R^2} \int_{\R^2} \big| \eta_1 (M x - M z) - \eta_1 (M x - M y) \big| \rho (\eps^{-1} z - \eps^{-1} y) \jbb{z}_{L, \mu} dz dy \\
&= \eps^2 M^2 \jbb{x}_{L, \mu}^{-1} \int_{\R^2} \int_{\R^2} \big| \eta_1 (M x - \eps M z) - \eta_1 (M x - \eps M y) \big| \rho (z - y) \jbb{\eps z}_{L, \mu} dz dy .
\end{split}
\label{rhoe1-1}
\end{align}

\noi
From the support property of $\rho$, we know that $|z - y| \leq \frac 14$, so that from the mean value theorem and the fact that $\eps M < 1$, we have
\begin{align}
\big| \eta_1 (M x - \eps M z) - \eta_1 (M x - \eps M y) \big| \leq \eps M \sup_{a \in B (M x - \eps M z, 1)} |\nb \eta_1 (a)| .
\label{rhoe1-2}
\end{align}

\noi
Thus, from \eqref{rhoe1-1}, \eqref{rhoe1-2}, and a change of variable, we have
\begin{align*}
\begin{split}
&\jbb{x}_{L, \mu}^{-1} \int_{\R^2} \int_{\R^2} \big| \eta_M (x - z) - \eta_M (x - y) \big| \rho_\eps (z - y) \jb{z}_L^\mu dz dy \\
&\leq \eps^3 M^3 \jbb{x}_{L, \mu}^{-1} \int_{\R^2} \sup_{a \in B(M x - \eps M z, 1)} |\nb \eta_1 (a)| \jbb{\eps z}_{L, \mu} dz \\
&= \eps M \jbb{x}_{L, \mu}^{-1} \int_{\R^2} \sup_{a \in B(z, 1)} |\nb \eta_1 (a)| \jbb{x - M^{-1} z}_{L, \mu} dz ,
\end{split}
\end{align*}

\noi
so that by using the fact that $\eta_1$ is a Schwartz function and
\begin{align*}
\jbb{x}_{L, \mu}^{-1} \jbb{x - M^{-1} z}_{L, \mu} \les \jbb{M^{-1} z}_{L, \mu} \les \jb{M^{-1} z}_L^\mu \leq \jb{M^{-1} z}^\mu \leq \jb{z}^\mu
\end{align*}

\noi
from Lemma~\ref{LEM:Lwei}~(v) and (ii), we get
\begin{align}
\jbb{x}_{L, \mu}^{-1} \int_{\R^2} \int_{\R^2} \big| \eta_M (x - z) - \eta_M (x - y) \big| \rho_\eps (z - y) \jb{z}_L^\mu dz dy \les \eps M .
\label{rhoe_est2}
\end{align}

\noi
Thus, from \eqref{rhoe_est1} and \eqref{rhoe_est2} along with the fact that $\eps N < \frac 12$, we obtain
\begin{align}
\| \Dl_N^L (\rho_{L, \eps} *_{\T_L^2} f_L) - \Dl_N^L f_L \|_{L_{- \mu}^\infty (\T_L^2)} \les \eps N \| \Dl_N^L f_L \|_{L_{-\mu}^\infty (\T_L^2)} < \eps^\dl N^\dl \| \Dl_N^L f_L \|_{L_{-\mu}^\infty (\T_L^2)} .
\label{eN2}
\end{align}

\noi
Combining \eqref{eN1} and \eqref{eN2}, we obtain the desired estimate.
\end{proof}

Lastly, we show the Strichartz estimates uniform in the size of the torus. Given $L \geq 1$ and $N \in \N$, let $P_N^L$ be the frequency projection onto frequencies $\{ n \in \Z^2 : \frac{N}{2} < |\frac{n}{L}| \leq N \}$ if $N \geq 2$ and $P_1^L$ be the projection onto $\{ n \in \Z^2 : |\frac{n}{L}| \leq 1 \}$.

\begin{lemma}
\label{LEM:L4}
Let $L \geq 1$ and let $f_L$ be an $L$-periodic function. Let $2 < q, r < \infty$ be such that $\frac{2}{q} + \frac{2}{r} = 1$. Let $N \in \N$ with $N \geq 10$ and let $I \subset \R_+$ be a closed interval with length $N^{-1}$. Then, we have
\begin{align*}
\| e^{- i t \Dl} P_N^L f_L \|_{L^q_I L^r (\T_L^2)} \les \| P_N^L f_L \|_{L^2 (\T_L^2)}.
\end{align*} 
\end{lemma}

\begin{proof}
We define $g (x) = f_L (Lx)$, so that $g$ is a periodic function on the torus $\T^2$ of size 1. Then, for any $t \in I$ and $x \in \T_L^2$, we have
\begin{align*}
e^{- i t \Dl} P_N^L f_L (x) &= \frac{1}{L} \sum_{\substack{n \in \Z^2 \\ \frac{N}{2} < |\frac{n}{L}| \leq N}} \ft{f_L} (\tfrac{n}{L}) e^{4 \pi^2 i t |\frac{n}{L}|^2} e^{2 \pi i \frac{n}{L} \cdot x} \\
&= \sum_{\substack{n \in \Z^2 \\ \frac{N}{2} < |\frac{n}{L}| \leq N}} \ft g (n) e^{4 \pi^2 i \frac{t}{L^2} |n|^2} e^{2 \pi i n \cdot \frac{x}{L}} \\
&= e^{- i \frac{t}{L^2} \Dl} Q_{LN} g (\tfrac{x}{L}),
\end{align*}

\noi
where $Q_{LN}$ denotes the frequency projection onto $\{ n \in \Z^2 : \frac{LN}{2} < |n| \leq LN \}$. Given $I = [a,b]$ for some $0 \leq a < b$, we denote $I_L = [L^{-2} a, L^{-2} b]$. Thus, by a change of variable, the Strichartz estimate in \cite{BGTz}, and a change of variable again, we obtain
\begin{align*}
\| e^{- i t \Dl} P_N^L f_L \|_{L^q_I L^r (\T_L^2)} &= L \| e^{- i t \Dl} Q_{LN} g \|_{L^q_{I_L} \! \! L^r (\T^2)} \les L \| g \|_{L^2 (\T^2)} = \| f_L \|_{L^2 (\T_L^2)} ,
\end{align*}

\noi
as desired.
\end{proof}

\section{Regularity and convergence of stochastic objects}
\label{SEC:sto}

In this section, we study the regularity and convergence properties of the stochastic objects introduced in Subsection~\ref{SUB:main}.

\subsection{$L$-periodic stochastic objects}
\label{SUB:sto}

In this subsection, we study the regularity properties of the $L$-periodic stochastic objects.

We recall that $\xi_L$ is the $L$-periodic spatial white noise defined in \eqref{defxiL}. Note that by applying $\xi_L$ to the orthonormal basis $\{ \frac{1}{L} e^{2 \pi i \frac{n}{L} \cdot x} \}_{n \in \Z^2}$ of $L^2 (\T_L^2)$ via \eqref{defxiL2}, we obtain the following Fourier series for $\xi_L$:
\begin{align}
\xi_L (x, \om) = \frac{1}{L} \sum_{n \in \Z^2} g_n^L (\om) e^{2 \pi i \frac{n}{L} \cdot x},
\label{xiL_fs}
\end{align}

\noi
where 
\begin{align*}
g_n^L \deff \frac{1}{L} ( \xi_L , e^{2 \pi i \frac{n}{L} \cdot x} ).
\end{align*}

\begin{lemma}
\label{LEM:gnL}
The set $\{ g_n^L \}_{n \in \Z^2}$ is a family of independent standard complex-valued centered Gaussian random variables conditioned such that $g_{-n}^L = \cj{g_n^L}$ for any $n \in \Z^2$ (in particular, $g_0^L$ is real-valued).
\end{lemma}

\begin{proof}
From the definition, it is easy to see that each $g_n^L$ is a complex-valued centered Gaussian random variable. Also, we have the isometry
\begin{align*}
\E [ |g_n^L|^2 ] = \frac{1}{L^2} \E [ | (\xi_L, e^{2 \pi i \frac{n}{L} \cdot x}) |^2 ] = \frac{1}{L^2} \int_{\T_L^2} |e^{2 \pi i \frac{n}{L} \cdot x}|^2 dx = 1.
\end{align*}

\noi
Moreover, for $n, m \in \Z^2$ with $n \neq \pm m$, $g_n^L$ is independent of $g_m^L$, which follows from the $L^2$-isometry, the polarization identity, and the fact that
\begin{align*}
\int_{\T_L^2} \phi_1 (2 \pi \tfrac{n}{L} \cdot x) \phi_2 (2 \pi \tfrac{m}{L} \cdot x) dx = 0 
\end{align*}

\noi
for $\phi_1, \phi_2 \in \{ \cos, \sin \}$. Lastly, the property that $g_{-n}^L = \cj{g_n^L}$ for any $n \in \Z^2$ follows from the fact that $(\xi_L, f)$ is real-valued if $f$ is real-valued.
\end{proof}

Let $\rho \in C_c^\infty (\R^2)$ be such that $\supp \rho \in \{ x \in \R^2: |x| < \frac 14 \}$, $\rho \geq 0$, and $\int_{\R^2} \rho = 1$. Given $0 < \eps < 1$, we define $\rho_\eps (x) = \eps^{-2} \rho (\eps^{-1} x)$. We also define the $L$-periodic extension $\rho_{L, \eps}$ as in \eqref{def_rhoLe}. We then define the regularized version of $\xi_L$ as
\begin{align}
\xi_{L, \eps} (x, \om) \deff (\rho_\eps *_{\R^2} \xi_L) (x, \om) = (\rho_{L, \eps} *_{\T_L^2} \xi_L) (x, \om) = \frac{1}{L} \sum_{n \in \Z^2} \mathcal{F}_{\R^2} ( \rho ) (\tfrac{\eps n}{L}) g_n^L (\om) e^{2 \pi i \frac{n}{L} \cdot x} ,
\label{xiLe_fs}
\end{align}

\noi
where in the second equality we used the identity in Lemma~\ref{LEM:rho_eps} and in the third equality we used \eqref{convo}, \eqref{xiL_fs}, and the fact that $\mathcal{F}_{\T_L^2} (\rho_{L, \eps}) = L^{-1} \mathcal{F}_{\R^2} (\rho_\eps)$. For simplicity below, we denote $\ft \rho = \mathcal{F}_{\R^2} (\rho)$. From \eqref{defY}, \eqref{eqY}, \eqref{xiL_fs}, Lemma~\ref{LEM:GL}, and similar steps in \eqref{xiLe_fs}, we have
\begin{align*}
(1 - \Dl) Y_L (x, \om) &= - \xi_L (x, \om) + (G *_{\R^2} \xi_L) (x, \om) - (\varphi *_{\R^2} \xi_L) (x, \om) \\
&= - \frac 1L \sum_{n \in \Z^2} (1 - \ft G (\tfrac{n}{L}) + \ft \varphi (\tfrac{n}{L})) g_n^L (\om) e^{2 \pi i \frac{n}{L} \cdot x} 
\end{align*}

\noi
with $\ft G = \mathcal{F}_{\R^2} (G)$ and $\ft \varphi = \mathcal{F}_{\R^2} (\varphi)$,  which gives
\begin{align}
Y_L (x, \om) = - \frac{1}{L} \sum_{n \in \Z^2} \frac{1 - \ft G (\frac{n}{L}) + \ft \varphi (\frac{n}{L})}{\jb{\frac{n}{L}}^2} g_n^L (\om) e^{2 \pi i \frac{n}{L} \cdot x}.
\label{YL_fs}
\end{align}

\noi
Since $G \in L^1 (\R^2)$ and $\varphi \in C^\infty_c (\R^2)$ with support in $\{ x \in \R^2 : |x| < \frac 14 \}$, we will simply bound $|\ft G (\xi)| \les 1$ and $|\ft \varphi (\xi)| \les 1$ in all the analysis below. 
We also define the regularized version $Y_L$ as
\begin{align}
Y_{L, \eps} (x, \om) \deff (\rho_\eps *_{\R^2} Y_L) (x, \om) = - \frac{1}{L} \sum_{n \in \Z^2} \frac{1 - \ft G (\frac{n}{L}) + \ft \varphi (\frac{n}{L})}{\jb{\frac{n}{L}}^2} \ft \rho ( \tfrac{\eps n}{L} ) g_n^L (\om) e^{2 \pi i \frac{n}{L} \cdot x}.
\label{YLe_fs}
\end{align}

\noi
We then define the Wick ordering $\wick{|\nb Y_{L, \eps}|^2}$ as
\begin{align}
\wick{|\nb Y_{L, \eps}|^2} \, \deff |\nb Y_{L, \eps}|^2 - C_{L, \eps},
\label{def_YLe2}
\end{align}

\noi
where
\begin{align*}
C_{L, \eps} = \E [ |\nb Y_{L, \eps}|^2 ] .
\end{align*}

\noi
Given $0 < \eps < \frac 12$, we have from \eqref{YLe_fs}, Lemma~\ref{LEM:gnL}, and the fact that $\ft \rho$ is a Schwartz function that
\begin{align}
\begin{split}
C_{L, \eps} &= \frac{1}{L^2} \sum_{n \in \Z^2} \jb{\tfrac{n}{L}}^{-2} | \ft \rho (\tfrac{\eps n}{L}) |^2 
\les \frac{1}{L^2} \sum_{n \in \Z^2} \jb{\tfrac{n}{L}}^{-2} \jb{\tfrac{\eps n}{L}}^{-1} \\
&\les \int_{\R^2} \jb{\xi}^{-2} \jb{\eps \xi}^{-1} d\xi 
\leq \int_{\{|\xi| \leq \eps^{-2}\}} \jb{\xi}^{-2} d \xi + \int_{\{ |\xi| > \eps^{-2} \}} \jb{\xi}^{- \frac 52} d \xi \\
&\les |\log \eps|.
\end{split}
\label{CLe}
\end{align}

\noi
We will show in Lemma~\ref{LEM:YregLe} below that $\wick{|\nb Y_{L, \eps}|^2}$ converges as $\eps \to 0$ to a limiting object $\wick{|\nb Y_L|^2}$ defined as
\begin{align}
\begin{split}
&\wick{ |\nb Y_L|^2 } (x, \om) \\
&\deff \frac{4 \pi^2}{L^2} \sum_{\substack{n, n_1, n_2 \in \Z^2 \\ n_1 - n_2 = n}} \frac{\big(1 - \ft G (\frac{n_1}{L}) + \ft \varphi (\frac{n_1}{L})\big) \big(1 - \cj{\ft G} (\frac{n_2}{L}) + \cj{\ft \varphi} (\frac{n_2}{L})\big) \frac{n_1}{L} \cdot \frac{n_2}{L}}{\jb{\frac{n_1}{L}}^2 \jb{\frac{n_2}{L}}^2} g_{n_1}^L (\om) \cj{g_{n_2}^L} (\om) e^{2 \pi i \frac{n}{L} \cdot x} \\
&\quad + \frac{4 \pi^2}{L^2} \sum_{n \in \Z^2 \setminus \{0\}} \frac{\big| 1 - \ft G (\frac{n}{L}) + \ft \varphi (\frac{n}{L}) \big|^2}{\jb{\frac{n}{L}}^2} (| g_n^L (\om) |^2 - 1).
\end{split}
\label{YL2_fs}
\end{align}

Let us now show the following lemma on the regularity of the aforementioned $L$-periodic stochastic objects and also the convergence of the regularized $L$-periodic stochastic objects.

\begin{lemma}
\label{LEM:YregLe}
Let $L \geq 1$, $0 < s < 1$, $0 < \mu \leq 1$, $0 < \dl < 1 - s$, $\varphi \in C_c^\infty (\R^2)$, and $a \in \R$. Then, there exists $\Om' \subset \Om$ with full probability measure such that for any $\om \in \Om'$, the following estimates hold.

\smallskip \noi
\textup{(i)} We have 
\begin{align*}
\| Y_L \|_{\C^s_{- \mu} (\T_L^2)} + \| \nb Y_L \|_{\C^{s - 1}_{- \mu} (\T_L^2)} + \sup_{\eps \in (0, 1)} \| Y_{L, \eps} \|_{\C^{s}_{- \mu} (\T_L^2)} + \sup_{\eps \in (0, 1)} \| \nb Y_{L, \eps} \|_{\C^{s - 1}_{- \mu} (\T_L^2)} \les_\om 1
\end{align*}

\noi
and, for any $0 < \eps < 1$,
\begin{align*}
\| Y_{L, \eps} - Y_{L} \|_{\C^{s}_{- \mu} (\T_L^2)} + \| \nb Y_{L, \eps} - \nb Y_{L} \|_{\C^{s - 1}_{- \mu} (\T_L^2)} \les_{\om} \eps^\dl .
\end{align*}

\smallskip \noi
\textup{(ii)} We have
\begin{align*}
\| \varphi *_{\R^2} \xi_L \|_{\C^s_{- \mu} (\T_L^2)} + \sup_{\eps \in (0, 1)} \| \varphi *_{\R^2} \xi_{L, \eps} \|_{\C^s_{- \mu} (\T_L^2)} \les_{\om} 1
\end{align*}

\noi
and, for any $0 < \eps < 1$,
\begin{align*}
\| \varphi *_{\R^2} \xi_{L, \eps} - \varphi *_{\R^2} \xi_L \|_{\C^s_{- \mu} (\T_L^2)} \les_{\om} \eps^\dl .
\end{align*}

\smallskip \noi
\textup{(iii)} We have 
\begin{align*}
\| e^{a Y_L} \|_{L^\infty_{- \mu} (\T_L^2)} + \sup_{\eps \in (0, 1)} \| e^{a Y_{L, \eps}} \|_{L^\infty_{- \mu} (\T_L^2)} + \| e^{a Y_L} \|_{\C^s_{- \mu} (\T_L^2)} + \sup_{\eps \in (0, 1)} \| e^{a Y_{L, \eps}} \|_{\C^s_{- \mu} (\T_L^2)} \les_{\om} 1
\end{align*}

\noi
and, for any $0 < \eps < 1$,
\begin{align*}
\| e^{a Y_{L, \eps}} - e^{a Y_L} \|_{L^\infty_{- \mu} (\T_L^2)} \les_{\om} \eps^\dl .
\end{align*}

\smallskip \noi
\textup{(iv)} We have
\begin{align*}
\| \wick{|\nb Y_L|^2} \|_{\C^{s - 1}_{- \mu} (\T_L^2)} + \sup_{\eps \in (0, 1)} \| \wick{|\nb Y_{L, \eps}|^2} \|_{\C^{s - 1}_{- \mu} (\T_L^2)} \les_{\om} 1
\end{align*}

\noi
and, for any $0 < \eps < 1$,
\begin{align*}
\big\| \wick{ |\nb Y_{L, \eps}|^2 } - \wick{ |\nb Y_{L}|^2 } \big\|_{\C^{s - 1}_{- \mu} (\T_L^2)} \les_{\om} \eps^\dl .
\end{align*}
\end{lemma}

\begin{proof}
(i) We first show the bound for $Y_L$. Let $N \geq 1$ be a dyadic number and $p \geq 2$ be such that $p \mu > 2$. Recall that $\Dl_N^L$ is as defined in \eqref{idDNL}. Then, by the Gaussian hypercontractivity, we get
\begin{align}
\E \Big[ \| \Dl_N^L Y_L \|_{L_{- \mu}^p (\T_L^2)}^p \Big] \leq p^{\frac{p}{2}} \int_{\T_L^2} \Big( \E \big[ | \Dl_N^L Y_L (x) |^2 \big] \Big)^{\frac{p}{2}} \jbb{x}_{L, - \mu}^p dx .
\label{Lreg1}
\end{align}

\noi
Note that for every $x \in \T_L^2$, by \eqref{YL_fs}, Lemma~\ref{LEM:gnL}, and the support property of $\psi_N$, we obtain
\begin{align}
\E \big[ | \Dl_N^L Y_L (x) |^2 \big] \les \frac{1}{L^2} \sum_{n \in \Z^2} \jb{\tfrac{n}{L}}^{-4} | \psi_N (\tfrac{n}{L}) |^2 \les \int_{\{ |\xi| \sim N \}} \jb{\xi}^{-4} d\xi \les N^{-2} .
\label{Lreg2}
\end{align}

\noi
Thus, combining \eqref{Lreg1} and \eqref{Lreg2} and using Lemma~\ref{LEM:Lwei}~(ii), we get
\begin{align}
\E \Big[ \| \Dl_N^L Y_L \|_{L_{- \mu}^p (\T_L^2)}^p \Big] \leq C^p p^{\frac{p}{2}} N^{-1} \int_{\T_L^2} \jbb{x}_{L, -\mu}^p dx \leq C^p p^{\frac{p}{2}} N^{-1} \int_{\T_L^2} \jb{x}_L^{- p \mu} dx \leq C^p p^{\frac{p}{2}} N^{-1}
\label{Lreg2-1}
\end{align}

\noi
for some constant $C > 0$. Then, we use the Besov embedding (Lemma~\ref{LEM:embL}~(iii)) to obtain
\begin{align}
\E \Big[ \| Y_L \|_{\C_{- \mu}^{s} (\T_L^2)}^p \Big] \leq C^p \sum_{\substack{N \geq 1 \\ \text{dyadic}}} N^{(s + \frac 2p) p} \E \Big[ \| \Dl_N^L Y_L \|_{L_{- \mu}^p (\T_L^2)}^p \Big] \leq C^p p^{\frac{p}{2}} 
\label{Lreg3}
\end{align}

\noi
for some constant $C > 0$ depending only on $s$, where we need to take $p \gg 1$ such that $s - 1 + \frac 2p < \frac{s - 1}{2} < 0$. 
This gives the desired almost sure bound for $\| Y_L \|_{\C^{s}_{- \mu} (\T_L^2)}$.
The bound for $\nb Y_L$ follows from a similar manner. 
The difference estimates for $Y_{L, \eps} - Y_L$ and $\nb Y_{L, \eps} - \nb Y_L$ then follows from Lemma~\ref{LEM:rho_eps}. 
These then imply the uniform bounds for $Y_{L, \eps}$ and $\nb Y_{L, \eps}$.

\smallskip \noi
(ii) From \eqref{xiL_fs} and Lemma~\ref{LEM:GL}, we have the following Fourier series for $\varphi *_{\R^2} \xi_L$:
\begin{align*}
\varphi *_{\R^2} \xi_L (x, \om) = \frac{1}{L} \sum_{n \in \Z^2} \ft \varphi (\tfrac{n}{L}) g_n^L (\om) e^{2 \pi i \frac{n}{L} \cdot x} ,
\end{align*}

\noi
where $\ft \varphi = \mathcal{F}_{\R^2} (\varphi)$. Since $\varphi \in C_c^\infty (\R^2)$, we know that $\ft \varphi$ is a Schwartz function, so that $| \ft \varphi (\frac{n}{L}) | \les_{M} \jb{\frac{n}{L}}^{-M}$ for any $M > 0$. Using this fact along with similar steps in part (i), we obtain the desired bound for $\varphi *_{\R^2} \xi_L$. The difference estimate for $\varphi *_{\R^2} \xi_{L, \eps} - \varphi *_{\R^2} \xi_L$ then follows from Lemma~\ref{LEM:rho_eps}. 
This then implies the uniform bounds for $\varphi *_{\R^2} \xi_{L, \eps}$.

\smallskip \noi
(iii) 
We first consider the bound for $e^{a Y_L}$. From the Taylor expansion and the product estimate in Lemma~\ref{LEM:prodL}~(ii), we have 
\begin{align}
\begin{split}
\| e^{a Y_{L}} \|_{\C^s_{- \mu} (\T_L^2)} \leq \sum_{k = 0}^\infty \frac{C^k}{k!} \| Y_{L} \|_{\C^s_{- \frac{\mu}{k}} (\T_L^2)}^k 
\end{split}
\label{eLreg1}
\end{align}

\noi
for some $C > 0$. Let $p \geq 2$ be large enough such that $p \mu > 2$ and $s - 1 + \frac{2}{p} < \frac{s - 1}{2} < 0$. By the Besov embedding (Lemma~\ref{LEM:embL}~(iii)) and \eqref{Lreg2-1}, we have
\begin{align}
\begin{split}
\E \Big[ \| Y_{L} \|_{\C^{s}_{- \frac{\mu}{k}} (\T_L^2)}^{p k} \Big] &\leq \E \Big[ \sup_{\substack{N \geq 1 \\ \text{dyadic}}} N^{s p k} \| \Dl_N^L Y_L \|_{L_{- \frac{\mu}{k}}^{\infty} (\T_L^2)}^{pk} \Big] \\
&\leq C^{pk} \sum_{\substack{N \geq 1 \\ \text{dyadic}}} N^{spk + 2} \E \Big[ \| \Dl_N^L Y_L \|_{L^{pk}_{- \frac{\mu}{k}} (\T_L^2)}^{pk} \Big] \\
&\les C^{p k} (p k)^{\frac{p k}{2}} .
\end{split}
\label{eLreg2}
\end{align}

\noi
Thus, combining \eqref{eLreg1} and \eqref{eLreg2} and using a ratio test, we have
\begin{align*}
\E \big[ \| e^{a Y_{L}} \|_{\C^\al_{- \mu} (\T_L^2)}^p \big]^{\frac{1}{p}} \leq \sum_{k = 0}^\infty \frac{C^k p^{\frac{k}{2}} k^{\frac{k}{2}}}{k!} < \infty ,
\end{align*}

\noi
which gives the desired almost sure bound for the $\C^s_{- \mu} (\T_L^2)$-norm. The bound for the $L^\infty_{-\mu} (\T_L^2)$-norm follows directly from the bound for the $\C^s_{- \mu} (\T_L^2)$-norm together with the embeddings in Lemma~\ref{LEM:embL}~(ii) and (i). The uniform bounds for $e^{a Y_{L, \eps}}$ follows from the same steps with an additional use of Lemma~\ref{LEM:rho_eps}.

For the difference $e^{a Y_{L, \eps}} - e^{a Y_L}$, we use the mean value theorem, Lemma~\ref{LEM:Lwei}~(iii), and the difference estimate in part (i) to obtain
\begin{align*}
\| e^{a Y_{L, \eps}} - e^{a Y_L} \|_{L_{- \mu}^\infty (\T_L^2)} 
\les \| Y_{L, \eps} - Y_L \|_{L_{- \frac{\mu}{2}}^\infty (\T_L^2)} \Big( \| e^{a Y_{L, \eps}} \|_{L_{- \frac{\mu}{2}}^\infty (\T_L^2)} + \| e^{a Y_L} \|_{L_{- \frac{\mu}{2}}^\infty (\T_L^2)} \Big) 
\les_{\om} \eps^\dl 
\end{align*}

\noi
for any $0 < \dl < 1 - s$, which gives the desired estimate.

\smallskip \noi
(iv) Let $N \geq 1$ be a dyadic number. By \eqref{YL2_fs}, Lemma~\ref{LEM:gnL}, the support property of $\psi_N$, and Lemma~\ref{LEM:convc}, we compute that for any $x \in \T_L^2$,
\begin{align*}
\E \big[ | \Dl_N^L (\wick{ |\nb Y_L|^2 }) (x) |^2 \big] &\les \frac{1}{L^4} \sum_{n \in \Z^2} \psi_N (\tfrac{n}{L}) \bigg( \sum_{\substack{n_1, n_2 \in \Z^2 \\ n_1 - n_2 = n}} \jb{\tfrac{n_1}{L}}^{-2} \jb{\tfrac{n_2}{L}}^{-2} + \jb{\tfrac{n}{L}}^{-4} \bigg) \\
&\les \int_{\{ |\xi| \sim N \}} \int_{\R^2} \jb{\xi + \xi'}^{-2} \jb{\xi'}^{-2} d \xi' d \xi + \int_{\R^2} \jb{\xi}^{-4} d\xi \\
&\les N^\kappa + 1 
\end{align*}

\noi
for any $\kappa > 0$ arbitrarily small. Thus, by using the fact that $\Dl_N^L (\wick{ |\nb Y_L|^2 })$ lies in the second Wiener chaos and the hypercontractivity of the second Wiener chaos, we can repeat the procedure in part (i) leading to \eqref{Lreg3} to obtain 
\begin{align*}
\E \big[ \| \wick{ | \nb Y_L |^2 } \|_{\C^{s - 1}_{- \mu} (\T_L^2)}^p \big] \les C^p p^{p} 
\end{align*}

\noi
for $p \gg 1$ sufficiently large and some constant $C > 0$. 
This gives the desired almost sure bound for $\| \wick{ | \nb Y_L |^2 } \|_{\C^{s - 1}_{- \mu} (\T_L^2)}$.

We now consider the difference $\wick{ |\nb Y_{L, \eps}|^2 } - \wick{ |\nb Y_{L}|^2 }\,$. From \eqref{def_YLe2} and \eqref{YL2_fs}, we write
\begin{align*}
&\wick{ | \nb Y_{L, \eps} |^2 } (x, \om) - \wick{ | \nb Y_L |^2 } (x, \om) \\
&= \frac{1}{L^2} \sum_{\substack{n, n_1, n_2 \in \Z^2 \\ n_1 - n_2 = n}} \frac{\big(1 - \ft G (\frac{n_1}{L}) + \ft \varphi (\frac{n_1}{L})\big) \big(1 - \cj{\ft G} (\frac{n_2}{L}) + \cj{\ft \varphi} (\frac{n_2}{L})\big) \frac{n_1}{L} \cdot \frac{n_2}{L}}{\jb{\frac{n_1}{L}}^2 \jb{\frac{n_2}{L}}^2} \\
&\qquad \qquad \qquad \times \big( \ft \rho (\tfrac{\eps n_1}{L}) \ft \rho (\tfrac{\eps n_2}{L}) - 1 \big) g_{n_1}^L (\om) \cj{g_{n_2}^L} (\om) e^{2 \pi i \frac{n}{L} \cdot x} \\
&\quad + \frac{1}{L^2} \sum_{n \in \Z^2} \frac{\big| 1 - \ft G (\frac{n}{L}) + \ft \varphi (\frac{n}{L}) \big|^2}{\jb{\frac{n}{L}}^2} \big( \ft \rho (\tfrac{\eps n}{L})^2 - 1 \big) (| g_n^L (\om) |^2 - 1)
\end{align*}

\noi
By the fact that $\int_{\R^2} \rho = 1$ and the mean value theorem, we have
\begin{align}
\big| \ft \rho (\tfrac{\eps n}{L}) - 1 \big| = \big| \ft \rho (\tfrac{\eps n}{L}) - \ft \rho (0) \big| \les \eps^{\dl_0} \jb{\tfrac{n}{L}}^{\dl_0}
\label{rho_diff}
\end{align}

\noi
for any $n \in \Z^2$ and $0 < \dl_0 < 1$. Thus, by repeating the procedure in the last paragraph with a slightly modified computation along with \eqref{rho_diff}, we get
\begin{align}
\E \Big[ \big\| \wick{|\nb Y_{L, \eps}|^2} - \wick{|\nb Y_L|^2} \big\|_{\C^{s - 1}_{- \mu} (\T_L^2)}^p \Big] \les C^p p^p \eps^{\dl_0 p}
\label{kolmo1}
\end{align}

\noi
for $p \gg 1$ sufficiently large and some constant $C > 0$. For any $0 < \eps_2 < \eps_1 \leq \frac 12$, using similar steps along with
\begin{align*}
\big| \ft \rho (\tfrac{\eps_1 n}{L}) - \ft \rho (\tfrac{\eps_2 n}{L}) \big| \les (\eps_1 - \eps_2)^{\dl_0} \jb{\tfrac{n}{L}}^{\dl_0}  ,
\end{align*}

\noi
we obtain
\begin{align}
\E \Big[ \big\| \wick{|\nb Y_{L, \eps_1}|^2} - \wick{|\nb Y_{L, \eps_2}|^2} \big\|_{\C^{s - 1}_{- \mu} (\T_L^2)}^p \Big] \les C^p p^p (\eps_1 - \eps_2)^{\dl_0 p} .
\label{kolmo2}
\end{align}

\noi
With \eqref{kolmo1} and \eqref{kolmo2}, we now take $p \gg 1$ to be large enough so that we can invoke the Kolmogorov continuity criterion (see \cite[Theorem~3.3]{DPZ}) to deduce that $\eps \mapsto \, \wick{ |\nb Y_{L, \eps}|^2 }$ from $[0, \frac 12]$ to $\C^{\al - 1} (\T_L^2)$ is almost surely H\"older continuous. This gives the desired difference estimate, which then implies the uniform bound for $\wick{ |\nb Y_{L, \eps}|^2 }$ .
\end{proof}

In the following lemma, we show that $\nb Y_L$ and $\wick{|\nb Y_L|^2}$ are only logarithmically away from being a function in any $L_{- \mu}^q (\T_L^2)$-space for any finite $q$ and $\mu > 0$.

\begin{lemma}
\label{LEM:YregLe2}
Let $L \geq 1$, $2 \leq q < \infty$, and $\mu > 0$ be such that $q \mu > 2$. Then, there exists $\Om' \subset \Om$ with full probability measure such that for any $\om \in \Om'$ and $0 < \eps < \frac 12$, we have
\begin{align}
\| \nb Y_{L, \eps} \|_{L_{- \mu}^q (\T_L^2)} \les_{\om} |\log \eps|
\label{loge1}
\end{align}

\noi
and
\begin{align}
\big\| \wick{ |\nb Y_{L, \eps}|^2 } \big\|_{L_{- \mu}^q (\T_L^2)} \les_{\om} |\log \eps|^2 .
\label{loge2}
\end{align}
\end{lemma}

\begin{proof}
The bound \eqref{loge2} follows direct from \eqref{loge1}, since from \eqref{def_YLe2}, \eqref{CLe}, and Lemma~\ref{LEM:Lwei}~(ii), we have
\begin{align*}
\big\| \wick{ |\nb Y_{L, \eps}|^2 } \big\|_{L_{- \mu}^q (\T_L^2)} &\leq \| \nb Y_{L, \eps} \|_{L_{- \frac{\mu}{2}}^{2q} (\T_L^2)}^2 + C_{L, \eps} \bigg( \int_{\T_L^2} \jbb{x}_{L, - \mu}^q dx \bigg)^{\frac{1}{q}} \\
&\les |\log \eps|^2 + |\log \eps| \bigg( \int_{\T_L^2} \jb{x}_{L}^{- q \mu} dx \bigg)^{\frac{1}{q}} \\
&\les |\log \eps|^2 .
\end{align*}

\noi
Thus, we focus on proving \eqref{loge1}.

We proceed as in \cite[Proposition~3.1]{DLTV}. From \eqref{YLe_fs}, the identity in Lemma~\ref{LEM:rho_eps}, the dyadic decomposition \eqref{idDNL}, and Young's convolution inequality along with Lemma~\ref{LEM:Lwei}~(v), we have
\begin{align*}
\| \nb Y_{L, \eps} \|_{L_{- \mu}^q (\T_L^2)} &\les \sum_{\substack{N \geq 1 \\ \text{dyadic}}} \| \Dl_N^L ( \rho_{L, \eps} *_{\T_L^2} \nb Y_L ) \|_{L_{- \mu}^q (\T_L^2)} \\
&\leq \sum_{\substack{N \geq 1 \\ \text{dyadic}}} \sum_{M = \frac{N}{2}, N, 2N} \| \Dl_M^L \rho_{L, \eps} *_{\T_L^2} \Dl_N^L \nb Y_L \|_{L_{- \mu}^q (\T_L^2)} \\
&\leq \sum_{\substack{N \geq 1 \\ \text{dyadic}}} \sum_{M = \frac{N}{2}, N, 2N} \| \Dl_M^L \rho_{L, \eps} \|_{L_\mu^1 (\T_L^2)} \| \Dl_N^L \nb Y_L \|_{L_{- \mu}^q (\T_L^2)} ,
\end{align*}

\noi
so that by using H\"older's inequality on the last expression and applying \eqref{rhoL2} in Lemma~\ref{LEM:rhoDNL}, we get
\begin{align}
\| \nb Y_{L, \eps} \|_{L_{- \mu}^q (\T_L^2)} 
\les \| \rho_{L, \eps} \|_{\B_{1, 1, \mu}^0 (\T_L^2)} \| \nb Y_L \|_{\B_{q, \infty, - \mu}^0 (\T_L^2)} 
\les |\log \eps| \| \nb Y_L \|_{\B_{q, \infty, - \mu}^0 (\T_L^2)} .
\label{YLq1}
\end{align}

\noi
To estimate $\| \nb Y_L \|_{\B_{q, \infty, - \mu}^0 (\T_L^2)}$, we first use similar steps as in \eqref{Lreg2-1} along with H\"older's inequality to obtain 
\begin{align}
\E \big[ \| \Dl_N^L \nb Y_L \|_{L_{- \mu}^q (\T_L^2)} \big] \les 1.
\label{YLq1-1}
\end{align}

\noi
Also, given $f_L \in L_{\mu}^{q'} (\T_L^2)$ with $\frac{1}{q} + \frac{1}{q'} = 1$, we use \eqref{YL_fs}, Lemma~\ref{LEM:gnL}, and Plancherel's identity \eqref{planL} to obtain
\begin{align*}
\E \big[ | ( \Dl_N^L \nb Y_L, f_L)_{\T_L^2} |^2 \big] &= \E \bigg[ \bigg| \sum_{n \in \Z^2} \frac{2 \pi (1 - \ft G (\frac{n}{L}) + \ft \varphi (\frac{n}{L})) \frac{n}{L}}{\jb{\frac{n}{L}}^2} \psi_N (\tfrac{n}{L}) g_n^L (\om) \cj{\ft{f_L} (\tfrac{n}{L})} \bigg|^2 \bigg] \\
&\les \sum_{n \in \Z^2} \jb{\tfrac{n}{L}}^{-2} \psi_N (\tfrac{n}{L})^2 \big| \ft{f_L} (\tfrac{n}{L}) \big|^2 \\
&\sim N^{-2} \| \Dl_N^L f_L \|_{L^2 (\T_L^2)}^2 ,
\end{align*}

\noi
so that by using the embedding in Lemma~\ref{LEM:embL}~(iii) on the last expression and applying Lemma~\ref{LEM:Lwei}~(i), we get
\begin{align}
\E \big[ | ( \Dl_N^L \nb Y_L, f_L)_{\T_L^2} |^2 \big] 
\les N^{- \frac{2 (2 q' - 2)}{q'}} \| f_L \|_{L^{q'} (\T_L^2)} 
\les N^{- \frac{2 (2 q' - 2)}{q'}} \| f_L \|_{L_{\mu}^{q'} (\T_L^2)} .
\label{YLq1-2}
\end{align}

\noi
Thus, from Lemma~\ref{LEM:Ver} along with \eqref{YLq1-1} and \eqref{YLq1-2}, we get
\begin{align}
\E \big[ \| \nb Y_L \|_{\B_{q, \infty, -\mu}^0 (\T_L^2)} \big] \les 1 .
\label{YLq2}
\end{align}

\noi
Combining \eqref{YLq1} and \eqref{YLq2}, we obtain the desired estimate \eqref{loge1}.
\end{proof}

\subsection{Large torus convergence of stochastic objects}
\label{SUB:sto_conv}

Again, we let $\rho \in C_c^\infty (\R^2)$ be such that $\supp \rho \in \{ x \in \R^2: |x| < \frac 14 \}$, $\rho \geq 0$, and $\int_{\R^2} \rho = 1$. Given $0 < \eps < 1$, we define $\rho_\eps (x) = \eps^{-2} \rho (\eps^{-1} x)$. We define the regularized version of $\xi$ and $Y$ in \eqref{defY} as
\begin{align}
\xi_\eps \deff \rho_\eps *_{\R^2} \xi
\label{xie}
\end{align}

\noi
and
\begin{align}
Y_\eps \deff \rho_\eps *_{\R^2} Y = \rho_\eps *_{\R^2} G *_{\R^2} \xi.
\label{Ye}
\end{align}

\noi
We then define
\begin{align}
\wick{| \nb Y_\eps |^2} \deff |\nb Y_\eps|^2 - C_\eps
\label{Ye2}
\end{align}

\noi
with (see \cite[Section~3]{HL15})
\begin{align}
C_\eps = \E [ |\nb Y_\eps|^2 ] = \| \rho_\eps *_{\R^2} \nb G \|_{L^2 (\R^2)}^2 \sim |\log \eps|.
\label{Ce}
\end{align}

We have the following estimates from \cite[Proposition~3.1 and Proposition~3.5]{DLTV}.
\begin{lemma}
\label{LEM:Ye_est}
\textup{(i)} Let $2 < q < \infty$ and $\mu > 0$ satisfying $\mu q > 2$. Then, there exists $\Om' \subset \Om$ with full probability measure such that for any $\om \in \Om'$ and $0 < \eps < \frac 12$, we have
\begin{align*}
\| \nb Y_\eps \|_{L^q_{- \mu} (\R^2)} \les_{\om} |\log \eps| .
\end{align*}

\smallskip \noi
\textup{(ii)} Let $\mu > 0$ and $a \in \R$. Then, there exists $\Om' \subset \Om$ with full probability measure such that for any $\om \in \Om'$, we have
\begin{align*}
\sup_{\eps \in (0, 1)} \| e^{a Y_\eps} \|_{L^\infty_{- \mu} (\R^2)} \les_\om 1.
\end{align*}
\end{lemma}

Given $L \geq 1$ and $0 < \eps < 1$, we recall $\xi_{L, \eps}$ in \eqref{xiLe_fs}, $Y_{L, \eps}$ in \eqref{YLe_fs}, and $\wick{|\nb Y_{L, \eps}|^2}$ in \eqref{def_YLe2}. We also recall the function $\s_R$ in \eqref{sigR}.

\begin{proposition}
\label{PROP:YconvL}
Let $0 < \mu \leq 4$, $2 < q < \infty$ satisfying $\mu q > 4$, $\varphi \in C_c^\infty (\R^2)$, and $a \in \R$. Then, for any $1 \leq R \leq L$, there exists $\Om' \subset \Om$ with full probability measure such that for any $\om \in \Om'$ and $0 < \eps < \frac 12$, the following estimates hold.

\smallskip \noi
\textup{(i)} We have
\begin{align*}
&\| \s_R^{-1} ( \nb Y_\eps - \nb Y_{L, \eps} ) \|_{L^q (\R^2)}  \les_{\omega} R^{\mu} L^{- \frac{\mu}{2}} |\log \eps| .
\end{align*}

\smallskip \noi
\textup{(ii)} We have
\begin{align*}
\| \s_R^{-1} ( \varphi *_{\R^2} \xi_\eps - \varphi *_{\R^2} \xi_{L, \eps} ) \|_{L^q (\R^2)} \les_{\omega} R^\mu L^{- \frac{\mu}{2}} .
\end{align*}

\smallskip \noi
\textup{(iii)} We have
\begin{align*}
\| \s_R^{-2} ( e^{a Y_\eps} - e^{a Y_{L, \eps}} ) \|_{L^q (\R^2)} \les_{\omega} R^{2 \mu} L^{- \frac{\mu}{4}} .
\end{align*}

\smallskip \noi
\textup{(iv)} We have
\begin{align*}
\big\| \s_R^{-2} \big( \wick{| \nb Y_\eps |^2} - \wick{| \nb Y_{L, \eps} |^2} \big) \big\|_{L^q (\R^2)} \les_{\omega} R^{2 \mu} L^{- \frac{\mu}{2} + \frac{2}{q}} |\log \eps|^2 .
\end{align*}
\end{proposition}

\begin{proof}
(i) From \eqref{sigR2}, \eqref{Ye}, \eqref{YLe_fs}, and the dyadic decomposition \eqref{DN}, we have
\begin{align}
\begin{split}
\| \s_R^{-1} ( \nb Y_\eps - \nb Y_{L, \eps} ) \|_{L^q (\R^2)} 
&\les R^{\mu} \| \nb Y_\eps - \nb Y_{L, \eps} \|_{L_{- \mu}^q (\R^2)} \\
&\leq R^{\mu} \sum_{\substack{N \geq 1 \\ \text{dyadic}}} \big\| \Dl_N \big( \rho_\eps *_{\R^2} (\nb Y - \nb Y_L) \big) \big\|_{L_{- \mu}^q (\R^2)} .
\end{split}
\label{Yc1-1}
\end{align}

\noi
By using Young's convolution inequality, we have
\begin{align}
\begin{split}
&\sum_{\substack{N \geq 1 \\ \text{dyadic}}} \big\| \Dl_N \big( \rho_\eps *_{\R^2} (\nb Y - \nb Y_L) \big) \big\|_{L_{- \mu}^q (\R^2)} \\
&\les \sum_{\substack{N \geq 1 \\ \text{dyadic}}} \sum_{M = \frac{N}{2}, N, 2N} \| \Dl_M \rho_\eps *_{\R^2} \Dl_N (\nb Y - \nb Y_L) \|_{L^{q}_{- \mu} (\R^2)} \\
&\les \sum_{\substack{N \geq 1 \\ \text{dyadic}}} \sum_{M = \frac{N}{2}, N, 2N} \| \Dl_M \rho_\eps \|_{L^1_{\mu} (\R^2)} \| \Dl_N (\nb Y - \nb Y_L) \|_{L_{- \mu}^q (\R^2)} 
\end{split}
\label{Yc1-2}
\end{align}

\noi
for any $\mu > 0$. Thus, from \eqref{Yc1-1}, \eqref{Yc1-2}, H\"older's inequality, Lemma~\ref{LEM:rhoB}, and Lemma~\ref{LEM:G} along with \eqref{defY}, we get
\begin{align}
\begin{split}
\| \s_R^{-1} ( \nb Y_\eps - \nb Y_{L, \eps} ) \|_{L^q (\R^2)} 
&\les R^{\mu} \| \rho_\eps \|_{\B_{1, 1, \mu}^0 (\R^2)} \| \nb Y - \nb Y_L \|_{\B_{q, \infty, - \mu}^0 (\R^2)} \\
&\les R^{\mu} |\log \eps| \| \xi - \xi_L \|_{\B_{q, \infty, - \mu}^{-1} (\R^2)} .
\end{split}
\label{Yconv1-0}
\end{align}

We now estimate $\| \xi - \xi_L \|_{\B_{q, \infty, - \mu}^{-1} (\R^2)}$ using Lemma~\ref{LEM:Ver}. Let $N \geq 1$ be a dyadic number and $\eta_N (\cdot) = N^2 \eta (N \cdot)$ be the kernel of $\Dl_N$ in \eqref{DN}. 
By the definition in \eqref{xiL}, we see that for any $x \in \R^2$,
\begin{align*}
\E &\big[ | \Dl_N (\xi - \xi_L) (x) |^2 \big] \\
&= N^4 \E \big[ \big| \big( \xi - \xi_L, \eta (N (x - \cdot)) \big)_{\R^2} \big|^2 \big] \\
&= N^4 \int_{\R^2} \Big| \eta (N (x - z)) - \sum_{k \in \Z^2} \eta (N (x - z + L k)) \ind_{[- \frac{L}{2}, \frac{L}{2})^2} (z) \Big|^2 dz ,
\end{align*}

\noi
so that we get
\begin{align}
\begin{split}
\E &\big[ | \Dl_N (\xi - \xi_L) (x) |^2 \big] \\
&\les N^4 \int_{\R^2 \setminus [- \frac{L}{2}, \frac{L}{2})^2} | \eta (N(x - z)) |^2 dz + N^4 \int_{[- \frac{L}{2}, \frac{L}{2})^2} \Big| \sum_{k \in \Z^2 \setminus \{0\}} \eta (N (x - z + Lk)) \Big|^2 dz \\
&\deff \textup{I}_1 + \textup{I}_2 .
\end{split}
\label{conv2-0}
\end{align}

\noi
For $\textup{I}_1$, we use the fact that $|z| \geq \frac{L}{2}$ and a change of variable to obtain
\begin{align}
\begin{split}
\textup{I}_1 &\les L^{- \mu} N^4 \int_{\R^2} \jb{z}^{\mu} | \eta (N (x - z)) |^2 dz \\
&\les L^{- \mu} N^4 \jb{x}^{\mu} \int_{\R^2} \jb{x - z}^{\mu} | \eta (N (x - z)) |^2 dz \\
&\les L^{- \mu} N^2 \jb{x}^{\mu} .
\end{split}
\label{conv2-1}
\end{align}

\noi
For $\textup{I}_2$, since $\eta$ is a Schwartz function, we have
\begin{align*}
\Big| \sum_{k \in \Z^2 \setminus \{0\}} \eta (N (x - z + Lk)) \Big| \les \sum_{k \in \Z^2} \frac{1}{\jb{N (x - z + Lk)}^{10}} \les 1
\end{align*}

\noi
uniform in $x$ and $z$, so that we use a similar treatment as that for $\textup{I}_1$ to obtain
\begin{align}
\begin{split}
\text{I}_2 &\les N^4 \int_{[- \frac{L}{2}, \frac{L}{2})^2} \sum_{k \in \Z^2 \setminus \{0\}} | \eta (N (x - z + L k)) | d z \\
&= N^4 \int_{\R^2 \setminus [ - \frac{L}{2}, \frac{L}{2} )^2} | \eta (N (x - z)) | dz \\
&\les L^{- \mu} N^2 \jb{x}^{\mu} .
\end{split}
\label{conv2-2}
\end{align}

\noi
Thus, by the Gaussian hypercontractivity, \eqref{conv2-0}, \eqref{conv2-1}, and \eqref{conv2-2}, we have
\begin{align*}
\begin{split}
\E \Big[ \| \Dl_N (\xi - \xi_L) \|_{L_{- \mu}^q (\R^2)}^q \Big] &= \int_{\R^2} \jb{x}^{- q \mu} \E \big[ | \Dl_N (\xi - \xi_L) (x) |^q \big] dx \\
&\les \int_{\R^2} \jb{x}^{- q \mu} \Big( \E \big[ |\Dl_N (\xi - \xi_L) (x)|^2 \big] \Big)^{\frac{q}{2}} dx \\
&\les L^{- \frac{q \mu}{2}} N^q \int_{\R^2} \jb{x}^{- \frac{q \mu}{2}} dx ,
\end{split}
\end{align*}

\noi
so that we get
\begin{align}
\E \Big[ \| \Dl_N (\xi - \xi_L) \|_{L_{- \mu}^q (\R^2)}^q \Big] 
\les_{\mu} L^{- \frac{q \mu}{2}} N^q ,
\label{conv2-3}
\end{align}

\noi
where we need to take $q$ large enough so that $\mu q > 4$. Thus, by \eqref{conv2-3} and H\"older's inequality, we obtain
\begin{align}
\E \big[ N^{-1} \| \Dl_N (\xi - \xi_L) \|_{L_{- \mu}^q (\R^2)} \big] \les L^{- \frac{\mu}{2}} .
\label{conv2i}
\end{align}

\noi
We now let $q'$ be the H\"older conjugate of $q$ and test $\Dl_N (\xi - \xi_L)$ against a function $f \in L^{q'}_{\mu} (\R^2)$. From the definition in \eqref{xiL}, we have
\begin{align*}
\E \big[ \big| \big( \Dl_N (\xi - \xi_L), f \big)_{\R^2} \big|^2 \big] &= \E \big[ | (\xi - \xi_L, \Dl_N f)_{\R^2} |^2 \big] \\
&= \int_{\R^2} \Big| \Dl_N f (z) - \sum_{k \in \Z^2} \Dl_N f (z + Lk) \ind_{[- \frac{L}{2}, \frac{L}{2})^2} (z) \Big|^2 dz ,
\end{align*}

\noi
so that we get
\begin{align}
\begin{split}
\E &\big[ \big| \big( \Dl_N (\xi - \xi_L), f \big)_{\R^2} \big|^2 \big] \\
&\les \int_{\R^2 \setminus [- \frac{L}{2}, \frac{L}{2})^2} |\Dl_N f (z)|^2 dz + \int_{[- \frac{L}{2}, \frac{L}{2})^2} \Big| \Dl_N f (z) - \sum_{k \in \Z^2} \Dl_N f (z + Lk) \Big|^2 dz \\
&\deff \textup{I}_3 + \textup{I}_4 .
\end{split}
\label{conv2-4}
\end{align}

\noi
For $\textup{I}_3$, we use the fact that $|z| \geq \frac{L}{2}$, Young's convolution inequality, and a change of variable to obtain
\begin{align}
\begin{split}
\textup{I}_3 &\les L^{- 2 \mu} \| \Dl_N f \|_{L^2_{\mu} (\R^2)}^2 \\
&\les L^{- 2 \mu} N^4 \| \eta (N \cdot) \|_{L^{\frac{2q}{q + 2}}_{\mu} (\R^2)}^2 \| f \|_{L^{q'}_{\mu} (\R^2)}^2 \\
&\les L^{-2 \mu} N^{\frac{2 (q - 2)}{q}} \| f \|_{L^{q'}_{\mu} (\R^2)}^2 .
\end{split}
\label{conv2-5}
\end{align}

\noi
For $\textup{I}_4$, we use \eqref{fLdiff} in Lemma~\ref{LEM:per} and similar steps in \eqref{conv2-5} to obtain
\begin{align}
\textup{I}_4 \les L^{- 2 \mu} \| \Dl_N f \|_{L_{\mu}^2 (\R^2)}^2 \les L^{-2 \mu} N^{\frac{2 (q - 2)}{q}} \| f \|_{L^{q'}_{\mu} (\R^2)}^2 .
\label{conv2-6}
\end{align}

\noi
Thus, combining \eqref{conv2-4}, \eqref{conv2-5}, and \eqref{conv2-6}, we get
\begin{align}
\E \big[ \big| \big( \Dl_N (\xi - \xi_L), f \big)_{\R^2} \big|^2 \big] \les L^{-2 \mu} N^{\frac{2 (q - 2)}{q}} \| f \|_{L^{q'}_{\mu} (\R^2)}^2 .
\label{conv2ii}
\end{align}

\noi
Therefore, by using Lemma~\ref{LEM:Ver} along with \eqref{conv2i} and \eqref{conv2ii}, we obtain
\begin{align}
\E \big[ \| \xi - \xi_L \|_{\B_{q, \infty, - \mu}^{-1} (\R^2)} \big] \les L^{- \frac{\mu}{2}} . 
\label{Yconv1-1}
\end{align}

\noi
Combining \eqref{Yconv1-0} and \eqref{Yconv1-1}, we obtain the desired estimate.

\smallskip \noi
(ii) By using \eqref{sigR2}, \eqref{xie}, \eqref{xiLe_fs}, and Young's convolution inequality, we have
\begin{align*}
\| &\s_R^{-1} ( \varphi *_{\R^2} \xi_\eps - \varphi *_{\R^2} \xi_{L, \eps} ) \|_{L^q (\R^2)} \\
&\les R^\mu \| \rho_\eps *_{\R^2} \varphi *_{\R^2} \xi - \rho_\eps *_{\R^2} \varphi *_{\R^2} \xi_L \|_{L^q_{-\mu} (\R^2)} \\
&\les R^\mu \| \rho_\eps \|_{L_\mu^1 (\R^2)} \| \varphi *_{\R^2} \xi - \varphi *_{\R^2} \xi_L \|_{L^q_{- \mu} (\R^2)} 
\end{align*}

\noi
for any $\mu > 0$, so that by further applying the embedding in Lemma~\ref{LEM:emb}~(ii) and Lemma~\ref{LEM:varphi}, we obtain
\begin{align}
\begin{split}
\| &\s_R^{-1} ( \varphi *_{\R^2} \xi_\eps - \varphi *_{\R^2} \xi_{L, \eps} ) \|_{L^q (\R^2)} \\
&\les R^{\mu} \| \varphi *_{\R^2} \xi - \varphi *_{\R^2} \xi_L \|_{\B^{0}_{q, 1, - \mu} (\R^2)} \\
&\les R^{\mu} \| \xi - \xi_L \|_{\B_{q, \infty, - \mu}^{-1} (\R^2)}
\end{split}
\label{Yconv2}
\end{align}

\noi
for any $\mu > 0$. The desired estimate follows directly from \eqref{Yconv2} and \eqref{Yconv1-1}.

\smallskip \noi
(iii) By the mean value theorem, H\"older's inequality, \eqref{sigR2}, Young's convolution inequality with \eqref{Ye} and \eqref{YLe_fs}, and \eqref{defY}, we have
\begin{align*}
\| &\s_R^{-2} (e^{a Y_\eps} - e^{a Y_{L, \eps}}) \|_{L^q (\R^2)} \\
&\les \big( \| \s_R^{-1} e^{a Y_\eps} \|_{L^\infty (\R^2)} + \| \s_R^{-1} e^{a Y_{L, \eps}} \|_{L^\infty (\R^2)} \big) \| \s_R^{-1} ( Y_\eps - Y_{L, \eps} ) \|_{L^q (\R^2)} \\
&\les R^\mu \big( R^\mu \| e^{a Y_\eps} \|_{L^\infty_{- \mu} (\R^2)} + \| e^{a Y_{L, \eps}} \|_{L^\infty (\T_L^2)} \big) \| G *_{\R^2} \xi - G *_{\R^2} \xi_L \|_{L^q_{- \mu} (\R^2)} \| \rho_\eps \|_{L_\mu^1 (\R^2)} 
\end{align*}

\noi
for any $\mu > 0$, so that by further applying Lemma~\ref{LEM:Ye_est}~(ii), Lemma~\ref{LEM:YregLe}~(iii), the embedding in Lemma~\ref{LEM:emb}~(ii), Lemma~\ref{LEM:G}, and the embedding in Lemma~\ref{LEM:emb}~(ii), we obtain
\begin{align}
\begin{split}
\| &\s_R^{-2} (e^{a Y_\eps} - e^{a Y_{L, \eps}}) \|_{L^q (\R^2)} \\
&\les_{\om} R^{2 \mu} L^{\frac{\mu}{4}} \| G *_{\R^2} \xi - G *_{\R^2} \xi_L \|_{\B^0_{q, 1, - \mu} (\R^2)} \\
&\les_{\om} R^{2 \mu} L^{\frac{\mu}{4}} \| \xi - \xi_L \|_{\B_{q, \infty, -\mu}^{-1} (\R^2)}
\end{split}
\label{Yconv3}
\end{align}

\noi
for any $0 < \mu \leq 4$. The desired estimate follows directly from \eqref{Yconv3} and \eqref{Yconv1-1}.

\smallskip \noi
(iv) From \eqref{Ye2} and \eqref{def_YLe2}, we have
\begin{align*}
\wick{| \nb Y_\eps |^2} - \wick{| \nb Y_{L, \eps} |^2} \, = ( | \nb Y_\eps |^2 - | \nb Y_{L, \eps} |^2 ) - \big( \E [| \nb Y_\eps |^2] - \E [| \nb Y_{L, \eps} |^2] \big) .
\end{align*}

\noi
By H\"older's inequality, part (i), \eqref{sigR2}, and Lemma~\ref{LEM:sig}, we have
\begin{align*}
\big\| &\s_R^{-2} ( | \nb Y_\eps |^2 - | \nb Y_{L, \eps} |^2 ) \big\|_{L^q (\R^2)} \\
&\leq \| \s_R^{-1} (\nb Y_\eps - \nb Y_{L, \eps}) \|_{L^{2q} (\R^2)} \big( \| \s_R^{-1} \nb Y_\eps \|_{L^{2q} (\R^2)} + \| \s_R^{-1} \nb Y_{L, \eps} \|_{L^{2q} (\R^2)} \big) \\
&\les_{\om} R^{2 \mu} L^{- \frac{\mu}{2}} |\log \eps| \big( \| \nb Y_\eps \|_{L_{- \mu}^{2q} (\R^2)} + \| \nb Y_{L, \eps} \|_{L^{2q} (\T_L^2)} \big) 
\end{align*}

\noi
for any $\mu > 0$, so that by applying Lemma~\ref{LEM:Ye_est}~(i) and Lemma~\ref{LEM:YregLe2} on the last expression, we obtain
\begin{align}
\big\| &\s_R^{-2} ( | \nb Y_\eps |^2 - | \nb Y_{L, \eps} |^2 ) \big\|_{L^q (\R^2)} \les_{\om} R^{2 \mu} L^{- \frac{\mu}{2} + \frac{2}{q}} |\log \eps|^2
\label{Yconv4-1}
\end{align}

\noi
as long as $q$ is large enough satisfying $2 \mu q > 4$. 
Also, by denoting $G_\eps = \rho_\eps *_{\R^2} G$, we have from \eqref{defY} and the Cauchy-Schwarz inequality that
\begin{align*}
\E [| \nb Y_\eps |^2] - \E [| \nb Y_{L, \eps} |^2] &= \int_{\R^2} |\nb G_\eps (x)|^2 dx - \int_{[- \frac{L}{2}, \frac{L}{2})^2} \Big| \sum_{k \in \Z^2} \nb G_\eps (x + Lk) \Big|^2 dx \\
&\les \bigg(\int_{[- \frac{L}{2}, \frac{L}{2})^2} \Big| \nb G_\eps (x) - \sum_{k \in \Z^2} \nb G_\eps (x + Lk) \Big|^2 dx \bigg)^{\frac 12} \\
&\quad \times \bigg( \int_{\R^2} |\nb G_\eps (x)|^2 dx + \int_{[- \frac{L}{2}, \frac{L}{2})^2} \Big| \sum_{k \in \Z^2} \nb G_\eps (x + Lk) \Big|^2 dx \bigg)^{\frac 12} ,
\end{align*}

\noi
so that by further applying Lemma~\ref{LEM:per}, the fact that $\jb{x}^\mu \les 1$ on the support of $G_\eps$, and \eqref{Ce}, we obtain
\begin{align}
\E [| \nb Y_\eps |^2] - \E [| \nb Y_{L, \eps} |^2] 
\les L^{- \mu} \| \nb G_\eps \|_{L^2_\mu (\R^2)}^2 \les L^{- \mu} \| \nb G_\eps \|_{L^2 (\R^2)}^2 
\les L^{- \mu} |\log \eps|.
\label{Yconv4-2}
\end{align}

\noi
The desired estimate then follows from \eqref{Yconv4-1}, \eqref{Yconv4-2}, and the fact that $\| \s_R^{-2} \|_{L^q (\R^2)} \les_q R^{\frac{2}{q}} < R^{2 \mu}$ given $\mu q > 4$.
\end{proof}

\section{Large torus limit of the linear dispersive Anderson model}
\label{SEC:convL1}

In this section, we show Theorem~\ref{THM:GWP} and Theorem~\ref{THM:convL} in the case when $\ld = 0$, which correspond to global well-posedness of the $L$-periodic linear DAM and the large torus limit of the linear DAM, respectively. In order to distinguish this from the nonlinear problem, we write $w$ as the unknown for the equation \eqref{vNLSAnd} on $\R^2$ with $\ld = 0$:
\begin{align}
\begin{cases}
i \dt w = \Dl w - 2 \nb Y \cdot \nb w + \wt{\wick{|\nb Y|^2}} \, w \\
w |_{t = 0} = w_0 ,
\end{cases}
\label{wNLSAnd}
\end{align}

\noi
and $w_L$ as the unknown for the equation \eqref{vNLSAndL} on $\T_L^2$ with $\ld = 0$:
\begin{align}
\begin{cases}
i \dt w_L = \Dl w_L - 2 \nb Y_L \cdot \nb w_L + \wt{\wick{|\nb Y_L|^2}} \, w_L \\
w_L |_{t = 0} = w_{0, L} ,
\end{cases}
\label{wNLSAndL}
\end{align}

\noi
where $w_0$ belongs to $H_\mu^2 (\R^2)$ with $\mu > 1$ and $w_{0, L}$ is defined by
\begin{align*}
w_{0, L} = \sum_{k \in \Z^2} w_0 (\cdot + Lk) .
\end{align*}

\noi
From Lemma~\ref{LEM:per}, we know that $w_{0, L} \in H^2 (\T_L^2)$.

\subsection{A priori bounds for mollified $L$-periodic linear solutions}

Let us consider the following mollified $L$-periodic linear DAM:
\begin{align}
\begin{cases}
i \dt w_{L, \eps} = \Dl w_{L, \eps} - 2 \nb Y_{L, \eps} \cdot \nb w_{L, \eps} + \wt{\wick{|\nb Y_{L, \eps}|^2}} \, w_{L, \eps} \\
w_{L, \eps} |_{t = 0} = w_{0, L} ,
\end{cases}
\label{wNLSLe}
\end{align}

\noi
where $Y_{L, \eps}$ is defined in \eqref{YLe_fs} and
\begin{align}
\wt{\wick{|\nb Y_{L, \eps}|^2}} \deff \, \wick{|\nb Y_{L, \eps}|^2} - \, \varphi *_{\R^2} \xi_{L, \eps}.
\label{YLe2}
\end{align}

\noi
with $\wick{|\nb Y_{L, \eps}|^2}$ being defined in \eqref{def_YLe2} and $\xi_{L, \eps}$ being defined in \eqref{xiLe_fs}. To see that \eqref{wNLSLe} is globally well-posed for $w_{0, L} \in H^2 (\T_L^2)$, we see that $u_{L, \eps} \deff e^{- Y_{L, \eps}} w_{L, \eps}$ satisfies
\begin{align*}
\begin{cases}
i \dt u_{L, \eps} = \Dl u_{L, \eps} + \xi_{L, \eps} u_{L, \eps} - C_{L, \eps} u_{L, \eps} \\
u_{L, \eps} |_{t = 0} = e^{- Y_{L, \eps}} w_{0, L} ,
\end{cases}
\end{align*}

\noi
whose unique global solution exists almost surely in $C(\R_+; H^2 (\T_L^2))$ (see \cite{Tsu, BGTz}). This shows that \eqref{wNLSLe} has a unique global solution in $C (\R_+; H^2 (\T_L^2))$.
For later convenience, we define $S_{L, \eps} (t)$ as the flow map of the equation \eqref{wNLSLe}.

Let us define the mass
\begin{align}
\mathcal{M}_{L, \eps} [w] \deff \int_{\T_L^2} e^{- 2 Y_{L, \eps}} |w|^2 dx ,
\label{defM1}
\end{align}

\noi
the energy
\begin{align}
\mathcal{E}_{L, \eps} [w] \deff \int_{\T_L^2} \Big( \frac 12 e^{-2 Y_{L, \eps}} |\nb w|^2 - \frac 12 e^{- 2 Y_{L, \eps}} \wt{\wick{|\nb Y_{L, \eps}|^2}} \, |w|^2 \Big) dx ,
\label{defE1}
\end{align}

\noi
and the modified energy (introduced in \cite{TzV23})
\begin{align}
\begin{split}
\mathcal{F}_{L, \eps} [w] 
&\deff \int_{\T_L^2} e^{-2 Y_{L, \eps}} |\Dl w|^2 dx - 4 \Re \int_{\T_L^2} e^{- 2 Y_{L, \eps}} \Dl w \nb Y_{L, \eps} \cdot \nb \cj{w} dx \\
&\quad - 4 \int_{\T_L^2} e^{-2 Y_{L, \eps}} ( \nb Y_{L, \eps} \cdot \nb w )^2 dx + 2 \Re \int_{\T_L^2} \wt{\wick{|\nb Y_{L, \eps}|^2}} \, w \nb (e^{- 2 Y_{L, \eps}}) \cdot \nb \cj{w} dx  \\
&\quad + 2 \Re \int_{\T_L^2} e^{- 2 Y_{L, \eps}} \wt{\wick{|\nb Y_{L, \eps}|^2}} \, \Dl w \cj{w} dx + \int_{\T_L^2} e^{- 2 Y_{L, \eps}} \big( \wt{\wick{|\nb Y_{L, \eps}|^2}} \big)^2 |w|^2 dx .
\end{split}
\label{defF1}
\end{align}

\noi
A direct computation yields that the above three quantities are conserved under the flow of \eqref{wNLSLe}:
\begin{align*}
\frac{d}{dt} \mathcal{M}_{L, \eps} [S_{L, \eps} (t) f_L] &= 0, \\
\frac{d}{dt} \mathcal{E}_{L, \eps} [S_{L, \eps} (t) f_L] &= 0, \\
\frac{d}{dt} \mathcal{F}_{L, \eps} [S_{L, \eps} (t) f_L] &= 0,
\end{align*}

\noi
for any $f_L \in H^2 (\T_L^2)$.

As a warm up, we first prove the following weighted $L^2$-bound for $S_{L, \eps} (t)$.
\begin{lemma}
\label{LEM:St_L2}
Let $L \geq 1$. Then, for any $0 < \mu \leq 1$, there exists $\Om' \subset \Om$ with full probability measure such that for any $\om \in \Om'$, we have
\begin{align*}
\sup_{\eps \in (0, 1)} \sup_{t \in \R_+} \| S_{L, \eps} (t) f_L \|_{L^2_{- \mu} (\T_L^2)} \les_\om \| f_L \|_{L^2_\mu (\T_L^2)} .
\end{align*}
\end{lemma}

\begin{proof}
Let $w_{L, \eps} (t) = S_{L, \eps} (t) f_L$ and $t \in \R_+$. By the conservation of mass $\mathcal{M}_{L, \eps}$ in \eqref{defM1}, we have
\begin{align}
\int_{\T_L^2} e^{- 2 Y_{L, \eps}} |w_{L, \eps} (t)|^2 dx = \int_{\T_L^2} e^{- 2 Y_{L, \eps}} |f_L|^2 dx .
\label{mass1}
\end{align}

\noi
Therefore, by H\"older's inequality, \eqref{mass1}, and H\"older's inequality along with Lemma~\ref{LEM:Lwei}~(iii), we have
\begin{align*}
\| w_{L, \eps} (t) \|_{L^2_{- \mu} (\T_L^2)}^2 &\leq \| e^{Y_{L, \eps}} \|_{L_{- \mu}^\infty (\T_L^2)}^2 \int_{\T_L^2} e^{- 2 Y_{L, \eps}} |w_{L, \eps} (t)|^2 dx \\
&\les \| e^{Y_{L, \eps}} \|_{L_{- \mu}^\infty (\T_L^2)}^2 \int_{\T_L^2} e^{- 2 Y_{L, \eps}} |f_{L}|^2 dx \\
&\les \| e^{Y_{L, \eps}} \|_{L_{- \mu}^\infty (\T_L^2)}^2 \| e^{- Y_{L, \eps}} \|_{L_{- \mu}^\infty (\T_L^2)}^2 \| f_L \|_{L^2_{\mu} (\T_L^2)}^2 ,
\end{align*}

\noi
so that we conclude the desired estimate using Lemma~\ref{LEM:YregLe}~(iii).
\end{proof}

Let us now prove the following weighted $L^2$-bounds for $S_{L, \eps} (t)$ and $\nb S_{L, \eps} (t)$.
\begin{lemma}
\label{LEM:St_H1}
Let $L \geq 1$. Then, for any $0 < \mu \leq \frac 17$ and $T \geq 1$, there exists $\Om' \subset \Om$ with full probability measure such that for any $\om \in \Om'$, we have the bounds
\begin{align}
\sup_{\eps \in (0, 1)} \| S_{L, \eps} (t) f_L \|_{C_T L_{3\mu}^2 (\T_L^2)} \les_{\om} \| f_L \|_{H_{5\mu}^1 (\T_L^2)} 
\label{L2bd1}
\end{align}

\noi
and
\begin{align}
\sup_{\eps \in (0, 1)} \| \nb S_{L, \eps} (t) f_L \|_{C_T L^2_{- \mu} (\T_L^2)} \les_{\om} \| f_L \|_{H_{5\mu}^1 (\T_L^2)}  .
\label{L2bd2}
\end{align}
\end{lemma}

\begin{proof}
We incorporate the $L$-periodic weights into the computation that goes back to \cite[Lemma~3.1]{DM} (see also \cite[Proposition~4.1 and Proposition~5.1]{DLTV}). Let $w_{L, \eps} (t) = S_{L, \eps} (t) f_L$ and fix $t \in \R_+$. By the conservation of energy $\mathcal{E}_{L, \eps}$ in \eqref{defE1}, we have
\begin{align}
\begin{split}
\frac 12 &\int_{\T_L^2} e^{-2 Y_{L, \eps}} |\nb w_{L, \eps} (t)|^2 dx - \frac 12 \int_{\T_L^2} e^{-2 Y_{L, \eps}} \wt{\wick{|\nb Y_{L, \eps}|^2}} \, |w_{L, \eps} (t)|^2 dx \\
&= \frac 12 \int_{\T_L^2} e^{-2 Y_{L, \eps}} |\nb f_L|^2 dx - \frac 12 \int_{\T_L^2} e^{-2 Y_{L, \eps}} \wt{\wick{|\nb Y_{L, \eps}|^2}} \, |f_L|^2 dx .
\end{split}
\label{L2b1-1}
\end{align}

\noi
By H\"older's inequality and Lemma~\ref{LEM:YregLe}~(iii), we have
\begin{align}
\begin{split}
\| \nb w_{L, \eps} (t) \|_{L_{- \mu}^2 (\T_L^2)}^2 &\leq \| e^{Y_{L, \eps}} \|_{L_{- \mu}^\infty (\T_L^2)}^2 \int_{\T_L^2} e^{- 2 Y_{L, \eps}} |\nb w_{L, \eps} (t)|^2 dx \\
&\les_\om \int_{\T_L^2} e^{- 2 Y_{L, \eps}} |\nb w_{L, \eps} (t)|^2 dx.
\end{split}
\label{L2b1-2}
\end{align}

\noi
Also, by the duality estimate in Lemma~\ref{LEM:dualL} and the product estimate in Lemma~\ref{LEM:prodL}~(ii), we have
\begin{align*}
\bigg| &\int_{\T_L^2} e^{-2 Y_{L, \eps}} \wt{\wick{|\nb Y_{L, \eps}|^2}} \, |w_{L, \eps} (t)|^2 dx \bigg| \\
&\les \big\| \wt{\wick{|\nb Y_{L, \eps}|^2}} \big\|_{\C^{-s}_{- \mu} (\T_L^2)} \big\| e^{-2 Y_{L, \eps}} |w_{L, \eps} (t)|^2 \big\|_{\B^s_{1, 1, \mu} (\T_L^2)} \\
&\les \big\| \wt{\wick{|\nb Y_{L, \eps}|^2}} \big\|_{\C^{-s}_{- \mu} (\T_L^2)} \| e^{-2 Y_{L, \eps}} \|_{\C^{2s}_{- \mu} (\T_L^2)} \| w_{L, \eps} (t) \|_{\B^{3s}_{2, 2, \mu} (\T_L^2)}^2 
\end{align*}

\noi
with $0 < s \leq \frac 16$, so that by further applying Lemma~\ref{LEM:YregLe}~(iii) and (iv) on the last expression, we get
\begin{align}
\bigg| \int_{\T_L^2} e^{-2 Y_{L, \eps}} \wt{\wick{|\nb Y_{L, \eps}|^2}} \, |w_{L, \eps} (t)|^2 dx \bigg| 
\les_\om \| w_{L, \eps} (t) \|_{\B^{3s}_{2, 2, \mu} (\T_L^2)}^2 
\label{L2b1-3}
\end{align}

\noi
To deal with the $\| w_{L, \eps} (t) \|_{\B^{3s}_{2, 2, \mu} (\T_L^2)}$ term, by using the interpolation estimate in Lemma~\ref{LEM:interpL}, Lemma~\ref{LEM:embL}~(v), Lemma~\ref{LEM:L2equiv}, the embedding in Lemma~\ref{LEM:embL}~(iv), and Cauchy's inequality, we have
\begin{align}
\begin{split}
\| w_{L, \eps} (t) \|_{\B^{3s}_{2, 2, \mu} (\T_L^2)}^2 
&\les_\om \| w_{L, \eps} (t) \|_{\B^0_{2, 2, 3 \mu} (\T_L^2)} \| w_{L, \eps} (t) \|_{\B^1_{2, 2, - \mu} (\T_L^2)} \\
&\sim \| w_{L, \eps} (t) \|_{L_{3 \mu}^2 (\T_L^2)} \big( \| w_{L, \eps} (t) \|_{L_{- \mu}^2 (\T_L^2)} + \| \nb w_{L, \eps} (t) \|_{L_{- \mu}^2 (\T_L^2)} \big) \\
&\les \dl^{-1} \| w_{L, \eps} (t) \|_{L_{3 \mu}^2 (\T_L^2)}^2 + \dl \| \nb w_{L, \eps} (t) \|_{L_{- \mu}^2 (\T_L^2)}^2
\end{split}
\label{L2b1-4}
\end{align}

\noi
for any $\dl > 0$. 
Thus, from \eqref{L2b1-1}, \eqref{L2b1-2}, \eqref{L2b1-3}, \eqref{L2b1-4} with $\dl > 0$ sufficiently small, H\"older's inequality along with Lemma~\ref{LEM:Lwei}~(iii), Lemma~\ref{LEM:YregLe}~(iii), and the embedding in Lemma~\ref{LEM:embL}~(iv), we get
\begin{align}
\begin{split}
\| \nb w_{L, \eps} (t) \|_{L_{- \mu}^2 (\T_L^2)}^2 &\les_\om \| w_{L, \eps} (t) \|_{L_{3 \mu}^2 (\T_L^2)}^2 + \| e^{- Y_{L, \eps}} \|_{L_{- \mu}^\infty (\T_L^2)}^2 \| \nb f_L \|_{L_\mu^2 (\T_L^2)}^2 \\
&\quad + \| f_L \|_{L_{3 \mu}^2 (\T_L^2)}^2 + \| \nb f_L \|_{L_{- \mu}^2 (\T_L^2)}^2 \\
&\les_\om \| w_{L, \eps} (t) \|_{L_{3 \mu}^2 (\T_L^2)}^2 + \| f_L \|_{H_{3 \mu}^1 (\T_L^2)}^2 .
\end{split}
\label{L2b1}
\end{align}

We now compute using the equation \eqref{wNLSLe} to get
\begin{align*}
\frac{d}{dt} &\int_{\T_L^2} \jbb{\cdot}_{L, 4\mu}^2 e^{-2 Y_{L, \eps}} |w_{L, \eps} (t)|^2 dx \\
&= 2 \Re \int_{\T_L^2} \jbb{\cdot}_{L, 4\mu}^2 e^{-2 Y_{L, \eps}} \dt w_{L, \eps} (t) \cj{w_{L, \eps} (t)} dx \\
&= 2 \Im \int_{\T_L^2} \jbb{\cdot}_{L, 4\mu}^2 e^{-2 Y_{L, \eps}} \big( \Dl w_{L, \eps} (t) - 2 \nb Y_{L, \eps} \cdot \nb w_{L, \eps} (t) \big) \cj{w_{L, \eps} (t)} dx ,
\end{align*}

\noi
so that by further using an integration by parts and Lemma~\ref{LEM:Lwei}~(vi) and (iii), we obtain
\begin{align}
\begin{split}
\frac{d}{dt} &\int_{\T_L^2} \jbb{\cdot}_{L, 4\mu}^2 e^{-2 Y_{L, \eps}} |w_{L, \eps} (t)|^2 dx \\
&= - 2 \Im \int_{\T_L^2} \nb (\jbb{\cdot}_{L, 4\mu}^2) e^{-2 Y_{L, \eps}} \cdot \nb w_{L, \eps} (t) \cj{w_{L, \eps} (t)} dx \\
&\les \int_{\T_L^2} \jbb{\cdot}_{L, 8 \mu - 1} e^{-2 Y_{L, \eps}} |\nb w_{L, \eps} (t)| |w_{L, \eps} (t)| dx .
\end{split}
\label{L2b2-0}
\end{align}

\noi
By integrating \eqref{L2b2-0} in time, we have
\begin{align*}
&\| e^{- Y_{L, \eps}} w_{L, \eps} (t) \|_{L^2_{4 \mu} (\T_L^2)}^2 \\
&\leq \| e^{- Y_{L, \eps}} f_{L} \|_{L^2_{4 \mu} (\T_L^2)}^2 + \| e^{-2 Y_{L, \eps}} \|_{L^\infty_{- \mu} (\T_L^2)} \int_0^t \| \nb w_{L, \eps} (t') \|_{L^2_{6\mu - 1} (\T_L^2)} \| w_{L, \eps} (t') \|_{L^2_{3 \mu} (\T_L^2)} dt' ,
\end{align*}

\noi
so that by further applying Lemma~\ref{LEM:Lwei}~(iii), H\"older's inequalities, and the embedding in Lemma~\ref{LEM:embL}~(iv) along with $0 < \mu \leq \frac 17$, we get
\begin{align}
\begin{split}
\| e^{- Y_{L, \eps}} w_{L, \eps} (t) \|_{L^2_{4 \mu} (\T_L^2)}^2 
&\les \| e^{- Y_{L, \eps}} \|_{L_{- \mu}^\infty (\T_L^2)}^2 \| f_{L} \|_{L^2_{5 \mu} (\T_L^2)}^2 \\
&\quad + \| e^{-2 Y_{L, \eps}} \|_{L^\infty_{- \mu} (\T_L^2)} \int_0^t  \| \nb w_{L, \eps} (t') \|_{L^2_{- \mu} (\T_L^2)} \| w_{L, \eps} (t') \|_{L^2_{3 \mu} (\T_L^2)} dt'
\end{split}
\label{L2b2-1}
\end{align}

\noi
Thus, from H\"older's inequality, \eqref{L2b2-1}, Lemma~\ref{LEM:YregLe}~(iii), Cauchy-Schwarz inequality, and \eqref{L2b1}, we obtain
\begin{align}
\begin{split}
\| w_{L, \eps} (t) \|_{L_{3 \mu}^2 (\T_L^2)}^2 &\leq \| e^{Y_{L, \eps}} \|_{L^\infty_{- \mu} (\T_L^2)}^2 \| e^{- Y_{L, \eps}} w_{L, \eps} \|_{L_{4 \mu}^2 (\T_L^2)}^2 \\
&\les_\om \| f_L \|_{L_{5 \mu}^2 (\T_L^2)}^2 + \int_0^t \| \nb w_{L, \eps} (t') \|_{L_{- \mu}^2 (\T_L^2)} \| w_{L, \eps} (t') \|_{L_{3 \mu}^2 (\T_L^2)} dt' \\
&\les_\om \| f_L \|_{L_{5 \mu}^2 (\T_L^2)}^2 + t \| f_L \|_{H_{3 \mu}^1 (\T_L^2)}^2 + \int_0^t \| w_{L, \eps} (t') \|_{L_{3 \mu}^2 (\T_L^2)}^2 dt' .
\end{split}
\label{L2b2-2}
\end{align}

\noi
By using Gronwall's inequality and the embedding in Lemma~\ref{LEM:embL}~(iv), we obtain the desired estimate \eqref{L2bd1}. The estimate \eqref{L2bd2} then follows from \eqref{L2b1}, \eqref{L2bd1}, and the embedding in Lemma~\ref{LEM:embL}~(iv).
\end{proof}

We also need a bound for $\Dl S_{L, \eps} (t)$. To achieve this, we first prove the following estimate on the modified energy $\mathcal{F}_{L, \eps}$ defined in \eqref{defF1}.
\begin{lemma}
\label{LEM:Fest}
Let $L \geq 1$. Then, for any $0 < \mu \leq \frac 12$ and $0 < \dl < 1$, there exists $\Om' \subset \Om$ with full probability measure such that for any $\om \in \Om'$ and $0 < \eps < \frac 12$, we have the bound
\begin{align*}
\bigg| \mathcal{F}_{L, \eps} [w] - \int_{\T_L^2} e^{-2 Y_{L, \eps}} |\Dl w|^2 dx \bigg| \les_\om \dl^{-7} |\log \eps|^{17} \| w \|_{L_{8\mu}^2 (\T_L^2)}^2 + \dl \| \Dl w \|_{L_{- \mu}^2 (\T_L^2)}^2 .
\end{align*}
\end{lemma}

\begin{proof}
By H\"older's inequalities along with Lemma~\ref{LEM:Lwei}~(iii), Lemma~\ref{LEM:YregLe}, Lemma~\ref{LEM:YregLe2}, Lemma~\ref{LEM:L2equiv}, and the embeddings in Lemma~\ref{LEM:embL}~(ii), (iv), and (v), we have
\begin{align*}
\bigg| &\int_{\T_L^2} e^{- 2 Y_{L, \eps}} \Dl w \nb Y_{L, \eps} \cdot \nb \cj{w} dx \bigg| \\
&\les \| e^{-2 Y_{L, \eps}} \|_{L_{- \mu}^\infty (\T_L^2)} \| \nb Y_{L, \eps} \|_{L_{- \mu}^r (\T_L^2)} \| \Dl w \|_{L_{- \mu}^2 (\T_L^2)} \| \nb w \|_{L_{3 \mu}^{\frac{2r}{r - 2}} (\T_L^2)} \\
&\les_\om |\log \eps| \| w \|_{\B_{2, 2, - \mu}^2 (\T_L^2)} \| w \|_{\B_{\frac{2r}{r - 2}, 1, 3 \mu}^1 (\T_L^2)} , \\
\bigg| &\int_{\T_L^2} e^{-2 Y_{L, \eps}} ( \nb Y_{L, \eps} \cdot \nb w )^2 dx \bigg| \\
&\les \| e^{-2 Y_{L, \eps}} \|_{L_{- 2 \mu}^\infty (\T_L^2)} \| \nb Y_{L, \eps} \|_{L_{- \mu}^{2r} (\T_L^2)}^2 \| \nb w \|_{L_{2 \mu}^{\frac{2r}{r - 1}} (\T_L^2)}^2 \\
&\les_\om |\log \eps|^2 \| w \|_{\B_{\frac{2r}{r - 1}, 1, 3 \mu}^1 (\T_L^2)}^2 , \\
\bigg| &\int_{\T_L^2} \wt{\wick{|\nb Y_{L, \eps}|^2}} \, w \nb (e^{- 2 Y_{L, \eps}}) \cdot \nb \cj{w} dx \bigg| \\
&\les \| e^{-2 Y_{L, \eps}} \|_{L_{- \mu}^\infty (\T_L^2)} \big\| \wt{\wick{|\nb Y_{L, \eps}|^2}} \big\|_{L_{- \mu}^{2r} (\T_L^2)} \| \nb Y_{L, \eps} \|_{L_{- \mu}^{2r} (\T_L^2)} \| w \|_{L_{\mu}^{\frac{2r}{r - 1}} (\T_L^2)} \| \nb w \|_{L_{2 \mu}^{\frac{2r}{r - 1}} (\T_L^2)} \\
&\les_\om |\log \eps|^3 \| w \|_{\B_{\frac{2r}{r - 1}, 1, 3 \mu}^1 (\T_L^2)}^2 , \\
\bigg| &\int_{\T_L^2} e^{- 2 Y_{L, \eps}} \wt{\wick{|\nb Y_{L, \eps}|^2}} \, \Dl w \cj{w} dx \bigg| \\
&\les \| e^{-2 Y_{L, \eps}} \|_{L_{- \mu}^\infty (\T_L^2)} \big\| \wt{\wick{|\nb Y_{L, \eps}|^2}} \big\|_{L_{- \mu}^{r} (\T_L^2)} \| \Dl w \|_{L_{- \mu}^2 (\T_L^2)} \| w \|_{L_{3 \mu}^{\frac{2r}{r - 2}} (\T_L^2)} \\
&\les_\om |\log \eps|^2 \| w \|_{\B_{2, 2, - \mu}^2 (\T_L^2)} \| w \|_{\B_{\frac{2r}{r - 2}, 1, 3 \mu}^1 (\T_L^2)} , \\
\bigg| &\int_{\T_L^2} e^{- 2 Y_{L, \eps}} \big( \wt{\wick{|\nb Y_{L, \eps}|^2}} \big)^2 |w|^2 dx \bigg| \\
&\les \| e^{-2 Y_{L, \eps}} \|_{L_{- \mu}^\infty (\T_L^2)} \big\| \wt{\wick{|\nb Y_{L, \eps}|^2}} \big\|_{L_{- \mu}^{2r} (\T_L^2)}^2 \| w \|_{L_{2 \mu}^{\frac{2r}{r - 1}} (\T_L^2)}^2 \\
&\les_\om |\log \eps|^4 \| w \|_{\B_{\frac{2r}{r - 1}, 1, 3 \mu}^1 (\T_L^2)}^2 
\end{align*}

\noi
for any $r > 2$ large satisfying $\mu r > 2$. Thus, the above estimates yield
\begin{align}
\begin{split}
\bigg| &\mathcal{F}_{L, \eps} [w] - \int_{T_L^2} e^{- 2 Y_{L, \eps}} |\Dl w|^2 dx \bigg| \\
&\les_{\om} |\log \eps|^4 \Big( \| w \|_{\B_{2, 2, - \mu}^2 (\T_L^2)} \| w \|_{\B_{\frac{2r}{r - 2}, 1, 3 \mu}^1 (\T_L^2)} + \| w \|_{\B_{\frac{2r}{r - 1}, 1, 3 \mu}^1 (\T_L^2)}^2 \Big) .
\end{split}
\label{Fest1}
\end{align}

\noi
From the embeddings in Lemma~\ref{LEM:embL}~(iii) and (ii) and the interpolation estimate in Lemma~\ref{LEM:interpL}, we have
\begin{align}
\| w \|_{\B_{\frac{2r}{r - 2}, 1, 3 \mu}^1 (\T_L^2)} + \| w \|_{\B_{\frac{2r}{r - 1}, 1, 3 \mu}^1 (\T_L^2)} &\les \| w \|_{\B_{2, 2, 3 \mu}^{1 + \frac{4}{r}} (\T_L^2)} \les \| w \|_{\B_{2, 2, \mu'}^0 (\T_L^2)}^{\frac 12 - \frac{2}{r}} \| w \|_{\B_{2, 2, - \mu}^2 (\T_L^2)}^{\frac 12 + \frac{2}{r}} ,
\label{Fest2}
\end{align}

\noi
where
\begin{align*}
\mu' = \frac{\frac 72 + \frac 2r}{\frac 12 - \frac 2r} \mu .
\end{align*}

\noi
Thus, the desired estimate follows from combining \eqref{Fest1} and \eqref{Fest2}, taking $r$ to be sufficiently large, and using Young's inequality along with the fact that
\begin{align*}
\| w \|_{\B_{2, 2, \mu'}^0 (\T_L^2)} \les \| w \|_{L^2_{8 \mu} (\T_L^2)}
\end{align*}

\noi
and
\begin{align*}
\| w \|_{\B_{2, 2, - \mu}^2 (\T_L^2)} \sim \| w \|_{L^2_{- \mu} (\T_L^2)} + \| \Dl w \|_{L^2_{- \mu} (\T_L^2)} \les \| w \|_{L^2_{8 \mu} (\T_L^2)} + \| \Dl w \|_{L^2_{- \mu} (\T_L^2)},
\end{align*}

\noi
which follow from Lemma~\ref{LEM:L2equiv}, Lemma~\ref{LEM:embL}~(v), and Lemma~\ref{LEM:embL}~(iv).
\end{proof}

We now show the following bound for $\Dl S_{L, \eps} (t)$.
\begin{lemma}
\label{LEM:St_H2}
Let $L \geq 1$. Then, for any $T \geq 1$ and $\mu > 0$ sufficiently small, there exists $\Om' \subset \Om$ with full probability measure such that for any $\om \in \Om'$ and $0 < \eps < \frac 12$, we have the bound
\begin{align*}
\| \Dl S_{L, \eps} (t) f_L \|_{C_T L_{- \mu}^2 (\T_L^2)} \les_{\om} |\log \eps|^{5} \| f_L \|_{H_{15 \mu}^2 (\T_L^2)} .
\end{align*}
\end{lemma}

\begin{proof}
Let $w_{L, \eps} (t) = S_{L, \eps} (t) f_L$ and fix $t \in \R_+$. From H\"older's inequality, Lemma~\ref{LEM:YregLe}~(iii), and the conservation of $\mathcal{F}_{L, \eps}$ in \eqref{defF1}, we have
\begin{align}
\begin{split}
\| \Dl w_{L, \eps} (t) \|_{L_{- \mu}^2 (\T_L^2)}^2 
&\leq \| e^{Y_{L, \eps}} \|_{L_{- \mu}^\infty (\T_L^2)}^2 \int_{\T_L^2} e^{- 2 Y_{L, \eps}} |\Dl w_{L, \eps} (t)|^2 dx \\
&\les_\om  \bigg| \int_{\T_L^2} e^{- 2 Y_{L, \eps}} |\Dl w_{L, \eps} (t)|^2 dx - \mathcal{F}_{L, \eps} [w_{L, \eps} (t)] \bigg| \\
&\quad + \bigg| \int_{\T_L^2} e^{- 2 Y_{L, \eps}} |\Dl f_{L}|^2 dx - \mathcal{F}_{L, \eps} [f_L] \bigg| + \int_{\T_L^2} e^{- 2 Y_{L, \eps}} |\Dl f_{L}|^2 dx .
\end{split}
\label{H2b1}
\end{align}

\noi
Then, from Lemma~\ref{LEM:Fest}, we get
\begin{align}
\begin{split}
\bigg| &\int_{\T_L^2} e^{- 2 Y_{L, \eps}} |\Dl w_{L, \eps} (t)|^2 dx - \mathcal{F}_{L, \eps} [w_{L, \eps} (t)] \bigg| \\
&\les_\om \dl^{-7} |\log \eps|^{17} \| w_{L, \eps} (t) \|_{L_{8 \mu}^2 (\T_L^2)}^2 + \dl \| \Dl w_{L, \eps} (t) \|_{L_{- \mu}^2 (\T_L^2)}^2 
\end{split}
\label{H2b2}
\end{align}

\noi
for any $0 < \dl < 1$, and
\begin{align}
\begin{split}
\bigg| &\int_{\T_L^2} e^{- 2 Y_{L, \eps}} |\Dl f_{L}|^2 dx - \mathcal{F}_{L, \eps} [f_L] \bigg| \\
&\les_\om |\log \eps|^{17} \| f_L \|_{L_{8 \mu}^2 (\T_L^2)}^2 + \| \Dl f_L \|_{L_{- \mu}^2 (\T_L^2)}^2 .
\end{split}
\label{H2b3}
\end{align}

\noi
Moreover, by H\"older's inequality along with Lemma~\ref{LEM:Lwei}~(iii), we have
\begin{align}
\int_{\T_L^2} e^{- 2 Y_{L, \eps}} |\Dl f_{L}|^2 dx 
\les \| e^{-Y_{L, \eps}} \|_{L^\infty_{- \mu} (\T_L^2)}^2 \| \Dl f_L \|_{L^2_{\mu} (\T_L^2)}^2
\les_\om \| \Dl f_L \|_{L^2_{\mu} (\T_L^2)}^2 .
\label{H2b4}
\end{align}

\noi
Combining \eqref{H2b1}, \eqref{H2b2}, \eqref{H2b3}, and \eqref{H2b4}, we obtain
\begin{align*}
\begin{split}
\| \Dl w_{L, \eps} (t) \|_{L_{- \mu}^2 (\T_L^2)}^2 
&\les_\om \dl^{-7} |\log \eps|^{17} \| w_{L, \eps} (t) \|_{L_{8 \mu}^2 (\T_L^2)}^2 + \dl \| \Dl w_{L, \eps} (t) \|_{L_{- \mu}^2 (\T_L^2)}^2  \\
&\quad + |\log \eps|^{17} \| f_L \|_{L_{8 \mu}^2 (\T_L^2)}^2 + \| \Dl f_L \|_{L_{\mu}^2 (\T_L^2)}^2 
\end{split}
\end{align*}

\noi
for any $0 < \dl < 1$. Thus, by taking $\dl = \dl (\om) > 0$ to be sufficiently small and using the $L^2$-bound in Lemma~\ref{LEM:St_H1} along with the embedding in Lemma~\ref{LEM:embL}~(iv), we obtain the desired estimate.
\end{proof}

As a consequence of the bounds established above, we obtain the following useful bound to be used in an essential way in the next subsection.
\begin{lemma}
\label{LEM:St_Bs}
Let $L \geq 1$ and $\mu_0 > 0$. Then, for any $T \geq 1$, $0 \leq s < 2$, and $\gamma > 0$ satisfying $s + \gamma \leq 2$, there exists small $\mu_1 = \mu_1 (\mu_0, s, \gamma) > 0$ such that for any $0 < \mu \leq \mu_1$, there exists $\Om' \subset \Om$ with full probability measure such that for any $\om \in \Om'$ and $0 < \eps < \frac 12$, we have the bound
\begin{align*}
\| S_{L, \eps} (t) f_L \|_{C_T \B_{2, 2, \mu}^{s} (\T_L^2)} \les_\om |\log \eps|^5 \| f_L \|_{\B_{2, 2, \mu_0}^{s + \gamma} (\T_L^2)} \les |\log \eps|^5 \| f_L \|_{H^2_{\mu_0} (\T_L^2)} .
\end{align*}
\end{lemma}

\begin{proof} 
The second bound follows directly from Lemma~\ref{LEM:embL}~(v) and Lemma~\ref{LEM:L2equiv}, and so we focus on proving the first bound. 
Given a dyadic $N \geq 1$, using the interpolation estimate in Lemma~\ref{LEM:interpL}, and Lemma~\ref{LEM:L2equiv}, we have
\begin{align*}
\| S_{L, \eps} (t) \Dl_N^L f_L \|_{C_T \B_{2, 2, \mu}^0 (\T_L^2)} 
&\les \| S_{L, \eps} (t) \Dl_N^L f_L \|_{C_T \B_{2, 2, - \mu}^0 (\T_L^2)}^{1 - \frac{\gamma}{2}} \| S_{L, \eps} (t) \Dl_N^L f_L \|_{C_T \B_{2, 2, \frac{4 - \gamma}{\gamma} \mu}^0 (\T_L^2)}^{\frac{\gamma}{2}} \\
&\sim \| S_{L, \eps} (t) \Dl_N^L f_L \|_{C_T L^2_{- \mu} (\T_L^2)}^{1 - \frac{\gamma}{2}} \| S_{L, \eps} (t) \Dl_N^L f_L \|_{C_T L^2_{\frac{4 - \gamma}{\gamma} \mu} (\T_L^2)}^{\frac{\gamma}{2}} 
\end{align*}

\noi
for any $\gamma > 0$ and $\mu = \mu (\mu_0, \gamma) > 0$ sufficiently small, so that by further using Lemma~\ref{LEM:St_L2}, Lemma~\ref{LEM:St_H1}, Lemma~\ref{LEM:embL}~(v) and (iv), and Young's inequality, we get
\begin{align}
\begin{split}
\| S_{L, \eps} (t) \Dl_N^L f_L \|_{C_T \B_{2, 2, \mu}^0 (\T_L^2)} 
&\les_\om \sum_{M = \frac{N}{2}, N, 2N} M^{\frac{\gamma}{2}} \| \Dl_N^L f_L \|_{L_\mu^2 (\T_L^2)}^{1 - \frac{\gamma}{2}} \| \Dl_M^L f_L \|_{L_{\frac{5 (4 - \gamma)}{3 \gamma} \mu}^2 (\T_L^2)}^{\frac{\gamma}{2}} \\
&\les \sum_{M = \frac{N}{2}, N, 2N} M^{\frac{\gamma}{2}} \| \Dl_M^L f_L \|_{L^2_{\mu_0} (\T_L^2)} .
\end{split}
\label{Bs1}
\end{align}

\noi 
Then, from the Littlewood-Paley decomposition \eqref{idDNL} and the interpolation estimate in Lemma~\ref{LEM:interpL}, we have
\begin{align}
\begin{split}
\| S_{L, \eps} (t) f_L \|_{C_T \B_{2, 2, \mu}^s (\T_L^2)} &\leq \sum_{\substack{N \geq 1 \\ \text{dyadic}}} \| S_{L, \eps} (t) \Dl_N^L f_L \|_{C_T \B_{2, 2, \mu}^s (\T_L^2)} \\
&\les \sum_{\substack{N \geq 1 \\ \text{dyadic}}} \| S_{L, \eps} (t) \Dl_N^L f_L \|_{C_T \B_{2, 2, \mu'}^0 (\T_L^2)}^{1 - \frac{s}{2}} \| S_{L, \eps} (t) \Dl_N^L f_L \|_{C_T \B_{2, 2, - \mu}^2 (\T_L^2)}^{\frac{s}{2}} ,
\end{split}
\label{Bs2}
\end{align}

\noi
where
\begin{align*}
\mu' = \frac{1 + \frac{s}{2}}{1 - \frac{s}{2}} \mu .
\end{align*}

\noi
From the equivalence in Lemma~\ref{LEM:embL}~(v),  Lemma~\ref{LEM:L2equiv}, Lemma~\ref{LEM:St_L2}, Lemma~\ref{LEM:St_H2}, and the estimate in Lemma~\ref{LEM:embL}~(v), we have
\begin{align}
\begin{split}
\| S_{L, \eps} (t) \Dl_N^L f_L \|_{C_T \B_{2, 2, - \mu}^2 (\T_L^2)}
&\les \| S_{L, \eps} (t) \Dl_N^L f_L \|_{C_T L^2_{- \mu} (\T_L^2)} + \| \Dl S_{L, \eps} (t) \Dl_N^L f_L \|_{C_T L^2_{- \mu} (\T_L^2)} \\
&\les |\log \eps|^5 \big( \| \Dl_N^L f_L \|_{L_{\mu_0}^2 (\T_L^2)} + \| \Dl_N^L \Dl f_L \|_{L_{\mu_0}^2 (\T_L^2)} \big) \\
&\les |\log \eps|^5 \sum_{M = \frac{N}{2}, N, 2N} M^2 \| \Dl_M^L f_L \|_{L_{\mu_0}^2 (\T_L^2)} .
\end{split}
\label{Bs3}
\end{align}

\noi
Thus, from \eqref{Bs2}, \eqref{Bs1} with $\mu = \mu (\mu_0, s, \gamma) > 0$ sufficiently small, \eqref{Bs3}, and Young's inequality, we obtain
\begin{align*}
\| &S_{L, \eps} (t) f_L \|_{C_T \B_{2, 2, \mu}^s (\T_L^2)} \\
&\les_\om |\log \eps|^5 \sum_{\substack{N \geq 1 \\ \text{dyadic}}} \bigg( \sum_{M = \frac{N}{2}, N, 2N} M^{\frac{\gamma}{2}} \| \Dl_M^L f_L \|_{L_{\mu_0}^2 (\T_L^2)} \bigg)^{1 - \frac{s}{2}} \bigg( \sum_{M = \frac{N}{2}, N, 2N} M^2 \| \Dl_M^L f_L \|_{L_{\mu_0}^2 (\T_L^2)} \bigg)^{\frac{s}{2}} \\
&\les |\log \eps|^5 \sum_{\substack{N \geq 1 \\ \text{dyadic}}} \sum_{M = \frac{N}{2}, N, 2N} M^{s + (1 - \frac{s}{2}) \frac{\gamma}{2}} \| \Dl_M^L f_L \|_{L_{\mu_0}^2 (\T_L^2)} ,
\end{align*}

\noi
so that the desired estimate follows from the Cauchy-Schwarz inequality in dyadic $N \geq 1$.
\end{proof}

\subsection{Convergence estimate for the $L$-periodic linear solution}
\label{SUB:GWP1}

In this subsection, we Theorem~\ref{THM:GWP} in the case when $\ld = 0$, which corresponds to global well-posedness of the $L$-periodic linear DAM \eqref{wNLSAndL}. In fact, this follows directly from the following proposition.
\begin{proposition}
\label{PROP:wLe_conv}
Let $L \geq 1$, $\mu_0 > 0$, and $w_{0, L} \in H_{\mu_0}^2 (\T_L^2)$. Given $0 < \eps < \frac 12$, let $w_{L, \eps}$ be the global-in-time solution to the mollified $L$-periodic linear DAM \eqref{wNLSLe} with $w_{L, \eps}|_{t = 0} = w_{0, L}$. Then, there exists a process $w_L$ such that, for any $T \geq 1$ and $0 \leq s < 2$, there exist $\dl > 0$ and $\Om' \subset \Om$ with full probability measure such that for any $\om \in \Om'$ and $0 < \eps < \frac 12$, we have the 
difference estimate
\begin{align}
\| w_{L, \eps} - w_L \|_{C_T H^s (\T_L^2)} \les_{\om, T} \eps^\dl \| w_{0, L} \|_{H_{\mu_0}^2 (\T_L^2)}.
\label{wLHs_diff}
\end{align}

\noi
The process $w_L$ is the unique global-in-time solution to the $L$-periodic linear DAM \eqref{wNLSAndL} with $w_L |_{t = 0} = w_{0, L}$ in $C(\R_+; H^s (\T_L^2))$ for any $1 < s < 2$. 
\end{proposition}

\begin{proof}
Let $\Om' \subset \Om$ be the event with full probability measure such that all the estimates in Section~\ref{SEC:sto} and Section~\ref{SEC:convL1} hold (with some parameters to be chosen below) and we fix $\om \in \Om'$.

Let $0 < \eps_2 < \eps_1 < \frac 12$ and $r_{L, \eps_1, \eps_2} = w_{L, \eps_1} - w_{L, \eps_2}$. Then, from \eqref{wNLSLe}, we see that $r_{L, \eps_1, \eps_2}$ satisfies the following equation with zero initial data:
\begin{align*}
i \dt r_{L, \eps_1, \eps_2} &= \Dl r_{L, \eps_1, \eps_2} - 2 \nb Y_{L, \eps_1} \cdot \nb r_{L, \eps_1, \eps_2} + \wt{\wick{ |\nb Y_{L, \eps_1}|^2 }} \, r_{L, \eps_1, \eps_2} \\
&\quad + \big( \wt{\wick{ |\nb Y_{L, \eps_1}|^2 }} - \wt{\wick{ |\nb Y_{L, \eps_2}|^2 }} \big) w_{L, \eps_2} - 2 (\nb Y_{L, \eps_1} - \nb Y_{L, \eps_2}) \cdot \nb w_{L, \eps_2} .
\end{align*}

\noi
Let $t \in [0, T]$. By using the equation for $r_{L, \eps_1, \eps_2}$, we compute that
\begin{align}
\frac 12 \frac{d}{dt} \int_{\T_L^2} e^{- 2 Y_{L, \eps_1}} |r_{L, \eps_1, \eps_2} (t)|^2 dx = \textup{I}_1 + \textup{I}_2 ,
\label{wLe0}
\end{align}

\noi
where
\begin{align*}
\textup{I}_1 &= -2 \Im \int_{\T_L^2} e^{- 2 Y_{L, \eps_1}} (\nb Y_{L, \eps_1} - \nb Y_{L, \eps_2}) \cdot \nb w_{L, \eps_2} (t) \cj{r_{L, \eps_1, \eps_2} (t)} dx , \\
\textup{I}_2 &= \Im \int_{\T_L^2} e^{- 2 Y_{L, \eps_1}} \big( \wt{\wick{ |\nb Y_{L, \eps_1}|^2 }} - \wt{\wick{ |\nb Y_{L, \eps_2}|^2 }} \big) w_{L, \eps_2} (t) \cj{r_{L, \eps_1, \eps_2} (t)} dx .
\end{align*}

\noi
Let $0 < s_0 < \frac{1}{100}$, $0 < \dl_0 < s_0$, and small $\mu > 0$ be fixed. For $\textup{I}_1$, using the duality estimate in  Lemma~\ref{LEM:dualL}, Lemma~\ref{LEM:YregLe}~(i), the product estimates (Lemma~\ref{LEM:prodL}~(ii)), Lemma~\ref{LEM:YregLe}~(iii), and Lemma~\ref{LEM:embL}~(v), we have
\begin{align}
\begin{split}
|\textup{I}_1| &\les \| \nb Y_{L, \eps_1} - \nb Y_{L, \eps_2} \|_{\C_{- \mu}^{- s_0} (\T_L^2)} \|  e^{- 2 Y_{L, \eps_1}} \nb w_{L, \eps_2} (t) \cj{r_{L, \eps_1, \eps_2} (t)} \|_{\B^{s_0}_{1, 1, \mu} (\T_L^2)} \\
&\les_{\om} \eps_1^{\dl_0} \| e^{- 2 Y_{L, \eps_1}} \|_{\C_{- \mu}^{3 s_0} (\T_L^2)} \| \nb w_{L, \eps_2} (t) \|_{\B_{2, 2, \mu}^{3 s_0} (\T_L^2)} \| r_{L, \eps_1, \eps_2} (t) \|_{\B_{2, 2, \mu}^{2 s_0} (\T_L^2)} \\
&\les_{\om} \eps_1^{\dl_0} \| w_{L, \eps_2} (t) \|_{\B_{2, 2, \mu}^{1 + 3 s_0} (\T_L^2)} \big( \| w_{L, \eps_1} (t) \|_{\B_{2, 2, \mu}^{2 s_0} (\T_L^2)} + \| w_{L, \eps_2} (t) \|_{\B_{2, 2, \mu}^{2 s_0} (\T_L^2)} \big) .
\end{split}
\label{wLe1}
\end{align}

\noi
For $\textup{I}_2$, using the duality estimate in Lemma~\ref{LEM:dualL}, Lemma~\ref{LEM:YregLe}~(ii) and (iv), the product estimates (Lemma~\ref{LEM:prodL}~(ii)), and Lemma~\ref{LEM:YregLe}~(iii), we have
\begin{align}
\begin{split}
|\textup{I}_2| &\les \big\| \wt{ \wick{|\nb Y_{L, \eps_1}|^2} } - \wt{ \wick{|\nb Y_{L, \eps_2}|^2} } \big\|_{\C_{- \mu}^{- s_0} (\T_L^2)} \|  e^{- 2 Y_{L, \eps_1}} v_{L, \eps_2} (t) \cj{r_{L, \eps_1, \eps_2} (t)} \|_{\B^{s_0}_{1, 1, \mu} (\T_L^2)} \\
&\les_{\om} \eps_1^{\dl_0} \| e^{- 2 Y_{L, \eps_1}} \|_{\C_{- \mu}^{3 s_0} (\T_L^2)} \| v_{L, \eps_2} (t) \|_{\B_{2, 2, \mu}^{3 s_0} (\T_L^2)} \| r_{L, \eps_1, \eps_2} (t) \|_{\B_{2, 2, \mu}^{2 s_0} (\T_L^2)} \\
&\les_{\om} \eps_1^{\dl_0} \| v_{L, \eps_2} (t) \|_{\B_{2, 2, \mu}^{3 s_0} (\T_L^2)} \big( \| v_{L, \eps_1} (t) \|_{\B_{2, 2, \mu}^{2 s_0} (\T_L^2)} + \| v_{L, \eps_2} (t) \|_{\B_{2, 2, \mu}^{2 s_0} (\T_L^2)} \big) .
\end{split}
\label{wLe2}
\end{align}

\noi
Combining \eqref{wLe0}, \eqref{wLe1}, and \eqref{wLe2} and using Lemma~\ref{LEM:St_Bs} (with $\mu > 0$ sufficiently small), we obtain
\begin{align*}
\sup_{t \in [0, T]} \int_{\T_L^2} e^{- 2 Y_{L, \eps_1}} |r_{L, \eps_1, \eps_2} (t)|^2 dx \les_{\om} \eps_1^{\dl_0} |\log \eps_2|^{10} T \| w_{0, L} \|_{H_{\mu_0}^2 (\T_L^2)}^{2} ,
\end{align*}

\noi
so that by H\"older's inequality and Lemma~\ref{LEM:YregLe}~(iii), we obtain
\begin{align}
\| w_{L, \eps_1} - w_{L, \eps_2} \|_{C_T L_{- \mu}^2 (\T_L^2)}^2 \les_{\om} \eps_1^{\dl_0} |\log \eps_2|^{10} T \| w_{0, L} \|_{H_{\mu_0}^2 (\T_L^2)}^{2} .
\label{wLe_diff}
\end{align}

Let us now take $\eps_1 = 2^{-k}$ and $\eps_2 = 2^{- (k + 1)}$ for $k \in \N$. Then, by using \eqref{wLe_diff}, we get
\begin{align}
\| w_{L, 2^{-k}} - w_{L, 2^{- (k + 1)}} \|_{C_T L_{- \mu}^2 (\T_L^2)}^2 \les_{\om} 2^{- \dl_0 k} (k + 1)^{10} T \| w_{0, L} \|_{H_{\mu_0}^2 (\T_L^2)}^{2} .
\label{wLe_diff2}
\end{align}

\noi
Let $0 \leq s < s_1 < 2$. 
By \eqref{Hs_equi} and the interpolation estimate in Lemma~\ref{LEM:interpL}, we have
\begin{align*}
\| &w_{L, 2^{-k}} - w_{L, 2^{- (k + 1)}} \|_{C_T H^s (\T_L^2)} \\
&\sim \| w_{L, 2^{-k}} - w_{L, 2^{- (k + 1)}} \|_{C_T \B_{2, 2}^s (\T_L^2)} \\
&\les \| w_{L, 2^{-k}} - w_{L, 2^{- (k + 1)}} \|_{C_T \B_{2, 2, - \mu}^0 (\T_L^2)}^{1 - \frac{s}{s_1}} \| w_{L, 2^{-k}} - w_{L, 2^{- (k + 1)}} \|_{C_T \B_{2, 2, \mu'}^{s_1} (\T_L^2)}^{\frac{s}{s_1}} 
\end{align*}

\noi
with
\begin{align*}
\mu' = \frac{1 - \frac{s}{s_1}}{\frac{s}{s_1}} \mu ,
\end{align*}

\noi
so that by further using Lemma~\ref{LEM:L2equiv}, Lemma~\ref{LEM:St_Bs} (with $\mu > 0$ sufficiently small), and \eqref{wLe_diff2}, we get
\begin{align}
\begin{split}
\| &w_{L, 2^{-k}} - w_{L, 2^{- (k + 1)}} \|_{C_T H^s (\T_L^2)} \\
&\les_{\om, T} \| w_{L, 2^{-k}} - w_{L, 2^{- (k + 1)}} \|_{C_T L_{- \mu}^2 (\T_L^2)}^{1 - \frac{s}{s_1}} (k + 1)^{\frac{5 s}{s_1}} \| w_{0, L} \|_{H_{\mu_0}^2 (\T_L^2)}^{\frac{s}{s_1}} \\
&\les_{\om, T} 2^{- \dl_1 k} \| w_{0, L} \|_{H_{\mu_0}^2 (\T_L^2)} 
\end{split}
\label{wLe_conv1}
\end{align}

\noi
for some $\dl_1 > 0$.
This shows that $\{ w_{L, 2^{-k}} \}_{k \in \N}$ is almost surely a Cauchy sequence in $C([0, T]; H^s (\T_L^2))$ and so converges to some process $w_L$ with values in $C([0, T]; H^s (\T_L^2))$. Similarly, we have
\begin{align}
\sup_{\eps \in (2^{- (k + 1)}, 2^{-k}]} \| w_{L, \eps} - w_{L, 2^{-k}} \|_{C_T H^s (\T_L^2)} \les_{\om, T} 2^{- \dl_2 k} \| w_{0, L} \|_{H_{\mu_0}^2 (\T_L^2)} 
\label{wLe_conv2}
\end{align}

\noi
for some $\dl_2 > 0$, so that the whole sequence $\{ w_{L, \eps} \}_{\eps \in (0, \frac 12)}$ converges to $w_L$ almost surely in $C([0, T]; H^s (\T_L^2))$ as $\eps \to 0$. The 
difference estimate \eqref{wLHs_diff} then follows immediately from \eqref{wLe_conv1} and \eqref{wLe_conv2}. The uniqueness part follows from the same analysis as in previous works \cite{DW, TzV23} and so we omit details.
\end{proof}

\subsection{Large torus convergence of mollified $L$-periodic linear solutions}

In this subsection, we show that the solution of the mollified $L$-periodic linear DAM \eqref{wNLSLe} converges to that of the following mollified version of the linear DAM \eqref{wNLSAnd} on $\R^2$:
\begin{align}
\begin{cases}
i \dt w_\eps = \Dl w_\eps - 2 \nb Y_\eps \cdot \nb w_\eps + \wt{\wick{ |\nb Y_\eps|^2 }} \, w_\eps \\
w_\eps |_{t = 0} = w_0 ,
\end{cases}
\label{wNLSe}
\end{align}

\noi
where $Y_\eps$ is as defined in \eqref{Ye} and 
\begin{align}
\wt{\wick{ |\nb Y_\eps|^2 }} \deff \, \wick{ |\nb Y_\eps|^2 } - \, \varphi *_{\R^2} \xi_\eps
\label{Ye22}
\end{align}

\noi
with $\wick{|\nb Y_\eps|^2}$ being defined in \eqref{Ye2} and $\xi_\eps$ being defined in \eqref{xie}. From \cite[Theorem~1.1]{DLTV}, we know that if $w_0 \in H^2_{\mu} (\R^2)$ for any $\mu > 0$, then \eqref{wNLSe} admits a unique global-in-time solution.

We first mention the following result regarding the weighted $H^1 (\R^2)$-bound of $w_\eps$, which follows from \cite[Lemma~4.3]{DLTV}.
\begin{lemma}
\label{LEM:SeH1}
Let $\mu_0 > 0$ and $w_0 \in H^2_{\mu_0} (\R^2)$. Given $0 < \eps < \frac 12$, let $w_\eps$ be the global-in-time solution to the mollified linear DAM \eqref{wNLSe} with $w_\eps |_{t = 0} = w_0$. Then, for any $T \geq 1$ and $\mu > 0$ sufficiently small, there exists $\Om' \subset \Om$ will full probability measure such that for any $\om \in \Om'$ and $0 < \eps < \frac 12$, we have the bound
\begin{align*}
\| w_\eps \|_{C_T H_\mu^1 (\R^2)} \les_{\om, T} \| w_0 \|_{H^2_{\mu_0} (\R^2)}.
\end{align*}
\end{lemma}

We now show the following proposition regarding the convergence of the solution $w_{L, \eps}$ of \eqref{wNLSLe} to the solution $w_\eps$ of \eqref{wNLSe}. 
We recall the function $\s_R$ in \eqref{sigR}.
\begin{proposition}
\label{PROP:wLewe}
Let $L \geq 1$, $\mu_0 > 0$, $w_0 \in H_{\mu_0}^2 (\R^2)$, and $w_{0, L} \in H_{\mu_0}^2 (\T_L^2)$. Given $0 < \eps < \frac 12$, let $w_\eps$ be the global-in-time solution to the mollified linear DAM \eqref{wNLSe} with $w_\eps |_{t = 0} = w_0$ and let $w_{L, \eps}$ be the global-in-time solution to the mollified $L$-periodic linear DAM \eqref{wNLSLe} with $w_{L, \eps} |_{t = 0} = w_{0, L}$. Then, for any $T \geq 1$, $1 \leq R \leq L$, $\mu > 0$ sufficiently small, and $0 < b \leq 4$, there exists $\Om' \subset \Om$ with full probability measure such that for any $\om \in \Om'$ and $0 < \eps < \frac 12$, we have the estimate
\begin{align*}
\| &\s_R^{-2} (w_{L, \eps} - w_\eps) \|_{C_T L_{- \mu}^2 (\R^2)} \\
&\les_{\om, T} \| \s_R^{-2} (w_{0, L} - w_0) \|_{L_\mu^2 (\R^2)} + (R^{-1} + R^{b} L^{- \frac{b}{8}} |\log \eps|) \big( \| w_{0, L} \|_{H^2_{\mu_0} (\T_L^2)} + \| w_0 \|_{H^2_{\mu_0} (\R^2)} \big) .
\end{align*}
\end{proposition}

\begin{proof}
Given $L \geq 1$ and $0 < \eps < \frac 12$, we let $r_{L, \eps} = w_{\eps} - w_{L, \eps}$. Then, we see that $r_{L, \eps}$ satisfies the following equation:
\begin{align*}
i \dt r_{L, \eps} &= \Dl r_{L, \eps} - 2 \nb Y_{L, \eps} \cdot \nb r_{L, \eps} + \wt{\wick{ |\nb Y_{L, \eps}|^2 }} \, r_{L, \eps}  - 2 (\nb Y_\eps - \nb Y_{L, \eps}) \cdot \nb w_{\eps} \\
&\quad + \big( \wt{\wick{|\nb Y_\eps|^2}} - \wt{\wick{|\nb Y_{L, \eps}|^2}} \big) w_{\eps} .
\end{align*}

\noi
Let $\Omega' \subset \Omega$ be the event with full probability measure such that all (already established) estimates in Section~\ref{SEC:sto} and Section~\ref{SEC:convL1} hold and we fix $\om \in \Omega'$.
Note that from \eqref{sigR1}, Lemma~\ref{LEM:Lwei}~(vi), (iii), and (i), Lemma~\ref{LEM:sig}, Lemma~\ref{LEM:L2equiv}, and Lemma~\ref{LEM:embL}~(v), we have
\begin{align*}
\big\| \s_R^{-1} \jbb{\cdot}_{L, \mu} w_{L, \eps} \big\|_{C_T H^1 (\R^2)} 
&\les \big\| \s_R^{-1} \jbb{\cdot}_{L, \mu} w_{L, \eps} \big\|_{C_T L^2 (\R^2)} + \big\| \s_R^{-1} \jbb{\cdot}_{L, \mu} \nb w_{L, \eps} \big\|_{C_T L^2 (\R^2)} \\
&\les \| w_{L, \eps} \|_{C_T L^2_\mu (\T_L^2)} + \| \nb w_{L, \eps} \|_{C_T L^2_\mu (\T_L^2)} \\
&\les \| w_{L, \eps} \|_{C_T \B_{2, 2, \mu}^1 (\T_L^2)} 
\end{align*}

\noi
for any $\mu > 0$ sufficiently small, so that by applying Lemma~\ref{LEM:St_Bs} on the last expression, we get
\begin{align}
\big\| \s_R^{-1} \jbb{\cdot}_{L, \mu} w_{L, \eps} \big\|_{C_T H^1 (\R^2)} 
\les_{\om, T} \| w_{0, L} \|_{H_{\mu_0}^2 (\T_L^2)} .
\label{RwLe_bdd}
\end{align}

Let $t \in [0, T]$. By using the equation and integration by parts, we obtain
\begin{align}
\frac 12 \frac{d}{dt} \int_{\R^2} \s_R^{-4} e^{- 2 Y_{L, \eps}} |r_{L, \eps} (t)|^2 dx 
= \textup{I}_1 + \textup{I}_2 + \textup{I}_3 ,
\label{wL0}
\end{align}

\noi  
where
\begin{align*}
\textup{I}_1 &= - \Im \int_{\R^2} e^{-2 Y_{L, \eps}} \nb (\s_R^{-4}) \cdot \nb r_{L, \eps} (t) \cj{r_{L, \eps} (t)} dx , \\
\textup{I}_2 &= - 2 \Im \int_{\R^2} \s_R^{-4} e^{- 2 Y_{L, \eps}} (\nb Y_\eps - \nb Y_{L, \eps}) \cdot \nb w_{\eps} (t) \cj{r_{L, \eps} (t)} dx , \\
\textup{I}_3 &= \Im \int_{\R^2} \s_R^{-4} e^{- 2 Y_{L, \eps}} \big( \wt{\wick{|\nb Y_\eps|^2}} - \wt{\wick{|\nb Y_{L, \eps}|^2}} \big) w_{\eps} (t) \cj{r_{L, \eps} (t)} dx .
\end{align*}

\noi
For $\textup{I}_1$, by \eqref{sigR1}, Cauchy's inequality, H\"older's inequality along with Lemma~\ref{LEM:Lwei}~(iii), Lemma~\ref{LEM:YregLe}~(iii), \eqref{RwLe_bdd}, and Lemma~\ref{LEM:SeH1}, we have
\begin{align}
\begin{split}
|\textup{I}_1| &\les R^{-2} \int_{\R^2} \s_R^{-4} e^{-2 Y_{L, \eps}} |\nb r_{L, \eps} (t)|^2 dx + \int_{\R^2} \s_R^{-4} e^{- 2 Y_{L, \eps}} |r_{L, \eps} (t)|^2 dx \\
&\leq R^{-2} \big\| \jbb{\cdot}_{L, -\mu} e^{- Y_{L, \eps}} \big\|_{L^\infty (\R^2)}^2 \big\| \s_R^{-1} \jbb{\cdot}_{L, \mu} \nb r_{L, \eps} (t) \big\|_{L^2 (\R^2)}^2 + \int_{\R^2} \s_R^{-4} e^{- 2 Y_{L, \eps}} |r_{L, \eps} (t)|^2 dx \\
&\les_{\om, T} R^{-2} \big( \| w_{0, L} \|_{H^2_{\mu_0} (\T_L^2)} + \| w_0 \|_{H_{\mu_0}^2 (\R^2)} \big)^2 + \int_{\R^2} \s_R^{-4} e^{- 2 Y_{L, \eps}} |r_{L, \eps} (t)|^2 dx 
\end{split}
\label{wL1}
\end{align}

\noi
for $\mu > 0$ sufficiently small. For $\textup{I}_2$, by H\"older's inequality along with Lemma~\ref{LEM:Lwei}~(iii), Lemma~\ref{LEM:YregLe}~(iii), Proposition~\ref{PROP:YconvL}~(i), Lemma~\ref{LEM:SeH1}, and Sobolev's inequality, we have
\begin{align}
\begin{split}
|\textup{I}_2| &\les \big\| \jbb{\cdot}_{L, - \mu} e^{-2 Y_{L, \eps}} \big\|_{L^\infty (\R^2)} \| \s_R^{-1} (\nb Y_\eps - \nb Y_{L, \eps}) \|_{L^q (\R^2)} \\
&\quad \times \| \s_R^{-1} \nb w_\eps (t) \|_{L^2 (\R^2)} \big\| \s_R^{-1} \jbb{\cdot}_{L, \mu} r_{L, \eps} (t) \big\|_{L^{\frac{2 q}{q - 2}} (\R^2)} \\
&\les_{\om, T} R^{b} L^{- \frac{b}{2}} |\log \eps| \| w_0 \|_{H_{\mu_0}^2 (\R^2)} \big\| \s_R^{-1} \jbb{\cdot}_{L, \mu} r_{L, \eps} (t) \big\|_{H^1 (\R^2)} 
\end{split}
\label{wL2}
\end{align}

\noi 
for any $0 < b \leq 4$ and $2 < q < \infty$ satisfying $b q > 4$.
For $\textup{I}_3$, by H\"older's inequality along with Lemma~\ref{LEM:Lwei}~(iii), Lemma~\ref{LEM:YregLe}~(iii), Proposition~\ref{PROP:YconvL}~(iv), Lemma~\ref{LEM:SeH1}, and Sobolev's inequality, we have
\begin{align}
\begin{split}
|\textup{I}_3| &\les \big\| \jbb{\cdot}_{L, - \mu} e^{-2 Y_{L, \eps}} \big\|_{L^\infty (\R^2)} \Big\| \s_R^{-2} \big( \wt{\wick{|\nb Y_\eps|^2}} - \wt{\wick{|\nb Y_{L, \eps}|^2}} \big) \Big\|_{L^q (\R^2)} \\
&\quad \times \| \s_R^{-1} w_\eps (t) \|_{L^2 (\R^2)} \big\| \s_R^{-1} \jbb{\cdot}_{L, \mu} r_{L, \eps} (t) \big\|_{L^{\frac{2q}{q - 2}} (\R^2)} \\
&\les_{\om, T} R^{2b} L^{- \frac{b}{2} + \frac{2}{q}} |\log \eps|^2 \| w_0 \|_{H^2_{\mu_0} (\R^2)} \big\| \s_R^{-1} \jbb{\cdot}_{L, \mu} r_{L, \eps} (t) \big\|_{H^1 (\R^2)} .
\end{split}
\label{wL3}
\end{align}

\noi
Combining \eqref{wL0}, \eqref{wL1},  \eqref{wL2}, and \eqref{wL3} and using \eqref{RwLe_bdd} and Lemma~\ref{LEM:SeH1}, we obtain
\begin{align*}
&\frac{d}{dt} \int_{\R^2} \s_R^{-4} e^{- 2 Y_{L, \eps}} |r_{L, \eps} (t)|^2 dx \\
&\les_{\om, T} \big( R^{-2} + R^{2b} L^{- \frac{b}{2} + \frac{2}{q}} |\log \eps|^2 \big) \big( \| w_{0, L} \|_{H^2_{\mu_0} (\T_L^2)} + \| w_0 \|_{H^2_{\mu_0} (\R^2)} \big)^2 + \int_{\R^2} \s_R^{-4} e^{- 2 Y_{L, \eps}} |r_{L, \eps} (t)|^2 dx.
\end{align*}

\noi
Using Gronwall's inequality and taking $q > 2$ to be large enough, we get 
\begin{align}
\begin{split}
&\sup_{t \in [0, T]} \int_{\R^2} \s_R^{-4} e^{- 2 Y_{L, \eps}} |r_{L, \eps} (t)|^2 dx \\
&\les_{\om, T} \int_{\R^2} \s_R^{-4} e^{- 2 Y_{L, \eps}} |w_0 - w_{0, L}|^2 dx \\
&\quad + (R^{-2} + R^{2b} L^{- \frac{b}{4}} |\log \eps|^2) \big( \| w_{0, L} \|_{H^2_{\mu_0} (\T_L^2)} + \| w_0 \|_{H^2_{\mu_0} (\R^2)} \big)^2
\end{split}
\label{wL4}
\end{align}

\noi
for any $0 < b \leq 4$. 
Therefore, by H\"older's inequality, \eqref{wL4}, and Lemma~\ref{LEM:YregLe}~(iii) along with the fact that $\jb{\cdot}^{- \mu} \les \jbb{\cdot}_{L, \mu}^{-1} \les \jbb{\cdot}_{L, - \mu}$ from Lemma~\ref{LEM:Lwei}~(ii) and (iii), we obtain
\begin{align*}
\| \s_R^{-2} r_{L, \eps} \|_{C_T L_{- \mu}^2 (\R^2)} &\leq \| \s_R^{-2} e^{- Y_{L, \eps}} r_{L, \eps} \|_{C_T L^2 (\R^2)} \| e^{Y_{L, \eps}} \|_{L_{- \mu}^\infty (\R^2)} \\
&\les_{\om, T} \| \s_R^{-2} e^{- Y_{L, \eps}} (w_0 - w_{0, L}) \|_{C_T L^2 (\R^2)} \\
&\quad + (R^{-1} + R^{b} L^{- \frac{b}{8}} |\log \eps|) \big( \| w_{0, L} \|_{H^2_{\mu_0} (\T_L^2)} + \| w_0 \|_{H^2_{\mu_0} (\R^2)} \big) ,
\end{align*}

\noi
so that by repeating the above procedure again, we get
\begin{align*}
\| \s_R^{-2} r_{L, \eps} \|_{C_T L_{- \mu}^2 (\R^2)} 
&\les_{\om, T} \| \s_R^{-2} (w_0 - w_{0, L}) \|_{C_T L_\mu^2 (\R^2)} \\
&\quad + (R^{-1} + R^{b} L^{- \frac{b}{8}} |\log \eps|) \big( \| w_{0, L} \|_{H^2_{\mu_0} (\T_L^2)} + \| w_0 \|_{H^2_{\mu_0} (\R^2)} \big) ,
\end{align*}

\noi
which gives the desired estimate.
\end{proof}

\subsection{Large torus convergence of the $L$-periodic linear solutions}

In this subsection, we show the proof of Theorem~\ref{THM:convL} in the case when $\ld = 0$, which corresponds to the convergence of the solution $w_L$ of the $L$-periodic linear DAM \eqref{wNLSAndL} to the solution $w$ of the linear DAM \eqref{wNLSAnd} on $\R^2$.

We first mention the following result regarding the convergence of the solution $w_\eps$ of the mollified linear DAM \eqref{wNLSe} on $\R^2$ to the solution $w$ of the linear DAM \eqref{wNLSAnd} on $\R^2$. This is covered by \cite[Theorem~1.2]{DLTV}.
\begin{lemma}
\label{LEM:wGWP}
Let $\mu_0 > 0$ and $w_0 \in H_{\mu_0}^2 (\R^2)$. Given $0 < \eps < \frac 12$, let $w_\eps$ be the global-in-time solution to the mollified linear DAM \eqref{wNLSe} with  $w_\eps|_{t = 0} = w_0$ and let $w$ be the global-in-time solution to the linear DAM \eqref{wNLSAnd} with $w|_{t = 0} = w_0$. Then, for any $T \geq 1$, $0 \leq s < 2$, and $\mu > 0$ sufficiently small, there exist $\dl > 0$ and $\Om' \subset \Om$ with full probability measure such that for any $\om \in \Om'$ and $0 < \eps < \frac 12$, we have the bound
\begin{align*}
\| w_\eps \|_{C_T H^s_\mu (\R^2)} \les_{\om, T} \| w_0 \|_{H^2_{\mu_0} (\R^2)}
\end{align*}

\noi
and the difference estimate
\begin{align*}
\| w_\eps - w \|_{C_T H^s_\mu (\R^2)} \les_{\om, T} \eps^\dl \| w_0 \|_{H^2_{\mu_0} (\R^2)}.
\end{align*}
\end{lemma}

We are now ready to show the large torus limit of the linear DAM.
\begin{proof}[Proof of Theorem~\ref{THM:convL} in the case $\ld = 0$]
Given $L \geq 1$, we let $w_0 = v_0 \in H_{\mu_0}^2 (\R^2)$ and $w_{0, L} = v_{0, L} \in H_{\mu_1}^2 (\T_L^2)$ be given by (see Lemma~\ref{LEM:per})
\begin{align*}
w_{0, L} = \sum_{k \in \Z^2} w_0 (\cdot + Lk) ,
\end{align*}

\noi
where $0 < \mu_1 < \mu_0 - 1$. From Lemma~\ref{LEM:perw}, we have
\begin{align}
\| \s_R^{-1} (w_{0, L} - w_0) \|_{L_{\mu_1}^2 (\R^2)} \les (L^{- \mu_0 + \mu_1} + R^{1 + 2 \mu_1} L^{-1 - \mu_1}) \| w_0 \|_{L_{\mu_0}^2 (\R^2)} .
\label{w0L_diff}
\end{align}

Given $0 < \eps < \frac 12$, we let $w_{L, \eps}$ and $w_L$ be the global-in-time solutions to \eqref{wNLSLe} and \eqref{wNLSAndL}, respectively, both with initial data $w_{0, L}$. We also let $w_\eps$ and $w$ be the global-in-time solutions to \eqref{wNLSe} and \eqref{wNLSAnd}, respectively, both with initial data $w_0$. Let $\Om' \subset \Om$ be the event with full probability measure such that all estimates in Section~\ref{SEC:sto} and Section~\ref{SEC:convL1} hold for some parameters to be fixed later and we fix $\om \in \Om'$.

Given any $T \geq 1$, $0 \leq s < 2$, and a bounded open set $U \subset \R^2$, we have
\begin{align}
\| w_L - w \|_{C_T H^s (U)} \leq \| w_L - w_{L, \eps} \|_{C_T H^s (U)} + \| w_{L, \eps} - w_\eps \|_{C_T H^s (U)} + \| w_\eps - w \|_{C_T H^s (U)}.
\label{wLlim0}
\end{align}

\noi
From Lemma~\ref{LEM:Hsloc}, Proposition~\ref{PROP:wLe_conv}, and Lemma~\ref{LEM:per}, we have
\begin{align}
\| w_L - w_{L, \eps} \|_{C_T H^s (U)} \les_{U} \| w_L - w_{L, \eps} \|_{C_T H^s (\T_L^2)} \les_{\om, T} \eps^\dl \| w_{0, L} \|_{H_{\mu_1}^2 (\T_L^2)} \les \eps^\dl \| w_{0} \|_{H_{\mu_0}^2 (\T_L^2)}
\label{wLlim1}
\end{align}

\noi
for some $\dl > 0$. Also, from Lemma~\ref{LEM:wGWP}, we have
\begin{align}
\| w_\eps - w \|_{C_T H^s (U)} \leq \| w_\eps - w \|_{C_T H^s_\mu (\R^2)} \les_{\om, T} \eps^\dl \| w_0 \|_{H_{\mu_0}^2 (\R^2)}
\label{wLlim2}
\end{align}

\noi
for some $\dl > 0$. 

We now estimate $\| w_{L, \eps} - w_\eps \|_{C_T H^s (U)}$. Let $\phi \in C_c^\infty (\R^2)$ be such that $\phi \equiv 1$ on $U$ and let $U' \subset \R^2$ be a bounded open set such that $\supp \phi \subset U'$. Let $0 \leq s < s' < s'' < 2$. Note that for any $f \in H^{s''} (U')$, we let $f'$ be an extension of $f$ on $\R^2$ such that $\| f' \|_{H^{s''} (\R^2)} \leq 2 \| f \|_{H^{s''} (U')}$, so that from the product estimate in Lemma~\ref{LEM:prod}, we have 
\begin{align}
\| \phi f \|_{H^{s'} (\R^2)} = \| \phi f' \|_{H^{s'} (\R^2)} \les \| f' \|_{H^{s''} (\R^2)} \les \| f \|_{H^{s''} (U')} .
\label{prodU}
\end{align}

\noi
Thus, by interpolation and \eqref{prodU}, we have
\begin{align}
\begin{split}
\| w_{L, \eps} - w_\eps \|_{C_T H^s (U)} &\leq \| \phi (w_{L, \eps} - w_\eps) \|_{C_T H^s (\R^2)} \\
&\les \| \phi (w_{L, \eps} - w_\eps) \|_{C_T L^2 (\R^2)}^{1 - \frac{s}{s'}} \| \phi (w_{L, \eps} - w_\eps) \|_{L_T^\infty H^{s'} (\R^2)}^{\frac{s}{s'}} \\
&\les \| w_{L, \eps} - w_\eps \|_{C_T L^2 (U')}^{1 - \frac{s}{s'}} \| w_{L, \eps} - w_\eps \|_{L_T^\infty H^{s''} (U')}^{\frac{s}{s'}} .
\end{split}
\label{wLlim3-0}
\end{align}

\noi
From Lemma~\ref{LEM:Hsloc}, \eqref{Hs_equi},  Lemma~\ref{LEM:St_Bs} along with Lemma~\ref{LEM:embL}~(iv), and Lemma~\ref{LEM:wGWP}, we get
\begin{align}
\begin{split}
\| w_{L, \eps} - w_\eps \|_{L_T^\infty H^{s''} (U')}
&\les \| w_{L, \eps} \|_{L_T^\infty H^{s''} (\T_L^2)} + \| w_\eps \|_{L_T^\infty H^{s''} (\R^2)} \\
&\sim \| w_{L, \eps} \|_{L_T^\infty \B_{2, 2}^{s''} (\T_L^2)} + \| w_\eps \|_{L_T^\infty H^{s''} (\R^2)} \\
&\les |\log \eps|^5 \big( \| w_{0, L} \|_{H^2_{\mu_1} (\T_L^2)} + \| w_{0} \|_{H^2_{\mu_0} (\R^2)} \big) .
\end{split}
\label{wLlim3-1}
\end{align}

\noi
Also, from Proposition~\ref{PROP:wLewe}, we have
\begin{align}
\begin{split}
\| w_{L, \eps} - w_\eps \|_{C_T L^2 (U')} &\les_{U'} \| \s_R^{-2} (w_{L, \eps} - w_\eps) \|_{C_T L_{- \mu}^2 (\R^2)} \\
&\les_{\om, T} \| \s_R^{-2} (w_{0, L} - w_0) \|_{C_T L_\mu^2 (\R^2)} \\
&\quad + (R^{-1} + R^{b} L^{- \frac{b}{8}} |\log \eps|) \big( \| w_{0, L} \|_{H^2_{\mu_1} (\T_L^2)} + \| w_0 \|_{H_{\mu_1}^2 (\R^2)} \big) 
\end{split}
\label{wLlim3-2}
\end{align}

\noi
for any $\mu > 0$ sufficiently small and $0 < b \leq 4$. Thus, from \eqref{wLlim3-0}, \eqref{wLlim3-1}, and \eqref{wLlim3-2} along with \eqref{w0L_diff}, Lemma~\ref{LEM:per}, and the choice $R \sim L^{\frac{1}{16}}$, we get
\begin{align}
\| w_{L, \eps} - w_\eps \|_{C_T H^s (U)} \les L^{- \kappa} |\log \eps|^{5} \| w_0 \|_{H_{\mu_0}^2 (\R^2)}
\label{wLlim3}
\end{align}

\noi
for some $\kappa > 0$.

Therefore, combining \eqref{wLlim0}, \eqref{wLlim1}, \eqref{wLlim2}, and \eqref{wLlim3}, we get
\begin{align*}
\| w_L - w \|_{C_T H^s (U)} &\les_{\om, T, U} (\eps^\dl + L^{- \kappa} |\log \eps|^{5} ) \| w_0 \|_{H_{\mu_0}^2 (\R^2)} .
\end{align*}

\noi
for some $\dl > 0$ and some $\kappa > 0$. Therefore, by choosing
\begin{align*}
 \eps \sim L^{- \frac{1}{\dl}} ,
\end{align*}

\noi
we obtain the desired convergence result.
\end{proof}

\section{Large torus limit of the nonlinear dispersive Anderson model}
\label{SEC:convL2}

In this section, we show Theorem~\ref{THM:GWP} and Theorem~\ref{THM:convL} in the case when $\ld > 0$, which correspond to global well-posedness of the $L$-periodic nonlinear DAM \eqref{vNLSAndL} and the large torus limit of the nonlinear DAM, respectively.

\subsection{Energy estimates for mollified $L$-periodic nonlinear solutions}
\label{SUB:mod}

We consider the following mollified $L$-periodic nonlinear DAM:
\begin{align}
\begin{cases}
i \dt v_{L, \eps} = \Dl v_{L, \eps} - 2 \nb Y_{L, \eps} \cdot \nb v_{L, \eps} + \wt{\wick{|\nb Y_{L, \eps}|^2}} \, v_{L, \eps} - \ld e^{- (p - 1) Y_{L, \eps}} |v_{L, \eps}|^{p - 1} v_{L, \eps} \\
v_{L, \eps} |_{t = 0} = v_{0, L} ,
\end{cases}
\label{vNLSLe}
\end{align}

\noi
where $Y_{L, \eps}$ is defined in \eqref{YLe_fs} and $\wt{\wick{|\nb Y_{L, \eps}|^2}}$ is defined in \eqref{YLe2}. To see that \eqref{vNLSLe} is globally well-posed for $v_{0, L} \in H^2 (\T_L^2)$, we see that $u_{L, \eps} \deff e^{- Y_{L, \eps}} v_{L, \eps}$ satisfies
\begin{align*}
\begin{cases}
i \dt u_{L, \eps} = \Dl u_{L, \eps} + \xi_{L, \eps} u_{L, \eps} - C_{L, \eps} u_{L, \eps} - \ld |u_{L, \eps}|^{p - 1} u_{L, \eps} \\
u_{L, \eps} |_{t = 0} = e^{- Y_{L, \eps}} v_{0, L} ,
\end{cases}
\end{align*}

\noi
whose unique global solution exists almost surely in $C(\R_+; H^2 (\T_L^2))$ (see \cite{Tsu, BGTz}). This shows that \eqref{vNLSLe} has a unique global solution in $C (\R_+; H^2 (\T_L^2))$.

Let us recall the mass $\mathcal{M}_{L, \eps}$ defined in \eqref{defM1}
and define the energy
\begin{align}
\mathcal{E}_{L, \ld, \eps} [v] \deff \int_{\T_L^2} \Big( \frac 12 e^{-2 Y_{L, \eps}} |\nb v|^2 - \frac 12 e^{- 2 Y_{L, \eps}} \wt{\wick{|\nb Y_{L, \eps}|^2}} \, |v|^2 + \frac{\ld}{p + 1} e^{- (p + 1) Y_{L, \eps}} |v|^{p + 1} \Big) dx .
\label{defE}
\end{align}

\noi
It is not difficult to compute that if $v_{L, \eps}$ satisfies \eqref{vNLSLe} then we have the conservation of mass and energy:
\begin{align*}
\frac{d}{dt} \mathcal{M}_{L, \eps} [v_{L, \eps} (t)] &= 0, \\
\frac{d}{dt} \mathcal{E}_{L, \ld, \eps} [v_{L, \eps} (t)] &= 0.
\end{align*}

We also recall the modified energy $\mathcal{F}_{L, \eps}$ defined in \eqref{defF1} and define the following modified energies (introduced in \cite{TzV23}):
\begin{align*}
\mathcal{G}_{L, \eps} [v] &\deff - \int_{\T_L^2} e^{- (p + 1) Y_{L, \eps}} |\nb v|^2 |v|^{p - 1} dx - 2 \Re \int_{\T_L^2} e^{- (p + 1) Y_{L, \eps}} v \nb ( |v|^{p - 1} ) \cdot \nb \cj{v} dx \\
&\quad + \frac{p - 1}{4} \int_{\T_L^2} e^{- (p + 1) Y_{L, \eps}} | \nb ( |v|^2 ) |^2 |v|^{p - 3} dx \\
&\quad + \frac{2}{p + 1} \int_{\T_L^2} e^{- (p + 1) Y_{L, \eps}} \wt{\wick{|\nb Y_{L, \eps}|^2}} \, |v|^{p + 1} dx \\
&\quad + 2 (p - 1) \Re \int_{\T_L^2} e^{- (p + 1) Y_{L, \eps}} |v|^{p - 1} v \nb Y_{L, \eps} \cdot \nb \cj{v} dx,
\end{align*}
\begin{align*}
\mathcal{H}_{L, \eps} [v] &\deff - \int_{\T_L^2} e^{- (p + 1) Y_{L, \eps}} \dt (|v|^{p - 1}) |\nb v|^2 dx - 2 \Re \int_{\T_L^2} e^{- (p + 1) Y_{L, \eps}} \dt v \nb (|v|^{p - 1}) \cdot \nb \cj{v} dx \\
&\quad - \frac{p - 1}{4} \int_{\T_L^2} e^{- (p + 1) Y_{L, \eps}} \dt ( |v|^{p - 3} ) | \nb ( |v|^2 ) |^2 dx \\
&\quad + 2 (p - 1) \Re \int_{\T_L^2} e^{- (p + 1) Y_{L, \eps}} \dt ( |v|^{p - 1} v ) \nb Y_{L, \eps} \cdot \nb \cj{v} dx. 
\end{align*}

\noi
From \cite[Proposition~4.1]{TzV23}, we have the following conservation of modified energies under the flow of \eqref{vNLSLe}:
\begin{align}
\frac{d}{dt} \big( \mathcal{F}_{L, \eps} [v_{L, \eps} (t)] - \ld \mathcal{G}_{L, \eps} [v_{L, \eps} (t)] \big) = - \ld \mathcal{H}_{L, \eps} [v_{L, \eps} (t)].
\label{FGH}
\end{align}

Let us first show the following estimates for $\mathcal{G}$ and $\mathcal{H}$.

\begin{lemma}
\label{LEM:Gest}
Let $p \geq 2$, $\ld > 0$, and $L \geq 1$. Then, for any $\mu > 0$ sufficiently small and $\nu > 0$ sufficiently small, there exists $\Om' \subset \Om$ with full probability measure such that for any $\om \in \Om'$ and $0 < \eps < \frac 12$, we have the bound
\begin{align*}
|\mathcal{G}_{L, \eps} [v]| \les_{\om} |\log \eps|^2 \big( \| v \|_{L_{- \mu}^2 (\T_L^2)} + \| \Dl v \|_{L_{- \mu}^2 (\T_L^2)} \big)^{1 + 2 \nu} \| v \|_{\B_{2, 2, 4 \mu}^{1 - \nu} (\T_L^2)}^{p - 2 \nu} .
\end{align*}
\end{lemma}

\begin{proof}
By H\"older's inequalities along with Lemma~\ref{LEM:Lwei}~(iii) and (iv), Lemma~\ref{LEM:YregLe}, Lemma~\ref{LEM:YregLe2}, Lemma~\ref{LEM:L2equiv}, and the embeddings in Lemma~\ref{LEM:embL}~(ii) and (v), we have
\begin{align*}
\bigg| &\int_{\T_L^2} e^{- (p + 1) Y_{L, \eps}} |\nb v|^2 |v|^{p - 1} dx \bigg| \\
&\les \| e^{- (p + 1) Y_{L, \eps}} \|_{L_{- (p - 1) \mu_1}^\infty (\T_L^2)} \| \nb v \|_{L^4 (\T_L^2)}^2 \| v \|_{L_{\mu_1}^{2 p - 2} (\T_L^2)}^{p - 1} \\
&\les_\om \| v \|_{\B_{4, 1}^{1} (\T_L^2)}^2  \| v \|_{\B_{2p - 2, 1, \mu_1}^{0} (\T_L^2)}^{p - 1} , \\
\bigg| &\int_{\T_L^2} e^{- (p + 1) Y_{L, \eps}} v \nb ( |v|^{p - 1} ) \cdot \nb \cj{v} dx \bigg| \\
&\les \| e^{- (p + 1) Y_{L, \eps}} \|_{L_{- (p - 1) \mu_1}^\infty (\T_L^2)} \| \nb v \|_{L^4 (\T_L^2)}^2 \| v \|_{L_{\mu_1}^{2 p - 2} (\T_L^2)}^{p - 1} \\
&\les_\om \| v \|_{\B_{4, 1}^{1} (\T_L^2)}^2  \| v \|_{\B_{2p - 2, 1, \mu_1}^{0} (\T_L^2)}^{p - 1} , \\
\bigg| &\int_{\T_L^2} e^{- (p + 1) Y_{L, \eps}} | \nb ( |v|^2 ) |^2 |v|^{p - 3} dx \bigg| \\
&\les \| e^{- (p + 1) Y_{L, \eps}} \|_{L_{- (p - 1) \mu_1}^\infty (\T_L^2)} \| \nb v \|_{L^4 (\T_L^2)}^2 \| v \|_{L_{\mu_1}^{2 p - 2} (\T_L^2)}^{p - 1} \\
&\les_\om \| v \|_{\B_{4, 1}^{1} (\T_L^2)}^2  \| v \|_{\B_{2p - 2, 1, \mu_1}^{0} (\T_L^2)}^{p - 1} , \\
\bigg| &\int_{\T_L^2} e^{- (p + 1) Y_{L, \eps}} \wt{\wick{|\nb Y_{L, \eps}|^2}} \, |v|^{p + 1} dx \bigg| \\
&\les \| e^{- (p + 1) Y_{L, \eps}} \|_{L_{- (p - \frac 32) \mu_1}^\infty (\T_L^2)} \big\| \wt{\wick{|\nb Y_{L, \eps}|^2}} \big\|_{L_{- \frac{\mu_1}{2}}^r (\T_L^2)} \| v \|_{L^4 (\T_L^2)}^2 \| v \|_{L_{\mu_1}^{\frac{2r (p - 1)}{r - 2}} (\T_L^2)}^{p - 1} \\
&\les_\om |\log \eps|^2 \| v \|_{\B_{4, 1}^{1} (\T_L^2)}^2 \| v \|_{\B_{\frac{2r (p - 1)}{r - 2}, 1, \mu_1}^{0} (\T_L^2)}^{p - 1} , \\
\bigg| &\int_{\T_L^2} e^{- (p + 1) Y_{L, \eps}} |v|^{p - 1} v \nb Y_{L, \eps} \cdot \nb \cj{v} dx \bigg| \\
&\les \| e^{- (p + 1) Y_{L, \eps}} \|_{L_{- (p - \frac 32) \mu_1}^\infty (\T_L^2)} \| \nb Y_{L, \eps} \|_{L_{- \frac{\mu_1}{2}}^r (\T_L^2)} \| \nb v \|_{L^4 (\T_L^2)} \| v \|_{L^4 (\T_L^2)} \| v \|_{L_{\mu_1}^{\frac{2r (p - 1)}{r - 2}} (\T_L^2)}^{p - 1} \\
&\les_\om |\log \eps| \| v \|_{\B_{4, 1}^{1} (\T_L^2)}^2 \| v \|_{\B_{\frac{2r (p - 1)}{r - 2}, 1, \mu_1}^{0} (\T_L^2)}^{p - 1}
\end{align*}

\noi
for any $\mu_1 > 0$ sufficiently small and $r > 2$ satisfying $\mu_1 r > 4$. Thus, the above estimates yield
\begin{align}
|\mathcal{G}_{L, \eps} [v]| 
\les_{\om} |\log \eps|^2 \| v \|_{\B_{4, 1}^{1} (\T_L^2)}^2 \Big( \| v \|_{\B_{2p - 2, 1, \mu_1}^{0} (\T_L^2)}^{p - 1} + \| v \|_{\B_{\frac{2r (p - 1)}{r - 2}, 1, \mu_1}^{0} (\T_L^2)}^{p - 1} \Big) .
\label{GLe1}
\end{align}

\noi
By the embeddings in Lemma~\ref{LEM:embL}~(i) and (iii), the interpolation estimate in Lemma~\ref{LEM:interpL}, Lemma~\ref{LEM:embL}~(v), and Lemma~\ref{LEM:L2equiv}, we have
\begin{align}
\begin{split}
\| v \|_{\B_{4, 1}^1 (\T_L^2)}
&\les \| v \|_{\B_{2, 2}^{\frac 32 + \frac{\nu}{2}} (\T_L^2)} \\
&\les \| v \|_{\B_{2, 2, -\mu_2}^2 (\T_L^2)}^{\frac 12 + \nu} \| v \|_{\B_{2, 2, \mu_1}^{(\frac 12 - \frac{3\nu}{2}) / (\frac 12 - \nu)} (\T_L^2)}^{\frac 12 - \nu}  \\
&\les \big( \| v \|_{L^2_{- \mu_2} (\T_L^2)} + \| \Dl v \|_{L^2_{- \mu_2} (\T_L^2)} \big)^{\frac 12 + \nu} \| v \|_{\B_{2, 2, \mu_1}^{1 - \nu} (\T_L^2)}^{\frac 12 - \nu}
\end{split}
\label{GLe2}
\end{align}

\noi
for any $0 < \nu < \frac 12$, where
\begin{align*}
\mu_2 = \frac{\frac 12 - \nu}{\frac 12 + \nu} \mu_1 .
\end{align*}

\noi
Also, with $q = 2p - 2$ or $q = \frac{2r (p - 1)}{r - 2}$, by the embeddings in Lemma~\ref{LEM:embL} (i) and (iii), we get
\begin{align}
\| v \|_{\B_{q, 1, \mu_1}^0 (\T_L^2)} \les \| v \|_{\B_{q, 2, \mu_1}^{\nu} (\T_L^2)} \les \| v \|_{\B_{2, 2, \mu_1}^{1 - \nu} (\T_L^2)}
\label{GLe3}
\end{align}

\noi
as long as $\nu > 0$ is sufficiently small. We make the restriction $0 < \nu < \frac 14$, so that $\mu_2 > \frac 14 \mu_1$.
With $\mu = \frac 14 \mu_1$, the desired estimate then follows from \eqref{GLe1}, \eqref{GLe2} along with the embedding in Lemma~\ref{LEM:embL}~(iv), and \eqref{GLe3}.
\end{proof}

\begin{lemma}
\label{LEM:Hest}
Let $p \geq 2$, $\ld > 0$, $L \geq 1$, and $T \geq 1$. Given $0 < \eps < \frac 12$, let $v_{L, \eps}$ be a solution to the mollified $L$-periodic nonlinear DAM \eqref{vNLSLe}. Then, for any $\mu > 0$ sufficiently small and $\nu > 0$ sufficiently small, there exists $\Om' \subset \Om$ with full probability measure such that for any $\om \in \Om'$ and $0 < \eps < \frac 12$, we have the bound
\begin{align*}
\int_0^T |\mathcal{H}_{L, \eps} [v_{L, \eps} (t)]| dt &\les_{\om} |\log \eps|^3 \Big( 1 + \| \Dl v_{L, \eps} \|_{L_T^\infty L_{- \mu}^2 (\T_L^2)}^{1 + 2 (p - 2) \nu} \Big) \\
&\quad \times \| v_{L, \eps} \|_{L_T^2 \B_{4, 1, 2  \mu}^1 (\T_L^2)}^2 \Big( 1 + \| v_{L, \eps} \|_{L_T^\infty \B_{2, 2, \mu}^{1 - \nu} (\T_L^2)}^{2p - 2} \Big) .
\end{align*}
\end{lemma}

\begin{proof}
By the diamagnetic inequality $|\partial |f|| \leq |\partial f|$, H\"older's inequalities along with Lemma~\ref{LEM:Lwei}~(iii), Lemma~\ref{LEM:YregLe}~(iii), Lemma~\ref{LEM:YregLe2}, and the embeddings in Lemma~\ref{LEM:embL}~(ii) and (v), we obtain
\begin{align*}
\int_0^T &\bigg| \int_{\T_L^2} e^{- (p + 1) Y_{L, \eps}} \dt (|v|^{p - 1}) |\nb v|^2 dx \bigg| dt \\
&\les \| e^{- (p + 1) Y_{L, \eps}} \|_{L_{- 2 \mu}^\infty (\T_L^2)} \| \dt v_{L, \eps} \|_{L^\infty_T L_{- 2 \mu}^2 (\T_L^2)} \| \nb v_{L, \eps} \|_{L^2_T L_{2 \mu}^4 (\T_L^2)}^2 \| v_{L, \eps} \|_{L^\infty_T L^\infty (\T_L^2)}^{p - 2} \\
&\les_\om \| \dt v_{L, \eps} \|_{L^\infty_T L_{- 2 \mu}^2 (\T_L^2)} \| v_{L, \eps} \|_{L^2_T \B_{4, 1, 2 \mu}^1 (\T_L^2)}^2 \| v_{L, \eps} \|_{L^\infty_T \B^0_{\infty, 1} (\T_L^2)}^{p - 2} , \\
\int_0^T &\bigg| \int_{\T_L^2} e^{- (p + 1) Y_{L, \eps}} \dt v \nb (|v|^{p - 1}) \cdot \nb \cj{v} dx \bigg| dt \\
&\les \| e^{- (p + 1) Y_{L, \eps}} \|_{L_{- 2 \mu}^\infty (\T_L^2)} \| \dt v_{L, \eps} \|_{L^\infty_T L_{- 2 \mu}^2 (\T_L^2)} \| \nb v_{L, \eps} \|_{L^2_T L_{2 \mu}^4 (\T_L^2)}^2 \| v_{L, \eps} \|_{L^\infty_T L^\infty (\T_L^2)}^{p - 2} \\
&\les_\om \| \dt v_{L, \eps} \|_{L^\infty_T L_{- 2 \mu}^2 (\T_L^2)} \| v_{L, \eps} \|_{L^2_T \B_{4, 1, 2 \mu}^1 (\T_L^2)}^2 \| v_{L, \eps} \|_{L^\infty_T \B^0_{\infty, 1} (\T_L^2)}^{p - 2} , \\
\int_0^T &\bigg| \int_{\T_L^2} e^{- (p + 1) Y_{L, \eps}} \dt ( |v|^{p - 3} ) | \nb ( |v|^2 ) |^2 dx \bigg| dt \\
&\les \| e^{- (p + 1) Y_{L, \eps}} \|_{L_{- 2 \mu}^\infty (\T_L^2)} \| \dt v_{L, \eps} \|_{L^\infty_T L_{- 2 \mu}^2 (\T_L^2)} \| \nb v_{L, \eps} \|_{L^2_T L_{2 \mu}^4 (\T_L^2)}^2 \| v_{L, \eps} \|_{L^\infty_T L^\infty (\T_L^2)}^{p - 2} \\
&\les_\om \| \dt v_{L, \eps} \|_{L^\infty_T L_{- 2 \mu}^2 (\T_L^2)} \| v_{L, \eps} \|_{L^2_T \B_{4, 1, 2 \mu}^1 (\T_L^2)}^2 \| v_{L, \eps} \|_{L^\infty_T \B^0_{\infty, 1} (\T_L^2)}^{p - 2} , \\
\int_0^T &\bigg| \int_{\T_L^2} e^{- (p + 1) Y_{L, \eps}} \dt ( |v|^{p - 1} v ) \nb Y_{L, \eps} \cdot \nb \cj{v} dx \bigg| dt \\
&\les \| e^{- (p + 1) Y_{L, \eps}} \|_{L_{- \mu}^\infty (\T_L^2)} \| \nb Y_{L, \eps} \|_{L_{- \mu}^r (\T_L^2)} \| \dt v_{L, \eps} \|_{L^\infty_T L_{- 2 \mu}^2 (\T_L^2)} \\
&\quad \times \| \nb v_{L, \eps} \|_{L^2_T L_{2 \mu}^4 (\T_L^2)} \| v_{L, \eps} \|_{L^2_T L^\infty_{2 \mu} (\T_L^2)} \| v_{L, \eps} \|_{L^\infty_T L^{\frac{4r (p - 2)}{r - 4}} (\T_L^2)}^{p - 2} \\
&\les_\om |\log \eps| \| \dt v_{L, \eps} \|_{L^\infty_T L_{- 2 \mu}^2 (\T_L^2)} \| v_{L, \eps} \|_{L^2_T \B_{4, 1, 2 \mu}^1 (\T_L^2)} \| v_{L, \eps} \|_{L^2_T \B_{\infty, 1, 2 \mu}^0 (\T_L^2)} \| v_{L, \eps} \|_{L^\infty_T \B_{\frac{4r (p - 2)}{r - 4}, 1}^0 (\T_L^2)}^{p - 2}
\end{align*}

\noi
for any $\mu > 0$ sufficiently small and $r > 4$ satisfying $\mu r > 2$. Thus, the above estimates yield
\begin{align}
\begin{split}
\int_0^T |\mathcal{H}_{L, \eps} [v_{L, \eps} (t)]| dt 
&\les_{\om} |\log \eps| \| \dt v_{L, \eps} \|_{L^\infty_T L_{- 2 \mu}^2 (\T_L^2)} \Big( \| v_{L, \eps} \|_{L^2_T \B_{4, 1, 2 \mu}^1 (\T_L^2)}^2 \| v_{L, \eps} \|_{L^\infty_T \B^0_{\infty, 1} (\T_L^2)}^{p - 2} \\
&\quad + \| v_{L, \eps} \|_{L^2_T \B_{4, 1, 2 \mu}^1 (\T_L^2)} \| v_{L, \eps} \|_{L^2_T \B_{\infty, 1, 2 \mu}^0 (\T_L^2)} \| v_{L, \eps} \|_{L^\infty_T \B_{\frac{4r (p - 2)}{r - 4}, 1}^0 (\T_L^2)}^{p - 2} \Big)
\end{split}
\label{Hbdd1}
\end{align}

\noi
Using the equation \eqref{vNLSLe} and H\"older's inequalities along with Lemma~\ref{LEM:Lwei}~(iii), we have
\begin{align*}
\| &\dt v_{L, \eps} \|_{L_{-2 \mu}^2 (\T_L^2)} \\
&\leq \| \Dl v_{L, \eps} \|_{L^2_{- 2\mu} (\T_L^2)} + \| \nb Y_{L, \eps} \|_{L_{-  \mu}^r (\T_L^2)} \| \nb v_{L, \eps} \|_{L_{- \mu}^{\frac{2r}{r - 2}} (\T_L^2)} \\
&\quad + \big\| \wt{\wick{|\nb Y_{L, \eps}|^2}} \big\|_{L_{- \mu}^r (\T_L^2)} \| v_{L, \eps} \|_{L_{- \mu}^{\frac{2r}{r - 2}} (\T_L^2)} + \| e^{- (p - 1) Y_{L, \eps}} \|_{L_{-2 \mu}^\infty (\T_L^2)} \| v_{L, \eps} \|_{L^{2 p} (\T_L^2)}^{p} ,
\end{align*}

\noi
so that by further applying Lemma~\ref{LEM:YregLe}, Lemma~\ref{LEM:YregLe2}, and the embeddings in Lemma~\ref{LEM:embL}~(ii), (iv), and (v), we get
\begin{align}
\| \dt v_{L, \eps} \|_{L_{-2 \mu}^2 (\T_L^2)} 
\les_{\om} |\log \eps|^2 \Big( \| \Dl v_{L, \eps} \|_{L^2_{- \mu} (\T_L^2)} + \|  v_{L, \eps} \|_{\B_{\frac{2r}{r - 2}, 1, -\mu}^1 (\T_L^2)} + \| v_{L, \eps} \|_{\B_{2 p, 1}^0 (\T_L^2)}^{p} \Big) .
\label{Hbdd2}
\end{align}

\noi
We still need to deal with the Besov norms in \eqref{Hbdd1} and \eqref{Hbdd2}. From the embeddings in Lemma~\ref{LEM:embL}~(i) and (iii) and the interpolation estimate in Lemma~\ref{LEM:interpL}, with $q = \infty$, $q = \frac{4 r (p - 2)}{r - 4}$, or $q = 2p$, we have
\begin{align}
\| v_{L, \eps} \|_{\B_{q, 1}^0 (\T_L^2)} \les \| v_{L, \eps} \|_{\B_{2, 2}^{1 + \nu} (\T_L^2)} \les \| v_{L, \eps} \|_{\B_{2, 2, \frac{2 \nu}{1 - \nu} \mu}^{1 - \nu} (\T_L^2)}^{\frac{1 - \nu_0}{1 + \nu}} \| v_{L, \eps} \|_{\B_{2, 2, - \mu}^{2} (\T_L^2)}^{\frac{2 \nu}{1 + \nu}} 
\label{Hbdd3}
\end{align}

\noi
for any $\nu > 0$. Also, from the Besov embedding in Lemma~\ref{LEM:embL}~(iii), we have
\begin{align}
\| v_{L, \eps} \|_{\B_{\infty, 1, 2 \mu}^0 (\T_L^2)} \les \| v_{L, \eps} \|_{\B_{4, 1, 2 \mu}^1 (\T_L^2)} .
\label{Hbdd4}
\end{align}

\noi
Moreover, to deal with the $\B_{\frac{2r}{r - 2}, 1, - \mu}^1$-term in \eqref{Hbdd2} and the $\B_{2, 2, -\mu}^2$-term in \eqref{Hbdd3}, we use the embeddings in Lemma~\ref{LEM:embL}~(i), (iii), and (v) and Lemma~\ref{LEM:L2equiv} to obtain
\begin{align}
\|  v_{L, \eps} \|_{\B_{\frac{2r}{r - 2}, 1, -\mu}^1 (\T_L^2)} \les \|  v_{L, \eps} \|_{\B_{2, 2, -\mu}^2 (\T_L^2)} \sim \| v_{L, \eps} \|_{\B_{2, 2, - \mu}^0 (\T_L^2)} + \| \Dl v_{L, \eps} \|_{L^2_{- \mu} (\T_L^2)} .
\label{Hbdd5}
\end{align}

\noi
The desired estimate then follows from combining \eqref{Hbdd1}, \eqref{Hbdd2}, \eqref{Hbdd3}, \eqref{Hbdd4}, and \eqref{Hbdd5} along with the embedding in Lemma~\ref{LEM:embL}~(iv).
\end{proof}

\subsection{A priori bounds for mollified $L$-periodic nonlinear solutions}
\label{SUB:H2L}

We now establish some a priori bounds for the solution $v_{L, \eps}$ to the mollified $L$-periodic nonlinear DAM \eqref{vNLSLe}.
We start with the following $L^2$ and $H^1$ a priori bounds for $v_{L, \eps}$.
\begin{lemma}
\label{LEM:H1bdd}
Let $p \geq 2$, $\ld > 0$, $L \geq 1$, $\mu_0 > 0$, and $v_{0, L} \in H_{\mu_0}^2 (\T_L^2)$. Given $0 < \eps < 1$, let $v_{L, \eps}$ be the global-in-time solution to the mollified $L$-periodic nonlinear DAM \eqref{vNLSLe} with $v_{L, \eps} |_{t = 0} = v_{0, L}$. Then, for any $T \geq 1$ and $\mu > 0$ sufficiently small, there exists $\Om' \subset \Om$ with full probability measure such that for any $\om \in \Om'$, we have the bounds
\begin{align*}
\sup_{\eps \in (0, 1)} \| v_{L, \eps} (t) \|_{C_T L_\mu^2 (\T_L^2)} \les_{\om} \| v_{0, L} \|_{H_{\mu_0}^1 (\T_L^2)}
\end{align*}

\noi
and
\begin{align*}
\sup_{\eps \in (0, 1)} \| \nb v_{L, \eps} (t) \|_{C_T L^2_{- \mu} (\T_L^2)} \les_{\om} \| v_{0, L} \|_{H_{\mu_0}^1 (\T_L^2)} + \| v_{0, L} \|_{H_{\mu_0}^1 (\T_L^2)}^{\frac{p + 1}{2}} .
\end{align*}
\end{lemma}

\begin{proof}
Due to the conservation of mass $\mathcal{M}_{L, \eps}$ in \eqref{defM1} and the energy $\mathcal{E}_{L, \ld, \eps}$ in \eqref{defE} along with the fact that $\ld > 0$, the proof follows similarly from that of Lemma~\ref{LEM:St_H1}, and so we will be brief. 

We fix $t \in \R_+$. By the conservation of energy $\mathcal{E}_{L, \ld, \eps}$ in \eqref{defE} and the fact that $\ld > 0$, we have
\begin{align}
\begin{split}
\frac 12 &\int_{\T_L^2} e^{-2 Y_{L, \eps}} |\nb v_{L, \eps} (t)|^2 dx - \frac 12 \int_{\T_L^2} e^{-2 Y_{L, \eps}} \wt{\wick{|\nb Y_{L, \eps}|^2}} \, |v_{L, \eps} (t)|^2 dx \\
&\leq \frac 12 \int_{\T_L^2} e^{-2 Y_{L, \eps}} |\nb v_{0, L}|^2 dx - \frac 12 \int_{\T_L^2} e^{-2 Y_{L, \eps}} \wt{\wick{|\nb Y_{L, \eps}|^2}} \, |v_{0, L}|^2 dx \\
&\quad + \frac{\ld}{p + 1} \int_{\T_L^2} e^{- (p + 1) Y_{L, \eps}} |v_{0, L}|^{p + 1} dx.
\end{split}
\label{H1-1}
\end{align}

\noi
By H\"older's inequality along with Lemma~\ref{LEM:Lwei}~(iii) and (iv), Lemma~\ref{LEM:YregLe}~(iii), the embeddings in Lemma~\ref{LEM:embL}~(ii), (i), and (v), and Lemma~\ref{LEM:L2equiv}, we have
\begin{align}
\begin{split}
\bigg| \int_{\T_L^2} e^{- (p + 1) Y_{L, \eps}} |v_{0, L}|^{p + 1} dx \bigg| &\les \| e^{- (p + 1) Y_{L, \eps}} \|_{L_{- (p + 1) \mu_0}^\infty (\T_L^2)} \| v_{0, L} \|_{L_{\mu_0}^{p + 1} (\T_L^2)}^{p + 1} \\
&\les_{\om} \| v_{0, L} \|_{\B_{2, 2, \mu_0}^{1} (\T_L^2)}^{p + 1} \\
&\sim \| v_{0, L} \|_{H_{\mu_0}^1 (\T_L^2)}^{p + 1} .
\end{split}
\label{H1-2}
\end{align}

\noi
Thus, from \eqref{H1-1}, \eqref{L2b1-2}, \eqref{L2b1-3}, and \eqref{H1-2}, we use similar steps in \eqref{L2b1} to obtain
\begin{align}
\| \nb v_{L, \eps} (t) \|_{L_{- \mu}^2 (\T_L^2)}^2 \les_\om \| v_{L, \eps} (t) \|_{L_{3 \mu}^2 (\T_L^2)}^2 + \| v_{0, L} \|_{H_{\mu_0}^1 (\T_L^2)} + \| v_{0, L} \|_{H_{\mu_0}^1 (\T_L^2)}^{p + 1} 
\label{H1-3}
\end{align}

\noi
for any $\mu > 0$ sufficiently small. Then, from \eqref{L2b2-0}, \eqref{L2b2-1}, \eqref{L2b2-2}, and \eqref{H1-3}, we obtain
\begin{align*}
\| v_{L, \eps} (t) \|_{L_{3 \mu}^2 (\T_L^2)}^2 &\les_\om \| v_{0, L} \|_{L_{\mu_0}^2 (\T_L^2)}^2 + \int_0^t \| \nb v_{L, \eps} (t') \|_{L_{- \mu}^2 (\T_L^2)} \| v_{L, \eps} (t') \|_{L_{3 \mu}^2 (\T_L^2)} dt' \\
&\les_{\om, T} \| v_{0, L} \|_{H_{\mu_0}^1 (\T_L^2)}^2 + \| v_{0, L} \|_{H_{\mu_0}^1 (\T_L^2)}^{p + 1} + \int_0^t \| v_{L, \eps} (t') \|_{L_{3 \mu}^2 (\T_L^2)}^2 dt' .
\end{align*}

\noi
The desired estimates then follow from Gronwall's inequality, the embedding in Lemma~\ref{LEM:embL}~(iv), and \eqref{H1-3}.
\end{proof}

We also need a bound for $\Dl v_{L, \eps}$ by using the conservation of modified energies \eqref{FGH}. To achieve this, we need to establish Strichartz estimates as in \cite{TzV23-2}. To avoid introducing more notations and preliminary embedding estimates, we choose to perform the analysis on Besov spaces. 

We first establish the following Strichartz estimates for the linear flow $S_{L, \eps} (t)$ of the mollified $L$-periodic linear DAM \eqref{wNLSLe}.

\begin{lemma}
\label{LEM:St_est}
Let $L \geq 1$ and $\mu_0 > 0$. Then, for any $T \geq 1$ and $\gamma > 0$, there exists $\Om' \subset \Om$ with full probability measure such that for any $\om \in \Om'$ and $0 < \eps < \frac 12$, we have the bound
\begin{align*}
\| S_{L, \eps} (t) f_L \|_{L_T^4 \B_{4, \infty}^{\frac 34 - \gamma} (\T_L^2)} \les_{\om, T} |\log \eps|^{7} \| f_L \|_{H_{\mu_0}^1 (\T_L^2)} .
\end{align*}
\end{lemma}

\begin{proof}
Let us first show that for any $\gamma > 0$ and $2 < q, r < \infty$ such that $\frac{2}{q} + \frac{2}{r} = 1$, we have on a full probability set that for any $0 < \eps < \frac 12$,
\begin{align}
\| S_{L, \eps} (t) f_L \|_{L_T^q L^r (\T_L^2)} \les_{\om, T} |\log \eps|^{7} \| f_L \|_{\B_{2, 2, \mu_0}^{\frac 1q + \gamma} (\T_L^2)} . 
\label{SLe_qr}
\end{align}

\noi
We use the Littlewood-Paley decomposition on $S_{L, \eps} (t) f_L$ to obtain 
\begin{align*}
S_{L, \eps} (t) f_L = \sum_{\substack{N_1, N_2 \geq 1 \\ \text{dyadic}}} \Dl_{N_1}^L S_{L, \eps} (t) \Dl_{N_2}^L f_L .
\end{align*}

\noi
In the case $N_1 \leq 8$, we easily obtain the desired estimate using the Besov embedding in Lemma~\ref{LEM:embL}~(iii) and the $L^2$-bound in Lemma~\ref{LEM:St_H1} along with Lemma~\ref{LEM:embL}~(iv), and so we only consider $N_1 \geq 16$.
We split the interval $[0, T]$ into an essentially disjoint union of intervals of size $\min (N_1^{-1}, N_2^{-1})$. Given one such interval $I = [a, b]$ with $b - a = \min (N_1^{-1}, N_2^{-1})$ and $t \in [a, b]$, we write
\begin{align}
\begin{split}
\Dl_{N_1}^L S_{L, \eps} (t) \Dl_{N_2}^L f_L &= \Dl_{N_1}^L e^{-i (t - a) \Dl} S_{L, \eps} (a) \Dl_{N_2}^L f_L \\
&\quad + i \int_a^t \Dl_{N_1}^L e^{- i (t - t') \Dl} \big( 2 \nb Y_{L, \eps} \cdot \nb (S_{L, \eps} (t') \Dl_{N_2}^L f_L) \\
&\qquad \quad - \wt{\wick{ |\nb Y_{L, \eps}|^2 }} \, S_{L, \eps} (t') \Dl_{N_2}^L f_L \big) dt' .
\end{split}
\label{L4bdd1}
\end{align}

\noi
Using the Strichartz estimate in Lemma~\ref{LEM:L4} and  Lemma~\ref{LEM:St_Bs} along with Lemma~\ref{LEM:embL}~(iv), we have
\begin{align}
\begin{split}
\big\| \Dl_{N_1}^L e^{-i (t - a) \Dl} S_{L, \eps} (a) \Dl_{N_2}^L f_L \big\|_{L^q_I L^r (\T_L^2)} &\les \| \Dl_{N_1}^L S_{L, \eps} (a) \Dl_{N_2}^L f_L \|_{L^2 (\T_L^2)} \\
&\les N_1^{- \frac 1q - \gamma_0} \| S_{L, \eps} (a) \Dl_{N_2}^L f_L \|_{\B_{2, 2}^{\frac{1}{q} + \gamma_0} (\T_L^2)} \\
&\les_{\om} |\log \eps|^5 N_1^{- \frac 1q - \gamma_0} \| \Dl_{N_2}^L f_L \|_{\B_{2, 2, \mu_0}^{\frac 1q + 2\gamma_0} (\T_L^2)}
\end{split}
\label{L4bdd2}
\end{align}

\noi
and, also from \eqref{Hs_equi},
\begin{align}
\begin{split}
\big\| \Dl_{N_1}^L e^{-i (t - a) \Dl} S_{L, \eps} (a) \Dl_{N_2}^L f_L \big\|_{L^q_I L^r (\T_L^2)} &\les \| S_{L, \eps} (a) \Dl_{N_2}^L f_L \|_{L^2 (\T_L^2)} \\
&\les_{\om} |\log \eps|^5 N_2^{- \frac 1q - \gamma_0} \| \Dl_{N_2}^L f_L \|_{\B_{2, 2, \mu_0}^{\frac 1q + 2 \gamma_0} (\T_L^2)}
\end{split}
\label{L4bdd3}
\end{align}

\noi
for any $\gamma_0 > 0$. By H\"older's inequalities along with Lemma~\ref{LEM:Lwei}~(iii), Lemma~\ref{LEM:YregLe2}, the embeddings in Lemma~\ref{LEM:embL}~(ii), (i), and (iii), and Lemma~\ref{LEM:St_Bs}, we have for any $t' \geq 0$ that
\begin{align*}
\big\| \nb Y_{L, \eps} \cdot \nb ( S_{L, \eps} (t') \Dl_{N_2}^L f_L ) \big\|_{L^2 (\T_L^2)} 
&\les \| \nb Y_{L, \eps} \|_{L_{- \mu}^{q_0} (\T_L^2)} \| \nb ( S_{L, \eps} (t') \Dl_{N_2}^L f_L ) \|_{L^{\frac{2q_0}{q_0 - 2}}_\mu (\T_L^2)} \\
&\les_\om |\log \eps| \| S_{L, \eps} (t') \Dl_{N_2}^L f_L \|_{\B_{2, 2, \mu}^{1 + \frac{2}{q_0} + \gamma_0} (\T_L^2)} \\
&\les |\log \eps|^6 N_2^{1 + \frac{2}{q_0} + 2 \gamma_0} \| \Dl_{N_2}^L f_L \|_{L_{\mu_0}^2 (\T_L^2)} 
\end{align*}

\noi
and
\begin{align*}
\big\| \wt{\wick{ |\nb Y_{L, \eps}|^2 }} \, S_{L, \eps} (t') \Dl_{N_2}^L f_L \big) \big\|_{L^2 (\T_L^2)} 
&\les \big\| \wt{\wick{ |\nb Y_{L, \eps}|^2 }} \big\|_{L^{q_0}_{- \mu} (\T_L^2)} \| S_{L, \eps} (t') \Dl_{N_2}^L f_L \|_{L^{\frac{2q_0}{q_0 - 2}}_{\mu} (\T_L^2)} \\
&\les_\om |\log \eps|^2 \| S_{L, \eps} (t') \Dl_{N_2}^L f_L \|_{\B_{2, 2, \mu}^{\frac{2}{q_0} + \gamma_0} (\T_L^2)} \\
&\les |\log \eps|^7 N_2^{\frac{2}{q_0} + 2 \gamma_0} \| \Dl_{N_2}^L f_L \|_{L_{\mu_0}^2 (\T_L^2)} 
\end{align*}

\noi
for any $\mu > 0$ sufficiently small and $q_0 > 2$ satisfying $\mu q_0 > 2$, so that by the Strichartz estimate in Lemma~\ref{LEM:L4} and the above two estimates, we get
\begin{align}
\begin{split}
\bigg\| &\int_a^t \Dl_{N_1}^L e^{- i (t - t') \Dl} \big( 2 \nb Y_{L, \eps} \cdot \nb (S_{L, \eps} (t') \Dl_{N_2}^L f_L ) - \wt{\wick{ |\nb Y_{L, \eps}|^2 }} \, S_{L, \eps} (t') \Dl_{N_2}^L f_L \big) dt' \bigg\|_{L_I^q L^r (\T_L^2)} \\
&\les \int_a^b \Big( \big\| \nb Y_{L, \eps} \cdot \nb ( S_{L, \eps} (t') \Dl_{N_2}^L f_L ) \big\|_{L^2 (\T_L^2)} + \big\| \wt{\wick{ |\nb Y_{L, \eps}|^2 }} \, S_{L, \eps} (t') \Dl_{N_2}^L f_L \big) \big\|_{L^2 (\T_L^2)} \Big) dt' \\
&\les_\om |\log \eps|^7 \max (N_1, N_2)^{-1} N_2^{1 + \frac{2}{q_0} + 2 \gamma_0} \| \Dl_{N_2}^L f_L \|_{L_{\mu_0}^2 (\T_L^2)} .
\end{split}
\label{L4bdd4}
\end{align}

\noi
Combining \eqref{L4bdd1}, \eqref{L4bdd2}, \eqref{L4bdd3}, and \eqref{L4bdd4}, we obtain
\begin{align*}
\| \Dl_{N_1}^L S_{L, \eps} (t) \Dl_{N_2}^L f_L \|_{L_I^q L^r (\T_L^2)} &\les_{\om} |\log \eps|^5 \max (N_1, N_2)^{- \frac 1q - \gamma_0} \| f_L \|_{\B_{2, 2, \mu_0}^{\frac 1q + 2 \gamma_0} (\T_L^2)} \\
&\quad + |\log \eps|^{7} \max (N_1, N_2)^{- \frac 1q - \gamma_0} \| f_L \|_{\B_{2, 2, \mu_0}^{\frac{1}{q} + \frac{2}{q_0} + 3 \gamma_0} (\T_L^2)}.
\end{align*}

\noi
The desired estimate \eqref{SLe_qr} then follows by taking the $q$th power, concatenating the time intervals, adding up the dyadic numbers $N_1 , N_2 \geq 1$, and taking $q_0$ to be sufficiently large and $\gamma_0$ to be sufficiently small.

We now show the $L_T^4 \B_{4, \infty}^{\frac 34 - \gamma}$-estimate. It suffices to show the estimate when $\gamma > 0$ is sufficiently small, say $0 < \gamma < \frac{1}{10}$. We write
\begin{align*}
S_{L, \eps} (t) f_L = \sum_{\substack{N \geq 1 \\ \text{dyadic}}} S_{L, \eps} (t) \Dl_N^L f_L.
\end{align*}

\noi
From the interpolation in Lemma~\ref{LEM:interpL}, the interpolation in time, and the embeddings in Lemma~\ref{LEM:embL}~(i) and (ii), we obtain
\begin{align}
\begin{split}
\| S_{L, \eps} (t) \Dl_N^L f \|_{L_T^4 \B_{4, \infty}^{\frac 34 - \gamma} (\T_L^2)} &\les \| S_{L, \eps} (t) \Dl_N^L f \|_{L_T^{4 (1 - \ta')} \B_{r, \infty}^0 (\T_L^2)}^{1 - \ta'} \| S_{L, \eps} (t) \Dl_N^L f \|_{L_T^\infty \B_{2, \infty}^{\frac 32} (\T_L^2)}^{\ta'} \\
&\les \| S_{L, \eps} (t) \Dl_N^L f_L \|_{L_T^{q} L^r (\T_L^2)}^{1 - \ta'} \| S_{L, \eps} (t) \Dl_N^L f_L \|_{L_T^\infty \B_{2, 2}^{\frac 32} (\T_L^2)}^{\ta'} ,
\end{split}
\label{L4bdd5}
\end{align}

\noi
where $0 < \ta' < \frac 12$ and $2 < r < \infty$ satisfy
\begin{align}
\frac{1}{4} = \frac{\ta'}{2} + \frac{1 - \ta'}{r} \quad \text{and} \quad \frac{3}{4} - \gamma = \frac{3 \ta'}{2}
\label{cond}
\end{align}

\noi
and $q = 4 (1 - \ta')$ satisfies $\frac{2}{q} + \frac{2}{r} = 1$. From \eqref{SLe_qr}, we have
\begin{align}
\| S_{L, \eps} (t) \Dl_N^L f_L \|_{L_T^{q} L^r (\T_L^2)} \les_{\om, T} |\log \eps|^{7} N^{\frac{1}{q} + \frac{\gamma}{2}} \| \Dl_N^L f_L \|_{L_{\mu_0}^2 (\T_L^2)} .
\label{L4bdd6}
\end{align}

\noi
Also, from Lemma~\ref{LEM:St_Bs} along with Lemma~\ref{LEM:Lwei}~(iii), we have
\begin{align}
\begin{split}
\| S_{L, \eps} (t) \Dl_N^L f_L \|_{L_T^\infty \B_{2, 2}^{\frac 32} (\T_L^2)} \les_{\om} |\log \eps|^5 N^{\frac 32 + \frac{\gamma}{2}} \| \Dl_N^L f_L \|_{L_{\mu_0}^2 (\T_L^2)}
\end{split}
\label{L4bdd7}
\end{align}

\noi
Combining \eqref{L4bdd5}, \eqref{L4bdd6}, and \eqref{L4bdd7}, we obtain
\begin{align*}
\| S_{L, \eps} (t) \Dl_N^L f_L \|_{L_T^4 \B_{4, \infty}^{\frac 34 - \gamma} (\T_L^2)} \les_{\om, T} |\log \eps|^{7} N^{\frac{1 - \ta'}{q} + \frac{3 \ta'}{2} + \frac{\gamma}{2}} \| \Dl_N^L f_L \|_{L_{\mu_0}^2 (\T_L^2)} .
\end{align*}

\noi
Since we have from \eqref{cond} that $\frac{1 - \ta'}{q} + \frac{3 \ta'}{2} + \frac{\gamma}{2} < \frac{1}{4} + \frac{3 \ta'}{2} + \gamma = 1$, we conclude the desired estimate by summing up dyadic $N \geq 1$, the embedding in Lemma~\ref{LEM:embL}~(i), and the fact that  
\begin{align*}
\| f_L \|_{\B_{2, 2, \mu_0}^1 (\T_L^2)} \sim \| f_L \|_{L_{\mu_0}^2 (\T_L^2)} + \| \nb f_L \|_{L_{\mu_0}^2 (\T_L^2)} = \| f_L \|_{H^1_{\mu_0} (\T_L^2)}, 
\end{align*}

\noi
which follows from Lemma~\ref{LEM:embL}~(v) and Lemma~\ref{LEM:L2equiv}.
\end{proof}

We now use the Strichartz estimate in Lemma~\ref{LEM:St_est} to deal with the $\| v_{L, \eps} \|_{L_T^2 \B_{4, 1, \mu}^1 (\T_L^2)}$ term appearing in Lemma~\ref{LEM:Hest}.
\begin{lemma}
\label{LEM:vL2L4}
Let $p \geq 2$, $\ld > 0$, $L \geq 1$, $\mu_0 > 0$, and $v_{0, L} \in H_{\mu_0}^2 (\T_L^2)$. Given $0 < \eps < \frac 12$, let $v_{L, \eps}$ be the solution to the mollified $L$-periodic nonlinear DAM \eqref{vNLSLe} with $v_{L, \eps} |_{t = 0} = v_{0, L}$.  Then, for any $T \geq 1$, $\nu > 0$ sufficiently small, and $\mu > 0$ sufficiently small, there exists $\Om' \subset \Om$ with full probability measure such that for any $\om \in \Om'$ and $0 < \eps < \frac 12$, we have
\begin{align*}
\| v_{L, \eps} \|_{L_T^2 \B_{4, 1, \mu}^1 (\T_L^2)} &\les_{\om, T} |\log \eps|^{4} \Big( 1 + \| v_{0, L} \|_{H_{\mu_0}^1 (\T_L^2)} + \| v_{L, \eps} \|_{L_T^\infty \B_{2, 2, 12 \mu}^{1 - \nu} (\T_L^2)}^{\frac{p}{2} + 1} \Big) \\
&\quad \times \Big( 1 + \| \Dl v_{L, \eps} \|_{L_T^\infty L^2_{- \mu} (\T_L^2)}^{\frac 38 + (p + 1) \nu} \Big) .
\end{align*}
\end{lemma}

\begin{proof}
From the interpolation estimate in Lemma~\ref{LEM:interpL}, the embeddings in Lemma~\ref{LEM:embL}~(i) and (iii), and H\"older's inequalities in time, we have
\begin{align}
\begin{split}
\| v_{L, \eps} \|_{L_T^2 \B_{4, 1, \mu}^1 (\T_L^2)} &\les \| v_{L, \eps} \|_{L_T^2 \B_{4, 1}^{\frac 34 - 2 \gamma} (\T_L^2)}^{\ta} \| v_{L, \eps} \|_{L_T^2 \B_{4, 1, \frac{\mu}{1 - \ta}}^{\frac 54 - 2 \gamma} (\T_L^2)}^{1 - \ta} \\
&\les_T \| v_{L, \eps} \|_{L_T^4 \B_{4, \infty}^{\frac 34 - \gamma} (\T_L^2)}^{\ta} \| v_{L, \eps} \|_{L_T^\infty \B_{2, 2, \frac{\mu}{1 - \ta}}^{\frac 74 - \gamma} (\T_L^2)}^{1 - \ta}
\end{split}
\label{L2L4bdd1}
\end{align}

\noi
for any $\gamma > 0$ small enough and $0 < \ta < 1$ satisfying
\begin{align*}
1 = \Big( \frac 34 - 2 \gamma \Big) \ta + \Big( \frac 54 - 2 \gamma \Big) (1 - \ta) \quad \Longleftrightarrow \quad \ta = \frac 12 - 4 \gamma.
\end{align*}

\noi
To deal with the $\B_{4, \infty}^{\frac 34 - \gamma}$-term in \eqref{L2L4bdd1}, we note that we have the following Duhamel formulation for $v_{L, \eps}$:
\begin{align*}
v_{L, \eps} (t) = S_{L, \eps} (t) v_{0, L} + i \ld \int_0^t S_{L, \eps} (t - t') \big( e^{- (p - 1) Y_{L, \eps}} |v_{L, \eps}|^{p - 1} v_{L, \eps} \big) (t') dt' .
\end{align*}

\noi
From Lemma~\ref{LEM:St_est}, we get
\begin{align}
\| v_{L, \eps} \|_{L_T^4 \B_{4, \infty}^{\frac 34 - \gamma} (\T_L^2)} 
\les_{\om, T} |\log \eps|^{7} \big( \| v_{0, L} \|_{H_{\mu_0}^1 (\T_L^2)} + \big\| e^{- (p - 1) Y_{L, \eps}} |v_{L, \eps}|^{p - 1} v_{L, \eps} \big\|_{L_T^\infty H_\mu^1 (\T_L^2)} \big) .
\label{L2L4bdd2-1}
\end{align}

\noi
For the second term on the right-hand-side of \eqref{L2L4bdd2-1}, we use H\"older's inequalities along with Lemma~\ref{LEM:Lwei}~(iii) to get
\begin{align*}
\big\| &e^{- (p - 1) Y_{L, \eps}} |v_{L, \eps}|^{p - 1} v_{L, \eps} \big\|_{L_T^\infty H_\mu^1 (\T_L^2)} \\
&\les \| e^{- (p - 1) Y_{L, \eps}} \|_{L_{- p \mu}^\infty (\T_L^2)} \| v_{L, \eps} \|_{L_T^\infty H_\mu^1 (\T_L^2)} \| v_{L, \eps} \|_{L_T^\infty L_\mu^\infty (\T_L^2)}^{p - 1} \\
&\quad + \| \nb Y_{L, \eps} \|_{L_{- \mu}^r (\T_L^2)} \| e^{- (p - 1) Y_{L, \eps}} \|_{L_{- (p - 1) \mu}^\infty (\T_L^2)} \| v_{L, \eps} \|_{L_T^\infty L_{\mu}^{\frac{2rp}{r - 2}} (\T_L^2)}^p ,
\end{align*}

\noi
for any $r > 2$ satisfying $\mu r > 2$, so that by further applying  Lemma~\ref{LEM:YregLe}~(iii), Lemma~\ref{LEM:YregLe2}, Lemma~\ref{LEM:L2equiv}, and the embeddings in Lemma~\ref{LEM:embL}~(v), (ii), (i), and (iii), we obtain
\begin{align}
\begin{split}
\big\| &e^{- (p - 1) Y_{L, \eps}} |v_{L, \eps}|^{p - 1} v_{L, \eps} \big\|_{L_T^\infty H_\mu^1 (\T_L^2)} \\
&\les_\om |\log \eps| \Big( \| v_{L, \eps} \|_{L_T^\infty \B_{2, 2, \mu}^{1} (\T_L^2)} \| v_{L, \eps} \|_{L_T^\infty \B_{\infty, 1, \mu}^0 (\T_L^2)}^{p - 1} 
+ \| v_{L, \eps} \|_{L_T^\infty \B_{\frac{2 r p}{r - 2}, 1, \mu}^0 (\T_L^2)}^p \Big) \\
&\les |\log \eps| \| v_{L, \eps} \|_{L_T^\infty \B_{2, 2, \mu}^{1 + \nu} (\T_L^2)}^p 
\end{split}
\label{L2L4bdd2-2}
\end{align}

\noi
for any $\nu > 0$. Combining \eqref{L2L4bdd2-1} and \eqref{L2L4bdd2-2}, we obtain
\begin{align}
\| v_{L, \eps} \|_{L_T^4 \B_{4, \infty}^{\frac 34 - \gamma} (\T_L^2)} 
\les_\om |\log \eps|^{8} \Big( \| v_{0, L} \|_{H_{\mu_0}^1 (\T_L^2)} + \| v_{L, \eps} \|_{L_T^\infty \B_{2, 2, \mu}^{1 + \nu} (\T_L^2)}^p \Big) .
\label{L2L4bdd2}
\end{align}

\noi
To deal with the $\B_{2, 2, \frac{\mu}{1 - \ta}}^{\frac 74 - \gamma}$-term in \eqref{L2L4bdd1}, we use the embedding in Lemma~\ref{LEM:embL}~(iv) and the interpolation estimate in Lemma~\ref{LEM:interpL} to obtain
\begin{align}
\| v_{L, \eps} \|_{L_T^\infty \B_{2, 2, \frac{\mu}{1 - \ta}}^{\frac 74 - \gamma} (\T_L^2)} \les \| v_{L, \eps} \|_{L_T^\infty \B_{2, 2, 2 \mu}^{\frac 74 + \frac{\nu}{4}} (\T_L^2)} \les \| v_{L, \eps} \|_{L_T^\infty \B_{2, 2, 11 \mu}^{1 + \nu} (\T_L^2)}^{\frac 14} \| v_{L, \eps} \|_{L_T^\infty \B_{2, 2, - \mu}^{2} (\T_L^2)}^{\frac 34} 
\label{L2L4bdd3}
\end{align}

\noi
Moreover, to bound the $\B_{2, 2, \mu}^{1 + \nu}$-term in \eqref{L2L4bdd2} and the $\B_{2, 2, 11 \mu}^{1 + \nu}$-term in \eqref{L2L4bdd3}, we use the embedding in Lemma~\ref{LEM:embL}~(iv), the interpolation estimate in Lemma~\ref{LEM:interpL}, Lemma~\ref{LEM:embL}~(v), and Lemma~\ref{LEM:L2equiv} to get
\begin{align}
\begin{split}
\| &v_{L, \eps} \|_{L_T^\infty \B_{2, 2, \mu}^{1 + \nu} (\T_L^2)} + \| v_{L, \eps} \|_{L_T^\infty \B_{2, 2, 11 \mu}^{1 + \nu} (\T_L^2)} \\
&\les \| v_{L, \eps} \|_{L_T^\infty \B_{2, 2, \frac{11 + 13 \nu}{1 - \nu} \mu}^{1 - \nu} (\T_L^2)}^{\frac{1 - \nu}{1 + \nu}} \| v_{L, \eps} \|_{L_T^\infty \B_{2, 2, - \mu}^{2} (\T_L^2)}^{\frac{2 \nu}{1 + \nu}} \\
&\sim \| v_{L, \eps} \|_{L_T^\infty \B_{2, 2, \frac{11 + 13 \nu}{1 - \nu} \mu}^{1 - \nu} (\T_L^2)}^{\frac{1 - \nu}{1 + \nu}} \big( \| v_{L, \eps} \|_{L_T^\infty L_{- \mu}^2 (\T_L^2)} + \| \Dl v_{L, \eps} \|_{L_T^\infty L_{- \mu}^2 (\T_L^2)} \big)^{\frac{2 \nu}{1 + \nu}}
\end{split}
\label{L2L4bdd4}
\end{align}

\noi
Thus, the desired estimate follows from combining \eqref{L2L4bdd1}, \eqref{L2L4bdd2}, \eqref{L2L4bdd3}, and \eqref{L2L4bdd4}, using Lemma~\ref{LEM:embL}~(ii), (i), and (iv) to deal with the $L_{- \mu}^2$-term in \eqref{L2L4bdd4}, and taking $\nu$ and $\gamma$ to be sufficiently small.
\end{proof}

We are now ready to prove a bound for $\Dl v_{L, \eps}$. Note that we do not claim any sharpness on the concrete numbers appearing below, but we choose to keep these numbers to make sure that they are independent of all parameters.

\begin{lemma}
\label{LEM:H2bdd}
Let $p \geq 2$, $\ld > 0$, $L \geq 1$, $\mu_0 > 0$, and $v_{0, L} \in H_{\mu_0}^2 (\T_L^2)$. Given $0 < \eps < \frac 12$, let $v_{L, \eps}$ be the solution to the mollified $L$-periodic nonlinear DAM \eqref{vNLSLe} with $v_{L, \eps} |_{t = 0} = v_{0, L}$. Then, for any $T \geq 1$ and $\mu > 0$ sufficiently small, there exists $\Om' \subset \Om$ with full probability measure such that for any $\om \in \Om'$ and $0 < \eps < \frac 12$, we have
\begin{align*}
\| \Dl v_{L, \eps} \|_{C_T L^2_{- \mu} (\T_L^2)} \les_{\om, T} |\log \eps|^{136} \Big( 1 + \| v_{0, L} \|_{H_{\mu_0}^2 (\T_L^2)}^{20 (p + 1)^2} \Big) .
\end{align*}
\end{lemma}

\begin{proof}
Let $t \in [0, T]$. From the conservation \eqref{FGH}, we have
\begin{align}
\begin{split}
\| &e^{- Y_{L, \eps}} \Dl v_{L, \eps} (t) \|_{L^2 (\T_L^2)}^2 \\
&\leq  \bigg| \mathcal{F}_{L, \eps} [v_{L, \eps} (t)] - \int_{\T_L^2} e^{-2 Y_{L, \eps}} |\Dl v_{L, \eps} (t)|^2 dx \bigg| + \bigg| \mathcal{F}_{L, \eps} [v_{0, L}] - \int_{\T_L^2} e^{-2 Y_{L, \eps}} |\Dl v_{0, L}|^2 dx \bigg| \\
&\quad + \| e^{- Y_{L, \eps}} \Dl v_{0, L} \|_{L^2 (\T_L^2)}^2 + \ld |\mathcal{G}_{L, \eps} [v_{L, \eps} (t)]| + \ld |\mathcal{G}_{L, \eps} [v_{0, L}]| + \ld \int_0^T |\mathcal{H}_{L, \eps} [v_{L, \eps} (t')]| dt' .
\end{split}
\label{vH2bdd1-1}
\end{align}

\noi
From Lemma~\ref{LEM:Fest}, we get
\begin{align}
\begin{split}
\bigg| &\mathcal{F}_{L, \eps} [v_{L, \eps} (t)] - \int_{\T_L^2} e^{-2 Y_{L, \eps}} |\Dl v_{L, \eps} (t)|^2 dx \bigg| \\
&\les_\om \dl^{-7} |\log \eps|^{17} \| v_{L, \eps} (t) \|_{L_{8 \mu}^2 (\T_L^2)}^{2} + \dl \| \Dl v_{L, \eps} (t) \|_{L^2_{- \mu} (\T_L^2)}^{2}
\end{split}
\label{vH2bdd1-2}
\end{align}

\noi
and
\begin{align}
\begin{split}
\bigg| &\mathcal{F}_{L, \eps} [v_{0, L}] - \int_{\T_L^2} e^{-2 Y_{L, \eps}} |\Dl v_{0, L}|^2 dx \bigg| \\
&\les_\om |\log \eps|^{17} \| v_{0, L} \|_{L_{8 \mu}^2 (\T_L^2)}^2 + \| \Dl v_{0, L} \|_{L_{- \mu}^2 (\T_L^2)}^{2} 
\end{split}
\label{vH2bdd1-3}
\end{align}

\noi
for any $0 < \dl < 1$ and $\mu > 0$. By H\"older's inequality along with Lemma~\ref{LEM:Lwei}~(iii) and Lemma~\ref{LEM:YregLe}~(iii), we have
\begin{align}
\| e^{- Y_{L, \eps}} \Dl v_{0, L} \|_{L^2 (\T_L^2)}^2 \les \| e^{- Y_{L, \eps}} \|_{L_{- \mu_0}^\infty (\T_L^2)}^2 \| \Dl v_{0, L} \|_{L_{\mu_0}^2 (\T_L^2)}^2 
\les_\om \| \Dl v_{0, L} \|_{L_{\mu_0}^2 (\T_L^2)}^2 .
\label{vH2bdd1-4}
\end{align}

\noi
Also, from Lemma~\ref{LEM:Gest}, we get
\begin{align}
|\mathcal{G}_{L, \eps} [v_{L, \eps} (t)]| \les_\om |\log \eps|^2 \big( \| v_{L, \eps} (t) \|_{L_{- \mu}^2 (\T_L^2)} + \| \Dl v_{L, \eps} (t) \|_{L_{- \mu}^2 (\T_L^2)} \big)^{1 + 2 \nu} \| v_{L, \eps} (t) \|_{\B_{2, 2, 4 \mu}^{1 - \nu} (\T_L^2)}^{p - 2 \nu}
\label{vH2bdd1-5}
\end{align}

\noi
and
\begin{align}
|\mathcal{G}_{L, \eps} [v_{0, L}]| \les_\om |\log \eps|^2 \big( \| v_{0, L} \|_{L_{- \mu}^2 (\T_L^2)} + \| \Dl v_{0, L} \|_{L_{- \mu}^2 (\T_L^2)} \big)^{1 + 2 \nu} \| v_{0, L} \|_{\B_{2, 2, 4 \mu}^{1 - \nu} (\T_L^2)}^{p - 2 \nu} 
\label{vH2bdd1-6}
\end{align}

\noi
for any $\mu > 0$ sufficiently small and $\nu > 0$ sufficiently small. Moreover, using Lemma~\ref{LEM:Hest} and Lemma~\ref{LEM:vL2L4} with $\nu > 0$ sufficiently small along with the embedding in Lemma~\ref{LEM:embL}~(iv), we have
\begin{align}
\begin{split}
\int_0^T |\mathcal{H}_{L, \eps} [v_{L, \eps} (t')]| dt' &\les_{\om, T} |\log \eps|^{11} \Big( 1 + \| v_{0, L} \|_{H_{\mu_0}^1 (\T_L^2)}^{2} \Big) \Big( 1 + \| v_{L, \eps} \|_{L_T^\infty \B_{2, 2, 24 \mu}^{1 - \nu} (\T_L^2)}^{3p} \Big) \\
&\quad \times \Big( 1 + \| \Dl v_{L, \eps} \|_{L_T^\infty L_{- \mu}^2 (\T_L^2)}^{\frac 74 + 4 p \nu} \Big) .
\end{split}
\label{vH2bdd1-7}
\end{align}

\noi
Thus, combining \eqref{vH2bdd1-1}, \eqref{vH2bdd1-2}, \eqref{vH2bdd1-3}, \eqref{vH2bdd1-4}, \eqref{vH2bdd1-5}, \eqref{vH2bdd1-6}, and \eqref{vH2bdd1-7} and using Lemma~\ref{LEM:embL}~(v) and (iv), Lemma~\ref{LEM:L2equiv}, and \eqref{H1H2}, we obtain
\begin{align}
\begin{split}
\| e^{- Y_{L, \eps}} \Dl v_{L, \eps} \|_{L_T^\infty L^2 (\T_L^2)}^2 &\les_{\om, T} \dl \| \Dl v_{L, \eps} \|_{L_T^\infty L_{- \mu}^2 (\T_L^2)}^2 +  \dl^{-7} |\log \eps|^{17} \Big( 1 + \| v_{0, L} \|_{H^2_{\mu_0} (\T_L^2)}^{p + 1} \Big) \\
&\quad \times \Big( 1 + \| v_{L, \eps} \|_{L_T^\infty \B_{2, 2, 24 \mu}^{1 - \nu} (\T_L^2)}^{3p} \Big)  \Big( 1 + \| \Dl v_{L, \eps} \|_{L_T^\infty L_{- \mu}^2 (\T_L^2)}^{\frac{15}{8}} \Big) .
\end{split}
\label{vH2bdd3}
\end{align}

\noi
Thus, by H\"older's inequality, Lemma~\ref{LEM:YregLe}~(iii), \eqref{vH2bdd3}, and Young's inequality, we get
\begin{align}
\begin{split}
\| \Dl v_{L, \eps} \|_{L_T^\infty L^2_{- \mu} (\T_L^2)}^2 &\leq \| e^{Y_{L, \eps}} \|_{L^\infty_{- \mu} (\T_L^2)}^2 \| e^{- Y_{L, \eps}} \Dl v_{L, \eps} \|_{L_T^\infty L^2 (\T_L^2)}^2 \\
&\les_{\om, T} \dl \| \Dl v_{L, \eps} \|_{L_T^\infty L_{- \mu}^2 (\T_L^2)}^2 \\
&\quad + \dl^{-127} |\log \eps|^{272} \Big( 1 + \| v_{0, L} \|_{H^2_{\mu_0} (\T_L^2)}^{16 (p + 1)} \Big) \Big( 1 + \| v_{L, \eps} \|_{L_T^\infty \B_{2, 2, 24 \mu}^{1 - \nu} (\T_L^2)}^{48 p} \Big)
\end{split}
\label{vH2bdd4} 
\end{align}

\noi
for any $0 < \dl < 1$. From the interpolation estimate in Lemma~\ref{LEM:interpL}, Lemma~\ref{LEM:embL}~(v) and (iv), Lemma~\ref{LEM:L2equiv}, Young's inequality, and Lemma~\ref{LEM:H1bdd}, we have
\begin{align}
\begin{split}
\| v_{L, \eps} \|_{L_T^\infty \B_{2, 2, 24 \mu}^{1 - \nu} (\T_L^2)} &\les \| v_{L, \eps} \|_{L_T^\infty \B_{2, 2, \frac{25 - \nu}{\nu} \mu}^{0} (\T_L^2)}^\nu \| v_{L, \eps} \|_{L_T^\infty \B_{2, 2, - \mu}^{1} (\T_L^2)}^{1 - \nu} \\
&\les \| v_{L, \eps} \|_{L_T^\infty L_{\frac{25 - \nu}{\nu} \mu}^{2} (\T_L^2)} + \| \nb v_{L, \eps} \|_{L_T^\infty L_{- \mu}^{2} (\T_L^2)} \\
&\les 1 + \| v_{0, L} \|_{H^1_{\mu_0} (\T_L^2)}^{\frac{p + 1}{2}} ,
\end{split}
\label{vH2bdd5}
\end{align}

\noi
where we need to take $\mu > 0$ to be sufficiently small.
Thus, from \eqref{vH2bdd4} and \eqref{vH2bdd5}, we take $\dl > 0$ to be sufficiently small to obtain the desired estimate.
\end{proof}

As a consequence of Lemma~\ref{LEM:H2bdd}, we have the following useful bound for $v_{L, \eps}$ with regularity $1 \leq s < 2$.

\begin{lemma}
\label{LEM:Hsbdd}
Let $p \geq 2$, $\ld > 0$, $L \geq 1$, $\mu_0 > 0$, and $v_{0, L} \in H_{\mu_0}^2 (\T_L^2)$. Given $0 < \eps < \frac 12$, let $v_{L, \eps}$ be the solution to the mollified $L$-periodic nonlinear DAM \eqref{vNLSLe} with $v_{L, \eps} |_{t = 0} = v_{0, L}$. Then, for any $T \geq 1$, $1 \leq s < 2$, $\nu > 0$, and $\mu > 0$ sufficiently small, there exists $\Om' \subset \Om$ with full probability measure such that for any $\om \in \Om'$ and $0 < \eps < \frac 12$, we have
\begin{align*}
\| v_{L, \eps} \|_{C_T \B_{2, 2, \mu}^s (\T_L^2)} \les_{\om, T} |\log \eps|^{136 (s - 1 + \nu)} \Big( 1 + \| v_{0, L} \|_{H_{\mu_0}^2 (\T_L^2)}^{21 (p + 1)^2} \Big) .
\end{align*}
\end{lemma}

\begin{proof}
From the interpolation estimate in Lemma~\ref{LEM:interpL}, Lemma~\ref{LEM:embL}~(v), Lemma~\ref{LEM:L2equiv}, and Lemma~\ref{LEM:H1bdd} along with Lemma~\ref{LEM:embL}~(iv), we have
\begin{align}
\begin{split}
\| v_{L, \eps} \|_{C_T \B_{2, 2, \mu'}^{1 - \nu} (\T_L^2)} &\les \| v_{L, \eps} \|_{C_T \B_{2, 2, \frac{2 - \nu}{\nu} \mu'}^0 (\T_L^2)}^\nu \| v_{L, \eps} \|_{C_T \B_{2, 2, - \mu'}^1 (\T_L^2)}^{1 - \nu} \\
&\les \| v_{L, \eps} \|_{C_T L_{\frac{2 - \nu}{\nu} \mu'}^2 (\T_L^2)}^\nu \Big( \| v_{L, \eps} \|_{C_T L^2_{- \mu'} (\T_L^2)} + \| \nb v_{L, \eps} \|_{C_T L^2_{- \mu'} (\T_L^2)} \Big)^{1 - \nu} \\
&\les_{\om, T} 1 + \| v_{0, L} \|_{H^1_{\mu_0} (\T_L^2)}^{\frac{p + 1}{2}} 
\end{split}
\label{B1v_bdd}
\end{align}

\noi
with $\mu' > 0$ sufficiently small. Thus, from the interpolation estimate in Lemma~\ref{LEM:interpL}, \eqref{B1v_bdd}, Lemma~\ref{LEM:embL}~(v), Lemma~\ref{LEM:L2equiv}, Lemma~\ref{LEM:H1bdd} along with Lemma~\ref{LEM:embL}~(iv), and Lemma~\ref{LEM:H2bdd}, we have
\begin{align*}
\| v_{L, \eps} \|_{C_T \B_{2, 2, \mu}^s (\T_L^2)} &\les \| v_{L, \eps} \|_{C_T \B_{2, 2, \mu'}^{1 - \nu} (\T_L^2)}^{\frac{2 - s}{1 + \nu}} \| v_{L, \eps} \|_{C_T \B_{2, 2, - \mu}^2 (\T_L^2)}^{\frac{s - 1 + \nu}{1 + \nu}} \\
&\les \Big( 1 + \| v_{0, L} \|_{H^1_{\mu_0} (\T_L^2)}^{\frac{p + 1}{2}} \Big)^{\frac{2 - s}{1 + \nu}} \Big( \| v_{L, \eps} \|_{C_T L^2_{- \mu'} (\T_L^2)} + \| \Dl v_{L, \eps} \|_{C_T L^2_{- \mu'} (\T_L^2)} \Big)^{\frac{s - 1 + \nu}{1 + \nu}} \\
&\les |\log \eps|^{136 (s - 1 + \nu)} \Big( 1 + \| v_{0, L} \|_{H^1_{\mu_0} (\T_L^2)}^{p + 1} \Big) \Big( 1 + \| v_{0, L} \|_{H^2_{\mu_0} (\T_L^2)}^{20 (p + 1)^2} \Big) ,
\end{align*}

\noi
where 
\begin{align*}
\mu' = \frac{s + 2 \nu}{2 - s} \mu 
\end{align*}

\noi
and we need to take $\mu > 0$ sufficiently small so that \eqref{B1v_bdd} holds. In view of \eqref{H1H2}, this gives the desired estimate.
\end{proof}

\subsection{Convergence estimate for the $L$-periodic nonlinear solution}
\label{SUB:GWP2}

In this subsection, we show Theorem~\ref{THM:GWP} in the case when $\ld > 0$, which corresponds to global well-posedness of the $L$-periodic nonlinear DAM \eqref{vNLSAndL}. In fact, this follows directly from the following proposition.
\begin{proposition}
\label{PROP:vLe_conv}
Let $p \geq 2$, $\ld > 0$, $L \geq 1$, $\mu_0 > 0$, and $v_{0, L} \in H_{\mu_0}^2 (\T_L^2)$. Given $0 < \eps < \frac 12$, let $v_{L, \eps}$ be the global-in-time solution to the mollified $L$-periodic nonlinear DAM \eqref{vNLSLe} with $v_{L, \eps} |_{t = 0} = v_{0, L}$. 
Then, there exists a process $v_L$ such that, for any $T \geq 1$ and $0 \leq s < 2$, there exist $\dl > 0$ and $\Om' \subset \Om$ with full probability measure such that for any $\om \in \Om'$ and $0 < \eps < \frac 12$, we have the difference estimate
\begin{align}
\| v_{L, \eps} - v_L \|_{C_T H^s (\T_L^2)} \les_{\om, T} \eps^\dl C \big( \| v_{0, L} \|_{H_{\mu_0}^2 (\T_L^2)} \big) 
\label{vLHs_diff}
\end{align}

\noi
for some constant $C ( \| v_{0, L} \|_{H_{\mu_0}^2 (\T_L^2)} ) > 0$. The process $v_L$ is a unique global-in-time solution to the $L$-periodic nonlinear DAM \eqref{vNLSAndL} with $v_L |_{t = 0} = v_{0, L}$ in $C(\R_+; H^s (\T_L^2))$ for any $1 < s < 2$. 
\end{proposition}

\begin{proof}
Let $\Om' \subset \Om$ be the event with full probability measure such that all the estimates in Section~\ref{SEC:sto} and Section~\ref{SEC:convL1} hold (with some parameters to be chosen below) and we fix $\om \in \Om'$.

Let $0 < \eps_2 < \eps_1 < \frac 12$ and $r_{L, \eps_1, \eps_2} \deff v_{L, \eps_1} - v_{L, \eps_2}$. Then, from \eqref{vNLSLe}, we see that $r_{L, \eps_1, \eps_2}$ satisfies the following equation with zero initial data:
\begin{align*}
i \dt r_{L, \eps_1, \eps_2} &= \Dl r_{L, \eps_1, \eps_2} - 2 \nb Y_{L, \eps_1} \cdot \nb r_{L, \eps_1, \eps_2} + \wt{\wick{ |\nb Y_{L, \eps_1}|^2 }} \, r_{L, \eps_1, \eps_2} \\
&\quad + \big( \wt{\wick{ |\nb Y_{L, \eps_1}|^2 }} - \wt{\wick{ |\nb Y_{L, \eps_2}|^2 }} \big) v_{L, \eps_2} - 2 (\nb Y_{L, \eps_1} - \nb Y_{L, \eps_2}) \cdot \nb v_{L, \eps_2} \\
&\quad - \ld \big( e^{- (p - 1) Y_{L, \eps_1}} |v_{L, \eps_1}|^{p - 1} v_{L, \eps_1} - e^{- (p - 1) Y_{L, \eps_2}} |v_{L, \eps_2}|^{p - 1} v_{L, \eps_2} \big) .
\end{align*}

\noi
Let $t \in [0, T]$. By using the equation for $r_{L, \eps_1, \eps_2}$, we compute that
\begin{align}
\frac 12 \frac{d}{dt} \int_{\T_L^2} e^{- 2 Y_{L, \eps_1}} |r_{L, \eps_1, \eps_2} (t)|^2 dx = \textup{I}_1 + \textup{I}_2 + \textup{I}_3 + \textup{I}_4 ,
\label{rLe0}
\end{align}

\noi
where
\begin{align*}
\textup{I}_1 &= -2 \Im \int_{\T_L^2} e^{- 2 Y_{L, \eps_1}} (\nb Y_{L, \eps_1} - \nb Y_{L, \eps_2}) \cdot \nb v_{L, \eps_2} (t) \cj{r_{L, \eps_1, \eps_2} (t)} dx , \\
\textup{I}_2 &= \Im \int_{\T_L^2} e^{- 2 Y_{L, \eps_1}} \big( \wt{\wick{ |\nb Y_{L, \eps_1}|^2 }} - \wt{\wick{ |\nb Y_{L, \eps_2}|^2 }} \big) v_{L, \eps_2} (t) \cj{r_{L, \eps_1, \eps_2} (t)} dx , \\
\textup{I}_3 &= - \ld \Im \int_{\T_L^2} e^{- 2 Y_{L, \eps_1}} \big( e^{- (p - 1) Y_{L, \eps_1}} - e^{- (p - 1) Y_{L, \eps_2}} \big) |v_{L, \eps_2} (t)|^{p - 1} v_{L, \eps_2} (t) \cj{r_{L, \eps_1, \eps_2} (t)} dx , \\
\textup{I}_4 &= - \ld \Im \int_{\T_L^2} e^{- 2 Y_{L, \eps_1}} e^{- (p - 1) Y_{L, \eps_1}} \big( |v_{L, \eps_1} (t)|^{p - 1} v_{L, \eps_1} (t) - |v_{L, \eps_2} (t)|^{p - 1} v_{L, \eps_2} (t) \big) \cj{r_{L, \eps_1, \eps_2} (t)} dx .
\end{align*}

\noi
Let $0 < s_0 < \frac{1}{100}$, $0 < \dl_0 < s_0$, and small $\mu > 0$ be fixed.  For $\textup{I}_1$, using the duality estimate in Lemma~\ref{LEM:dualL}, Lemma~\ref{LEM:YregLe}~(i), the product estimates (Lemma~\ref{LEM:prodL}~(ii)), Lemma~\ref{LEM:YregLe}~(iii), and Lemma~\ref{LEM:embL}~(v), we have
\begin{align}
\begin{split}
|\textup{I}_1| &\les \| \nb Y_{L, \eps_1} - \nb Y_{L, \eps_2} \|_{\C_{- \mu}^{- s_0} (\T_L^2)} \|  e^{- 2 Y_{L, \eps_1}} \nb v_{L, \eps_2} (t) \cj{r_{L, \eps_1, \eps_2} (t)} \|_{\B^{s_0}_{1, 1, \mu} (\T_L^2)} \\
&\les_{\om} \eps_1^{\dl_0} \| e^{- 2 Y_{L, \eps_1}} \|_{\C_{- \mu}^{3 s_0} (\T_L^2)} \| \nb v_{L, \eps_2} (t) \|_{\B_{2, 2, \mu}^{3 s_0} (\T_L^2)} \| r_{L, \eps_1, \eps_2} (t) \|_{\B_{2, 2, \mu}^{2 s_0} (\T_L^2)} \\
&\les_{\om} \eps_1^{\dl_0} \| v_{L, \eps_2} (t) \|_{\B_{2, 2, \mu}^{1 + 3 s_0} (\T_L^2)} \big( \| v_{L, \eps_1} (t) \|_{\B_{2, 2, \mu}^{2 s_0} (\T_L^2)} + \| v_{L, \eps_2} (t) \|_{\B_{2, 2, \mu}^{2 s_0} (\T_L^2)} \big) .
\end{split}
\label{rLe1}
\end{align}

\noi
For $\textup{I}_2$, using the duality estimate in Lemma~\ref{LEM:dualL}, Lemma~\ref{LEM:YregLe}~(ii) and (iv), the product estimates (Lemma~\ref{LEM:prodL}~(ii)), and Lemma~\ref{LEM:YregLe}~(iii), we have
\begin{align}
\begin{split}
|\textup{I}_2| &\les \big\| \wt{ \wick{|\nb Y_{L, \eps_1}|^2} } - \wt{ \wick{|\nb Y_{L, \eps_2}|^2} } \big\|_{\C_{- \mu}^{- s_0} (\T_L^2)} \|  e^{- 2 Y_{L, \eps_1}} v_{L, \eps_2} (t) \cj{r_{L, \eps_1, \eps_2} (t)} \|_{\B^{s_0}_{1, 1, \mu} (\T_L^2)} \\
&\les_{\om} \eps_1^\dl \| e^{- 2 Y_{L, \eps_1}} \|_{\C_{- \mu}^{3 s_0} (\T_L^2)} \| v_{L, \eps_2} (t) \|_{\B_{2, 2, \mu}^{3 s_0} (\T_L^2)} \| r_{L, \eps_1, \eps_2} (t) \|_{\B_{2, 2, \mu}^{2 s_0} (\T_L^2)} \\
&\les_{\om} \eps_1^\dl \| v_{L, \eps_2} (t) \|_{\B_{2, 2, \mu}^{3 s_0} (\T_L^2)} \big( \| v_{L, \eps_1} (t) \|_{\B^{2 s_0}_{2, 2, \mu} (\T_L^2)} + \| v_{L, \eps_2} (t) \|_{\B^{2 s_0}_{2, 2, \mu} (\T_L^2)} \big) .
\end{split}
\label{rLe2}
\end{align}

\noi
For $\textup{I}_3$, using H\"older's inequality along with Lemma~\ref{LEM:Lwei}~(iii) and (iv), Lemma~\ref{LEM:YregLe}~(iii), and the embeddings in Lemma~\ref{LEM:embL}~(ii), (i), and (iii), we have
\begin{align}
\begin{split}
|\textup{I}_3| &\les \| e^{- 2 Y_{L, \eps_1}} \|_{L_{- \mu}^\infty (\T_L^2)} \big\| e^{- (p - 1) Y_{L, \eps_1}} - e^{- (p - 1) Y_{L, \eps_2}} \big\|_{L_{- p \mu}^\infty (\T_L^2)} \\
&\qquad \times \| v_{L, \eps_2} (t) \|_{L_{\mu}^{2 p} (\T_L^2)}^p \| r_{L, \eps_1, \eps_2} (t) \|_{L_\mu^2 (\T_L^2)} \\
&\les_{\om} \eps_1^\dl \| v_{L, \eps_2} (t) \|_{\B^1_{2, 2, \mu} (\T_L^2)}^p \big( \| v_{L, \eps_1} (t) \|_{L_\mu^2 (\T_L^2)} + \| v_{L, \eps_2} (t) \|_{L_\mu^2 (\T_L^2)} \big) .
\end{split}
\label{rLe3}
\end{align}

\noi
For $\textup{I}_4$, using the mean value theorem, H\"older's inequality along with Lemma~\ref{LEM:Lwei}~(iii) and (iv), Lemma~\ref{LEM:YregLe}~(iii), and Lemma~\ref{LEM:embL}~(ii), (i), and (iii), we have
\begin{align}
\begin{split}
|\textup{I}_4| &\les \| e^{- (p - 1) Y_{L, \eps}} \|_{L_{- (p - 1) \mu}^\infty (\T_L^2)} \Big( \| v_{L, \eps_1} (t) \|_{L_\mu^\infty (\T_L^2)}^{p - 1} + \| v_{L, \eps_2} (t) \|_{L_\mu^\infty (\T_L^2)}^{p - 1} \Big) \\
&\quad \times \| e^{- Y_{L, \eps_1}} r_{L, \eps_1, \eps_2} (t) \|_{L^2 (\T_L^2)}^2 \\
&\les_{\om} \Big( \| v_{L, \eps_1} (t) \|_{\B_{2, 2, \mu}^{1 + s_0} (\T_L^2)}^{p - 1} + \| v_{L, \eps_2} (t) \|_{\B_{2, 2, \mu}^{1 + s_0} (\T_L^2)}^{p - 1} \Big) \| e^{- Y_{L, \eps_1}} r_{L, \eps_1, \eps_2} (t) \|_{L^2 (\T_L^2)}^2 .
\end{split}
\label{rLe4}
\end{align}

\noi
Combining \eqref{rLe0}, \eqref{rLe1}, \eqref{rLe2}, \eqref{rLe3}, and \eqref{rLe4} and using Lemma~\ref{LEM:Hsbdd} (with $\mu > 0$ sufficiently small), we obtain
\begin{align*}
\frac 12 &\frac{d}{dt} \int_{\T_L^2} e^{-2 Y_{L, \eps_1}} |r_{L, \eps_1, \eps_2} (t)|^2 dx \\
&\les_{\om} \eps_1^\dl |\log \eps_2| C \big(\| v_{0, L} \|_{H_{\mu_0}^2 (\T_L^2)}\big) + |\log \eps_2|^{272 (p - 1) s_0} C \big( \| v_{0, L} \|_{H_{\mu_0}^2 (\T_L^2)} \big) \int_{\T_L^2} e^{-2 Y_{L, \eps_1}} |r_{L, \eps_1, \eps_2} (t)|^2 dx 
\end{align*}

\noi
for some constant $C ( \| v_{0, L} \|_{H_{\mu_0}^2 (\T_L^2)} ) > 0$. By Gronwall's inequality, we get
\begin{align*}
\begin{split}
\sup_{t \in [0, T]} \int_{\T_L^2} e^{-2 Y_{L, \eps_1}} |r_{L, \eps_1, \eps_2} (t)|^2 dx &\les_{\om} \eps_1^\dl |\log \eps_2| C\big( \| v_{0, L} \|_{H_{\mu_0}^2 (\T_L^2)} \big) \\
&\quad \times \exp \big( |\log \eps_2|^{272 (p - 1) s_0} C\big(\| v_{0, L} \|_{H_{\mu_0}^2 (\T_L^2)}\big) T \big) ,
\end{split}
\end{align*}

\noi
so that by H\"older's inequality and Lemma~\ref{LEM:YregLe}~(iii), we obtain
\begin{align}
\begin{split}
\| v_{L, \eps_1} - v_{L, \eps_2} \|_{C_T L_{- \mu}^2 (\T_L^2)}^2 &\les_{\om} \eps_1^\dl |\log \eps_2| C\big( \| v_{0, L} \|_{H_{\mu_0}^2 (\T_L^2)} \big) \\
&\quad \times \exp \big( |\log \eps_2|^{272 (p - 1) s_0} C\big(\| v_{0, L} \|_{H_{\mu_0}^2 (\T_L^2)}\big) T \big) .
\end{split}
\label{rLe_gron}
\end{align}

Let us now take $\eps_1 = 2^{-k}$ and $\eps_2 = 2^{- (k + 1)}$ for $k \in \N$. Then, by using \eqref{rLe_gron} and taking $s_0$ to be sufficiently small so that $272 (p - 1) s_0  < 1$, we get
\begin{align}
\begin{split}
\| v_{L, 2^{-k}} - v_{L, 2^{- (k + 1)}} \|_{C_T L_{- \mu}^2 (\T_L^2)}^2 &\les_{\om} 2^{- \dl k} (k + 1) C\big(\| v_{0, L} \|_{H_{\mu_0}^2 (\T_L^2)}\big) \\
&\quad \times \exp \big( (k + 1)^{272 (p - 1) s_0} C\big(\| v_{0, L} \|_{H_{\mu_0}^2 (\T_L^2)}\big) T \big) \\
&\les_{\om, T} 2^{- \frac{\dl k}{4}} C\big(\| v_{0, L} \|_{H_{\mu_0}^2 (\T_L^2)}\big) .
\end{split}
\label{vLe_diff2}
\end{align}

\noi
Let $0 \leq s < s_1 < 2$. By \eqref{Hs_equi} and the interpolation estimate in Lemma~\ref{LEM:interpL}, we have
\begin{align*}
\| &v_{L, 2^{-k}} - v_{L, 2^{- (k + 1)}} \|_{C_T H^s (\T_L^2)} \\
&\sim \| v_{L, 2^{-k}} - v_{L, 2^{- (k + 1)}} \|_{C_T \B_{2, 2}^s (\T_L^2)} \\
&\les \| v_{L, 2^{-k}} - v_{L, 2^{- (k + 1)}} \|_{C_T \B_{2, 2, - \mu}^0 (\T_L^2)}^{1 - \frac{s}{s_1}}  \| v_{L, 2^{-k}} - v_{L, 2^{- (k + 1)}} \|_{C_T \B_{2, 2, \mu'}^{s_1} (\T_L^2)}^{\frac{s}{s_1}} 
\end{align*}

\noi
with
\begin{align*}
\mu' = \frac{1 - \frac{s}{s_1}}{\frac{s}{s_1}} \mu ,
\end{align*}

\noi
so that by further using Lemma~\ref{LEM:L2equiv}, Lemma~\ref{LEM:Hsbdd} (with $\mu > 0$ sufficiently small), and \eqref{vLe_diff2}, we get
\begin{align}
\begin{split}
\| &v_{L, 2^{-k}} - v_{L, 2^{- (k + 1)}} \|_{C_T H^s (\T_L^2)} \\
&\les_{\om, T} \| v_{L, 2^{-k}} - v_{L, 2^{- (k + 1)}} \|_{C_T L_{- \mu}^2 (\T_L^2)}^{1 - \frac{s}{s_1}} (k + 1)^{\frac{136 s}{s_1}} C\big(\| v_{0, L} \|_{H_{\mu_0}^2 (\T_L^2)}\big) \\
&\les_{\om, T} 2^{- \dl_1 k} C\big(\| v_{0, L} \|_{H_{\mu_0}^2 (\T_L^2)}\big)
\end{split}
\label{vLe_conv1}
\end{align}
 
\noi
for some $\dl_1 > 0$. This shows that $\{ v_{L, 2^{-k}} \}_{k \in \N}$ is almost surely a Cauchy sequence in $C([0, T]; H^s (\T_L^2))$ and so converges to some process $v_L$ with values in $C([0, T]; H^s (\T_L^2))$. Similarly, we have
\begin{align}
\sup_{\eps \in (2^{- (k + 1)}, 2^{-k}]} \| v_{L, \eps} - v_{L, 2^{-k}} \|_{C_T H^s (\T_L^2)} \les_{\om, T} 2^{- \dl_2 k} C\big( \| v_{0, L} \|_{H_{\mu_0}^2 (\T_L^2)} \big)
\label{vLe_conv2}
\end{align}

\noi
for some $\dl_2 > 0$, so that the whole sequence $\{ v_{L, \eps} \}_{\eps \in (0, \frac 12)}$ converges to $v_L$ almost surely in $C([0, T]; H^s (\T_L^2))$ as $\eps \to 0$. The difference estimate \eqref{vLHs_diff} then follows immediately from \eqref{vLe_conv1} and \eqref{vLe_conv2}. 
The uniqueness part follows from the same analysis as in previous works \cite{DW, TzV23} and so we omit details.
\end{proof}

\subsection{Large torus convergence of mollified $L$-periodic nonlinear solutions}

In this subsection, we show that the solution $v_{L, \eps}$ to the mollified $L$-periodic nonlinear DAM \eqref{vNLSLe} converges to that of the following mollified version of the nonlinear DAM \eqref{vNLSAnd} on $\R^2$:
\begin{align}
\begin{cases}
i \dt v_\eps = \Dl v_\eps - 2 \nb Y_\eps \cdot \nb v_\eps + \wt{\wick{|\nb Y_\eps|^2}} \, v_\eps - \ld e^{- (p - 1) Y_\eps} |v_\eps|^{p - 1} v_\eps \\
v_\eps |_{t = 0} = v_0 ,
\end{cases}
\label{vNLSe}
\end{align}

\noi
where $Y_\eps$ is as defined in \eqref{Ye} and $\wt{\wick{|\nb Y_\eps|^2}}$ is as defined in \eqref{Ye22}.
From \cite[Theorem~1.1]{DLTV}, we know that if $v_0 \in H^2_\mu (\R^2)$ for any $\mu > 0$, then \eqref{vNLSe} admits a unique global-in-time solution.

We first mention the following result regarding the weighted $H^s (\R^2)$-bound of $v_\eps$, which follows from \cite[Theorem~1.1 and Theorem~1.2]{DLTV}.
\begin{lemma}
\label{LEM:veHs}
Let $\mu_0 > 0$ and $v_0 \in H_{\mu_0}^2 (\R^2)$. Given $0 < \eps < \frac 12$, let $w_\eps$ be the global-in-time solution to the mollified linear DAM \eqref{vNLSe} with $v_\eps |_{t = 0} = v_0$. Then, for any $T \geq 1$, $0 \leq s < 2$, and $\mu > 0$ sufficiently small, there exists $\Om' \subset \Om$ will full probability measure such that for any $\om \in \Om'$ and $0 < \eps < \frac 12$, we have the bound
\begin{align*}
\| v_\eps \|_{C_T H_\mu^s (\R^2)} \les_{\om, T} C \big( \| v_0 \|_{H^2_{\mu_0} (\R^2)} \big)
\end{align*}

\noi
for some constant $C \big( \| v_0 \|_{H^2_{\mu_0} (\R^2)} \big) > 0$.
\end{lemma}

We now show the following proposition regarding the convergence of the solution $v_{L, \eps}$ of \eqref{vNLSLe} to the solution $v_\eps$ of \eqref{vNLSe}. We recall the function $\s_R$ in \eqref{sigR}.

\begin{proposition}
\label{PROP:vLeve}
Let $L \geq 1$, $p \geq 2$, $\ld > 0$, $\mu_0 > 0$, $v_0 \in H_{\mu_0}^2 (\R^2)$, and $v_{0, L} \in H_{\mu_0}^2 (\T_L^2)$. Given $0 < \eps < \frac 12$, let $v_{\eps}$ be the global-in-time solution to the mollified $L$-periodic nonlinear DAM \eqref{vNLSLe} with $v_{\eps} |_{t = 0} = v_{0}$ and let $v_{L, \eps}$ be the global-in-time solution to the mollified nonlinear DAM \eqref{vNLSe} with $v_{L, \eps}|_{t = 0} = v_{0, L}$. Then, for any $T \geq 1$, $1 \leq R \leq L$, $\mu > 0$ sufficiently small, $\kappa > 0$, and $0 < b \leq 4$, there exists $\Om' \subset \Om$ with full probability measure such that for any $\om \in \Om'$ and $0 < \eps < \frac 12$, we have the estimate
\begin{align*}
\| &\s_R^{-2} (v_{L, \eps} - v_\eps) \|_{C_T L^2_{- \mu} (\R^2)} \\
&\les_{\om, T} \eps^{- \kappa} \big( \| \s_R^{-2} (v_0 - v_{0, L}) \|_{C_T L_\mu^2 (\R^2)} + R^{-1} + R^{b} L^{- \frac{b}{8}} |\log \eps| \big) C \big( \| v_{0, L} \|_{H^2_{\mu_0} (\T_L^2)}, \| v_0 \|_{H_{\mu_0}^2 (\R^2)} \big)
\end{align*}

\noi
for some constant $C (\| v_{0, L} \|_{H^2_{\mu_0} (\T_L^2)}, \| v_0 \|_{H_{\mu_0}^2 (\R^2)}) > 0$.
\end{proposition}

\begin{proof}
Given $L \geq 1$ and $0 < \eps < \frac 12$, we let $r_{L, \eps} = v_\eps - v_{L, \eps}$. Then, we see that $r_{L, \eps}$ satisfies the following equation:
\begin{align*}
i \dt r_{L, \eps} &= \Dl r_{L, \eps} - 2 \nb Y_{L, \eps} \cdot \nb r_{L, \eps} + \wt{\wick{ |\nb Y_{L, \eps}|^2 }} \, r_{L, \eps}  - 2 (\nb Y_\eps - \nb Y_{L, \eps}) \cdot \nb v_{\eps} \\
&\quad + \big( \wt{\wick{|\nb Y_\eps|^2}} - \wt{\wick{|\nb Y_{L, \eps}|^2}} \big) v_{\eps} - \ld \big( e^{- (p - 1) Y_\eps} |v_\eps|^{p - 1} v_\eps - e^{- (p - 1) Y_{L, \eps}} |v_{L, \eps}|^{p - 1} v_{L, \eps} \big) .
\end{align*}

\noi
Let $\Omega' \subset \Omega$ be the event with full probability measure such that all (already established) estimates in Section~\ref{SEC:sto} and Section~\ref{SEC:convL1} hold and we fix $\om \in \Omega'$.
Note that from \eqref{sigR1}, Lemma~\ref{LEM:Lwei}~(vi), (iii), and (i), Lemma~\ref{LEM:sig}, Lemma~\ref{LEM:L2equiv}, and Lemma~\ref{LEM:embL}~(v), we have
\begin{align*}
\big\| \s_R^{-1} \jbb{\cdot}_{L, \mu} v_{L, \eps} \big\|_{C_T H^1 (\R^2)} 
&\les \big\| \s_R^{-1} \jbb{\cdot}_{L, \mu} v_{L, \eps} \big\|_{C_T L^2 (\R^2)} + \big\| \s_R^{-1} \jbb{\cdot}_{L, \mu} \nb v_{L, \eps} \big\|_{C_T L^2 (\R^2)} \\
&\les \| v_{L, \eps} \|_{C_T L^2_\mu (\T_L^2)} + \| \nb v_{L, \eps} \|_{C_T L^2_\mu (\T_L^2)} \\
&\les \| v_{L, \eps} \|_{C_T \B_{2, 2, \mu}^1 (\T_L^2)}
\end{align*}

\noi
for any $\mu > 0$ sufficiently small,
so that by applying Lemma~\ref{LEM:Hsbdd} on the last expression, we get
\begin{align}
\big\| \s_R^{-1} \jbb{\cdot}_{L, \mu} v_{L, \eps} \big\|_{C_T H^1 (\R^2)} 
\les_{\om, T} C \big( \| v_{0, L} \|_{H_{\mu_0}^2 (\T_L^2)} \big)
\label{RvLe_bdd}
\end{align}

\noi
for some constant $C ( \| v_{0, L} \|_{H_{\mu_0}^2 (\T_L^2)} ) > 0$.

Let $t \in [0, T]$. By using the equation and integration by parts, we obtain
\begin{align}
\frac 12 \frac{d}{dt} \int_{\R^2} \s_R^{-4} e^{- 2 Y_{L, \eps}} |r_{L, \eps} (t)|^2 dx 
= \textup{I}_1 + \textup{I}_2 + \textup{I}_3 + \textup{I}_4 + \textup{I}_5 ,
\label{vL0}
\end{align}

\noi
where
\begin{align*}
\textup{I}_1 &= - \Im \int_{\R^2} e^{-2 Y_{L, \eps}} \nb (\s_R^{-4}) \cdot \nb r_{L, \eps} (t) \cj{r_{L, \eps} (t)} dx , \\
\textup{I}_2 &= - 2 \Im \int_{\R^2} \s_R^{-4} e^{- 2 Y_{L, \eps}} (\nb Y_\eps - \nb Y_{L, \eps}) \cdot \nb v_{\eps} (t) \cj{r_{L, \eps} (t)} dx , \\
\textup{I}_3 &= \Im \int_{\R^2} \s_R^{-4} e^{- 2 Y_{L, \eps}} \big( \wt{\wick{|\nb Y_\eps|^2}} - \wt{\wick{|\nb Y_{L, \eps}|^2}} \big) v_{\eps} (t) \cj{r_{L, \eps} (t)} dx , \\
\textup{I}_4 &= - \ld \Im \int_{\R^2} \s_R^{-4} e^{- 2 Y_{L, \eps}} \big( e^{- (p - 1) Y_\eps} - e^{- (p - 1) Y_{L, \eps}} \big) |v_\eps (t)|^{p - 1} v_\eps (t) \cj{r_{L, \eps} (t)} dx , \\
\textup{I}_5 &= - \ld \Im \int_{\R^2} \s_R^{-4} e^{- 2 Y_{L, \eps}} e^{- (p - 1) Y_{L, \eps}} \big( |v_\eps (t)|^{p - 1} v_\eps (t) - |v_{L, \eps} (t)|^{p - 1} v_{L, \eps} (t) \big) \cj{r_{L, \eps} (t)} dx .
\end{align*}

\noi  
For $\textup{I}_1$, by \eqref{sigR1}, Cauchy's inequality, H\"older's inequality along with Lemma~\ref{LEM:Lwei}~(iii), Lemma~\ref{LEM:YregLe}~(iii), \eqref{RvLe_bdd}, and Lemma~\ref{LEM:veHs}, we have
\begin{align}
\begin{split}
|\textup{I}_1| &\les R^{-2} \int_{\R^2} \s_R^{-4} e^{-2 Y_{L, \eps}} |\nb r_{L, \eps} (t)|^2 dx + \int_{\R^2} \s_R^{-4} e^{- 2 Y_{L, \eps}} |r_{L, \eps} (t)|^2 dx \\
&\leq R^{-2} \big\| \jbb{\cdot}_{L, -\mu} e^{- 2 Y_{L, \eps}} \big\|_{L^\infty (\R^2)} \big\| \s_R^{-1} \jbb{\cdot}_{L, \mu} \nb r_{L, \eps} (t) \big\|_{L^2 (\R^2)}^2 + \int_{\R^2} \s_R^{-4} e^{- 2 Y_{L, \eps}} |r_{L, \eps} (t)|^2 dx \\
&\les_{\om, T} R^{-2} C \big( \| v_{0, L} \|_{H^2_{\mu_0} (\T_L^2)}, \| v_0 \|_{H_{\mu_0}^2 (\R^2)} \big) + \int_{\R^2} \s_R^{-4} e^{- 2 Y_{L, \eps}} |r_{L, \eps} (t)|^2 dx 
\end{split}
\label{vL1}
\end{align}

\noi
for $\mu > 0$ sufficiently small and some constant $C ( \| v_{0, L} \|_{H^2_{\mu_0} (\T_L^2)}, \| v_0 \|_{H_{\mu_0}^2 (\R^2)} ) > 0$. For $\textup{I}_2$, by H\"older's inequality along with Lemma~\ref{LEM:Lwei}~(iii), Lemma~\ref{LEM:YregLe}~(iii), Proposition~\ref{PROP:YconvL}~(i), Lemma~\ref{LEM:veHs}, and Sobolev's inequality, we have
\begin{align}
\begin{split}
|\textup{I}_2| &\les \big\| \jbb{\cdot}_{L, - \mu} e^{-2 Y_{L, \eps}} \big\|_{L^\infty (\R^2)} \| \s_R^{-1} (\nb Y_\eps - \nb Y_{L, \eps}) \|_{L^q (\R^2)} \\
&\quad \times \| \nb v_\eps (t) \|_{L^2 (\R^2)} \big\| \s_R^{-1} \jbb{\cdot}_{L, \mu} r_{L, \eps} (t) \big\|_{L^{\frac{2 q}{q - 2}} (\R^2)} \\
&\les_{\om, T} R^{b} L^{- \frac{b}{2}} |\log \eps| C \big( \| v_0 \|_{H_{\mu_0}^2 (\R^2)} \big) \big\| \s_R^{-1} \jbb{\cdot}_{L, \mu} r_{L, \eps} (t) \big\|_{H^1 (\R^2)} 
\end{split}
\label{vL2}
\end{align}

\noi 
for any $0 < b \leq 4$ and $2 < q < \infty$ satisfying $b q > 4$.
For $\textup{I}_3$, by H\"older's inequality along with Lemma~\ref{LEM:Lwei}~(iii), Lemma~\ref{LEM:YregLe}~(iii), Proposition~\ref{PROP:YconvL}~(iv), Lemma~\ref{LEM:veHs}, and Sobolev's inequality, we have
\begin{align}
\begin{split}
|\textup{I}_3| &\les \big\| \jbb{\cdot}_{L, - \mu} e^{-2 Y_{L, \eps}} \big\|_{L^\infty (\R^2)} \Big\| \s_R^{-2} \big( \wt{\wick{|\nb Y_\eps|^2}} - \wt{\wick{|\nb Y_{L, \eps}|^2}} \big) \Big\|_{L^q (\R^2)} \\
&\quad \times \| v_\eps (t) \|_{L^2 (\R^2)} \big\| \s_R^{-1} \jbb{\cdot}_{L, \mu} r_{L, \eps} (t) \big\|_{L^{\frac{2q}{q - 2}} (\R^2)} \\
&\les_{\om, T} R^{2b} L^{- \frac{b}{2} + \frac{2}{q}} |\log \eps|^2 C \big( \| v_0 \|_{H^2_{\mu_0} (\R^2)} \big) \big\| \s_R^{-1} \jbb{\cdot}_{L, \mu} r_{L, \eps} (t) \big\|_{H^1 (\R^2)} .
\end{split}
\label{vL3}
\end{align}

\noi
For $\textup{I}_4$, by H\"older's inequality along with Lemma~\ref{LEM:Lwei}~(iii), Lemma~\ref{LEM:YregLe}~(iii), Proposition~\ref{PROP:YconvL}~(iii), and Sobolev's inequalities, we have
\begin{align}
\begin{split}
|\textup{I}_4| &\les \big\| \jbb{\cdot}_{L, - \mu} e^{-2 Y_{L, \eps}} \big\|_{L^\infty (\R^2)} \| \s_R^{-2} (e^{a Y_\eps} - e^{a Y_{L, \eps}}) \|_{L^q (\R^2)} \\
&\quad \times \| v_\eps (t) \|_{L^{2p} (\R^2)}^p \big\| \s_R^{-1} \jbb{\cdot}_{L, \mu} r_{L, \eps} (t) \big\|_{L^{\frac{2q}{q - 2}} (\R^2)} \\
&\les_{\om, T} R^{2 b} L^{- \frac{b}{4}} \| v_\eps (t) \|_{H^1 (\R^2)}^p \big\| \s_R^{-1} \jbb{\cdot}_{L, \mu} r_{L, \eps} (t) \big\|_{H^1 (\R^2)} .
\end{split}
\label{vL4}
\end{align}

\noi
For $\textup{I}_5$, by the mean value theorem, H\"older's inequality along with Lemma~\ref{LEM:Lwei}~(iii) and (iv), Lemma~\ref{LEM:YregLe}~(iii), Sobolev's embedding, and the embeddings in Lemma~\ref{LEM:embL}~(ii), (i), and (iii), we have
\begin{align}
\begin{split}
|\textup{I}_5| &\les \big\| \jbb{\cdot}_{L, - (p - 1) \mu} e^{- (p - 1) Y_{L, \eps}} \big\|_{L^\infty (\R^2)} \Big( \big\| \jbb{\cdot}_{L, \mu} v_\eps (t) \big\|_{L^\infty (\R^2)}^{p - 1} + \big\| \jbb{\cdot}_{L, \mu} v_{L, \eps} (t) \big\|_{L^\infty (\R^2)}^{p - 1} \Big) \\
&\quad \times \int_{\R^2} \s_R^{-4} e^{- 2 Y_{L, \eps}} |r_{L, \eps} (t)|^2 dx \\
&\les_\om \Big( \| v_\eps (t) \|_{H^{1 + s_0}_\mu (\R^2)}^{p - 1} + \| v_{L, \eps} (t) \|_{\B_{2, 2, \mu}^{1 + s_0} (\T_L^2)}^{p - 1} \Big) \int_{\R^2} \s_R^{-4} e^{- 2 Y_{L, \eps}} |r_{L, \eps} (t)|^2 dx 
\end{split}
\label{vL5}
\end{align}

\noi
for any $s_0 > 0$. Combining \eqref{vL0}, \eqref{vL1}, \eqref{vL2}, \eqref{vL3} (with $q$ large enough), \eqref{vL4}, and \eqref{vL5} and using \eqref{RvLe_bdd}, Lemma~\ref{LEM:veHs}, and Lemma~\ref{LEM:Hsbdd}, we obtain
\begin{align*}
&\frac{d}{dt} \int_{\R^2} \s_R^{-4} e^{- 2 Y_{L, \eps}} |r_{L, \eps} (t)|^2 dx \\
&\les_{\om, T} C_v \big( R^{-2} + R^{2b} L^{- \frac{b}{4}} |\log \eps|^2 \big) + C_v |\log \eps|^{272 (p - 1) s_0} \int_{\R^2} \s_R^{-4} e^{- 2 Y_{L, \eps}} |r_{L, \eps} (t)|^2 dx 
\end{align*}

\noi
for some constant constant $C_v = C ( \| v_{0, L} \|_{H^2_{\mu_0} (\T_L^2)}, \| v_0 \|_{H^2_{\mu_0} (\R^2)} ) > 0$. Using Gronwall's inequality, we get
\begin{align}
\begin{split}
&\sup_{t \in [0, T]} \int_{\R^2} \s_R^{-4} e^{- 2 Y_{L, \eps}} |r_{L, \eps} (t)|^2 dx \\
&\les_{\om, T} \! \bigg( \int_{\R^2} \s_R^{-4} e^{- 2 Y_{L, \eps}} |v_0 - v_{0, L}|^2 dx +  R^{-2} + R^{2b} L^{- \frac{b}{4}} |\log \eps|^2 \bigg) \\
&\quad \times \exp \big( C_v |\log \eps|^{272 (p - 1) s_0} T \big) .
\end{split}
\label{vL6}
\end{align}

\noi
We take $s_0 > 0$ to be small enough such that $272 (p - 1) s_0 < 1$, so that we have $$\exp (|\log \eps|^{272(p - 1) s_0}) \les_{\kappa_0} \eps^{-\kappa_0}$$ 

\noi
for any $\kappa > 0$. Therefore, by  H\"older's inequality, \eqref{vL6}, and Lemma~\ref{LEM:YregLe}~(iii) along with the fact that $\jb{\cdot}^{- \mu} \les \jbb{\cdot}_{L, \mu}^{-1} \les \jbb{\cdot}_{L, - \mu}$ from Lemma~\ref{LEM:Lwei}~(ii) and (iii), we obtain
\begin{align*}
\| &\s_R^{-2} r_{L, \eps} \|_{C_T L_{- \mu}^2 (\R^2)} \\
&\leq \| \s_R^{-2} e^{- Y_{L, \eps}} r_{L, \eps} \|_{C_T L^2 (\R^2)} \| e^{Y_{L, \eps}} \|_{L_{- \mu}^\infty (\R^2)} \\
&\les_{\om, T, \kappa_0} \eps^{- C_v T \kappa_0} \big( \| \s_R^{-2} e^{- Y_{L, \eps}} (v_0 - v_{0, L}) \|_{C_T L^2 (\R^2)} + R^{-1} + R^{b} L^{- \frac{b}{8}} |\log \eps| \big) 
\end{align*}

\noi
for any $\kappa_0 > 0$ and $0 < b \leq 4$, so that by repeating the above procedure again, we get
\begin{align}
\| \s_R^{-2} r_{L, \eps} \|_{C_T L_{- \mu}^2 (\R^2)} \les_{\om, T, \kappa_0} \eps^{- C_v T \kappa_0} \big( \| \s_R^{-2} (w_0 - w_{0, L}) \|_{C_T L_\mu^2 (\R^2)} + R^{-1} + R^{b} L^{- \frac{b}{8}} |\log \eps| \big) .
\label{vL7}
\end{align}

\noi
The desired estimate then follows from \eqref{vL7} with $\kappa_0 = (C_v T)^{-1} \kappa$.
\end{proof}

\subsection{Large torus convergence of the $L$-periodic nonlinear solutions}

In this subsection, we prove Theorem~\ref{THM:convL} in the case when $\ld > 0$, the convergence of $v_L$ satisfying the $L$-periodic nonlinear DAM \eqref{vNLSAndL} to $v$ satisfying the nonlinear DAM \eqref{vNLSAnd} on $\R^2$.

We first mention the following result regarding the convergence of the solution $v_\eps$ of the mollified nonlinear DAM \eqref{vNLSe} on $\R^2$ to the solution $v$ of the nonlinear DAM \eqref{vNLSAnd} on $\R^2$. This is covered by \cite[Theorem~1.2]{DLTV}.

\begin{lemma}
\label{LEM:vGWP}
Let $p \geq 2$, $\ld > 0$, $\mu_0 > 0$, and $v_0 \in H_{\mu_0}^2 (\R^2)$. Given $0 < \eps < \frac 12$, let $v_\eps$ be the global-in-time solution to the mollified nonlinear DAM \eqref{vNLSe} with $v_\eps |_{t = 0} = v_0$ and let $v$ be the global-in-time solution to the nonlinear DAM \eqref{vNLSAnd} with $v |_{t = 0} = v_0$. Then, for any $T \geq 1$, $0 \leq s < 2$, and $\mu > 0$ sufficiently small, there exists $\dl > 0$ and $\Om' \subset \Om$ with full probability measure such that for any $\om \in \Om'$ and $0 < \eps < \frac 12$, we have the bound
\begin{align*}
\| v \|_{C_T H_\mu^s (\R^2)} \les_{\om, T} C \big( \| v_0 \|_{H_{\mu_0}^2 (\R^2)} \big)
\end{align*}

\noi
and the difference estimate
\begin{align*}
\| v_\eps - v \|_{C_T H_\mu^s (\R^2)} \les_{\om, T} C \big( \| v_0 \|_{H_{\mu_0}^2 (\R^2)} \big)
\end{align*}

\noi
for some constant $C (\| v_0 \|_{H_{\mu_0}^2 (\R^2)}) > 0$. 
\end{lemma}

We are now ready to show the large torus limit of the nonlinear DAM.

\begin{proof}[Proof of Theorem~\ref{THM:convL} for $\ld > 0$]
From Lemma~\ref{LEM:perw}, we have
\begin{align}
\| \s_R^{-1} (v_{0, L} - v_0) \|_{L_{\mu_1}^2 (\R^2)} \les (L^{- \mu_0 + \mu_1} + R^{1 + 2 \mu_1} L^{-1 - \mu_1}) \| v_0 \|_{L_{\mu_0}^2 (\R^2)} 
\label{v0L_diff}
\end{align}

\noi
for any $0 < \mu_1 < \mu_0 - 1$.
Given $0 < \eps < \frac 12$, we let $v_{L, \eps}$ and $v_L$ be the global-in-time solutions to \eqref{vNLSLe} and \eqref{vNLSAndL}, respectively, both with initial data $v_{0, L}$. We also let $v_\eps$ and $v$ be the global-in-time solutions to \eqref{vNLSe} and \eqref{vNLSAnd}, respectively, both with initial data $w_0$. Let $\Om' \subset \Om$ be the event with full probability measure such that all estimates in Section~\ref{SEC:sto} and Section~\ref{SEC:convL1} hold for some parameters to be fixed later and we fix $\om \in \Om'$.

Given any $T \geq 1$, $0 \leq s < 2$, and a bounded open set $U \subset \R^2$, we have
\begin{align}
\| v_L - v \|_{C_T H^s (U)} \leq \| v_L - v_{L, \eps} \|_{C_T H^s (U)} + \| v_{L, \eps} - v_\eps \|_{C_T H^s (U)} + \| v_\eps - v \|_{C_T H^s (U)}.
\label{vLlim0}
\end{align}

\noi
From Lemma~\ref{LEM:Hsloc}, Proposition~\ref{PROP:vLe_conv}, and Lemma~\ref{LEM:per}, we have
\begin{align}
\| v_L - v_{L, \eps} \|_{C_T H^s (U)} \les_{U} \| v_L - v_{L, \eps} \|_{C_T H^s (\T_L^2)} \les_{\om, T} \eps^\dl C \big( \| v_{0} \|_{H_{\mu_0}^2 (\R^2)} \big)
\label{vLlim1}
\end{align}

\noi
for some $\dl > 0$ and some constant $C ( \| v_{0} \|_{H_{\mu_0}^2 (\T_L^2)} ) > 0$. Also, from Lemma~\ref{LEM:vGWP}, we have
\begin{align}
\| v_\eps - v \|_{C_T H^s (U)} \leq \| v_\eps - v \|_{C_T H^s_\mu (\R^2)} \les_{\om, T} \eps^\dl C \big( \| v_0 \|_{H_{\mu_0}^2 (\R^2)} \big)
\label{vLlim2}
\end{align}

\noi
for some $\dl > 0$. 

We now estimate $\| v_{L, \eps} - v_\eps \|_{C_T H^s (U)}$. Let $\phi \in C_c^\infty (\R^2)$ be such that $\phi \equiv 1$ on $U$ and let $U' \subset \R^2$ be a bounded open set such that $\supp \phi \subset U'$. Let $0 \leq s < s' < s'' < 2$. By interpolation and \eqref{prodU}, we have
\begin{align}
\begin{split}
\| v_{L, \eps} - v_\eps \|_{C_T H^s (U)} &\leq \| \phi (v_{L, \eps} - v_\eps) \|_{C_T H^s (\R^2)} \\
&\les \| \phi (v_{L, \eps} - v_\eps) \|_{C_T L^2 (\R^2)}^{1 - \frac{s}{s'}} \| \phi (v_{L, \eps} - v_\eps) \|_{L_T^\infty H^{s'} (\R^2)}^{\frac{s}{s'}} \\
&\les \| v_{L, \eps} - v_\eps \|_{C_T L^2 (U')}^{1 - \frac{s}{s'}} \| v_{L, \eps} - v_\eps \|_{L_T^\infty H^{s''} (U')}^{\frac{s}{s'}} .
\end{split}
\label{vLlim3-0}
\end{align}

\noi
From Lemma~\ref{LEM:Hsloc}, \eqref{Hs_equi},  Lemma~\ref{LEM:Hsbdd} along with Lemma~\ref{LEM:embL}~(iv), and Lemma~\ref{LEM:vGWP}, we get
\begin{align}
\begin{split}
\| v_{L, \eps} - v_\eps \|_{L_T^\infty H^{s''} (U')} 
&\les \| v_{L, \eps} \|_{L_T^\infty H^{s''} (\T_L^2)} + \| v_\eps \|_{L_T^\infty H^{s''} (\R^2)} \\
&\sim \| v_{L, \eps} \|_{L_T^\infty \B_{2, 2}^{s''} (\T_L^2)} + \| v_\eps \|_{L_T^\infty H^{s''} (\R^2)} \\
&\les |\log \eps|^{136} C \big( \| v_{0, L} \|_{H^2_{\mu_1} (\T_L^2)}, \| v_{0} \|_{H^2_{\mu_0} (\R^2)} \big) 
\end{split}
\label{vLlim3-1}
\end{align}

\noi
for some constant $C ( \| v_{0, L} \|_{H^2_{\mu_1} (\T_L^2)}, \| v_{0} \|_{H^2_{\mu_0} (\R^2)} ) > 0$. 
Also, from Proposition~\ref{PROP:vLeve}, we have
\begin{align}
\begin{split}
\| v_{L, \eps} - v_\eps \|_{C_T L^2 (U')} 
&\les_{U'} \| \s_R^{-2} (v_{L, \eps} - v_\eps) \|_{C_T L_{- \mu}^2 (\R^2)} \\
&\les_{\om, T} \eps^{- \kappa} \big( \| \s_R^{-2} (v_{0, L} - v_0) \|_{C_T L_\mu^2 (\R^2)} + R^{-1} + R^{b} L^{- \frac{b}{8}} |\log \eps| \big) \\
&\quad \times C \big( \| v_{0, L} \|_{H^2_{\mu_1} (\T_L^2)}, \| v_0 \|_{H_{\mu_1}^2 (\R^2)} \big)
\end{split}
\label{vLlim3-2}
\end{align}

\noi
for any $\kappa > 0$, $\mu > 0$ sufficiently small, and $0 < b \leq 4$.
Thus, from \eqref{vLlim3-0}, \eqref{vLlim3-1}, and \eqref{vLlim3-2} along with \eqref{v0L_diff}, Lemma~\ref{LEM:perw}, Lemma~\ref{LEM:per}, and the choice $R \sim L^{\frac{1}{16}}$, we get
\begin{align}
\| v_{L, \eps} - v_\eps \|_{C_T H^s (U)} \les \eps^{- \kappa'} L^{- \kappa_0} |\log \eps|^{136} C \big( \| v_0 \|_{H_{\mu_0}^2 (\R^2)} \big)
\label{vLlim3}
\end{align}

\noi
for any $\kappa' > 0$ and some $\kappa_0 > 0$.

Therefore, combining \eqref{vLlim0}, \eqref{vLlim1}, \eqref{vLlim2}, and \eqref{vLlim3}, we get
\begin{align*}
\| v_L - v \|_{C_T H^s (U)} &\les_{\om, T, U} (\eps^\dl + \eps^{- \kappa'} L^{- \kappa_0} |\log \eps|^{136} ) C \big( \| v_0 \|_{H_{\mu_0}^2 (\R^2)} \big) 
\end{align*}

\noi
for any $\kappa' > 0$ and some $\dl > 0$ and  $\kappa_0 > 0$. Therefore, by choosing
\begin{align*}
\eps \sim L^{- \frac{\kappa_0}{2 \kappa'}} ,
\end{align*}

\noi
we obtain the desired convergence result.
\end{proof}

\appendix

\section{Large torus limit of the deterministic nonlinear Schr\"odinger equation}
\label{SEC:NLS}

In this section of the appendix, we study the large torus limit problem for the deterministic NLS in the two-dimensional setting. To be more specific, we consider the following Cauchy problem on $\R^2$:
\begin{align}
\begin{cases}
i \dt u = \Dl u - \lambda |u|^{p - 1} u \\
u|_{t = 0} = u_0,
\end{cases}
\label{dNLSu}
\end{align}

\noi
where $p \geq 2$ and $\lambda \geq 0$. Our goal is to show that the global-in-time solution to \eqref{dNLSu} can be realized as a limit of the solution to the following $L$-periodic NLS as $L$ tends to infinity:
\begin{align}
\begin{cases}
i \dt u_L = \Dl u_L - \lambda |u_L|^{p - 1} u_L \\
u_L|_{t = 0} = u_{0, L} ,
\end{cases}
\label{dNLSLu}
\end{align}

\noi
where $u_{0, L}$ is the following $L$-periodized version of $u_0$:
\begin{align}
u_{0, L} \deff \sum_{k \in \Z^2} u_0 (\cdot + Lk) .
\label{u0L2}
\end{align}

We assume that the initial data $u_0$ for \eqref{dNLSu} lies in $H^2_{\mu_0} (\R^2)$ with $\mu_0 > 1$. From \eqref{u0L2} and Lemma~\ref{LEM:per}, we know that $u_{0, L} \in H^2 (\T_L^2)$. We know from \cite{GV79, Kat87, Tsu, Caz03, BGTz} that \eqref{dNLSu} has a unique global-in-time solution in $C(\R_+; H^2 (\R^2))$ and that \eqref{dNLSLu} has a unique global-in-time solution in $C(\R_+; H^2 (\T_L^2))$. Alternatively, we can use the steps in Subsection~\ref{SUB:H2L} to obtain global well-posedness of \eqref{dNLSLu} in $H^2 (\T_L^2)$. Note that for any $T \geq 1$, an analog to Lemma~\ref{LEM:H2bdd} yields
\begin{align}
\| u_L \|_{C_T H^2 (\T_L^2)} \les_T C (\| u_{0, L} \|_{H^2 (\T_L^2)})
\label{uL_H2bdd}
\end{align}

\noi
for some constant $C (\| u_{0, L} \|_{H^2 (\T_L^2)}) > 0$ uniformly in $L$ (we can put all the weight parameters to be zero since there is no stochastic terms in this deterministic setting).

We now prove the following statement.
\begin{theorem}
\label{THM:NLS}
Let $p \geq 2$, $\ld \geq 0$, $\mu_0 > 1$, $u_0 \in H^2_{\mu_0} (\R^2)$, and $u_{0, L} \in H^2 (\T_L^2)$ be as defined in \eqref{u0L2} given any $L \geq 1$. Let $u$ be the global-in-time solution to NLS \eqref{dNLSu} with initial data $u_0$ and let $u_L$ be the global-in-time solution to NLS \eqref{dNLSLu} with initial data $u_{0, L}$. Then, for any $0 \leq s < 2$, $T \geq 1$, and bounded Lipschitz domain $U \subset \R^2$, we have the estimate
\begin{align*}
\| u_L - u \|_{C_T H^s (U)} \leq L^{- \frac 12 + \frac{s}{4}} C \big( \| u_0 \|_{H_{\mu_0}^2 (\R^2)}, T \big)
\end{align*}

\noi
for some constant $C (\| u_0 \|_{H_{\mu_0}^2 (\R^2)}, T) > 0$, which implies that $u_L$ converges to $u$ in $C(\R_+; H^s (U))$ (endowed with the compact-open topology in time) as $L \to \infty$.
\end{theorem}

As in the case of Theorem~\ref{THM:convL}, we mention that the condition on the weight parameter $\mu_0$ in Theorem~\ref{THM:NLS} is sharp up to the endpoint in the following sense: if $\mu_0 < 1$, then there exists $u_0 \in H_{\mu_0}^2 (\R^2)$ such that the periodic function $u_{0, L}$ may not be constructed via the form \eqref{u0L2}; see Remark~\ref{RMK:mu0}.

\begin{proof}[Proof of Theorem~\ref{THM:NLS}]
Given $L \geq 1$, we let $r_L = u - u_L$. Then, we see that $r_L$ satisfies the following equation:
\begin{align*}
i \dt r_L = \Dl r_L - \ld \big( |u|^{p - 1} u - |u_L|^{p - 1} u_L \big) .
\end{align*}

\noi
Given $1 \leq R \leq L$, we recall the function $\s_R$ in \eqref{sigR}. From Lemma~\ref{LEM:sig}, \eqref{uL_H2bdd}, and Lemma~\ref{LEM:per}, we have
\begin{align}
\begin{split}
\| \s_R^{-1} u_L \|_{C_T L^2 (\R^2)} + \| \s_R^{-1} \nb u_L \|_{C_T L^2 (\R^2)} 
&\les \| u_L \|_{C_T L^2 (\T_L^2)} + \| \nb u_L \|_{C_T L^2 (\T_L^2)} \\
&\les_T C (\| u_{0, L} \|_{H^2 (\T_L^2)}) \\
&\les C \big( \| u_{0} \|_{H_{\mu_0}^2 (\R^2)} \big) .
\end{split}
\label{sRbdd}
\end{align}

Let $T \geq 1$ and $t \in [0, T]$. By using the equation and integration by parts, we obtain
\begin{align}
\begin{split}
\frac 12 \frac{d}{dt} \int_{\R^2} \s_R^{-2} |r_L (t)|^2 dx &= - \Im \int_{\R^2} \nb (\s_R^{-2}) \cdot \nb r_L (t) \cj{r_L (t)} dx \\
&\quad - \ld \Im \int_{\R^2} \s_R^{-2} \big( |u (t)|^{p - 1} u (t) - |u_L (t)|^{p - 1} u_L (t) \big) \cj{r_L (t)} dx \\
&\deff \textup{I}_1 + \textup{I}_2.
\end{split}
\label{uL1}
\end{align}

\noi
For $\textup{I}_1$, by \eqref{sigR1}, Cauchy's inequality, and \eqref{sRbdd}, we have
\begin{align}
\begin{split}
|\textup{I}_1| &\les R^{-2} \int_{\R^2} \s_R^{-2} |\nb r_L (t)|^2 dx + \int_{\R^2} \s_R^{-2} |r_L (t)|^2 dx \\
&\les_T R^{-2} C \big( \| u_0 \|_{H_{\mu_0}^2 (\R^2)} \big) + \int_{\R^2} \s_R^{-2} |r_L (t)|^2 dx
\end{split}
\label{uL2}
\end{align}

\noi
for some constant $C (\| u_0 \|_{H_{\mu_0}^2 (\R^2)}) > 0$. For $\textup{I}_2$, by the mean value theorem, Sobolev's embedding, the embeddings in Lemma~\ref{LEM:embL}~(ii), (i), and (iii) along with \eqref{Hs_equi}, and \eqref{sRbdd}, we have
\begin{align}
\begin{split}
|\textup{I}_2| &\les \big( \| u (t) \|_{L^\infty (\R^2)}^{p - 1} + \| u_L (t) \|_{L^\infty (\R^2)}^{p - 1} \big) \int_{\R^2} \s_R^{-2} |r_L (t)|^2 dx \\
&\les \big( \| u (t) \|_{H^{1 + s_0} (\R^2)}^{p - 1} + \| u_L (t) \|_{H^{1 + s_0} (\T^2)}^{p - 1} \big) \int_{\R^2} \s_R^{-2} |r_L (t)|^2 dx \\
&\les_T C \big( \| u_0 \|_{H_{\mu_0}^2 (\R^2)} \big) \int_{\R^2} \s_R^{-2} |r_L (t)|^2 dx
\end{split}
\label{uL3}
\end{align}

\noi
for any $s_0 > 0$. Combining \eqref{uL1}, \eqref{uL2}, and \eqref{uL3}, we obtain
\begin{align*}
\frac{d}{dt} \int_{\R^2} \s_R^{-2} |r_L (t)|^2 dx \les_T C \big( \| u_0 \|_{H_{\mu_0}^2 (\R^2)} \big) R^{-2} + C \big( \| u_0 \|_{H_{\mu_0}^2 (\R^2)} \big) \int_{\R^2} \s_R^{-2} |r_L (t)|^2 dx .
\end{align*}

\noi
Using Gronwall's inequality and Lemma~\ref{LEM:perw}, we obtain
\begin{align}
\begin{split}
\sup_{t \in [0, T]} \int_{\R^2} \s_R^{-2} |r_L (t)|^2 dx &\les \bigg( \int_{\R^2} \s_R^{-2} |u_0 - u_{0, L}|^2 dx + R^{-2} \bigg) \exp \big( C \big( \| u_0 \|_{H_{\mu_0}^2 (\R^2)} \big) T \big) \\
&\les (L^{- 2\mu_0 + 2 \mu} + R^{2 + 4 \mu} L^{-2 - 2 \mu} + R^{-2}) \exp \big( C \big( \| u_0 \|_{H_{\mu_0}^2 (\R^2)} \big) T \big) 
\end{split}
\label{uLd}
\end{align}

\noi
for any $0 < \mu < \mu_0 - 1$.

Let $U \subset \R^2$ be a bounded Lipschitz domain. Let $\phi \in C_c^\infty (\R^2)$ be such that $\phi \equiv 1$ on $U$ and let $U' \subset \R^2$ be a bounded Lipschitz domain such that $\supp \phi \subset U'$. By interpolation, we have
\begin{align*}
\| r_L \|_{C_T H^s (U)} \leq \| \phi r_L \|_{C_T H^s (\R^2)} 
\les \| \phi r_L \|_{C_T L^2 (\R^2)}^{1 - \frac{s}{2}} \| \phi r_L \|_{C_T H^2 (\R^2)}^{\frac{s}{2}} ,
\end{align*}

\noi
so that by the equivalence in \eqref{loc_equi}, Lemma~\ref{LEM:Hsloc}, and \eqref{sRbdd}, we get
\begin{align}
\begin{split}
\| r_L \|_{C_T H^s (U)} 
&\les \| r_L \|_{C_T L^2 (U')}^{1 - \frac{s}{2}} \| r_L \|_{C_T H^2 (U')}^{\frac{s}{2}} \\
&\les \| r_L \|_{C_T L^2 (U')}^{1 - \frac{s}{2}} \big( \| u \|_{C_T H^2 (\R^2)} + \| u_L \|_{C_T H^2 (\T_L^2)} \big)^{\frac{s}{2}} \\
&\les_T \| r_L \|_{C_T L^2 (U')}^{1 - \frac{s}{2}} C \big( \| u_0 \|_{H_{\mu_0}^2 (\R^2)} \big) .
\end{split}
\label{uLd1}
\end{align}

\noi
From \eqref{uLd}, we have
\begin{align}
\begin{split}
\| r_L \|_{C_T L^2 (U')} &\les_{U'} \| \s_R^{-1} r_L \|_{C_T L^2 (\R^2)} \\
&\les (L^{- \mu_0 + \mu} + R^{1 + 2 \mu} L^{-1 - \mu} + R^{-1}) \exp \big( C \big( \| u_0 \|_{H_{\mu_0}^2 (\R^2)} \big) T \big) .
\end{split}
\label{uLd2}
\end{align}

\noi
Thus, from \eqref{uLd1} and \eqref{uLd2} along with the choice $R = L^{\frac 12}$, we get
\begin{align*}
\| r_L \|_{C_T H^s (U)} \les L^{- \frac 12 + \frac{s}{4}} C' \big( \| u_0 \|_{H_{\mu_0}^2 (\R^2)}, T \big)
\end{align*}

\noi
for some constant $C' (\| u_0 \|_{H_{\mu_0}^2 (\R^2)}, T) > 0$. This gives the desired estimate and the convergence result.
\end{proof}

\begin{remark} \rm
We see from the proof of Theorem~\ref{THM:NLS} that we only need that the solutions $u$ and $u_L$ have regularities $1 + s_0$ for any $s_0 > 0$. It is possible to reduce the regularity requirement of the initial data $u_0$ in Theorem~\ref{THM:NLS}. Indeed, we know that \eqref{dNLSu} and \eqref{dNLSLu} are globally well-posed in $H^{1 + s_0} (\R^2)$  and $H^{1 + s_0} (\T_L^2)$, respectively. This can be proved by the classical $H^1$-global well-posedness of \eqref{dNLSu} and \eqref{dNLSLu} (see \cite{GV79, CW88, Bour93}) and a persistence of regularity result for showing that $u \in C(\R_+; H^{1 + s_0} (\R^2))$ and $u_L \in C(\R_+; H^{1 + s_0} (\T_L^2))$.
Nevertheless, we still require an extension of Lemma~\ref{LEM:per} so that it works for functions with fractional derivatives.
\end{remark}

\begin{remark} \rm
The large torus convergence in Theorem~\ref{THM:NLS} still holds if we only assume that $U$ is a bounded open set just as in Theorem~\ref{THM:convL}. In this case, instead of using the equivalence \eqref{loc_equi}, we need to invoke the produce estimate in Lemma~\ref{LEM:prod}, but we have to lose some slight amount of regularity. This will result in a less optimal convergence rate in $L$.
\end{remark}

\section{Large torus limit of the parabolic Anderson model}
\label{SEC:heat}

In this section of the appendix, we study the large torus limit problem for the parabolic Anderson model (PAM) in the two-dimension setting. It is known that the parabolic smoothing property holds in a wide class of weighted $L^p$-based spaces for any $1 \leq p \leq \infty$. This provides us with the advantage to treat the functions on the full space and the functions on the $L$-periodic domain in a unified manner. In other words, we can put periodic functions in function spaces on $\R^d$ (with $d \in \N$) with decaying weights. Such observation has already been used in \cite{MW17, GH19, DGR}, where the authors obtain uniform-in-$L$ bounds for $L$-periodic solutions in weighted function spaces on $\R^d$. We remark that this is in contrast to the situation of the dispersive Anderson model, for which, due to the lack of dispersive smoothing property on weighted function spaces, we need to use $L$-periodic weighted function spaces constructed in Section~\ref{SEC:BesL} to obtain uniform-in-$L$ bounds for $L$-periodic solutions.

Let us consider the following Cauchy problem for the linear heat equation with a multiplicative spatial white noise on $\R^2$:
\begin{align}
\begin{cases}
\dt u = \Dl u + \xi u \\
u |_{t = 0} = u_0 ,
\end{cases}
\label{PAMu}
\end{align}

\noi
where $\xi$ denotes the real-valued white noise in space. Our goal is to show that the global-in-time solution to \eqref{PAMu} can be realized as a limit of the solution to the following $L$-periodic PAM as $L$ tends to infinity:
\begin{align}
\begin{cases}
\dt u_L = \Dl u_L + \xi_L u_L \\
u_L |_{t = 0} = u_{0, L} ,
\end{cases}
\label{PAMLu}
\end{align}

\noi
where $\xi_L$ is the $L$-periodized version of $\xi$ defined in \eqref{defxiL} and $u_{0, L}$ is a suitable $L$-periodized version of $u_0$ to be specified below (see \eqref{PAMu0L}).

We proceed as in \eqref{gauge} by defining
\begin{align*}
v = e^Y u \quad \text{and} \quad v_L = e^{Y_L} u_L ,
\end{align*}

\noi
where $Y$ and $Y_L$ are defined in \eqref{defY}. Then, recalling Wick orderings $\wick{|\nb Y|^2}$ and $\wick{|\nb Y_L|^2}$ defined in \eqref{Y2_sto} and \eqref{YL2_fs}, respectively, and also the notations $\wt{\wick{|\nb Y|^2}}$ and $\wt{\wick{|\nb Y_L|^2}}$ in \eqref{wY2t}, we get the following equations for $v$ and $v_L$:
\begin{align}
\begin{cases}
\dt v = \Dl v - 2 \nb Y \cdot \nb v + \wt{\wick{ |\nb Y|^2 }} \, v \\
v |_{t = 0} = v_0 \deff e^{- Y} u_0 
\end{cases}
\label{PAM}
\end{align}

\noi
and
\begin{align}
\begin{cases}
\dt v_L = \Dl v_L - 2 \nb Y_L \cdot \nb v_L + \wt{\wick{ |\nb Y_L|^2 }} \, v_L \\
v_L |_{t = 0} = v_{0, L} \deff e^{- Y_L} u_{0, L} ,
\end{cases}
\label{PAML}
\end{align}

\noi
where $v_{0, L}$ is defined by
\begin{align}
v_{0, L} \deff \sum_{k \in \Z^2} v_0 (\cdot + Lk).
\label{v0L2}
\end{align}

\noi
This also provides the choice of the initial data $u_{0, L}$ for \eqref{PAMLu}:
\begin{align}
u_{0, L} = e^{Y_L} v_{0, L} = e^{Y_L} \sum_{k \in \Z^2} e^{- Y (\cdot + Lk)} u_0 (\cdot + Lk) .
\label{PAMu0L}
\end{align} 

\noi
For later convenience, we refer to the equation \eqref{PAM} as PAM and the equation \eqref{PAML} as the $L$-periodic PAM.

We will use the weighted $L^\infty$-based spaces to show the convergence of $v_L$ to $v$ as $L \to \infty$ in some suitable sense. The reason for the choice of this type of space is that one can then treat the functions on the full space and the functions on the $L$-periodic domain in a unified manner.

We assume that the initial data $v_0$ for \eqref{PAM} lies in $L^\infty_{\mu_0} (\R^2)$ with $\mu_0 > 2$. In view of the global well-posedness statement in \cite[Theorem~1.4]{HL15}, the regularity assumption of $v_0$ can be lowered, but we choose to make our presentation simpler and focus more on the large torus limit. The requirement $\mu_0 > 2$ is due to the following lemma, which is an $L^\infty$-version of Lemma~\ref{LEM:per}.
\begin{lemma}
\label{LEM:per3}
Let $L \geq 1$, $\mu_0 > 2$, and $f \in L^\infty_{\mu_0} (\R^2)$. Let $f_L$ be an $L$-periodic function defined by
\begin{align*}
f_L \deff \sum_{k \in \Z^2} f (\cdot + Lk).
\end{align*}

\noi
Then, $f_L \in L^\infty (\R^2)$ and, more precisely, 
\begin{align}
\| f_L \|_{L^\infty (\R^2)} \les \| f \|_{L^\infty_{\mu_0} (\R^2)}.
\label{fLbdd2}
\end{align}

\noi
Moreover, for any $\mu > 0$, we have
\begin{align}
\| f - f_L \|_{L^\infty_{- \mu} (\R^2)} \les L^{- \min (\mu, \mu_0)} \| f \|_{L^\infty_{\mu_0} (\R^2)}.
\label{fLdiff2}
\end{align}
\end{lemma}

\begin{proof}
We first show \eqref{fLbdd2}. For any $x \in \R^2$, we have
\begin{align*}
|f_{L} (x)| &= \Big| \sum_{k \in \Z^2} f (x + L k) \Big| \\
&\leq \Big| \sum_{k \in \Z^2} \jb{x + L k}^{- \mu_0} \Big| \sup_{k \in \Z^2} | \jb{x + L k}^{\mu_0} f (x + L k) | \\
&\les \| f \|_{L_{\mu_0}^\infty (\R^2)} ,
\end{align*}

\noi
which gives the desired estimate.

We now show \eqref{fLdiff2}. Given $x \in \R^2$, if $|x| \leq \frac{L}{2}$, we have
\begin{align*}
\jb{x}^{- \mu} |f (x) - f_L (x)| &\leq \Big| \sum_{k \in \Z^2 \setminus \{0\}} f (x + L k) \Big| \\
&\leq \Big| \sum_{k \in \Z^2 \setminus \{0\}} \jb{x + L k}^{- \mu_0} \Big| \sup_{k \in \Z^2} | \jb{x + L k}^{\mu_0} f (x + L k) | \\ 
&\les L^{- \mu_0} \| f \|_{L_{\mu_0}^\infty (\R^2)}.
\end{align*}

\noi
If $x > \frac{L}{2}$, we have
\begin{align*}
\jb{x}^{- \mu} |f (x) - f_L (x)| &\les L^{- \mu} |f (x) - f_L (x)| \\
&\leq L^{- \mu} \Big| \sum_{k \in \Z^2 \setminus \{0\}} \jb{x + L k}^{- \mu_0} \Big| \sup_{k \in \Z^2} | \jb{x + L k}^{\mu_0} f (x + L k) | \\ 
&\les L^{- \mu} \| f \|_{L_{\mu_0}^\infty (\R^2)}.
\end{align*}

\noi
Thus, we obtain the desired difference estimate.
\end{proof}

\subsection{Besov spaces with exponential weights}
\label{SUB:exp}

In this subsection, we introduce Besov spaces with exponential weights.

Let us first recall the Gevrey class. See also \cite[Subsection~2.1]{MW17}. Given $\ta \geq 1$, we define the Gevrey class $\mathcal{G}^\ta \subset C^\infty (\R^2)$ of order $\ta$ to be the set of functions $f$ such that, given any compact set $K \subset \R^2$, there exists $C_K > 0$ such that
\begin{align*}
\sup_{x \in K} |\partial^\al f (x)| \leq (C_K)^{|\al| + 1} (\al !)^\ta
\end{align*}

\noi
for any multi-index $\al$. We denote by $\mathcal{G}^\ta_c$ the set of compactly supported functions in $\mathcal{G}^\ta$. Note that when $\ta > 1$, the set $\mathcal{G}^\ta_c$ is not empty; see \cite[Subsection~2.1]{MW17} or \cite[Example~1.4.9]{Rod}.

From now on, we fix $\ta > 1$. Let $\chi \in \mathcal{G}_c^\ta (\R^2)$ be supported on $\{ \xi \in \R^2 : |\xi| \leq \frac 85 \}$ and $\chi \equiv 1$ on $\{ \xi \in \R^2 : |\xi| \leq \frac 54 \}$, which is possible according to \cite[Subsection~2.1]{MW17}.
Given a dyadic number $N \geq 1$, we define
\begin{align*}
\zeta_N (\xi) =
\begin{cases}
\chi (\xi)  & \text{if } N = 1 \\
\chi (\tfrac{\xi}{N}) - \chi (\tfrac{2 \xi}{N})  & \text{if } N \geq 2 .
\end{cases} 
\end{align*}

\noi
We then define the Littlewood-Paley projector $\wt \Dl_N$ as the Fourier multiplier operator with symbol $\zeta_N$, so that
\begin{align*}
\Id = \sum_{\substack{N \geq 1 \\ \text{dyadic}}} \wt \Dl_N .
\end{align*}


We now fix $\ta > 1$ and $0 < \vartheta < \frac{1}{\ta}$ throughout the rest of this section. Given $\nu \in \R$, we define $\wt L^p_\nu (\R^2)$ as the $L^p$ space with an exponential weight $e^{\nu \jb{x}^\vartheta}$ via the norm
\begin{align*}
\| f \|_{\wt L_\nu^p (\R^2)} \deff \bigg( \int_{\R^2} |f (x)|^p e^{p \nu \jb{x}^\vartheta} dx \bigg)^{\frac 1p} = \big\| e^{\nu \jb{\cdot}^\vartheta} f \big\|_{L^p (\R^2)} ,
\end{align*}

\noi
with the obvious interpretation if $p = \infty$. Given $s, \nu \in \R$ and $1 \leq p, q \leq \infty$, we define $\wt \B_{p, q, \nu}^s (\R^2)$ as the inhomogeneous Besov space with an exponential weight $e^{\nu \jb{x}^\vartheta}$ via the norm
\begin{align*}
\| f \|_{\wt \B_{p, q, \nu}^s (\R^2)} \deff \bigg( \sum_{\substack{N \geq 1 \\ \text{dyadic}}} N^{sq} \| \wt \Dl_N f \|_{\wt L_\nu^p (\R^2)}^q \bigg)^{\frac 1q} = \bigg( \sum_{\substack{N \geq 1 \\ \text{dyadic}}} N^{sq} \big\| e^{\nu \jb{\cdot}^\vartheta} \wt \Dl_N f \big\|_{L^p (\R^2)}^q \bigg)^{\frac 1q} .
\end{align*}

\noi
When $p = q = \infty$, we denote
\begin{align*}
\wt \C_\nu^s (\R^2) \deff \wt \B_{\infty, \infty, \nu}^s (\R^2) .
\end{align*}

\noi
When $\nu = 0$, the above weighted spaces coincide with the usual unweighted spaces: $\wt \B_{p, q, 0}^s (\R^2) = \B_{p, q}^s (\R^2)$ and $\wt \C_0^s (\R^2) = \C^s (\R^2)$.
Note that unlike \cite{MW17}, we put the weight inside the $L^p$ norm in order to be consistent with the case of polynomial weights in the previous subsection. Nevertheless, most results in \cite{MW17} on Besov spaces with exponential weights can be carried over to our setting with straightforward modifications.

Let us mention the following useful estimates for Besov norms with exponential weights.

\begin{lemma}
\label{LEM:embex}
\textup{(i)} Let $1 \leq p \leq \infty$ and $s_1, s_2 \in \R$ be such that $s_2 \geq s_1 + \frac{2}{p}$, and $\nu \in \R$. Then, we have
\begin{align*}
\| f \|_{\wt{\mathcal{C}}_{\nu}^{s_1} (\R^2)} \les \| f \|_{\wt{\B}_{p, \infty, \nu}^{s_2} (\R^2)} .
\end{align*}

\smallskip \noi
\textup{(ii)} Let $\nu_1, \nu_2 \in \R$ be such that $\nu_1 \leq \nu_2$. Then, we have
\begin{align*}
\| f \|_{\wt{\mathcal{C}}_{\nu_1}^s (\R^2)} \leq  \| f \|_{\wt{\mathcal{C}}_{\nu_2}^s (\R^2)} .
\end{align*} 

\smallskip \noi
\textup{(iii)} Let $s \in \R$ and $\nu_0 < 0$. Then, for any $\nu \geq \nu_0$, we have
\begin{align*}
\| \nb f \|_{\wt{\mathcal{C}}_\nu^s (\R^2)} \les \| f \|_{\wt{\mathcal{C}}_\nu^{s + 1} (\R^2)} ,
\end{align*}

\noi
where the underlying constant is independent of $\nu$.
\end{lemma}

\begin{proof}
Part (i) follows from \cite[Proposition~2]{MW17}. Part (ii) follows directly from the definition. Part (iii) is a special case of \cite[Proposition~3]{MW17}.
\end{proof}

We also have the following estimate on convolution with a smooth and compactly supported function.
\begin{lemma}
\label{LEM:varphiex}
For any $s_1, s_2 \in \R$, $\nu \geq 0$, and $\varphi \in \C_c^\infty (\R^2)$, we have
\begin{align*}
\| \varphi *_{\R^2} f \|_{\wt \C_{- \nu}^{s_1} (\R^2)} \les_\varphi \| f \|_{\wt \C_{- \nu}^{s_2} (\R^2)} .
\end{align*}
\end{lemma}

\begin{proof}
The estimate follows directly from the proof of \cite[Lemma~2.5]{DLTV} along with the additional observation that
\begin{align*}
e^{- \nu \jb{x + y}^\vartheta} \leq e^{\nu \jb{x}^\vartheta} e^{- \nu \jb{y}^{\vartheta}}
\end{align*}

\noi
for any $x, y \in \R^2$.
\end{proof}

We then record the following Schauder estimate for the linear heat kernel $e^{t \Dl}$. For a proof, see \cite[Proposition~5]{MW17}.
\begin{lemma}
\label{LEM:Sch}
Let $s_1 \geq s_2$ and $\nu_0 < 0$. Then, for any $\nu \geq \nu_0$ and $t > 0$, we have
\begin{align*}
\| e^{t \Dl} f \|_{\wt \C_\nu^{s_1} (\R^2)} \les t^{- \frac{s_1 - s_2}{2}} \| f \|_{\wt \C_\nu^{s_2} (\R^2)} ,
\end{align*}

\noi
where the underlying constant is independent of $\nu$ and $t$.
\end{lemma}

Lastly, we show the following product estimate, which involves both the polynomial weight (see Subsection~\ref{SUB:poly}) and the exponential weight.
\begin{lemma}
\label{LEM:prodex2}
Let $s, s_1, s_2 \in \R$ be such that $s_1 + s_2 > 0$ and $s = \min (s_1 + s_2, s_1, s_2)$, $\mu \in \R$, and $\nu_1 < \nu_2$. Then, we have
\begin{align*}
\| fg \|_{\wt \C_{\nu_1}^{s} (\R^2)} \les (\nu_2 - \nu_1)^{\frac{\mu}{\vartheta}} \| f \|_{\C_{\mu}^{s_1} (\R^2)} \| g \|_{\wt \C_{\nu_2}^{s_2} (\R^2)} ,
\end{align*}

\noi
where $0 < \vartheta < 1$ is a fixed parameter defined in Subsection~\ref{SUB:exp} and the underlying constant is independent of $\nu_1$ and $\nu_2$.
\end{lemma}

\begin{proof}
The proof of the estimate follows directly from \cite[Theorem~3.1]{MW17} along with the following elementary estimate:
\begin{align*}
\frac{e^{\nu_1 \jb{x}^\vartheta}}{\jb{x}^{\mu} e^{\nu_2 \jb{x}^\vartheta}} = \jb{x}^{-\mu} e^{- (\nu_2 - \nu_1) \jb{x}^\vartheta} \les (\nu_2 - \nu_1)^{\frac{\mu}{\vartheta}}
\end{align*}

\noi
uniformly in $x \in \R^2$, where we used the fact that $e^{- t} \les_a t^{a}$ for any $t > 0$ and $a \in \R$.
\end{proof}

\subsection{Large torus convergence of stochastic objects in weighted Besov spaces}

In this subsection, we establish the convergence of stochastic objects as $L \to \infty$ in Besov spaces with polynomial weights (see Subsection~\ref{SUB:poly}). This is in contrast with Proposition~\ref{PROP:YconvL}, where we showed the large torus convergence of mollified stochastic objects in Lebesgue spaces with some exponential decaying weight.

For $\wick{ |\nb Y_{L}|^2 }$ defined in \eqref{YL2_fs}, we want to use the form \eqref{YL2_sto}, but it is not obvious to see that \eqref{YL2_sto} is the same as \eqref{YL2_fs}. To see this equivalence, we take the regularized version of $Y_{L, \eps}$ in \eqref{YLe_fs} and write
\begin{align*}
Y_{L, \eps} &= \rho_\eps *_{\R^2} Y_L = (\rho_\eps *_{\R^2} G) *_{\R^2} \xi_L \\
&= \int_{\R^2} (\rho_\eps *_{\R^2} G) (\cdot - z) \xi_L (dz) = \int_{[ - \frac{L}{2}, \frac{L}{2} )^2} \sum_{k \in \Z^2} (\rho_\eps *_{\R^2} G) (\cdot - z + Lk) \xi (dz) .
\end{align*}

\noi
Then, we have
\begin{align*}
\nb Y_{L, \eps} = \int_{[ - \frac{L}{2}, \frac{L}{2} )^2} \sum_{k \in \Z^2}  (\rho_\eps *_{\R^2} \nb G) (\cdot - z + Lk) \xi (dz) ,
\end{align*}

\noi
which makes sense as a Wiener integral since the integrand is smooth thanks to the convolution with $\rho_\eps$. Then, in view of the product formula in \cite[Proposition~1.1.3]{Nua} along with the definition in \eqref{def_YLe2}, we have the following multiple Wiener-Ito integral representation for $\wick{ |\nb Y_{L, \eps}|^2 }$:
\begin{align*}
\wick{ |\nb Y_{L, \eps}|^2 } \, &= \int_{[ - \frac{L}{2}, \frac{L}{2} )^2} \int_{[ - \frac{L}{2}, \frac{L}{2} )^2} \sum_{k_1, k_2 \in \Z^2} (\rho_\eps *_{\R^2} \nb G) (\cdot - z_1 + L k_1) \\
&\qquad \cdot (\rho_\eps *_{\R^2} \nb G) (\cdot - z_2 + L k_2) \xi (dz_1) \xi (dz_2) .
\end{align*}

\noi
Given $f \in \mathcal{D} (\R^2)$, by using the stochastic Fubini theorem (see \cite[Lemma~B.2]{OWZ} in a more general setting), we have
\begin{align*}
( \wick{ |\nb Y_{L, \eps}|^2 } \, , f )_{\R^2} &= \int_{[ - \frac{L}{2}, \frac{L}{2} )^2} \int_{[ - \frac{L}{2}, \frac{L}{2} )^2} \bigg( \int_{\R^2} \sum_{k_1, k_2 \in \Z^2} (\rho_\eps *_{\R^2} \nb G) (x - z_1 + L k_1) \\
&\qquad \cdot (\rho_\eps *_{\R^2} \nb G) (x - z_2 + L k_2) \cj{f (x)} dx \bigg) \xi (dz_1) \xi (dz_2).
\end{align*}


\begin{lemma}
\label{LEM:YL2_sto}
Let $L \geq 1$. Then, for any $f \in \mathcal{D} (\R^2)$, we have
\begin{align}
\begin{split}
( \wick{ |\nb Y_L|^2 } \, , f )_{\R^2} &= \int_{[ - \frac{L}{2}, \frac{L}{2} )^2} \int_{[ - \frac{L}{2}, \frac{L}{2} )^2} \bigg( \int_{\R^2} \sum_{k_1, k_2 \in \Z^2} \nb G (x - z_1 + L k_1) \\
&\qquad \cdot \nb G (x - z_2 + L k_2) \cj{f (x)} dx \bigg) \xi (dz_1) \xi (dz_2).
\end{split}
\label{YL2_sto2}
\end{align}
\end{lemma}

\begin{proof}
We first show that $\wick{ |\nb Y_{L, \eps}|^2 }$ converges to $\wick{ |\nb Y_L|^2 }$ almost surely in $\mathcal{D}' (\R^2)$. Let $f \in \mathcal{D} (\R^2)$ and define $f_L$ to be the $L$-periodized version of $f$ given by
\begin{align*}
f_L \deff \sum_{k \in \Z^2} f (\cdot + L k).
\end{align*}

\noi
Then, from Lemma~\ref{LEM:Sfull}, the definition in \eqref{Sper}, the duality estimate in Lemma~\ref{LEM:dualL}, and Lemma~\ref{LEM:YregLe}~(iv), we have
\begin{align*}
( \wick{ |\nb Y_{L}|^2 } - \wick{ |\nb Y_{L, \eps}|^2 } \, , f )_{\R^2} &= ( \wick{ |\nb Y_{L}|^2 } - \wick{ |\nb Y_{L, \eps}|^2 } \, , f_L )_{\T_L^2} \\
&\les \big\| \wick{ |\nb Y_{L}|^2 } - \wick{ |\nb Y_{L, \eps}|^2 } \big\|_{\C^{-s} (\T_L^2)} \| f_L \|_{\B_{1, 1}^s (\T_L^2)} \\
&\les_{\om, L, \dl} \eps^\dl \| f_L \|_{\B_{1, 1}^s (\T_L^2)}
\end{align*}

\noi
for any $0 < s < 1$ and $\dl > 0$ small. This shows that $\wick{ |\nb Y_{L, \eps}|^2 }$ converges to $\wick{ |\nb Y_L|^2 }$ almost surely in $\mathcal{D}' (\R^2)$. 

It remains to that, given any $f \in \mathcal{D} (\R^2)$,
\begin{align*}
I_\eps (z_1, z_2) \deff \int_{\R^2} \sum_{k_1, k_2 \in \Z^2} (\rho_\eps *_{\R^2} \nb G) (x - z_1 + L k_1) \cdot (\rho_\eps *_{\R^2} \nb G) (x - z_2 + L k_2) \cj{f (x)} dx
\end{align*}

\noi
converges in $L^2 ( [-\frac{L}{2}, \frac{L}{2})^2 \times [-\frac{L}{2}, \frac{L}{2})^2 )$ to 
\begin{align*}
I (z_1, z_2) \deff \int_{\R^2} \sum_{k_1, k_2 \in \Z^2} \nb G (x - z_1 + L k_1) \cdot \nb G (x - z_2 + L k_2) \cj{f (x)} dx .
\end{align*}

\noi
This convergence then implies that $(\wick{ |\nb Y_{L, \eps}|^2 } \, f)$ converges in $L^2 (\Om)$ to the right-hand-side of \eqref{YL2_sto2}, which then implies the expression \eqref{YL2_sto2} in view of the almost sure convergence of $\wick{ |\nb Y_{L, \eps}|^2 }$ to $\wick{ |\nb Y_{L}|^2 }$ in $\mathcal{D}' (\R^2)$. 

To show the convergence of $I_\eps$ to $I$, we first observe that since $z_1, z_2 \in [ -\frac{L}{2}, \frac{L}{2} )^2$ and $\rho_\eps$, $G$, and $f$ have compact supports, we must have $|k_1|, |k_2| \leq K$ for some $K > 0$ in order to get non-zero contributions. We also note that $\nb G \in L^p (\R^2)$ for $1 \leq p < 2$. Thus, by Young's convolution inequality and H\"older's inequality, we have
\begin{align*}
&\| I (z_1, z_2) \|_{L_{z_1, z_2}^2 ( [ -\frac{L}{2}, \frac{L}{2} )^2 \times [ -\frac{L}{2}, \frac{L}{2} )^2 )}^2 \\
&\les_{L} \bigg\| \sum_{\substack{k_1 \in \Z^2 \\ |k_1| \leq K}} \nb G (\cdot + L k_1) \bigg\|_{L^{\frac{4}{3}} (\R^2)}^2 \bigg\| \sum_{\substack{k_2 \in \Z^2 \\ |k_2| \leq K}} \nb G (\cdot + L k_2) \bigg\|_{L^{\frac{4}{3}} (\R^2)}^2 \| f \|_{L^\infty (\R^2)}^2 \\
&\les_{K} \| \nb G \|_{L^{\frac 43} (\R^2)}^4  \| f \|_{L^\infty (\R^2)}^2 .
\end{align*}

\noi
Similarly, by Young's convolution inequality and the fact that $\| \rho_\eps \|_{L^1 (\R^2)} = 1$, we have
\begin{align*}
\| &I_\eps (z_1, z_2) - I (z_1, z_2) \|_{L_{z_1, z_2}^2 ( [ -\frac{L}{2}, \frac{L}{2} )^2 \times [ -\frac{L}{2}, \frac{L}{2} )^2 )}^2 \\
&\les_{L, K} \| \rho_\eps *_{\R^2} \nb G - \nb G \|_{L^{\frac 43} (\R^2)}^2 \big( \| \rho_\eps *_{\R^2} \nb G \|_{L^{\frac 43} (\R^2)} + \| \nb G \|_{L^{\frac 43} (\R^2)} \big)^2 \| f \|_{L^\infty (\R^2)}^2 \\
&\les \| \rho_\eps *_{\R^2} \nb G - \nb G \|_{L^{\frac 43} (\R^2)}^2 \| \nb G \|_{L^{\frac 43} (\R^2)}^2 \| f \|_{L^\infty (\R^2)}^2 ,
\end{align*}

\noi
which converges to 0 as $\eps \to 0$ since $\{ \rho_\eps \}_{\eps > 0}$ is an approximation identity. This finishes the proof of the lemma.
\end{proof}

We also need the following convolution lemma.

\begin{lemma}
\label{LEM:conv}
Let $0 < \al, \be < 2$. Then, for any $a \neq 0$, we have
\begin{align*}
\int_{\R^2} \frac{\ind_{\{ |x - a| \leq 1 \}} \ind_{\{ |x| \leq 1 \}}}{|x - a|^\al |x|^\be} dx \les 
\begin{cases}
1 + |a|^{2 - \al - \be} & \text{if } \al + \be \neq 2 \\ 
\log ( 2 + |a| ) & \text{if } \al + \be = 2 .
\end{cases}
\end{align*}
\end{lemma}

\begin{proof}
We write
\begin{align*}
\int_{\R^2} \frac{\ind_{\{ |x - a| \leq 1 \}} \ind_{\{ |x| \leq 1 \}}}{|x - a|^\al |x|^\be} dx = \textup{I}_1 + \textup{I}_2 + \textup{I}_3 ,
\end{align*}

\noi
where $\textup{I}_1$ denotes the contribution from $\{ |x - a| \leq \frac{|a|}{2} \}$, $\textup{I}_2$ denotes the contribution from $\{ |x| \leq \frac{|a|}{2} \}$, and $\textup{I}_3$ denotes the contribution from $\{ |x - a| > \frac{|a|}{2}, |x| > \frac{|a|}{2} \}$. For $\textup{I}_1$, by the triangle inequality, we have $|x| \geq |a| - |x - a| \geq \frac{|a|}{2}$, so that
\begin{align*}
\textup{I}_1 \les |a|^{- \be} \int_{\{ |x - a| \leq \frac{|a|}{2} \}} \frac{1}{|x - a|^\al} dx \les |a|^{2 - \al - \be} .
\end{align*}

\noi
For $\textup{I}_2$, by the triangle inequality, we have $|x - a| \geq |a| - |x| \geq \frac{|a|}{2}$, so that
\begin{align*}
\textup{I}_2 \les |a|^{- \al} \int_{\{ |x| \leq \frac{|a|}{2} \}} \frac{1}{|x|^\be} dx \les |a|^{2 - \al - \be} .
\end{align*}

\noi
For $\textup{I}_3$, by the condition $|x - a| > \frac{|a|}{2}$ and the triangle inequality, we have
\begin{align*}
|x - a| > \frac{3}{8}|a| + \frac 14 |x - a| \geq \frac 14 |x| + \frac 18 |a| > \frac 14 |x| .
\end{align*}

\noi
Thus, we have
\begin{align*}
\textup{I}_3 \les \int_{\{ \frac{|a|}{2} < |x| \leq 1 \}} \frac{1}{|x|^{\al + \be}} dx \les \begin{cases}
1 + |a|^{2 - \al - \be} & \text{if } \al + \be \neq 2 \\ 
\log ( 2 + |a| ) & \text{if } \al + \be = 2 .
\end{cases}
\end{align*}

\noi
Thus, we have finished the proof.
\end{proof}

We are now ready to prove the following proposition on the regularity and convergence property of stochastic objects.

\begin{proposition}
\label{PROP:YconvL2}
Let $0 < s < 1$, $\mu > 0$, and $\varphi \in C_c^\infty (\R^2)$. Then, for any $L \geq 1$ and $0 < \dl < \mu$, there exists $\Om' \subset \Om$ with full probability measure such that for any $\om \in \Om'$, the following estimates hold.

\smallskip \noi
\textup{(i)} We have 
\begin{align*}
&\| \varphi *_{\R^2} \xi \|_{\C_{- \mu}^s (\R^2)} + \| \varphi *_{\R^2} \xi_L \|_{\C_{- \mu}^s (\R^2)} \les_\om 1, \\
&\| \nb Y \|_{\C_{- \mu}^{s - 1} (\R^2)} +  \| \nb Y_L \|_{\C_{- \mu}^{s - 1} (\R^2)} \les_\om 1, \\
&\| \wick{|\nb Y|^2} \|_{\C_{- \mu}^{s - 1} (\R^2)} + \| \wick{|\nb Y_L|^2} \|_{\C_{- \mu}^{s - 1} (\R^2)} \les_\om 1.
\end{align*}

\smallskip \noi
\textup{(ii)} We have
\begin{align}
&\| \varphi *_{\R^2} \xi - \varphi *_{\R^2} \xi_L \|_{\C_{- \mu}^s (\R^2)} \les_{\om} L^{- \dl} , \label{YconvL1} \\
&\| \nb Y - \nb Y_L \|_{\C_{- \mu}^{s - 1} (\R^2)} \les_{\om} L^{- \dl} , \label{YconvL2} \\
&\big\| \wick{|\nb Y|^2} - \wick{|\nb Y_L|^2} \big\|_{\C_{- \mu}^{s - 1} (\R^2)} \les_{\om} L^{- \dl} . \label{YconvL3}
\end{align}
\end{proposition}

\begin{proof}
The bounds for $\varphi *_{\R^2} \xi$, $\nb Y$, and $\wick{|\nb Y|^2}$ in part (i) follow from \cite{HL15} along with Lemma~\ref{LEM:varphi}. The uniform-in-$L$ bounds for $\varphi *_{\R^2} \xi_L$, $\nb Y_L$, and $\wick{|\nb Y_L|^2}$ in part (i) will follow directly from part (ii). 
Thus, we focus on proving part (ii).

We first show \eqref{YconvL1} and \eqref{YconvL2} by first showing that
\begin{align}
\| \xi - \xi_L \|_{\C_{- \mu}^{s - 2} (\R^2)} \les_{\omega, \dl} L^{- \dl} 
\label{conv_goal1}
\end{align}

\noi
for almost sure $\om \in \Om$. Let $N \geq 1$ be a dyadic number and $\eta_N (\cdot) = N^2 \eta (N \cdot)$ be the kernel of $\Dl_N$ in \eqref{DN}. 
By the definition in \eqref{xiL}, we see that for any $x \in \R^2$, we have
\begin{align*}
\E \big[ | \Dl_N (\xi - \xi_L) (x) |^2 \big] 
&= N^4 \E \big[ \big| \big( \xi - \xi_L, \eta (N (x - \cdot)) \big)_{\R^2} \big|^2 \big] \\
&= N^4 \int_{\R^2} \Big| \eta (N (x - z)) - \sum_{k \in \Z^2} \eta (N (x - z + L k)) \ind_{[- \frac{L}{2}, \frac{L}{2})^2} (z) \Big|^2 dz ,
\end{align*}

\noi
so that we get
\begin{align}
\begin{split}
\E &\big[ | \Dl_N (\xi - \xi_L) (x) |^2 \big] \\
&\les N^4 \int_{\R^2 \setminus [- \frac{L}{2}, \frac{L}{2})^2} | \eta (N(x - z)) |^2 dz + N^4 \int_{[- \frac{L}{2}, \frac{L}{2})^2} \Big| \sum_{k \in \Z^2 \setminus \{0\}} \eta (N (x - z + Lk)) \Big|^2 dz \\
&\deff \textup{I}_1 + \textup{I}_2 .
\end{split}
\label{convL2-0}
\end{align}

\noi
For $\textup{I}_1$, we use the fact that $|z| \geq \frac{L}{2}$ to obtain
\begin{align}
\begin{split}
\textup{I}_1 &\les L^{-2 \dl} N^4 \int_{\R^2} \jb{z}^{2 \dl} | \eta (N (x - z)) |^2 dz \\
&\les L^{-2 \dl} N^4 \jb{x}^{2 \dl} \int_{\R^2} \jb{x - z}^{2 \dl} | \eta (N (x - z)) |^2 dz \\
&\les L^{-2 \dl} N^2 \jb{x}^{2 \dl} .
\end{split}
\label{convL2-1}
\end{align}

\noi
For $\textup{I}_2$, since $\eta$ is a Schwartz function, we have
\begin{align*}
\Big| \sum_{k \in \Z^2 \setminus \{0\}} \eta (N (x - z + Lk)) \Big| \les \sum_{k \in \Z^2} \frac{1}{\jb{N (x - z + Lk)}^{10}} \les 1
\end{align*}

\noi
uniform in $x$ and $z$, so that we use a similar treatment as that for $\textup{I}_1$ to obtain
\begin{align}
\begin{split}
\text{I}_2 &\les N^4 \int_{[- \frac{L}{2}, \frac{L}{2})^2} \sum_{k \in \Z^2 \setminus \{0\}} | \eta (N (x - z + L k)) | d z \\
&= N^4 \int_{\R^2 \setminus [ - \frac{L}{2}, \frac{L}{2} )^2} | \eta (N (x - z)) | dz \\
&\les L^{-2 \dl} N^2 \jb{x}^{2 \dl} .
\end{split}
\label{convL2-2}
\end{align}

\noi
Thus, for any $p \geq 2$, by the Gaussian hypercontractivity, \eqref{convL2-0}, \eqref{convL2-1}, and \eqref{convL2-2}, we have
\begin{align*}
\E \Big[ \| \Dl_N (\xi - \xi_L) \|_{L_{- \mu}^p (\R^2)}^p \Big] &= \int_{\R^2} \jb{x}^{- \mu p} \E \big[ | \Dl_N (\xi - \xi_L) (x) |^p \big] dx \\
&\les_p \int_{\R^2} \jb{x}^{- \mu p} \Big( \E \big[ |\Dl_N (\xi - \xi_L) (x)|^2 \big] \Big)^{\frac{p}{2}} dx \\
&\les L^{- \dl p} N^p \int_{\R^2} \jb{x}^{- (\mu - \dl) p} dx ,
\end{align*}

\noi
so that by taking $p$ large enough such that $(\mu - \dl) p > 2$, we get
\begin{align}
\E \Big[ \| \Dl_N (\xi - \xi_L) \|_{L_{- \mu}^p (\R^2)}^p \Big] 
\les L^{- \dl p} N^p .
\label{convL2-3}
\end{align}

\noi
We now use the Besov embedding (Lemma~\ref{LEM:embex}) and \eqref{convL2-3} to obtain
\begin{align}
\E \Big[ \| \xi - \xi_L \|_{\C_{- \mu}^{s - 2} (\R^2)}^p \Big] \les \sum_{\substack{N \geq 1 \\ \text{dyadic}}} N^{(s - 2 + \frac{2}{p}) p} \E \Big[ \| \Dl_N (\xi - \xi_L) \|_{L^p_{- \mu} (\R^2)}^p \Big] \les_{p} L^{- \dl p} ,
\label{convL2-4}
\end{align}

\noi
where we need to take $p$ sufficiently large such that $s - 1 + \frac{2}{p} < 0$. 
This gives the desired estimate \eqref{conv_goal1}.
By using Lemma~\ref{LEM:varphiex} and \eqref{conv_goal1}, we obtain
\begin{align*}
\| \varphi *_{\R^2} \xi - \varphi *_{\R^2} \xi_L \|_{\C_{- \mu}^s (\R^2)} \les \| \xi - \xi_L \|_{\C_{- \mu}^{s - 2} (\R^2)} \les_{\omega} L^{- \dl} ,
\end{align*}

\noi
which proves \eqref{YconvL1}. By using \eqref{defY}, Lemma~\ref{LEM:G}, and \eqref{conv_goal1}, we obtain
\begin{align*}
\| \nb Y - \nb Y_L \|_{\C_{- \mu}^{s - 1} (\R^2)} \les \| \xi - \xi_L \|_{\C_{- \mu}^{s - 2} (\R^2)} \les_{\omega} L^{- \dl} ,
\end{align*}

\noi
which proves \eqref{YconvL2}.

Finally, we consider \eqref{YconvL3}. Let $N \geq 1$ be a dyadic number. From \eqref{Y2_sto} and Lemma~\ref{LEM:YL2_sto}, we see that for any $x \in \R^2$, we have
\begin{align*}
\Dl_N ( \wick{ |\nb Y|^2 } - \wick{ |\nb Y_L|^2 } ) (x) &= N^2 \int_{\R^2} \int_{\R^2} \bigg( \int_{\R^2} \eta ( N (x - y) ) \Big( \nb G (y - z_1) \cdot \nb G (y - z_2) \\
&\quad - \sum_{k_1, k_2 \in \Z^2} \nb G (y - z_1 + L k_1) \ind_{[ - \frac{L}{2}, \frac{L}{2} )^2} (z_1)  \\
&\qquad \cdot \nb G (y - z_2 + L k_2)  \ind_{[ - \frac{L}{2}, \frac{L}{2} )^2} (z_2) \Big) dy \bigg) \xi (d z_1) \xi (d z_2) .
\end{align*}

\noi
Thus, we obtain
\begin{align}
\E \big[ | \Dl_N ( \wick{ |\nb Y|^2 } - \wick{ |\nb Y_L|^2 } ) (x) |^2 \big] \les \textup{I}_3 + \textup{I}_4 + \textup{I}_5 + \textup{I}_6 ,
\label{conv4-0}
\end{align}

\noi
where
\begin{align*}
\textup{I}_3 &\deff N^4 \int_{\R^2} \int_{\R^2 \setminus [ - \frac{L}{2}, \frac{L}{2} )^2} \bigg| \int_{\R^2} \eta ( N (x - y) ) \nb G (y - z_1) \cdot \nb G (y - z_2) dy \bigg|^2 d z_2 d z_1 , \\
\textup{I}_4 &\deff N^4 \int_{\R^2} \int_{[ - \frac{L}{2}, \frac{L}{2} )^2} \bigg| \int_{\R^2} \eta ( N (x - y) ) \sum_{k_2 \in \Z^2 \setminus \{ 0 \}} \nb G (y - z_1) \cdot \nb G (y - z_2 + L k_2) dy \bigg|^2 d z_2 d z_1 , \\
\textup{I}_5 &\deff N^4 \int_{\R^2 \setminus [ - \frac{L}{2}, \frac{L}{2} )^2} \int_{[ - \frac{L}{2}, \frac{L}{2} )^2} \bigg| \int_{\R^2} \eta ( N (x - y) ) \\
&\qquad \times \sum_{k_2 \in \Z^2} \nb G (y - z_1) \cdot \nb G (y - z_2 + L k_2) dy \bigg|^2 d z_2 d z_1 , \\
\textup{I}_6 &\deff N^4 \int_{[ - \frac{L}{2}, \frac{L}{2} )^2} \int_{[ - \frac{L}{2}, \frac{L}{2} )^2} \bigg| \int_{\R^2} \eta ( N (x - y) ) \\
&\qquad \times \sum_{k_1 \in \Z^2 \setminus \{0\}} \sum_{k_2 \in \Z^2} \nb G (y - z_1 + L k_1) \cdot \nb G (y - z_2 + L k_2) dy \bigg|^2 d z_2 d z_1 .
\end{align*}

\noi
For $\textup{I}_3$, since $|z_2| \geq \frac{L}{2}$, we have
\begin{align}
L^\dl \les \jb{z_2}^\dl \les \jb{y - z_2}^\dl \jb{y}^\dl \les \jb{y - z_2}^\dl \jb{x - y}^\dl \jb{x}^\dl .
\label{expL1}
\end{align}

\noi
Thus, by \eqref{expL1}, we compute that
\begin{align*}
\textup{I}_3 &\les L^{-2 \dl} N^4 \int_{\R^2} \int_{\R^2} \bigg( \int_{\R^2} | \eta (N (x - y)) | \frac{\ind_{\{ |y - z_1| \leq \frac 14 \}} \ind_{\{ |y - z_2| \leq \frac 14 \}} \jb{y - z_2}^\dl \jb{y}^\dl}{|y - z_1| |y - z_2|} dy \bigg)^2 d z_2 d z_1 \\
&\les L^{-2 \dl} N^4 \jb{x}^{2 \dl} \int_{\R^2} \int_{\R^2}  \bigg( \int_{\R^2} \frac{\ind_{\{ |y_1 - z_1| \leq \frac 14 \}} \ind_{ \{ |y_2 - z_1| \leq \frac 14 \} }}{|y_1 - z_1| |y_2 - z_1|} d z_1 \bigg) \\
&\qquad \qquad \times \bigg( \int_{\R^2} \frac{\ind_{\{ |y_1 - z_2| \leq \frac 14 \}} \ind_{ \{ |y_2 - z_2| \leq \frac 14 \} } \jb{y_1 - z_2}^{\dl} \jb{y_2 - z_2}^{\dl}}{|y_1 - z_2| |y_2 - z_2|} d z_2 \bigg)  \\
&\qquad \qquad \times \jb{x - y_1}^{\dl} |\eta (N (x - y_1))| \jb{x - y_2}^{\dl} |\eta (N (x - y_2))| d y_2 d y_1 ,
\end{align*}

\noi
so that by Lemma~\ref{LEM:conv} and the fact that $\log (2 + |y_1 - y_2|) \leq \log (2 + |y_1 - z_j| + |y_2 - z_j|) \les 1$ given $|y_1 - z_j| \leq \frac 14$ and $|y_2 - z_j| \leq \frac 14$ for $j = 1, 2$, we obtain
\begin{align}
\begin{split}
\textup{I}_3 &\les L^{-2 \dl} N^4 \jb{x}^{2 \dl} \int_{\R^2} \int_{\R^2} \jb{x - y_1}^\dl |\eta (N (x - y_1))| \jb{x - y_2}^\dl |\eta (N (x - y_2))| d y_2 d y_1 \\
&\les L^{-2 \dl} N^4 \jb{x}^{2 \dl} \int_{\R^2} \jb{N y_1}^\dl |\eta (N y_1)| dy_1 \int_{\R^2} \jb{N y_2}^\dl |\eta (N y_2)| dy_2 \\
&\les L^{-2 \dl} \jb{x}^{2 \dl} .
\end{split}
\label{conv4-1}
\end{align}

\noi
For $\textup{I}_4$, we proceed in the same way as that for $\textup{I}_3$ and then encounter the term
\begin{align*}
\textup{I}' = \int_{[- \frac{L}{2}, \frac{L}{2} )^2} \sum_{k, k' \in \Z^2 \setminus \{0\}} \frac{\ind_{\{ |y_1 - z + L k| \leq \frac 14 \}} \ind_{\{ |y_2 - z + L k'| \leq \frac 14 \}}}{|y_1 - z + L k| |y_2 - z + L k'|} d z .
\end{align*}

\noi
Note that since $z \in [ - \frac{L}{2}, \frac{L}{2} )^2$, $|y_1 - z + Lk| \leq \frac 14$, $|y_2 - z + Lk| \leq \frac 14$, and $k, k' \in \Z^2 \setminus \{0\}$, we must have $|y_1| \geq \frac{L}{4}$ and $|y_2| \geq \frac{L}{4}$ in order to get a nonzero contribution, which gives
\begin{align}
L^{2 \dl'} \leq \jb{y_1}^{\dl'} \jb{y_2}^{\dl'} \leq \jb{x - y_1}^{\dl'} \jb{x - y_2}^{\dl'} \jb{x}^{2 \dl'} 
\label{expL2}
\end{align} 

\noi
for $0 < \dl < \dl' < \mu$. Moreover, from the same conditions, for any fixed $y_1 \in \R^2$ (and $y_2 \in \R^2$), we must have $|y_1 + L k| < L$ (and $|y_2 + L k'| < L$) to get a nonzero contribution, which gives at most $O(1)$ choices of $k$ (and $k'$, respectively).
Thus, by \eqref{expL2} and Lemma~\ref{LEM:conv}, we compute that
\begin{align*}
\begin{split}
\textup{I}' &\leq L^{- 2 \dl'} \jb{x - y_1}^{\dl'} \jb{x - y_2}^{\dl'} \jb{x}^{2 \dl'} \textup{I}' \\
&\les L^{- 2 \dl'} \jb{x - y_1}^{\dl'} \jb{x - y_2}^{\dl'} \jb{x}^{2 \dl'} \\
&\quad \times \sum_{k, k' \in \Z^2 \setminus \{0\}} \log (2 + |y_1 + L k - y_2 - L k'|) \ind_{\{ |y_1 + L k| < L \}} \ind_{\{ |y_2 + L k'| < L \}} ,
\end{split}
\end{align*}

\noi
so that we get
\begin{align}
\textup{I}' \les_{\dl, \dl'} L^{- 2 \dl} \jb{x - y_1}^{\dl'} \jb{x - y_2}^{\dl'} \jb{x}^{2 \dl'} .
\label{Ip_comput}
\end{align}

\noi
Then, by using a similar treatment as that for $\textup{I}_3$ along with \eqref{Ip_comput}, we obtain
\begin{align}
\textup{I}_4 \les L^{-2 \dl} \jb{x}^{2 \dl'} .
\label{conv4-2}
\end{align}

\noi
For $\textup{I}_5$, we proceed in the same way as that for $\textup{I}_3$ and then encounter the term 
\begin{align*}
\textup{I}'' \deff \int_{[ - \frac{L}{2}, \frac{L}{2} )^2} \sum_{k, k' \in \Z^2} \frac{\ind_{ \{ |y_1 - z + L k| \leq \frac 14 \} } \ind_{ \{ |y_2 - z + L k'| \leq \frac 14 \} }}{|y_1 - z + L k| |y_2 - z + L k'|} dz .
\end{align*}

\noi
From the conditions that $z \in [- \frac{L}{2}, \frac{L}{2})^2$ and $|y_1 - z + L k| \leq \frac 14$ (and $|y_2 - z + L k'| \leq \frac 14$), for any fixed $y_1 \in \R^2$ (and $y_2 \in \R^2$), we must have $|y_1 + L k| < L$ (and $|y_2 + L k'| < L$) to get a nonzero contribution, which gives at most $O(1)$ choices of $k$ (and $k'$, respectively). Thus, by Lemma~\ref{LEM:conv}, we obtain
\begin{align}
\textup{I}'' \les \sum_{k, k' \in \Z^2} \log (2 + |y_1 + L k - y_2 - L k'|) \ind_{\{ |y_1 + L k| < L \}} \ind_{\{ |y_2 + L k'| < L \}} \les \log (2 + L).
\label{conv4h2}
\end{align}

\noi
Then, by using a similar treatment as that for $\textup{I}_3$ along with \eqref{conv4h2}, we obtain
\begin{align}
\textup{I}_5 \les L^{-2 \dl} \jb{x}^{2 \dl} .
\label{conv4-3}
\end{align}

\noi
For $\textup{I}_{6}$, we use similar steps as above along with \eqref{Ip_comput} and \eqref{conv4h2} to obtain
\begin{align}
\textup{I}_{6} \les L^{-2 \dl} \jb{x}^{2 \dl'} 
\label{conv4-4}
\end{align}

\noi
for $0 < \dl < \dl' < \mu$. Combining \eqref{conv4-0}, \eqref{conv4-1}, \eqref{conv4-2}, \eqref{conv4-3}, and \eqref{conv4-4}, we obtain
\begin{align*}
\E \big[ | \Dl_N ( \wick{ |\nb Y|^2 } - \wick{ |\nb Y_L|^2 } ) (x) |^2 \big] \les L^{-2 \dl} \jb{x}^{2 \dl'} 
\end{align*}

\noi
for $0 < \dl < \dl' < \mu$. Thus, by proceeding as in the proof of \eqref{YconvL1} using the Besov embedding and the Gaussian hypercontractivity,
we obtain the desired estimate \eqref{YconvL3}.
\end{proof}

\subsection{Large torus convergence of the $L$-periodic PAM}

Given $s > 0$, $\ta > 0$, $T > 0$, and $\nu \in \R$, we define the space $\mathcal{X}_{\nu, T}^{s, \ta}$ via the norm
\begin{align*}
\| u \|_{\mathcal{X}_{\nu, T}^{s, \ta}} \deff \sup_{t \in (0, T]} t^{\ta} \| u (t, \cdot) \|_{\wt \C_{\nu - t}^s (\R^2)} ,
\end{align*}

\noi
where we recall that $\C_\nu^s$ is the notation for the function space with an exponential weight defined in Subsection~\ref{SUB:exp}. Our goal is to prove the following statement.
\begin{theorem}
\label{THM:heat}
Let $\mu_0 > 2$, $T \geq 1$, $v_0 \in L_{\mu_0}^\infty (\R^2)$, and $v_{0, L}$ be as defined in \eqref{v0L2} given any $L \geq 1$.

\smallskip \noi
\textup{(i)} For any $1 < s < 2$, there exists a unique solution $v \in \mathcal{X}_{0, T}^{s, \frac{s}{2}}$ to PAM \eqref{PAM} with initial data $v_0$. Also, for any $1 < s < 2$, there exists a unique solution $v_L \in \mathcal{X}_{0, T}^{s, \frac{s}{2}}$ to the $L$-periodic PAM \eqref{PAML} with initial data $v_{0, L}$.

\smallskip \noi
\textup{(ii)} For any $1 < s < 2$ and $0 < \ta < \frac{2 - s}{2}$, $v_L$ converges to $v$ in $\mathcal{X}_{0, T}^{s, \frac{s}{2} + \ta}$ in probability as $L \to \infty$.
\end{theorem}

\begin{proof}
We mainly follow \cite{HL15}. Let $\Om' \subset \Om$ be the event with full probability measure such that Proposition~\ref{PROP:YconvL2} holds (with some parameters to be chosen later) and we fix $\om \in \Om'$. 

We first focus on the well-posedness statement in part (i). For this purpose, we work on the Duhamel formulation of \eqref{PAM}:
\begin{align*}
v (t) = \Gamma_{v_0} [v] (t) \deff e^{t \Dl} v_0 + \int_0^t e^{(t - t') \Dl} \Big( - 2 \nb Y \cdot \nb v (t') + \wt{\wick{|\nb Y|^2}} \, v (t') \Big) dt' .
\end{align*}

\noi
Let $\mu > 0$ and $s_0 > 0$ be small enough such that $\frac{s + s_0}{2} + \frac{\mu}{\vartheta} < 1$, where $0 < \vartheta < 1$ is a fixed parameter defined in Subsection~\ref{SUB:exp}. From Lemma~\ref{LEM:Sch}, Lemma~\ref{LEM:embex}~(ii), and Lemma~\ref{LEM:emb}~(ii), we have
\begin{align}
\| e^{\cdot \Dl} v_0 \|_{\mathcal{X}_{0, T}^{s, \frac{s}{2}}} = \sup_{t \in (0, T]} t^{\frac{s}{2}} \| e^{t \Dl} v_0 \|_{\wt \C_{-t}^s (\R^2)} \les \sup_{t \in (0, T]} \| v_0 \|_{\wt \C_{- t}^0 (\R^2)} \leq \| v_0 \|_{\C^0 (\R^2)} \les \| v_0 \|_{L^\infty (\R^2)}.
\label{heat1}
\end{align}

\noi
By Lemma~\ref{LEM:prodex2}, Proposition~\ref{PROP:YconvL2}~(i), and Lemma~\ref{LEM:embex}~(iii) along with the fact that $1 + 2 s_0 < s$, we have for any $t > 0$ and $0 \leq t' \leq t$ that
\begin{align*}
\begin{split}
\| \nb Y \cdot \nb v (t') \|_{\wt \C_{-t}^{- s_0} (\R^2)} 
&\les (t - t')^{- \frac{\mu}{\vartheta}} \| \nb Y \|_{\C_{- \mu}^{- s_0} (\R^2)} \| \nb v (t') \|_{\wt \C_{- t'}^{2 s_0} (\R^2)} \\
&\les_\om (t - t')^{- \frac{\mu}{\vartheta}} (t')^{- \frac{s}{2}} \| v \|_{\mathcal{X}_{0, T}^{s, \frac{s}{2}}}
\end{split}
\end{align*}

\noi
and
\begin{align*}
\begin{split}
\big\| \wt{\wick{|\nb Y|^2}} \, v (t') \big\|_{\wt \C_{-t}^{- s_0} (\R^2)} 
&\les (t - t')^{- \frac{\mu}{\vartheta}} \big\| \wt{\wick{|\nb Y|^2}} \big\|_{\C_{- \mu}^{- s_0} (\R^2)} \| v (t') \|_{\wt \C_{- t'}^{2 s_0} (\R^2)} \\
&\les_\om (t - t')^{- \frac{\mu}{\vartheta}} (t')^{- \frac{s}{2}} \| v \|_{\mathcal{X}_{0, T}^{s, \frac{s}{2}}} ,
\end{split}
\end{align*}

\noi
so that by Lemma~\ref{LEM:Sch} and the above two estimates, we get
\begin{align}
\begin{split}
\bigg\| &\int_0^\cdot e^{(\cdot - t') \Dl} \Big( - 2 \nb Y \cdot \nb v (t') + \wt{\wick{|\nb Y|^2}} \, v (t') \Big) dt' \bigg\|_{\mathcal{X}_{0, T}^{s, \frac{s}{2}}} \\
&\les  \sup_{t \in (0, T]} t^{\frac{s}{2}} \int_0^t (t - t')^{- \frac{s + s_0}{2}} \Big( \| \nb Y \cdot \nb v (t') \|_{\wt \C_{-t}^{- s_0} (\R^2)} + \big\| \wt{\wick{|\nb Y|^2}} \, v (t') \big\|_{\wt \C_{-t}^{- s_0} (\R^2)} \Big) dt' \\
&\les_\om T^{1 - \frac{s + s_0}{2} - \frac{\mu}{\vartheta}} \| v \|_{\mathcal{X}_{0, T}^{s, \frac{s}{2}}} .
\end{split}
\label{heat2}
\end{align}

\noi
Combining \eqref{heat1} and \eqref{heat2}, we get
\begin{align*}
\| \Gamma_{v_0} [v] \|_{\mathcal{X}_{0, T}^{s, \frac{s}{2}}} \les_\om \| v_0 \|_{L^\infty (\R^2)} + T^{1 - \frac{s + s_0}{2} - \frac{\mu}{\vartheta}} \| v \|_{\mathcal{X}_{0, T}^{s, \frac{s}{2}}} .
\end{align*}

\noi
Using similar steps, we obtain
\begin{align*}
\| \Gamma_{v_0} [v_1] - \Gamma_{v_0} [v_2] \|_{\mathcal{X}_{0, T}^{s, \frac{s}{2}}} \les_\om T^{1 - \frac{s + s_0}{2} - \frac{\mu}{\vartheta}} \| v_1 - v_2 \|_{\mathcal{X}_{0, T}^{s, \frac{s}{2}}} .
\end{align*}

\noi
Thus, by choosing $T_0 > 0$ to be sufficiently small, we use a standard contraction argument to get a solution $v \in \mathcal{X}_{0, T_0}^{s, \frac{s}{2}}$ to PAM \eqref{PAM}. Since $T_0$ can be chosen to be independent of the initial data and also the exponential weight factor, we can then use $v (T_0) \in \wt \C_{- T_0}^s (\R^2)$ as the new initial data and run the above contraction argument (with slight modifications of \eqref{heat1}) in the space $\mathcal{X}_{- T_0, T_0}^{s, \frac{s}{2}}$ on the time interval $[T_0, 2T_0]$. This gives a unique solution $v$ on the time interval $[0, 2T_0]$ in the space $\mathcal{X}_{0, 2T_0}^{s, \frac{s}{2}}$. By iterating this procedure, we obtain a unique solution $v \in \mathcal{X}_{0, T}^{s, \frac{s}{2}}$ for any given $T \geq 1$. Well-posedness for the $L$-periodic PAM \eqref{PAML} follows from the same steps thanks to the fact that $v_{0, L} \in L^\infty (\R^2)$ from Lemma~\ref{LEM:per3} and the uniform bounds of the $L$-periodic stochastic objects in Proposition~\ref{PROP:YconvL2}~(i). In particular, we have
\begin{align}
\| v \|_{\mathcal{X}_{0, T}^{s, \frac{s}{2}}} + \| v_L \|_{\mathcal{X}_{0, T}^{s, \frac{s}{2}}} \les_{\om, T} \| v_{0} \|_{L^\infty (\R^2)} + \| v_{0, L} \|_{L^\infty (\R^2)} \les \| v_0 \|_{L_{\mu_0}^\infty (\R^2)} .
\label{vvL}
\end{align}

We now show the large torus convergence in part (ii). Let $0 < \ta < \frac{2 - s}{2}$. 
From Lemma~\ref{LEM:Sch}, the fact that $e^{- t \jb{x}^\vartheta} \les t^{- \ta} \jb{x}^{- \ta \vartheta}$, Lemma~\ref{LEM:emb}~(ii), and Lemma~\ref{LEM:per3}, we have
\begin{align}
\begin{split}
\| e^{\cdot \Dl} (v_0 - v_{0, L}) \|_{\mathcal{X}_{0, T}^{s, \frac{s}{2} + \ta}} 
&\les \sup_{t \in (0, T]} t^{\ta} \| v_0 - v_{0, L} \|_{\wt \C_{- t}^0 (\R^2)} \\
&\les \| v_0 - v_{0, L} \|_{L_{- \ta \vartheta}^\infty (\R^2)} \\
&\les L^{- \ta \vartheta} \| v_0 \|_{L_{\mu_0}^\infty (\R^2)} .
\end{split}
\label{heat3-1}
\end{align}

\noi
Also, by using similar steps in \eqref{heat2} along with Proposition~\ref{PROP:YconvL2}~(ii), we have
\begin{align}
\begin{split}
\bigg\| &\int_0^\cdot e^{(\cdot - t') \Dl} \Big( - 2 \big( \nb Y \cdot \nb v (t') - \nb Y_L \cdot \nb v_L (t') \big) \\
&\quad + \big( \wt{\wick{|\nb Y|^2}} \, v (t') - \wt{\wick{|\nb Y_L|^2}} \, v_L (t') \big) \Big) dt' \bigg\|_{\mathcal{X}_{0, T}^{s, \frac{s}{2}}} \\
&\les_\om T^{1 - \frac{s + s_0}{2} - \frac{\mu}{\vartheta} + \ta} L^{- \frac{\mu}{2}} \Big( \| v \|_{\mathcal{X}_{0, T}^{s, \frac{s}{2}}} + \| v_L \|_{\mathcal{X}_{0, T}^{s, \frac{s}{2}}} \Big) 
+ T^{1 - \frac{s + s_0}{2} - \frac{\mu}{\vartheta}} \| v - v_L \|_{X_{0, T}^{s, \frac{s}{2} + \ta}} .
\end{split}
\label{heat3-2}
\end{align}

\noi
Thus, by taking the difference of the Duhamel formulations of \eqref{PAM} and \eqref{PAML} and using \eqref{heat3-1}, \eqref{heat3-2}, and \eqref{vvL}, we obtain
\begin{align}
\begin{split}
\| &v - v_L \|_{\mathcal{X}_{0, T}^{s, \frac{s}{2} + \ta}} \\
&\les_{\om, T} L^{- \ta \vartheta} \| v_0 \|_{L^\infty_\mu (\R^2)} + T^{1 - \frac{s + s_0}{2} - \frac{\mu}{\vartheta} + \ta} L^{- \frac{\mu}{2}} \| v_0 \|_{L_{\mu_0}^\infty (\R^2)} + T^{1 - \frac{s + s_0}{2} - \frac{\mu}{\vartheta}} \| v - v_L \|_{X_{0, T}^{s, \frac{s}{2} + \ta}}.
\end{split}
\label{heat3}
\end{align}

\noi
Thus, by choosing $T_1 > 0$ to be sufficiently small, we obtain
\begin{align*}
\| v - v_L \|_{\mathcal{X}_{0, T_1}^{s, \frac{s}{2} + \ta}} \les_{\om, T} L^{- \ta \vartheta} \| v_0 \|_{L^\infty_\mu (\R^2)} + L^{- \frac{\mu}{2}} \| v_0 \|_{L_{\mu_0}^\infty (\R^2)} .
\end{align*}

\noi
Since $T_1$ can be chosen to be independent of the initial data and also the exponential weight factor, we can then iterate the procedure in \eqref{heat3} and obtain
\begin{align*}
\| v - v_L \|_{\mathcal{X}_{0, T}^{s, \frac{s}{2} + \ta}} \les_{\om, T} L^{- \ta \vartheta} \| v_0 \|_{L^\infty_\mu (\R^2)} + L^{- \frac{\mu}{2}} \| v_0 \|_{L_{\mu_0}^\infty (\R^2)}
\end{align*}

\noi
for any $T \geq 1$. This gives the desired convergence result.
\end{proof}

\begin{remark} \rm
\label{RMK:as2}
By slightly modifying the proof of Proposition~\ref{PROP:YconvL2}, we can prove in Theorem~\ref{THM:heat} almost sure convergence of $v_L$ to $v$ if we make the restriction that $L$ goes along positive integers. The key point here is to make the full probability measure set $\Om' \subset \Om$ in Proposition~\ref{PROP:YconvL2} independent of $L \in \N$. This can be done via
\begin{align*}
\E \Big[ \sup_{L \in \N} L^{\frac{\dl p}{2}} \| \xi - \xi_L \|_{\C_{- \mu}^{s - 2} (\R^2)}^p \Big] \leq \E \Big[ \sum_{L \in \N} L^{\frac{\dl p}{2}} \| \xi - \xi_L \|_{\C_{- \mu}^{s - 2} (\R^2)}^p \Big] \les \sum_{L \in \N} L^{- \frac{\dl p}{2}} < \infty
\end{align*}

\noi
for $p$ large enough such that $\dl p > 2$, where we used \eqref{convL2-4}. The other difference estimates in Proposition~\ref{PROP:YconvL2}~(ii) can be done in a similar way and can then be used to obtain uniform-in-$L$ bounds in Proposition~\ref{PROP:YconvL2}~(i), all of which hold on a full probability measure set $\Om' \subset \Om$ independent of $L \in \N$. It would be of interest to investigate if the convergence in Theorem~\ref{THM:heat} holds almost surely for $L$ going along the continuum.
\end{remark}

\begin{ackno} \rm
The authors would like to thank Tadahiro Oh for helpful discussions. The authors are also grateful to Bjoern Bringmann for pointing out a serious mistake in their first draft. R.L. was supported by the DFG through the Hausdorff Center for Mathematics under Germany's Excellence Strategy - GZ 2047/1, Project-ID 390685813, SFB 1060, Project-ID 211504053  and SFB 1720, Project-ID 539309657.
N.T. was partially supported by the ANR project Smooth ANR-22-CE40-0017 and by the European Research Council (ERC) under the European Union Horizon 2020 research and innovation programme (Grant agreement 101097172 - GEOEDP).
\end{ackno}

\end{document}